\def\DateTime{March 27th, 2021}
\def\Version{Version $1$}
\def\FrontCoverColor{white} 

\documentclass[a4paper,leqno,10pt,english]{smfbook}
\usepackage{smfthm,appendix,CJK}
\SwapTheoremNumbers

\usepackage[T2A,T1]{fontenc}
\usepackage[english,francais]{babel}
\usepackage[utf8]{inputenc}
\usepackage{latexsym,amscd,color}
\usepackage{amsmath,amsfonts,amssymb,mathrsfs,mathtools}
\usepackage{qsymbols,euscript,enumitem}
\usepackage{multirow,bigdelim}
\usepackage{bm,comment}
\usepackage{url,xspace}
\usepackage{verbatim}
\usepackage{tikz}
\usepackage{tikz-cd}
\usetikzlibrary{positioning,arrows,shapes}
\usepackage[normalem]{ulem}
\usepackage{manfnt} 
\usepackage{etoolbox}

\usepackage{hyperref}
\usepackage{xstring}

\makeindex

\theoremstyle{plain}
\author{Huayi CHEN}
\address{Universit\'e de Paris and Sorbonne Universit\'e, CNRS, IMJ-PRG, F-75006 Paris, France}
\email{huayi.chen@imj-prg.fr}
\author{Atsushi MORIWAKI}
\address{University of Kyoto}
\email{moriwaki@math.kyoto-u.ac.jp}
\title[Arithmetic intersection theory over adelic curves]{Arithmetic intersection theory over adelic curves}
\date{\DateTime, \Version}
\input xy
\xyoption{all}
\numberwithin{equation}{chapter}

\newcommand{\emptyinnprod}{\langle\kern.15em,\kern-.02em\rangle}
\newcommand{\hooklongrightarrow}{\lhook\joinrel\longrightarrow}
\renewcommand{\emptyset}{\varnothing}

\newcommand{\ZZ}{{\mathbb{Z}}}
\newcommand{\QQ}{{\mathbb{Q}}}
\newcommand{\RR}{{\mathbb{R}}}

\newcommand{\CC}{{\mathbb{C}}}
\newcommand{\PP}{{\mathbb{P}}}
\newcommand{\NN}{{\mathbb{N}}}

\newcommand{\A}{{\mathbb{A}}}
\newcommand{\OO}{{\mathcal{O}}}
\newcommand{\oo}{{\mathfrak{o}}}

\newcommand{\UUU}{{\mathscr{U}}}
\newcommand{\XXX}{{\mathscr{X}}}

\newcommand{\Proj}{\operatorname{Proj}}

\newcommand{\Spec}{\operatorname{Spec}}

\newcommand{\Supp}{\operatorname{Supp}}
\newcommand{\codim}{\operatorname{codim}}

\newcommand{\rank}{\operatorname{rk}}
\newcommand{\ord}{\operatorname{ord}}

\newcommand{\adeg}{\widehat{\operatorname{deg}}}

\newcommand{\an}{\operatorname{an}}
\newcommand{\fin}{\operatorname{fin}}

\newcommand{\ndot}{\raisebox{.4ex}{.}}

\newcommand{\rest}[2]{\left.{#1}\right\vert_{{#2}}}

\newcommand{\IP}{{\mathcal{IP}}}
\newcommand{\length}{\operatorname{length}}

\newcommand{\indic}{1\hspace{-0.25em}\mathrm{l}}
\newcommand\redsout{\bgroup\markoverwith{\textcolor{mred}{\rule[0.5ex]{2pt}{0.7pt}}}\ULon}
\def\colorsout#1{\bgroup\markoverwith{\textcolor{#1}{\rule[0.5ex]{2pt}{0.7pt}}}\ULon} 
\def\coloruline#1{\bgroup\markoverwith{\textcolor{#1}{\rule[-0.5ex]{2pt}{0.7pt}}}\ULon} 

\DeclareSymbolFont{bbold}{U}{bbold}{m}{n}
\DeclareMathSymbol{\bbalpha}{\mathord}{bbold}{"0B}
\DeclareMathSymbol{\bbbeta}{\mathord}{bbold}{"0C}
\DeclareMathSymbol{\bbgamma}{\mathord}{bbold}{"0D}
\DeclareMathSymbol{\bbdelta}{\mathord}{bbold}{"0E}
\DeclareMathSymbol{\bbespilon}{\mathord}{bbold}{"0F}
\DeclareMathSymbol{\bbzeta}{\mathord}{bbold}{"10}
\DeclareMathSymbol{\bbeta}{\mathord}{bbold}{"11}
\DeclareMathSymbol{\bbtheta}{\mathord}{bbold}{"12}
\DeclareMathSymbol{\bbiota}{\mathord}{bbold}{"13}
\DeclareMathSymbol{\bbkappa}{\mathord}{bbold}{"14}
\DeclareMathSymbol{\bblambda}{\mathord}{bbold}{"15}
\DeclareMathSymbol{\bbmu}{\mathord}{bbold}{"16}
\DeclareMathSymbol{\bbnu}{\mathord}{bbold}{"17}
\DeclareMathSymbol{\bbxi}{\mathord}{bbold}{"18}
\DeclareMathSymbol{\bbpi}{\mathord}{bbold}{"19}
\DeclareMathSymbol{\bbrho}{\mathord}{bbold}{"1A}
\DeclareMathSymbol{\bbsigma}{\mathord}{bbold}{"1B}
\DeclareMathSymbol{\bbtau}{\mathord}{bbold}{"1C}
\DeclareMathSymbol{\bbupsilon}{\mathord}{bbold}{"1D}
\DeclareMathSymbol{\bbphi}{\mathord}{bbold}{"1E}
\DeclareMathSymbol{\bbchi}{\mathord}{bbold}{"1F}
\DeclareMathSymbol{\bbpsi}{\mathord}{bbold}{"20}

\definecolor{ruby}{rgb}{0.88, 0.07, 0.37}
\definecolor{coolblack}{rgb}{0.0, 0.18, 0.39}
\definecolor{darkspringgreen}{rgb}{0.09, 0.45, 0.27}
\definecolor{emerald}{rgb}{0.31, 0.78, 0.47}
\definecolor{lavenderindigo}{rgb}{0.58, 0.34, 0.92}
\definecolor{mred}{rgb}{0.83, 0.0, 0.0}

\begin{document}
\frontmatter 
\maketitle

\tableofcontents
\mainmatter%
\begin{abstract}
We establish an arithmetic intersection theory in the framework of Arakelov geometry over adelic curves. To each projective scheme over an adelic curve, we associate a multi-homogenous form on the group of adelic Cartier divisors, which can be written as an integral of local intersection numbers along the adelic curve. The integrability of the local intersection number is justified by using the theory of resultants. 
\end{abstract}


\chapter*{Introduction}

Since the seminal work of Dedekind and Weber \cite{MR1579901}, the similarity between number fields and fields of algebraic functions of one variable has been known  and has deeply influenced researches in algebraic geometry and number theory. Inspired by the discovery of Hensel and Hasse on embeddings of a number field into diverse local fields, Weil \cite{Weil39} considered all places, finite or infinite, of a number field, which made a decisive step toward the unification of number theory and algebraic geometry. Many works have then been done along this direction. On the one hand, the analogue of Diophantine problems (notably Mordell's conjecture) in the function field setting has been studied by Manin \cite{MR0154868}, Grauert \cite{MR222087} and Samuel \cite{MR204430}; on the other hand, through Weil's height machine \cite{MR42169} and the theory of N\'{e}ron-Tate's height \cite{MR179173}, methods of algebraic geometry have been systematically applied to the research of Diophantine problems, and it has been realized that the understanding of the arithmetic of algebraic varieties over a number field, which should be analogous to algebraic geometry over a smooth projective curve,  is indispensable in the geometrical approach of Diophantine problems. Under such a circonstance Arakelov \cite{MR0466150,MR0472815} has developed the arithmetic intersection theory for arithmetic surfaces (namely relative curves over $\Spec\mathbb Z$). Note that the transcription of the intersection theory to the arithmetic setting is by no means automatic. The key idea of Arakelov is to introduce transcendental objects, notably Hermitian metrics or Green functions, over the infinite place, in order to ``compactify'' arithmetic surfaces. To each pair of compactified arithmetic divisor, he attached a family of local intersection numbers parametrized by the set of places of the base number field. The global intersection number is obtained by taking the sum of local intersection numbers.   Arakelov's idea has soon led to spectacular advancements in Diophantine geometry, especially Faltings' proof \cite{MR740897} of Mordell's conjecture. 

The fundament of Arakelov geometry for higher dimensional arithmetic varieties has been established by Gillet and Soul\'{e}. The have introduced an arithmetic intersection theory \cite{MR770447,MR1087394} for general arithmetic varieties and proved an ``arithmetic Riemann-Roch theorem'' \cite{MR1055224}. They have introduced the arithmetic Chow group, which is an hybride construction of the classic Chow group in algebraic geometry and currents in complex analytic geometry. Applications of arithmetic intersection theory in Diophantine geometry have then been developed, notably to build up an intrinsic height theory for arithmetic projective varieties (see for example \cite{MR1109353,MR1260106}). Arakelov's height theory becomes now an important tool in arithmetic geometry. Upon the need of including several constructions of local heights (such as the canonical local height for subvarieties in an Abelian variety) in the setting of Arakelov geometry, Zhang \cite{MR1311351} has introduced the notion of adelic metrics for ample line bundles on a projective variety over a number field, which could be considered as uniform limit of Hermitian line bundles (with possibly different integral models). 

Inspired by the similarity between Diophantine analysis and Nevanlinna theory, Gubler \cite{MR1472498} has proposed a vast generalization of height theory in the framework of $M$-fields. Recall that a $M$-field is a field $K$ equipped with a measure space $M$ and a map from $K\times M$ to $\mathbb R_{\geqslant 0}$ which behaves almost everywhere like absolute values on $K$. Combining the intersection product of Green currents in the Archimedean case and the local height of Chow forms, he has introduced local heights (parametrized by the measure space $M$) for a projective variety over an $M$-field. Assuming the integrability of the function of local heights on the measure space $M$, he has defined the global height of the variety as the integral of local heights. Interesting examples have been discussed in the article, which show that in many cases the function of local heights is indeed integrable. 

In \cite{MR3783789}, we have developed an Arakelov geometry over adelic curves. Our framework is similar to $M$-field of Gubler, with a slightly different point of view: an adelic curve is a field equipped with a family of absolute values  parametrized by a measure space (in particular, we require the absolute values to be defined everywhere).  These absolute values play the role of places in algebraic number theory. Hence we can view an adelic curve as a measure space of ``places'' of a given field, except that we allow  possibly equivalent absolute values in the family, or even copies of the same absolute values. Natural examples of adelic curves contain global fields, countably generated fields over global fields (as we will show in the second chapter of the current article), field equipped with copies of the trivial absolute value, and also the amalgamation of different adelic structures of the same field.  Our motivation was to establish  a theory of adelic vector bundles (generalizing previous works of Stuhler \cite{MR424707}, Grayson \cite{MR780079}, Bost \cite{MR1423622} and Gaudron \cite{MR2431505}), which is  analogous to geometry of numbers and hence provides tools to consider Diophantine analysis in a general and flexible setting. By using the theory of adelic vector bundles, the arithmetic birational invariants are discussed in a systematic way.

The first contribution of the current article is to discuss transcendental coverings of adelic curves. Let $S=(K,(\Omega,\mathcal A,\nu),\phi)$ be an adelic curve, where $K$ is a countable field, $(\Omega,\mathcal A,\nu)$ is a measure space, and $\phi:\omega\mapsto|\ndot|_\omega$ is a map from $\Omega$ to the set of all absolute values of $K$, such that, for any $a\in K^{\times}$, the function $(\omega\in\Omega)\mapsto\ln|a|_\omega$ is measurable. In \cite[Chapter 3]{MR3783789}, for any algebraic extension $L/K$, we have constructed a measure space $(\Omega_{L},\mathcal A_L,\nu_L)$, which is fibered over $(\Omega,\mathcal A,\nu)$ and admits a family of disintegration probability measures. To each $\omega\in\Omega$, we correspond the fiber $\Omega_{L,\omega}$ to the family of all absolute values of $L$ extending $|\ndot|_\omega$. Thus we obtain a structure of adelic curve on $L$ which is called an \emph{algebraic covering} of $S$. 

In \cite[\S3.2.5]{MR3783789}, we have illustrated the construction of an adelic curve structure on $\mathbb Q(T)$, which takes into account the arithmetic of $\mathbb Q$ and the geometry of $\mathbb P^1$. In the current article, we generalizes and systemize such a construction on a purely transcendental and countably generated extension of the underlying field $K$ of the adelic curve $S$. For simplicity, we explain here the case of rational function of finitely many variables. Let $n$ be an integer such that $n\geqslant 1$ and $\boldsymbol{T}=(T_1,\ldots,T_n)$ be variables. Let $L$ be the rational function field $K(\boldsymbol{T})=K(T_1,\ldots,T_n)$, which is by definition the field of fractions of the polynomial ring $K[\boldsymbol{T}]=K[T_1,\ldots,T_n]$. 
To each $\omega\in\Omega$ such that the absolute value $|\ndot|_{\omega}$ is non-Archimedean, by Gauss's lemma, we extends $|\ndot|_{\omega}$ to be an absolute value on $L$ such that
\[\forall\,f=\sum_{\boldsymbol{d}\in\mathbb N^n}a_{\boldsymbol{d}}(f)\boldsymbol{T}^{\boldsymbol{d}}\in K[\boldsymbol{T}],\quad |f|_{\omega}=\max_{\boldsymbol{d}\in\mathbb N^n}|a_{\boldsymbol{d}}|_\omega.\]
We then take $\Omega_{L,\omega}$ to be the one point set $\{\omega\}$, which is equipped with the trivial probability measure.
In the case where the absolute value $|\ndot|_\omega$ is Archimedean, we fix an embedding $\iota_{\omega}:K\rightarrow\mathbb C$ such that $|\ndot|_{\omega}$ is the composition of the usual absolute value $|\ndot|$ on $\mathbb C$ with $\iota_\omega$ (by a measurable selection argument, we can arrange that the family of $\iota_{\omega}$ parametrized by Archimedean places is $\mathcal A$-measurable). We let 
\[ \Omega_{L,\omega}:=\left\{ (t_1,\ldots,t_n) \in  [0,1]^n  \left|
\begin{array}{l}
\text{\small $( e(t_1),\ldots,e(t_n))$ is algebraically} \\
\text{\small independent over $\iota_{\omega}(K)$} 
\end{array} 
\right\}\right., 
\]
where for each $t\in[0,1]$, $e(t)$ denotes $\mathrm{e}^{2\pi it}$. Note that, if we equip $[0,1]^n$ with the Borel $\sigma$-algebra and the uniform probability measure, then $\Omega_{L,\omega}$ is a Borel set of measure $1$. Moreover, each element $\boldsymbol{t}=(t_1,\ldots,t_n)\in\Omega_{L,\omega}$ gives rise to an absolute value $|\ndot|_{\boldsymbol{t}}$ on $L$ such that 
\[\forall\,f=\sum_{\boldsymbol{d}\in\mathbb N^n}a_{\boldsymbol{d}}(f)\boldsymbol{T}^{\boldsymbol{d}}\in K[\boldsymbol{T}],\quad |f|_{\boldsymbol{t}}=\bigg|\sum_{\boldsymbol{d}\in\mathbb N^n}\iota_\omega(a_{\boldsymbol{d}}(f))e(t_1)^{d_1}\cdots e(t_n)^{d_n}\bigg|.\]
It turns out that the disjoint union $\Omega_{L}$ of $(\Omega_{L,\omega})_{\omega\in\Omega}$ forms a structure of adelic curve on the field $L$, which is fibered over that of $S$, and admits a family of disintegration probability measures. We denote by $S_{L}=(L,(\Omega_{L},\mathcal A_{\Omega_{L}},\nu_{L}),\phi_{L})$ the corresponding adelic curve. 

In the case where the adelic curve $S$ is proper, namely the following equality holds for any $a\in K^{\times}$
\[\int_{\Omega}\ln|a|_\omega\,\nu(\mathrm{d}\omega)=0,\]
it is not true in general that the adelic curve $S_{L}$ is also proper. In the article, we propose several natural ``compactifications'' of the adelic curve. Here we explain one of them which has an ``arithmetic nature''. We say that two irreducible polynomials $P$ and $Q$ in $K[T_1,\ldots,T_n]$ are equivalent if they differ by a factor of non-zero element of $K$. This is an equivalence relation on the set of all irreducible polynomials In each equivalence class we pick a representative to form a family $\mathscr P$ of irreducible polynomials. Then every non-zero element $f$ of $K$ can be written in a unique way as 
\[f=a(f)\prod_{F\in\mathscr P}F^{\operatorname{ord}_F(f)},\]
where $a(f)\in K^{\times}$, and $\operatorname{ord}_F(\ndot):L\rightarrow\mathbb Z\cup\{+\infty\}$ is the discrete valuation associated with $F$, we denote by $|\ndot|_F=\mathrm{e}^{-\operatorname{ord}_F(\ndot)}$ the corresponding absolute value on $L$. Moreover, the degree function on $K[\boldsymbol{T}]$ extends naturally to $L$ so that $-\deg(\ndot)$ is a discrete valuation on $L$. Moreover, the following equality holds (see Proposition \ref{prop:admissible:fibration:K:T})
\[\forall\,f\in K(\boldsymbol{T}),\quad \sum_{F\in\mathscr P}\deg(F)\operatorname{ord}_F(f)=\deg(f).\] 
We let $|\ndot|_\infty$ be the absolute value on $L$ such that $|\ndot|_\infty=\mathrm{e}^{\deg(\ndot)}$. Note that, for any $F\in\mathscr P$, one has
\[h_{S_L}(F):=\int_{\Omega}\nu(\mathrm{d}\omega)\int_{\Omega_{L,\omega}}\ln|F|_x\,\nu_{L,\omega}(\mathrm{d}x)\geqslant 0.\]
We fix a positive real number $\lambda$. Let $(\Omega_L^\lambda,\mathcal A_{L}^{\lambda},\nu_L^\lambda)$ be the disjoint union of $(\Omega_L,\mathcal A_L,\nu_L)$ and $\mathscr P\cup\{\infty\}$, which is equipped with the measure $\nu_L^\lambda$ extending $\nu_L$ and such that $\nu_L^\lambda(\{\infty\})=\lambda$ and
\[\forall\, F\in\mathscr P,\quad \nu_L^\lambda(\{F\})=h_{S_L}(F)+\lambda\deg(F).\]
Let $\phi_L^\lambda$ be the map from $\Omega_L^\lambda$  to the set of absolute values on $L$, sending $x\in\Omega_L^\lambda$ to $|\ndot|_x$. 
Then we establish the following result (see \S\ref{sec:Arithmetic adelic structure}, notably Propositions \ref{prop:proper:L:lambda} and \ref{prop:Northcott:property}, see also Proposition \ref{Pro: compactification 2} for the general construction).

\begin{enonce*}{Theorem A} Assume that the adelic curve $S$ is proper.
\begin{enumerate}[label=\rm(\arabic*)]
\item For any $\lambda>0$, the adelic curve $S_L^\lambda=(L,(\Omega_L^\lambda,\mathcal A_L^\lambda,\nu_L^\lambda),\phi_L^\lambda)$ is proper.
\item If the adelic curve $S$ satisfies the Northcott property, namely, for any $C\geqslant 0$, the set
\[\bigg\{a\in K\,\bigg|\,\int_{\Omega}\max\{\ln|a|_\omega,0\}\,\nu(\mathrm{d}\omega)\leqslant C\bigg\}\]
is finite, then, for any $\lambda>0$, the adelic curve $S_L^\lambda$  satisfies the Northcott property.
\end{enumerate}
\end{enonce*}
Together with the algebraic covering of adelic curves mentioned above. This construction provides a large family of adelic structures for finitely generated extensions of $\mathbb Q$, which behaves well from the view of geometry of numbers. Note however that the compactification $S_L^\lambda$ is not fibered over $S$, but rather fibered over the amalgamation of $S$ with copies of the trivial absolute value on $K$. This phenomenon suggest that it is a need of dealing with the trivial absolute value in the consideration of the relative geometry of adelic curves.

To build up a more complete picture of Arakelov geometry over an adelic curve, it is important to develop an arithmetic intersection theory and relate it to the heights of  projective varieties over an adelic curve. Although the local intersection theory is now well understood, thanks to works such as \cite{MR1472498,MR1629925,MR2244803,MR3498148}, it remains a challenging problem to show that the local intersection numbers form an integrable function over the parametrizing measure space. In this article, we resolve this integrability problem and thus establish a global intersection theory in the framework of Arakelov geometry over adelic curves. Recall that the function of local heights for an adelic line bundle is only well defined up to the function of absolute values of a non-zero scalar. One way to make explicit the local height function is to fix a family of global sections of the line bundle which intersect properly. Note that each global sections determine a Cartier divisor on the projective variety, and the adelic metrics of the adelic line bundle determine a family of Green functions of the Cartier divisor parametrized by the measure space of ``places''. For this reason, we choose to work in the framework of adelic Cartier divisors.   

Let $S$ be an adelic curve, which consists of a field $K$, a measure space $(\Omega,\mathcal A,\nu)$ and a family $(|\ndot|_{\omega})_{\omega\in\Omega}$ of absolute values on $K$ parametrized by $\Omega$. Let $X$ be a projective scheme over $\Spec K$ and $d$ be the Krull dimension of $X$. By adelic Cartier divisor on $X$, we mean the datum $\overline D$ consisting of a Cartier divisor $D$ on $X$ together with a family $g=(g_\omega)_{\omega\in\Omega}$ parametrized by $\Omega$, where $g_\omega$ is a Green function of $D_\omega$, the pull back of $D$ on $X_\omega=X\otimes_K K_\omega$, with $K_\omega$ being the completion of $K$ with respect to $|\ndot|_\omega$. Conditions of measurability and dominancy (with respect to $\omega\in\Omega$) for the family $g$ are also required (see \S\S\ref{Sec: Reminder on adelic vector bundles}--\ref{Sec: Integrability of local intersection numbers} for more details). We first introduce the local intersection product for adelic Cartier divisors.  More precisely, if $\overline D_i=(D_i,g_i)$, $i\in\{0,\ldots,d\}$, form a family of integrable metrized Cartier divisors on $X$ (namely a Cartier divisor equipped with a Green function) such that $D_0,\ldots,D_d$ intersect properly, we define, for any $\omega\in\Omega$ a local intersection number 
\[(\overline D_0,\ldots,\overline D_d)_\omega\in\mathbb R\]
in a recursive way using Bedford-Taylor theory \cite{MR445006} and its non-Archimedean analogue \cite{MR2244803}. 
In the case where $|\ndot|_{\omega}$ is a trivial absolute value, we need a careful definition of
the local intersection number
(see Definition~\ref{Def: intersection in trivial valuation case}, for details).
Note the local intersection number is a multi-linear function on the set of $(d+1)$-uplets $(\overline D_0,\ldots,\overline D_d)$ such that $D_0,\ldots,D_d$ intersect properly. 

To establish a global intersection theory, we need to show that the function of local intersection numbers
\[(\omega\in\Omega)\longmapsto (\overline D_0,\ldots,\overline D_d)_\omega \]
is measurable and integrable with respect to $\nu$, where the measurability part is more subtle. Although the Green function families of $\overline D_0,\ldots,\overline D_d$ are supposed to be measurable, the corresponding products of Chern currents (or their non-Archimedean analogue) depend on the local analytic geometry relatively to the absolute values $|\ndot|_\omega$. It seems to be a difficult (but interesting) problem to precisely describe the measurability of the local geometry of the analytic spaces $X_\omega^{\mathrm{an}}$.  For places $\omega$ which are Archimedean, as we can embed all  local completions $K_{\omega}$ in the same field $\mathbb C$, by a measurable selection theorem one can show that the family of Monge-Amp\`{e}re measures is measurable with respect to $\omega$ (see Theorem \ref{thm:measurable:archimedean}). However, for non-Archimedean places, such embeddings in a common valued field do not exist in general, and the classic approach of taking a common integral model for all non-Archimedean places is not adequate in the setting of adelic curves, either.  

To overcome this difficulty, our approach consist in relating the local intersection number to the local length of the mixed resultant and hence reduce the problem to the measurability of the function of local length of the mixed resultant, which is known by the theory of adelic vector bundles developed in \cite{MR3783789}. This approach is inspired by previous results of Philippon \cite{MR876159} on height of algebraic cycles via the theory of Chow forms and the comparison \cite{MR1092175,MR1306549,MR1144338,MR1260106}  between Philippon's height and Faltings height (defined by the arithmetic intersection theory). Note that the similar idea has also been used in \cite{MR1472498} to construct the local height in the setting of $M$-fields.

Let us briefly recall the theory of mixed resultant. It is a multi-homogeneous generalization of Chow forms, which allows to describe the interactions of several embeddings of a variety in projective spaces by a multi-homogeneous polynomial. One of its original forms is the discriminant of a quadratic polynomial, or more generally the resultant of $n+1$ polynomials $P_0,\ldots,P_n$ in $n$ variables over an algebraically closed field, which is an irreducible polynomial in the coefficients of $P_0,\ldots,P_n$, which vanishes precisely when these polynomials have a common root. The modern algebraic approach of resultants goes back to the elimination theory of Cayley \cite{Cayley48}, where he related resultant to the determinant of Koszul complex. We use here a geometric reformulation as in the book \cite{MR1264417} of Gel'fand, Kapranov and Zelevinsky. In Diophantine geometry, mixed resultant has been used by R\'emond \cite{MR1837827} to study multi-projective heights.

We assume that the Cartier divisors $D_i$ are very ample and thus determine closed immersions $f_i$ from $X$ to the projective space of the linear system $E_i$ of the divisor $D_i$. By \emph{incidence variety} of $(f_0,\ldots,f_d)$, we mean the closed subscheme $I_X$ of $X\times_K\mathbb P(E_0^\vee)\times_K\cdots\times_K\mathbb P(E_{d}^\vee)$ parametrizing points $(x,\alpha_0,\ldots,\alpha_d)$ such that 
\[\alpha_0(x)=\cdots=\alpha_d(x)=0.\]
One can also consider $I_X$ as a multi-projective bundle over $X$ (of $E_i^\vee$ quotient by the tautological line subbundle).
Therefore, the projection of $I_X$ in $\mathbb P(E_0^\vee)\times_K\cdots\times_K\mathbb P(E_d^\vee)$ consists of a family of hyperplanes in $\mathbb P(E_0),\ldots,\mathbb P(E_d)$ respectively, which contain at least one common point of $X$. It turns out that this projection is actually a multi-homogeneous hypersurface of $\mathbb P(E_0^\vee)\times_K\cdots\times_K\mathbb P(E_d^\vee)$, which is defined by a multi-homogeneous polynomial $R_{f_0,\ldots,f_d}^X$, called a \emph{resultant} of $X$ with respect to the embeddings of $f_0,\ldots,f_d$. We refer the readers to \cite[\S3.3]{MR1264417} for more details, see also \cite{MR3098424} for applications in arithmetic Nullstellensatz. When $K$ is a number field, the height of the polynomial $R_{f_0,\ldots,f_d}^X$ can be viewed as a height of the arithmetic variety $X$, and, in the particular case where the image of $D_i$ in the Picard group are colinear, an explicit comparison between the height of resultant and the Faltings height of $X$ has been discussed in \cite[Theorem 4.3.2]{MR1260106} (see also \S4.3.4 of \emph{loc. cit.}).

Usually the resultant is well defined up to a factor in $K^{\times}$. In the classic setting of number field, this is anodyne for the study of the global height, thanks to the product formula. However, in our setting, this dependence on the choice of a non-zero scalar could be annoying, especially when the adelic curve does not satisfy a product formula.  In order to obtain a local height equality, we introduce, for each vector 
\[(s_0,\ldots,s_d)\in E_0\times\cdots \times E_d\] 
such that $\operatorname{div}(s_0),\ldots,\operatorname{div}(s_d)$ intersect properly on $X$, a specific resultant $R^{X,s_0,\ldots,s_d}_{f_0,\ldots,f_d}$ of $X$ with respect to the embeddings, which is the only resultant such that 
\[R^{X,s_0,\ldots,s_d}_{f_0,\ldots,f_d}(s_0,\ldots,s_d)=1.\] We then show that the local height for this resultant coincides with the local height of $X$ defined by the local intersection theory. By using this comparison of local height and properties of adelic vector bundles over an adelic curve (see \cite[\S4.1.4]{CMArakelovAdelic}), we prove the integrability of the local height function on non-Archimedean places. Moreover, the integral of the local height equalities leads to a equality between the global height of the resultant and the arithmetic intersection number (see Remark \ref{Rem: global heights equality}), which generalizes the height comparison results in \cite{MR1092175,MR1260106}. In resume, we obtain the following result (see  Theorems \ref{Thm: equality of local height} and \ref{Thm: integrability of local intersection}).

\begin{enonce*}{Theorem B}
Let $S=(K,(\Omega,\mathcal A,\nu),\phi)$ be an adelic curve, $X$ be a projective scheme over $S$, $d$ be the dimension of $X$, $D_0,\ldots, D_d$ be Cartier divisors on $X$, which are equipped with Green function families $g_0,\ldots,g_d$, respectively, such that $(D_{i,\omega},g_{i,\omega})$ is integrable for any $\omega\in\Omega$ and $i\in\{0,\ldots,d\}$. 
\begin{enumerate}[label=\rm(\arabic*)]
\item\label{Item: local height formula} Assume that the Cartier divisors $D_0,\ldots,D_d$ are very ample. For any $i\in\{0,\ldots,d\}$, let $E_i=H^0(X,\mathcal O_X(D_i))$, $f_i:X\rightarrow\mathbb P(E_i)$ be the closed embedding and $s_i\in E_i$ be the regular meromorphic section of $\mathcal O_X(D_i)$ corresponding to $D_i$. Assume that the continuous metric family $\varphi_{g_i}$ corresponding to the Green function family $g_i$ is consisting of the orthogonal quotient metrics induces by a Hermitian norm family $\xi_i=(\|\ndot\|_{i,\omega})_{\omega\in\Omega}$ on $E_i$. Then, for any $\omega\in\Omega$, then following equalities hold.
\begin{enumerate}[label=\rm(1.\alph*)]
\item In the case where $|\ndot|_\omega$ is non-Archimedean, one has
\[(\overline D_0\cdots\overline D_d)_\omega=\ln\Big\|R_{f_0,\ldots,f_d}^{X,s_0,\ldots,s_d}\Big\|_{\omega,\varepsilon},\]
where the norm $\|\ndot\|_{\omega,\varepsilon}$ on the space of multi-homogeneous polynomials is the $\varepsilon$-tensor product of $\varepsilon$-symmetric power of $\|\ndot\|_{i,\omega,*}$.
\item In the case where $|\ndot|_\omega$ is Archimedean, one has
\[\begin{split}(\overline D_0\cdots\overline D_d)_\omega=\int_{\mathbb S(E_{0,\omega})\times\cdots\times\mathbb S(E_{d,\omega})}&\ln\Big|R_{f_0,\ldots,f_d}^{X,s_0,\ldots s_d}(z_0,\ldots,z_d)\Big|_\omega\,\mathrm{d}z_0\cdots\mathrm{d}z_d\\
&+\frac 12\sum_{i=0}^d\delta_i\sum_{\ell=1}^{r_i}\frac{1}{\ell},
\end{split}\]
where $\mathbb{S}(E_{i,\omega})$ denotes the unit sphere of $(E_{i,\omega},\|\ndot\|_{i,\omega})$, $\mathrm{d}z_i$ is the Borel probability measure on $\mathbb{S}(E_{i,\omega})$ invariant by the unitary group, $r_i$ is the dimension of $E_i$, and $\delta_i$ is the intersection number
\[(D_0\cdots D_{i-1}D_{i+1}\cdots D_d).\]
\end{enumerate}
\item Assume that, either the $\sigma$-algebra $\mathcal A$ is discrete, or the field $K$ admits a countable subfield, which is dense in each $K_\omega$. If all couples $\overline D_i=(D_i,g_i)$ are integrable adelic Cartier divisors on $X$, the the function 
\[(\omega\in\Omega)\longrightarrow (\overline D_0\cdots\overline D_d)_\omega\]
is $\nu$-integrable.
\end{enumerate}
\end{enonce*} 
As an application, we can define the multi-height of the projective scheme $X$ with respect to $\overline D_0,\ldots,\overline D_d$ as
\[h_{\overline D_0\cdots\overline D_d}(X)=\int_{\Omega}(\overline D_0\cdots\overline D_d)_\omega\,\nu(\mathrm{d}\omega),\]
and, under the assumptions of the point \ref{Item: local height formula} in the above theorem, we can relate the multi-height with the height of the resultant, by taking the integral of the local height equality.

From the methodological point of view, the approach of \cite{MR1092175} works  within $\mathbb P^N(\mathbb C)$ and uses elimination theory and complex analysis of the Fubini-Study metric; that of \cite{MR1260106} relies on a choice of integral model and computations in the arithmetic Chow groups. In our setting, we need to deal with general non-Archimedean metrics. Hence these approaches do not fit well with the framework of adelic curves.  Our method consists in computing the local height of 
\[X\times_K\mathbb P(E_0^\vee)\times_K\cdots\times_K\mathbb P(E_{d}^\vee)\]
in two ways (see Lemma \ref{Lem: h s in h IX}
 for details). We first consider this scheme as a fibration of multi-projective space over $X$ and relate this local height to that of $X$ by taking the local intersection along the fibers. We then relate the height of this product scheme to that of the incidence subscheme $I_X$ and then use the identification of $I_X$ with a multi-projective bundle over $X$ to compute recursively the height of $I_X$. Our method allows to obtain a local height equality in considering the Archimedean case and the non-Archimedean case in a uniform way.

It is worth mentioning that an intersection theory of arithmetic cycles and a Riemann-Roch theory could be expected for the setting of adelic curves. However, new ideas are needed to establish a good formulation of the measurability for various arithmetic objects arising in such a theory.   

The rest of the article is organized as follows. In the first chapter, we remind several basic constructions used in the article, including multi-linear subset and multi-linear functions, Cartier divisors on general scheme,  proper intersection of Cartier divisors on a projective scheme, multi-homogeneous polynomials, incidence subscheme and resultants, and linear projections of closed subschemes in a projective space. The second chapter is devoted to the construction of adelic structures. After a brief reminder on the definition of adelic curves and their algebraic covers, we introduce transcendental fibrations of adelic curves and their compactifications. These constructions provide a large family of examples of adelic curves. In the third chapter, we consider the local intersection theory in the setting of projective schemes over a complete valued field. We first remind the notions of continuous metrics on an invertible sheaf and its semi-positivity. Then we explain the notion of Green functions of Cartier divisors and their relation  with continuous metrics. The construction of Monge-Amp\`{e}re mesures and local intersection numbers is then discussed. The last sections are devoted to establish the link between the local intersection number and the length (in the non-Archimedean case) or Mahler measure (in the Archimedean case) of the corresponding resultant, respectively. In the fourth and last chapter, we prove the integrability of the local height function and construct the global multi-height.


\chapter{Reminder and preliminaries}

\section{Symmetric and multi-linear subsets}

In this section, we fix a commutative and unitary ring $k$.
\begin{defi}\label{def:domain}
Let $d$ be a non-negative integer and $V$ be a $k$-module.
We say that a subset $S$ of $V^{d+1}$ is \emph{multi-linear} if, for any $j\in\{0,\ldots,d\}$ and for any $(x_0,\ldots,x_{j-1},x_{j+1},\ldots,x_d)\in V^d$, the subset
\[\{x_j\in V\,|\,(x_0,\ldots,x_{j-1},x_j,x_{j+1},\ldots,x_d)\in S\}\]
of $V$ is either empty or a sub-$k$-module. If in addition
\[(x_0,\ldots,x_d)\in S\;\Longrightarrow\;(x_{\sigma(0)},\ldots,x_{\sigma(d)})\in S\]
for any bijection $\sigma:\{0,\ldots,d\}\rightarrow\{0,\ldots,d\}$, we say that the multi-linear subset $S$ is \emph{symmetric}.
\end{defi}

\begin{prop}
Let $d\in{\mathbb Z}_{\geqslant 0}$, $V$ be a $k$-module and $S$ be a multi-linear subset of $V^{d+1}$. 
For any $j\in\{0,\ldots,d\}$, let $I_j$ be a non-empty finite set, $(x_{j,i})_{i\in I_j}$ be a family of elements of $V$, $(\lambda_{j,i})_{i\in I_j}$ be a family of elements of $k$, and $y_j=\sum_{i\in I_j}\lambda_{j,i}x_{j,i}$. Assume that, for any $(i_0,\ldots,i_d)\in I_0\times\cdots\times I_d$, one has $(x_{0,i_0},\ldots,x_{d,i_d})\in S$. Then $(y_0,\ldots,y_d)\in S$. 
\end{prop}

\begin{proof}
We reason by induction on $d$. In the case where $d=0$, $S$ is  a sub-$k$-module of $V$ when it is not empty. Since $y_0$ is a $k$-linear combination  of elements of $S$, we obtain that $y_0\in S$. 

We now assume that $d\geqslant 1$ and that the statement holds for multi-linear subsets of $V^d$. Let \[S'=\{(z_0,\ldots,z_{d-1})\in V^d\,|\,(z_0,\ldots,z_{d-1},y_d)\in S\}.\]
Since $S$ is a multi-linear subset of $V^{d+1}$, for any $(i_0,\ldots,i_{d-1})\in I_0\times\cdots\times I_{d-1}$, one has $(x_{i_0},\ldots,x_{i_{d-1}},y_d)\in S$ and hence $(x_{i_0},\ldots,x_{i_{d-1}})\in S'$. Moreover, $S'$ is a multi-linear subset of $V^d$. Hence the induction hypothesis leads to $(y_0,\ldots,y_{d-1})\in S'$ and thus $(y_0,\ldots,y_d)\in S$.
\end{proof}

\begin{defi}Let $d\in {\mathbb Z}_{\geqslant 0}$, $V$ and $W$ be two $k$-modules, and $S$ be a multi-linear subset of $V^{d+1}$. We say that a map $f:S\rightarrow W$ is \emph{multi-linear} if, for any $j\in\{0,\ldots,d\}$ and for any $(x_0,\ldots,x_{j-1},x_{j+1},\ldots,x_d)\in V^d$, the map
\[\{x_j\in V\,|\,(x_0,\ldots,x_{j-1},x_j,x_{j+1},\ldots,x_d)\in S\}\longrightarrow W,\quad x_j\mapsto f(x_0,\ldots,x_d),\]
is $k$-linear once \[\{x_j\in V\,|\,(x_0,\ldots,x_{j-1},x_j,x_{j+1},\ldots,x_d)\in S\}\] is not empty. If in addition $S$ is symmetric and $f(x_0,\ldots,x_d)=f(x_{\sigma(0)},\ldots,x_{\sigma(d)})$ for any $(x_0,\ldots,x_d)\in S$ and any bijection $\sigma:\{0,\ldots,d\}\rightarrow\{0,\ldots,d\}$, we say that $f$ is a \emph{symmetric multi-linear map}.
\end{defi}

\begin{prop}\label{pro:domain}
Let $d\in{\mathbb Z}_{\geqslant 0}$, $V$ and $W$ be two $k$-modules, $S$ be a multi-linear subset of $V^{d+1}$, and $f:S\rightarrow W$ be a multi-linear map.
Let $( x_{j,i} )_{(j,i)\in\{0,\ldots,d\}^2}$ be a matrix consisting of elements of $V$
such that $(x_{0,i_0}, \ldots, x_{d,i_d}) \in S$ for any $(i_0, \ldots, i_d) \in \{ 0 , \ldots, d\}^{d+1}$.
Then
\begin{equation}\label{Equ:somme alternee}\begin{split}
\sum_{\sigma \in \mathfrak{S}(\{0, \ldots, d\})} f&(x_{0,\sigma(0)}, \ldots, x_{d,\sigma(d)}) \\
&= \sum_{\emptyset \not= I \subseteq \{ 0, \ldots, d \}} (-1)^{d+1 - \#I} f\Big(\sum\nolimits_{i_0 \in I} x_{0,i_0} , \ldots,\sum\nolimits_{i_d \in I} x_{d,i_d} \Big),\end{split}
\end{equation}
where $\mathfrak{S}(\{0, \ldots, d\})$ is the permutation group of $\{ 0, \ldots, d\}$.
\end{prop}

\begin{proof} By the multi-linearity of $f$, we can rewrite the right-hand side of the equality \eqref{Equ:somme alternee} as
\[\begin{split}&\quad\;\sum_{\varnothing\neq I\subseteq\{0,\ldots,d\}}(-1)^{d+1-\#I}\sum_{(i_0,\ldots,i_d)\in I^{d+1}}f(x_{0,i_0},\ldots,x_{d,i_d})\\
&=\sum_{(i_0,\ldots,i_d)\in\{0,\ldots,d\}^{d+1}}\bigg(\sum_{\{i_0,\ldots,i_d\}\subseteq I\subseteq\{0,\ldots,d\}}(-1)^{d+1-\#I}\bigg)f(x_{0,i_0},\ldots,x_{d,i_d}).
\end{split}\]
Note that, for $(i_0,\ldots,i_d)\in\{0,\ldots,d\}^{d+1}$ such that $\{i_0,\ldots,i_d\}\subsetneq \{0,\ldots,d\}$, one has
\[\sum_{\{i_0,\ldots,i_d\}\subseteq I\subseteq\{0,\ldots,d\}}(-1)^{d+1-\#I}=(-1)^{d+1-\#\{i_0,\ldots,i_d\}}\sum_{J\subseteq\{0,\ldots,d\}\setminus\{i_0,\ldots,i_d\}}(-1)^{-\#J}=0\]
since
\[\sum_{J\subseteq\{0,\ldots,d\}\setminus\{i_0,\ldots,i_d\}}(-1)^{-\#J}=(1+(-1))^{d+1-\#\{i_0,\ldots,i_d\}}=0.\]
Therefore the equality \eqref{Equ:somme alternee} holds.

\end{proof}

\section{Cartier divisors}

In this section, let us recall the notion of Cartier divisor on a general scheme. The main references are \cite[$\text{IV}_4$, \S\S20-21]{EGA} and \cite{MR570309}.

\begin{defi}Let $X$ be a locally ringed space. We denote by $\mathcal O_X$ the structural sheaf of $X$. Let $\mathscr M_X$ be the sheaf of \emph{meromorphic functions on $X$}. Recall that $\mathscr M_X$ is the sheaf of commutative and unitary rings associated with the presheaf
\[U\longmapsto\mathcal O_X(U)[S_X(U)^{-1}],\]
where $S_X(U)$ denotes the multiplicative sub-monoid of $\mathcal O_X(U)$ consisting of local non-zero-divisors of $\mathcal O_X(U)$, that is,
$s\in\mathcal O_X(U)$ such that the homothety
\[\mathcal O_{X,x}\longrightarrow\mathcal O_{X,x},\quad a\longmapsto as_x\]
is injective for any $x\in U$
(here $s_x$ denotes the canonical image of $s$ in the local ring $\mathcal O_{X,x}$). We refer the readers to \cite{MR570309} for a clarification on the construction of the sheaf of meromorphic functions comparing to \cite[$\text{IV}_4$.(20.1.3)]{EGA}. 
\end{defi}

\begin{rema}Note that, for any $x\in X$, $\mathscr M_{X,x}$ identifies with $\mathcal O_{X,x}(S_{X,x}^{-1})$, where $S_{X,x}$ denotes the direct limit of $S_X(U)$ with $U$ running over the set of open neighbourhoods of $x$, viewed as a multiplicative submonoid of $\mathcal O_{X,x}$, which is contained in the sub-monoid of non-zero-divisors. Therefore, $\mathscr M_{X,x}$ could be considered as a sub-ring of the total fraction ring of $\mathcal O_{X,x}$, namely the localization of $\mathcal O_{X,x}$ with respect to the set of non-zero-divisors. In general the local ring $\mathscr M_{X,x}$ is different from the ring of total fractions of $\mathcal O_{X,x}$ even if $X$ is an affine scheme. The equality holds notably when $X$ is a locally Noetherian scheme or a reduced scheme whose set of  irreducible component is locally finite. We refer the readers to \cite{MR570309} for counter-examples and more details.
\end{rema}

\begin{defi}Let $X$ be a locally ringed space. We denote by $\mathscr M_X^{\times}$ the subsheaf of multiplicative monoids of $\mathscr M_X$ consisting of invertible elements. In other words, for any open subset $U$ of $X$, $\mathscr M_X^{\times}(U)$ is consisting of sections $s\in\mathscr M_X^{\times}(U)$ such that, for any $x\in U$, the homothety 
\[\mathscr M_{X,x}\longrightarrow\mathscr M_{X,x},\quad a\longmapsto as_x\]
is an isomorphism of $\mathscr M_{X,x}$-modules. An element of $\mathscr M_X^{\times}(U)$ is called a \emph{regular meromorphic function} on $X$. Similarly, let $\mathcal O_X^{\times}$ be the subsheaf of multiplicative monoids of $\mathcal O_X$ consisting of invertible elements~: for any open subset $U$ of $X$, $\mathcal O_X^{\times}(U)$ consists of sections $s\in\mathcal O_X(U)$ such that, for any $x\in U$, the homothety
\[\mathcal O_{X,x}\longrightarrow\mathcal O_{X,x},\quad a\longmapsto as_x\]
is an isomorphism of $\mathcal O_{X,x}$-modules.   Note that, for each $s\in\mathcal O_X(U)$, the homothety $s_x:\mathcal O_{X,x}\rightarrow\mathcal O_{X,x}$ induces by passing to localisation an homothety $\mathscr M_{X,x}\rightarrow\mathscr M_{X,x}$, which is an isomorphism of $\mathscr M_{X,x}$-modules if $s_x:\mathcal O_{X,x}\rightarrow\mathcal O_{X,x}$ is an isomorphism. Therefore, the canonical morphism $\mathcal O_X\rightarrow\mathscr M_X$ induces a morphism of sheaves of abelian groups $\mathcal O_X^{\times}\rightarrow\mathscr M_X^{\times}$.
\end{defi}

\begin{defi}
We call \emph{Cartier divisor} on $X$ any global section of the sheaf $\mathscr M_X^{\times}/\mathcal O_X^{\times}$. 
By definition, a Cartier divisor $D$ is represented by the following data: (i) an open covering $X = \bigcup_i U_i$ of $X$ and (ii) $f_i \in \mathscr M_X^{\times}(U_i)$ for each $i$ such that $f_i/f_j \in \mathcal O^{\times}_X$ on $U_i \cap U_j$ for all $i, j$.
The regular meromorphic function $f_i$ is called a \emph{local equation} of $D$ over $U_i$.
The group of Cartier divisors is denoted by $\operatorname{Div}(X)$ and the group law of $\operatorname{Div}(X)$ is written additively. Note that the exact sequence 
\[\xymatrix{1\ar[r]&\mathcal O_X^{\times}\ar[r]&\mathscr M_X^{\times}\ar[r]&\mathscr M_X^{\times}/\mathcal O_X^{\times}\ar[r]&0}\]
induces an exact sequence of cohomological groups
\begin{equation}\label{Equ: cohomological exact sequence}\xymatrix@C=1.17pc{1\ar[r]&\Gamma(X,\mathcal O_{X}^{\times})\ar[r]&\Gamma(X,\mathscr M_X^{\times})\ar[r]&\operatorname{Div}(X)\ar[r]&H^1(X,\mathcal O_X^{\times})\ar[r]&H^1(X,\mathscr M_X^{\times})}.\end{equation}
We denote by $\operatorname{div}(\ndot)$ the morphism $\Gamma(X,\mathscr M_X^{\times})\rightarrow\operatorname{Div}(X)$ in this exact sequence. Since the group law of $\operatorname{Div}(X)$ is written additively,  one has
\[\operatorname{div}(fg)=\operatorname{div}(f)+\operatorname{div}(g)\] 
for any couple of regular meromorphic functions $f$ and $g$ on $X$. A Cartier divisor belonging to the image of $\operatorname{div}(\ndot)$ is said to be \emph{principal}. If $D_1$ and $D_2$ are two Cartier divisors such that $D_1-D_2$ is principal, we say that $D_1$ and $D_2$ are \emph{linearly equivalent}, denoted by $D_1\sim D_2$.
\end{defi}

\begin{rema}Recall that $H^1(X,\mathcal O_X^{\times})$ identifies with the Picard group $\operatorname{Pic}(X)$ of $X$, namely the group of isomorphism classes of invertible $\mathcal O_X$-modules (see \cite[{\bf0}.(5.6.3)]{MR3075000}). Similarly, $H^1(X,\mathscr M_X^{\times})$ identifies with the group of isomorphism classes of invertible $\mathscr M_X$-modules. If $L$ is an invertible $\mathcal O_X$-module, then $\mathscr M_X\otimes_{\mathcal O_X}L$ is an invertible $\mathscr M_X$-module. The homomorphism $H^1(X,\mathcal O_X^{\times})\rightarrow H^1(X,\mathscr M_X^{\times})$ sends the isomorphism class of an invertible $\mathcal O_X$-module $L$ to that of the invertible $\mathscr M_X$-module $\mathscr M_X\otimes_{\mathcal O_X}L$. 
\end{rema}

\begin{defi}
Let $L$ be an invertible $\mathcal O_X$-module and $U$ be a non-empty open subset of $X$. We call \emph{regular meromorphic section} of $L$ on $U$ any element of $\Gamma(U,\mathscr M_X\otimes_{\mathcal O_X}L)$ which defines an isomorphism from $\mathscr M_U$ to $\mathscr M_U\otimes_{\mathcal O_U}L|_U$. Therefore, $\mathscr M_X\otimes_{\mathcal O_X}L$ is isomorphic as $\mathscr M_X$-module to $\mathscr M_X$ if and only if $L$ admits a regular meromorphic section on $X$.
\end{defi} 

\begin{rema}\label{Rem: regular meromorphic section}
Let $X$ be a locally Noetherian scheme or a reduced scheme whose set of  irreducible component is locally finite. For any $x\in X$, the local ring $\mathscr M_{X,x}$ identifies with the ring of total fractions of $\mathcal O_{X,x}$. Therefore, if $L$ is an invertible $\mathcal O_X$-module and if $U$ is an open subset of $X$, an element $s\in \Gamma(U,\mathscr M_X\otimes_{\mathcal O_X}L)$ is a regular meromorphic section of $L$ on $U$ if and only if it defines an injective homomorphism from $\mathcal O_U$ to $\mathscr M_U\otimes_{\mathcal O_U}L$. In particular, an element $s\in\Gamma(U,L)$ defines a regular meromorphic section of $L$ on $U$ if and only if, for any $x\in U$, $s_x\in\mathcal O_{X,x}\otimes_{\mathcal O_X}L$ is of the form $f_xs_{0,x}$, where $f_x$ is a non-zero-divisor of $\mathcal O_{X,x}$, and $s_{0,x}$ is a local trivialization of $L$ at $x$. This condition is also equivalent to $s(y)\neq 0$ for any associate point $y\in U$. Recall that a point $y\in X$ is called an \emph{associated point} if there exists $a\in\mathcal O_{X,y}$ such that the maximal ideal of $\mathcal O_{X,y}$ identifies with \[\operatorname{ann}(a):=\{f\in\mathcal O_{X,y}\,|\,af=0\}.\]
Let $x$ be a point of $X$. Assume that $s_x=f_xs_{0,x}$ where $f_x$ is a zero-divisor in $\mathcal O_{X,x}$, then $f_x$ belongs to an associated prime ideal of $\mathcal O_{X,x}$, which corresponds to an associated point $y\in X$ such that $x\in\overline{\{y\}}$ and  $s(y)=0$.

By \cite[$\text{IV}_4.(21.3.5)$]{EGA}, if $X$ is a Noetherian scheme, which admits an ample invertible $\mathcal O_X$-module, then the set of all associated points of $X$ is contained in an affine open subset of $X$, and any invertible $\mathcal O_X$-module admits a regular meromorphic section.   
\end{rema}

\begin{defi}\label{Def: divisor of a section}
Let $D$ be a Cartier divisor on $X$. The homomorphism $\operatorname{Div}(X)\rightarrow H^1(X,\mathcal O_X^{\times})$ in the exact sequence \eqref{Equ: cohomological exact sequence} sends $D$ to an isomorphism class of invertible $\mathcal O_X$-modules. One can actually construct explicitly an invertible $\mathcal O_X$-module $\mathcal O_X(D)$ in this class as follows. Let $(U_i)_{i\in I}$ be an open cover of the topological space such that $D$ is represented on each $U_i$ by a regular meromorphic function $f_i\in\Gamma(U_i,\mathscr M_{U_i}^{\times})$. For any couple $(i,j)\in I^2$, $f_i|_{U_i\cap U_j}f_j|_{U_i\cap U_j}^{-1}$ defines an isomorphism\[(f_i^{-1}\mathcal O_{U_i})|_{U_i\cap U_j}\longrightarrow(f_j^{-1}\mathcal O_{U_j})|_{U_i\cap U_j}.\]Moreover, these isomorphisms clearly satisfy the cocycle condition.  Thus the gluing of the sheaves $f_i^{-1}\mathcal O_{U_i}$ leads to an invertible sub-$\mathcal O_X$-module of $\mathscr M_X$ which we denote by $\mathcal O_X(D)$. Note that the gluing of meromorphic sections \[f_i\otimes f_i^{-1}\in\Gamma(U_i,\mathscr M_{U_i}\otimes\mathcal O_X(D))\] leads to a global regular meromorphic section of $\mathcal O_X(D)$, which we denote by $s_D$ and call \emph{canonical regular meromorphic section} of $\mathcal O_X(D)$. Hence $\mathscr M_X\otimes_{\mathcal O_X}\mathcal O_X(D)$ is canonically isomorphic to $\mathscr M_X$. Note that two Cartier divisors $D_1$ and $D_2$ are linearly equivalent if and only if the invertible $\mathcal O_X$-modules $\mathcal O_X(D_1)$ and $\mathcal O_X(D_2)$ are isomorphic.

 Conversely, the exactness of the diagram \eqref{Equ: cohomological exact sequence} shows that, an invertible $\mathcal O_X$-module $L$ is isomorphic to  an invertible $\mathcal O_X$-module of the form $\mathcal O_X(D)$ if and only if it admits a regular meromorphic section on $X$. One can also construct explicitly a Cartier divisor from a regular meromorphic section $s$ of $L$. In fact, let $(U_i)_{i\in I}$ be an open cover of $X$ such that each $L|_{U_i}$ is trivialized by a section $s_i\in L(U_i)$. For any $i\in I$, let $f_i$ be the unique regular meromorphic function on $U_i$ such that $s=f_is_i$.  Then the family $(f_i)_{i\in I}$ of regular meromorphic functions defines a Cartier divisor on $X$ which we denote by $\operatorname{div}(L;s)$, or by $\operatorname{div}(s)$ for simplicity.
\end{defi}

\begin{rema}
In the case where $X$ is a quasi-projective scheme over a field, any invertible $\mathcal O_X$-module admits a global regular meromorphic section and therefore is isomorphic to an invertible $\mathcal O_X$-module of the form $\mathcal O_X(D)$, where $D$ is a Cartier divisor. Hence one has an exact sequence
\[\xymatrix{1\ar[r]&\Gamma(X,\mathcal O_{X}^{\times})\ar[r]&\Gamma(X,\mathscr M_X^{\times})\ar[r]&\operatorname{Div}(X)\ar[r]&H^1(X,\mathcal O_X^{\times})\ar[r]&1}.\]   
\end{rema}

\begin{rema}\label{remark:Cartier:div:on:0:dim}
Let $X$ be a $0$-dimensional projective scheme over a field $k$.
Then there is a $k$-algebra $A$ which is finite-dimensional as a vector space over $k$, and such that
$X = \Spec(A)$.
Note that the canonical homomorphism
$A \to \bigoplus_{x \in X} A_x$ is an isomorphism. 
Let $f_x$ be a regular element of $A_x$. As the homotethy map $A_x \to A_x$, $a\mapsto f_x a$, is injective and $A_x$ is a finite-dimensional vector space over $k$, this homothety map is actually an isomorphism, that is,
$f_x \in A^{\times}_x$. Thus $\mathscr M_X^{\times} = \mathcal O_{X}^{\times}$.
Therefore, every Cartier divisor on $X$ can be represented by $1 \in A$.
\end{rema}

\begin{rema}\label{remark:expansion:Cartier:div:as:cycle}
Let $X$ be a Noetherian scheme. We denote by $X^{(1)}$ the set of all height 1 points of $X$, that is, 
$x \in X$ with $\dim (\mathcal O_{X, x}) = 1$. For $x \in X^{(1)}$ and a regular element $f$ of $\mathcal O_{X, x}$,
we set \[\ord_x(f) := \operatorname{length}_{\mathcal O_{X, x}}(\mathcal O_{X, x}/f\mathcal O_{X, x}).\]
Then $\ord_x(fg) = \ord_x(f) + \ord_x(g)$ for
all regular elements $f, g$ of $\mathcal O_{X, x}$ (cf \cite[the last paragraph of Section~1.3]{MArakelov}), so that
$\ord_x$ extends to a homomorphism $\mathscr M^{\times}_{X,x} \to \mathbb Z$.
Let $D$ be a Cartier divisor on $X$ and $f$ a local equation of $D$ at $x$.
Then it is easy to see that $\ord_{x}(f)$ does not depend on the choice of $f$, so that $\ord_x(f)$ is denoted by $\ord_x(D)$.
We call the cycle
\[\sum_{x \in X^{(1)}} \ord_x(D) \overline{\{ x \}}\]
\emph{the cycle associated with $D$}, which is denoted by $z(D)$.
Let $X_1, \ldots, X_{\ell}$ be the irreducible components of $X$ and $\eta_1, \ldots, \eta_{\ell}$
be the generic points of $X_1, \ldots, X_{\ell}$, respectively.
Then
\begin{equation}\label{eqn:remark:expansion:Cartier:div:as:cycle:01}
z(D) = \sum_{j=1}^{\ell} \length_{\OO_{X, \eta_j}}(\OO_{X, \eta_j}) z(\rest{D}{X_j}).
\end{equation}
Indeed, by \cite[(6) of Lemma~1.7]{MArakelov}, $\ord_x(D) = \sum_{j \in J_x} b_j \ord_x(\rest{D}{X_j})$, where $b_j =
\length_{\OO_{X, \eta_j}}(\OO_{X, \eta_j})$ and $J_x = \{ j \mid x \in X_j \}$. Thus if we set
\[
a_{x, j} = \begin{cases} \ord_x(\rest{D}{X_j}) & \text{if $x \in X_j$}, \\
0 & \text{if $x\not\in X_j$},
\end{cases}
\]
then $\ord_x(D) = \sum_{j=1}^{\ell} a_{x, j} b_j$. Thus
\begin{align*}
z(D) & = \sum_{x \in X^{(1)}} \ord_x(D) \overline{\{ x \}} = \sum_{x \in X^{(1)}}  \Big( \sum_{j=1}^{\ell} a_{x, j} b_j \Big)\overline{\{ x \}} \\
& = \sum_{j=1}^{\ell} b_j \sum_{x \in X^{(1)}} a_{x, j} \overline{\{ x \}} = \sum_{j=1}^{\ell} b_j \sum_{x \in X_j^{(1)}}  \ord_x(\rest{D}{X_j})  \overline{\{ x \}} =
 \sum_{j=1}^{\ell} b_j z(\rest{D}{Z_j}),
\end{align*}
as required.

Let $L$ be an invertible $\OO_X$-module and $s$ a regular meromorphic section of $L$ over $X$.
For $x \in X^{(1)}$, $\ord_x(s)$ is defined by $\ord_x(f)$, where $f$ is given by $s = f \omega$ for some local basis $\omega$ of $L$ around $x$.
Note that $\ord_x(s)$ does not depend on the choice of the local basis $\omega$ around $x$.
Then the cycle $z(L; s)$ associated  with $\operatorname{div}(L; s)$ is defined by
\[ z(L;s) := \sum_{x \in X^{(1)}} \ord_x(s) \overline{\{ x \}}.\]

\end{rema}

{
\begin{defi}
Let $\varphi:X\rightarrow Y$ be a morphism of locally ringed space. If $U$ is an open subset of $Y$, we denote by $S_\varphi(U)$ the preimage of $S_X(\varphi^{-1}(U))$ by the structural ring homomorphism 
\[\mathcal O_Y(U)\longrightarrow\mathcal O_X(\varphi^{-1}(U)).\]
We denote by $\mathscr M_\varphi$ the sheaf of commutative and unitary rings associated with the presheaf
\[U\longmapsto \mathcal O_Y(U)[S_\varphi(U)^{-1}].\]
It is a subsheaf of $\mathscr M_Y$. Moreover, the structural morphism of sheaves $\mathcal O_Y\rightarrow \varphi_*(\mathcal O_X)$ induces by localization a morphism $\mathscr M_\varphi\rightarrow \varphi_*(\mathscr M_X)$, which defines a morphism of locally ringed spaces $(X,\mathscr M_X)\rightarrow (Y,\mathscr M_\varphi)$.
\end{defi}

\begin{rema}\label{Rem: pull back}
There are several situations in which $\mathscr M_\varphi$ identifies with $\mathscr M_Y$, notably when one of the following conditions is satisfied (see \cite[$\text{IV}_4.(21.4.5)$]{EGA}):
\begin{enumerate}[label=\rm(\arabic*)]
\item $\varphi$ is flat, namely for any $x\in X$, the morphism of rings $\varphi_x:\mathcal O_{Y,\varphi(x)}\rightarrow\mathcal O_{X,x}$ defines a structure of flat $\mathcal O_{Y,\varphi(x)}$-algebra on $\mathcal O_{X,x}$, 
\item $X$ and $Y$ are locally Noetherian schemes, and $f$ sends any associated point of $X$ to an associated point of $Y$,
\item $X$ and $Y$ are schemes, the set of irreducible components of $Y$ is locally finite, $X$ is reduced, and any irreducible component of $X$ dominates an irreducible component of $Y$.
\end{enumerate}
\end{rema}

\begin{defi}\label{Def: pull back}
Let $\varphi:X\rightarrow Y$ be a morphism of locally ringed spaces, and $D$ be a Cartier divisor on $Y$.  Assume that both $D$ and $-D$ are global sections of $(\mathscr M_Y^{\times}\cap\mathscr M_\varphi)/\mathcal O_X^{\times}$, or equivalently, for any local equation $f$ of $D$ over an open subset $U$ of $Y$, one has $\{f,f^{-1}\}\subset\mathscr M_\varphi(U)$. Then the canonical regular meromorphic section $s_D$ of $\mathcal O_Y(D)$ actually defines an isomorphism 
\[\mathscr M_\varphi\longrightarrow \mathscr M_\varphi\otimes_{\mathcal O_Y}\mathcal O_Y(D).\]
which induces an isomorphisme
\[\varphi^*(s_D):\mathscr M_X\longrightarrow\mathscr M_{X}\otimes_{\mathcal O_X}\varphi^*(\mathcal O_Y(D)).\]
We denote by $\varphi^*(D)$ the Cartier divisor $\operatorname{div}(\varphi^*(\mathcal O_Y(D));\varphi^*(s_D))$ corresponding to this regular meromorphic section, and call it the \emph{pull-back} of $D$ by $\varphi$. In the case where $\varphi$ is an immersion, the Cartier divisor $\varphi^*(D)$ is also denoted by $D|_X$.
\end{defi}
}

Finally let us consider the following lemma.

\begin{lemm}\label{lem:eq:regular:mero:generic:total}
Let $\mathfrak o$ be an integral domain, $A$ be an $\mathfrak o$-algebra and $S := \mathfrak o \setminus \{ 0 \}$.
If $A$ is flat over $\mathfrak o$, then we have the following:
\begin{enumerate}[label=\rm(\arabic*)]
\item For $s \in S$, the homomorphism $s \cdot : A \to A$ given by $a \mapsto s \cdot a$ is injective.
In particular, the structure homomorphism $\mathfrak o \to A$ is injective, so that in the following,
$\mathfrak o$ is considered as a subring of $A$.
\item The natural homomorphism $A \to A_S$ is injective.
\item For $a \in A$, $a$ is a non-zero-divisor in $A$ if and only if $a/1$ is a non-zero-divisor in $A_S$.
In particular, a non-zero-divisor of $A_S$ can be written by a form $a/s$ where
$a$ is a non-zero-divisor of $A$ and $s \in S$.
\item Let $Q(A)$ and $Q(A_S)$ be the total quotient rings of $A$ and $A_S$, respectively.
The homomorphism $Q(A) \to Q(A_S)$ induced by $A \to A_S$ is well-defined and bijective.
In particular, $Q(A)^{\times} = Q(A_S)^{\times}$.
\end{enumerate}
\end{lemm}

\begin{proof}
(1) is obvious because $\mathfrak o$ is an integral domain and $A$ is flat over $\mathfrak o$.
(2) follows from (1).

(3) The assertion follows from (1) and the following commutative diagram:
\[
\begin{CD}
A @>>> A_S \\
@V{a\cdot}VV @VV{a\cdot}V \\
A @>>> A_S
\end{CD}
\]

(4) By (3), if $a \in A$ is a non-zero-divisor, then $a/1$ is a non-zero-divisor in $A_S$, so that 
$Q(A) \to Q(A_S)$ is well-defined. The injectivity of $Q(A) \to Q(A_S)$ follows from (2).
For its surjectivity, observe the following:
\[
\frac{b/t}{a/s} = \frac{(st/1)(b/t)}{(st/1)(a/s)} = \frac{sb/1}{ta/1}.
\]
\end{proof}

\section{Proper intersection}

Let $d$ be a non-negative integer and $X$ be a $d$-dimensional  scheme of finite type over a field $k$. 
Let $D$ be a Cartier divisor on $X$. We define the support of $D$ to be
\[
\operatorname{Supp}(D) := \{ x \in X \mid f_x \not\in \mathcal O_{X, x}^{\times} \},
\]
where $f_x$ is a local equation of $D$ at $x$. Note that the above definition does not depend on the choice of $f_x$ since two local equations of $D$ at $x$ differ by a factor in $\mathcal O_{X,x}^{\times}$.

\begin{prop}\label{prop:Supp:closed}
\begin{enumerate}[label=\rm(\arabic*)]
\item $\operatorname{Supp}(D)$ is a Zariski closed subset of $X$.

\item $\operatorname{Supp}(D + D') \subseteq \operatorname{Supp}(D) \cup \operatorname{Supp}(D')$.
\end{enumerate}
\end{prop}

\begin{proof}
(1) Clearly we may assume that $X$ is affine and $D$ is principal, that is,
$X = \operatorname{Spec}(A)$ and $D$ is defined by a regular meromorphic function $f$ on $X$, which could be considered as an element of the total fraction ring of $A$ (that is, the localization of $A$ with respect to the subset of non-zero-divisors). By \cite{MR570309}, for any prime ideal $\mathfrak p$ of $A$, there is a canonical ring homomorphism from the total fraction ring of $A$ to that of $A_{\mathfrak p}$.
We set $\mathfrak a = \{a \in A \mid af \in A \}$ and $\mathfrak b = \mathfrak af$. Then $\mathfrak a$ and $\mathfrak b$ are ideals of $A$.  Note that, for $\mathfrak p \in \operatorname{Spec}(A)$,
\[\mathfrak a_{\mathfrak p} = \{ u \in A_{\mathfrak p} \mid u f \in A_{\mathfrak p} \}.\] In fact, clearly one has $\mathfrak a_{\mathfrak p}\subseteq\{u\in A_{\mathfrak p}\,|\,uf\in A_{\mathfrak p}\}$. Conversely, if $u=a/s$ (with $a\in A$ and $s\in A\setminus\mathfrak p$) is an element of $A_{\mathfrak p}$ such that $uf\in A_{\mathfrak p}$, then there exists $t\in A\setminus\mathfrak p$ such that $at\in\mathfrak a$ and hence $u=at/st\in\mathfrak a_{\mathfrak p}$.
Thus 
\[
\mathfrak p \not\in \operatorname{Supp}(D) \Longleftrightarrow f \in A_\mathfrak p^{\times} \Longleftrightarrow \text{$\mathfrak a_{\mathfrak p} = A_{\mathfrak p}$ and $\mathfrak b_{\mathfrak p} = A_{\mathfrak p}$}
\Longleftrightarrow \mathfrak p \not\in V(\mathfrak a) \cup V(\mathfrak b),
\]
that is, $\operatorname{Supp}(D) = V(\mathfrak a) \cup V(\mathfrak b)$, as desired.

\medskip
(2) Let $f_x$ and $f'_x$ be local equations of $D$ and $D'$ at $x$, respectively.
Then
\begin{align*}
x \not\in \operatorname{Supp}(D) \cup \operatorname{Supp}(D') & \Longrightarrow f_x, f'_x \in \mathcal O^{\times}_{X, x} \Longrightarrow
f_xf'_x \in \mathcal O^{\times}_{X, x} \\
&  \Longrightarrow  x \not\in \operatorname{Supp}(D + D'),
\end{align*}
as required.
\end{proof}

\begin{defi}\label{def:meet:properly}
Let $n$ be an integer such that $0 \leqslant n \leqslant d$.
Let $D_0, \ldots, D_n$ be Cartier divisors on $X$.
We say that \emph{$D_0, \ldots, D_n$ intersect properly} if, for any non-empty subset $J$ of $\{ 0, \ldots, n\}$,
\[\dim \Big( \bigcap\nolimits_{j \in J} \operatorname{Supp}(D_j) \Big) \leqslant d - \operatorname{card}(J).\] By convention, $\dim(\emptyset)$ is defined to be $-1$.
We set
\[
\mathcal{IP}_X^{(n)} := \{ (D_0, \ldots, D_n) \in \operatorname{Div}(X)^{n+1} \mid \text{$D_0, \ldots, D_n$ intersect properly} \}.
\]
In the case where $n=d$, we often denote $\mathcal{IP}_X^{(n)}$ by $\mathcal{IP}_X$.
\end{defi}

\begin{lemm}\label{lemm:equality:supp:extension:field}
Let $k'/k$ be an extension of fields.
Let $A$ be a $k$-algebra and $A' := A \otimes_k k'$.
Let $\pi : \Spec(A') \to \Spec(A)$ be the morphism induced by the natural homomorphism $A \to A'$.
Let $Q(A)$ (resp. $Q(A')$) be the total quotient ring of $A$ (resp. $A'$).
Let $\alpha \in Q(A)^{\times}$ and $\alpha' := \alpha \otimes_k 1 \in Q(A) \otimes_k k'$.
If we set
\[
\begin{cases}
\Supp(\alpha) := \{ P \in \Spec(A) \mid \alpha \not\in A_P^{\times} \}, \\
\Supp(\alpha') := \{ P' \in \Spec(A') \mid \alpha' \not\in {A'}_{P'}^{\times} \},
\end{cases}
\]
then $\Supp(\alpha') = \pi^{-1}(\Supp(\alpha))$.
\end{lemm}

\begin{proof}
First of all, note that $Q(A) \otimes_k k' \subseteq Q(A')$ and $\alpha' \in (Q(A) \otimes_k k')^{\times} \subseteq Q(A')^{\times}$ because $\pi$ is flat. 
Let 
$I := \{ a \in A \mid a \alpha \in A \}$, $J:= I \alpha$, 
$I':= \{ a' \in A' \mid a' \alpha' \in A' \}$ and $J':= I' \alpha'$.
Then one has the following.

\begin{enonce}{Claim}\label{Claim:lemm:equality:supp:extension:field:01}
\begin{enumerate}[label=\rm(\arabic*)]
\item $\Supp(\alpha) = \Spec(A/I) \cup \Spec(A/J)$ and $\Supp(\alpha') = \Spec(A'/I') \cup \Spec(A'/J')$.
\item $I' = I \otimes_k k'$ and $J'=J \otimes_k k'$.
\item $\Spec(A'/I') = \pi^{-1}(\Spec(A/I))$ and $\Spec(A'/J') = \pi^{-1}(\Spec(A/J))$.
\end{enumerate}
\end{enonce}

\begin{proof}
Let $\{ x_\lambda \}_{\lambda\in\Lambda}$ be a basis of $k'$ over $k$.
Note that $V \otimes_k k' = \bigoplus_{\lambda \in \Lambda} V \otimes_k k x_{\lambda}$ for a $k$-module $V$.

(1) It is sufficient to prove the first equality. The second is similar to the first.
Note that $I_P = \{ a \in A_P \mid a \alpha \in A_P^{\times} \}$.
Thus, if $\alpha \in A_P^{\times}$, then $I_P = J_P = A_P$, so that $P \not\in \Spec(A/I) \cup \Spec(A/J)$.
Conversely, we assume that $P \not\in \Spec(A/I) \cup \Spec(A/J)$, that is, $I \not\subseteq P$ and
$J \not\subseteq P$. Thus $I_P = J_P = A_P$, and hence $\alpha \in A_P^{\times}$.

(2) 
Obviously $I \otimes_k k' \subseteq  I'$. We assume $a' \in I'$.
Then there are $a_{\lambda}$'s such that $a_{\lambda} \in A$ and $a' = \sum_{\lambda} a_{\lambda} \otimes x_{\lambda}$. By our assumption,
we can find $b_{\lambda}$'s such that $b_{\lambda} \in A$ and
\[
\sum\nolimits_{\lambda} a_{\lambda} \alpha \otimes x_{\lambda} = a' \alpha' = \sum\nolimits_{\lambda} b_{\lambda} \otimes x_{\lambda},
\]
so that $a_{\lambda}\alpha = b_{\lambda} \in A$ for all $\lambda$. Thus $a_{\lambda} \in I$.
Therefore the first assertion follows. The second is a consequence of the first.

(3) follows from (2).
\end{proof}
By using (1) and (3) of the above claim,
\begin{align*}
\pi^{-1}(\Supp(\alpha)) & = \pi^{-1}(\Spec(A/I) \cup \Spec(A/J)) \\
& = \pi^{-1}(\Spec(A/I)) \cup \pi^{-1}(\Spec(A/J)) \\
& = \Supp(A'/I') \cup \Supp(A'/J') = \Supp(\alpha'),
\end{align*}
as required.
\end{proof}

\begin{rema}\label{Rem: extension of scalars}
Let $k'/k$ be an extension of fields, $X_{k'}=X\times_{\Spec k}\Spec k'$ and $\pi:X_{k'}\rightarrow X$ be the morphism of projection. Since the canonical morphism $\Spec k'\rightarrow\Spec k$ is flat, also is the morphism of projection $\pi$ (see \cite[$\text{IV}_1.(2.1.4)$]{EGA}). Therefore, for any Cartier divisor $D$ on $X$, the pull-back $\pi^*(D)$ is well defined as a Cartier divisor on $X_{k'}$, which we denote by $D_{k'}$.

By Lemma~\ref{lemm:equality:supp:extension:field}, one has
\[\operatorname{Supp}(D_{k'})=\pi^{-1}(\operatorname{Supp}(D)).\]
In particular, if $D_0,\ldots,D_n$ are Cartier divisors on $X$, which intersect properly, then, for any subset $J$ of $\{0,\ldots,n\}$, one has (see for example \cite[Proposition 5.38]{MR2675155} for the equality in the middle)
\[\begin{split}\dim \Big( \bigcap\nolimits_{j \in J} \operatorname{Supp}(D_{j,k'}) \Big)&= \dim \Big( \pi^{-1}\Big(\bigcap\nolimits_{j \in J} \operatorname{Supp}(D_j) \Big)\Big)\\&=\dim \Big( \bigcap\nolimits_{j \in J} \operatorname{Supp}(D_j) \Big)\leqslant d-\operatorname{card}(J).
\end{split}\] 
Therefore, the Cartier divisors $D_{0,k'},\ldots,D_{n,k'}$ on $X_{k'}$ intersect properly.
\end{rema}

\begin{lemm}\label{lem:meet:properly}
The set $\mathcal{IP}_X^{(n)}$ forms a symmetric and multi-linear subset of $\operatorname{Div}(X)^{n+1}$ in the sense of Definition~\ref{def:domain}.
\end{lemm}

\begin{proof}
It is sufficient to show that if $(D_0, D_1, \ldots, D_n), (D'_0, D_1, \ldots, D_n) \in \mathcal{IP}_X^{(n)}$, then
$(D_0 + D'_0, D_1, \ldots, D_n) \in \mathcal{IP}_X^{(n)}$. We set 
\[
D''_i = \begin{cases}
D_0 + D'_0, & \text{if $i = 0$}, \\
D_i, & \text{if $i \geqslant 1$}.
\end{cases}
\]
If $(D''_0, D''_1, \ldots, D''_n) \not\in \mathcal{IP}_X^{(n)}$, then there is a non-empty subset $J$ of $\{ 0, \ldots, n \}$ such that
\[
\dim \Big( \bigcap\nolimits_{j \in J} \operatorname{Supp}\big(D''_{j} \big) \Big) >  d - \#(J).
\]
Clearly $0 \in J$.
We can find a schematic point $P \in X$ such that $\dim \overline{\{ P \}} > d - \#(J)$ and
$P \in \operatorname{Supp}\big(D''_{j} \big)$ for all $j \in J$, so that $P \in \operatorname{Supp}(D_0 + D'_0)$ and $P \in \operatorname{Supp}(D_j)$ for $j \in J \setminus \{ 0 \}$.
Thus, by Proposition~\ref{prop:Supp:closed}, $P \in \operatorname{Supp}(D_0)$ or $P \in \operatorname{Supp}(D'_0)$, which is a contradiction.
\end{proof}

\begin{lemm}\label{lem:IP:projective:ample}
We assume that $X$ is projective. Let $n$ be an integer such that $0 \leqslant n \leqslant d$.
\begin{enumerate}[label=\rm(\arabic*)]
\item
Let $L_0, \ldots, L_n$ be invertible $\mathcal O_X$-modules.
Then there are regular meromorphic sections $s_0, \ldots, s_n$ of $L_0, \ldots, L_n$,
respectively, such that, if we set $D_i = \operatorname{div}(L_i; s_i)$ for $i \in \{ 0, \ldots, n\}$,
then $D_0, \ldots, D_n$ intersect properly.

\item
If $(D_0, D_1, \ldots, D_n), (D'_0, D'_1, \ldots, D'_n)
\in \mathcal{IP}^{(n)}_X$ and $D_0 \sim D'_0$, then there is $D''_0$ such that
$D''_0 \sim D_0 \,(\sim D'_0)$ and $(D''_0, D_1, \ldots, D_n), (D''_0, D'_1, \ldots, D'_n)
\in \mathcal{IP}_X^{(n)}$.
\end{enumerate}
\end{lemm}

\begin{proof}
(1) We prove it by induction on $n$ in incorporating the proof of the initial case in the induction procedure. 
By the hypothesis of induction (when $n\geqslant 1$),
there are regular meromorphic sections $s_0, \ldots, s_{n-1}$ of $L_0, \ldots, L_{n-1}$,
respectively, such that if we set $D_i = \operatorname{div}(L_i; s_i)$ for $i \in \{ 0, \ldots, n-1\}$,
then $D_0, \ldots, D_{n-1}$ intersect properly.

We now introduce the following claim, which (in the case where $n=0$) also prove the initial case of induction.
\begin{enonce}{Claim}\label{Claim: decomposition}
There exist very ample invertible $\mathcal O_X$-modules $L'_n$ and $L''_n$, and global sections $s'_n$ and $s''_n$ of $L'_n$ and $L''_n$, which satisfy the following conditions~:
\begin{enumerate}[label=\rm(\roman*)]
\item $L_n = L'_n \otimes {L''_n}^{-1}$,
\item $s'_n$ and $s''_n$ define regular meromorphic sections of $L_n'$ and $L_n''$, respectively,
\item\label{Item intersection properly} if we set $D'_n = \operatorname{div}(L'_n; s'_n)$ and $D''_n = \operatorname{div}(L''_n; s''_n)$,
then both families of Cartier divisors $D_0, \ldots, D_{n-1}, D'_n$ and $D_0, \ldots, D_{n-1}, D''_n$
intersect properly.
\end{enumerate} 
\end{enonce}
\begin{proof}[Proof of Claim \ref{Claim: decomposition}]
Since $X$ is projective, there exists a very ample $\mathcal O_X$-module $L$. By \cite[II.(4.5.5)]{EGA}, there exists an integer $\ell_0\in\mathbb N_{\geqslant 1}$ such that both invertible $\mathcal O_X$-modules $L^{\otimes \ell_0}$ and $L^{\otimes\ell_0}\otimes L_n^{-1}$ are generated by global sections.  Let $\Sigma$ be the set of generic points of 
\[\bigcap_{i=0}^{n-1}\operatorname{Supp}(D_i).\]
We equip the set $\Sigma\cup\operatorname{Ass}(X)$ with the order $\succ$ of generalization, namely $x\succ y$ if and only if $y$ belongs to the Zariski closure of $\{x\}$. We denote by $\{y_1,\ldots,y_b\}$ the set of all minimal elements of the set $\Sigma\cup\operatorname{Ass}(X)$ .

For any $i\in\{1,\ldots,b\}$, one has
\[y_i\in X\setminus\bigcup_{\begin{subarray}{c}j\in\{1,\ldots,b\}\\
j\neq i\end{subarray}}\overline{\{y_j\}}.\]
By \cite[II.(4.5.4)]{EGA}, for any $i\in\{1,\ldots,b\}$, there exists $\ell_i\in\mathbb N_{\geqslant 1}$ and a section $t_i\in H^0(X,L^{\otimes n})$ such that $t_i(y_i)\neq 0$ and that $t_i(y_j)=0$ for any $j\in\{1,\ldots,b\}\setminus\{i\}$. Moreover, by replacing the global sections $t_1,\ldots,t_b$ by suitable powers, we may assume, without loss of generality, that all $\ell_1,\ldots,\ell_b$ are equal to a positive integer $\ell$. For any $i\in\{1,\ldots,b\}$, let $u_i\in H^0(X,L^{\otimes{\ell_0}})$ and $v_i\in H^0(X,L^{\otimes\ell_0}\otimes L_n^{-1})$ be such that $u_i(y_i)\neq 0$ and $v_i(y_i)\neq 0$. These sections exist since the invertible $\mathcal O_X$-modules $L^{\otimes\ell_0}$ and $L^{\otimes\ell_0}\otimes L_n^{-1}$ are generated by global sections. Now we take 
\[L_n'=L^{\otimes(\ell_0+\ell)},\qquad L_n''=L_n'\otimes L_n^{-1}=(L^{\ell_0}\otimes L_n^{-1})\otimes L^{\otimes\ell}, \]
and
\[s_n'=\sum_{i=1}^bu_it_i,\qquad s_n''=\sum_{i=1}^bv_it_i.\]
Then, for any $i\in\{1,\ldots,b\}$, one has $s_n'(y_i)\neq 0$ and $s_n''(y_i)\neq 0$. In particular, $s_n'$ and $s_n''$ do not vanish on any of the associated points of $X$ and hence are regular meromorphic sections (see Remark \ref{Rem: regular meromorphic section}). Moreover, since these sections do not vanish on any point of $\Sigma$, we obtain the condition \ref{Item intersection properly} above.
\end{proof}

Thus, by Lemma \ref{lem:meet:properly},
we can see that $D_1, \ldots, D_{n-1}, D_n$ intersect properly,
where $D_n = D'_n - D''_n = \operatorname{div}(L_n; s_n \otimes {s_n}^{-1})$,
as required.

\medskip
(2) We can find very ample Cartier divisors $A$ and $B$ on $X$ such that $D_0 = A - B$.
Then,  by the same argument as the induction procedure in the proof of (1), we obtain that there are $A'$ and $B'$ such that $A' \sim A$, $B' \sim B$ and
\[(A', D_1, \ldots, D_n), (A', D'_1, \ldots, D'_n),(B', D_1, \ldots, D_n), (B', D'_1, \ldots, D'_n)
\in \mathcal{IP}_X^{(n)}.\]
Thus if we set $D''_0 = A' - B'$, then, by Lemma~\ref{lem:meet:properly}, one has the conclusion.
\end{proof}

\begin{rema}\label{Rem: further explanation to Claim}
Claim \ref{Claim: decomposition} has its own interest and will be used in further chapters  in the following way. Let $X$ be a $d$-dimensional projective scheme over $\Spec k$ and $D_0,\ldots,D_d$ be Cartier divisors on $X$. We suppose that $D_0,\ldots, D_d$ intersect properly. Let $D_0=A_0-A_0'$ be a decomposition of $D_0$ into the difference of two very ample Cartier divisors. A priori $A_0,D_1,\ldots,D_d$ do not intersect properly. However, by Claim \ref{Claim: decomposition}, one can find a very ample invertible $\mathcal O_X$-module $L$ and a global section $s$ of $L\otimes\mathcal O_X(A_0)$ defining a regular meromorphic section, such that $\operatorname{div}(L;s),D_1,\ldots,D_d$ intersect properly. Let $B=\operatorname{div}(L;s)-A_0$. This is a very ample Cartier divisor since $\mathcal O_X(B)$ is isomorphic to $L$. Moreover, both $(A_0+B,D_1,\ldots,D_d)$ and $(A_0'+B,D_1,\ldots,D_d)$ belong to $\mathcal {IP}^{(d)}_X$ since the former one and their difference do.
\end{rema}

\section{Multi-homogeneous polynomials}\label{Sec: multi-homogeneous polynomial}

Let $k$ be a field and $(E_i)_{i=0}^d$ be a family of finite-dimensional vector spaces over $k$. Let $(\delta_0,\ldots,\delta_{d})$ be a multi-index in $\mathbb N^{d+1}$. 
\begin{defi}We call \emph{multi-homogeneous polynomial of multi-degree $(\delta_0,\ldots,\delta_d)$} on $E_0\times\cdots\times E_d$ any element of
\[S^{\delta_0}(E_0^\vee)\otimes_k\cdots\otimes_kS^{\delta_d}(E_d^\vee),\]
where $S^{\delta_i}(E_i^\vee)$ denotes the $\delta_i$-th symmetric power of the vector space $E_i^\vee$.
\end{defi}
Recall that one dual vector space of $S^{\delta_i}(E_i^\vee)$ is given by \[\Gamma^{\delta_i}(E_i):=(E_i^{\otimes \delta_i})^{\mathfrak S_{\delta_i}},\]
where $\mathfrak S_{\delta_i}$ is the symmetric group on $\{1,\ldots,\delta_i\}$, which acts on $E_i^{\otimes \delta_i}$ by permuting tensor factors (see \cite[Chapitre IV, \S5, no. 11, proposition 20]{MR0274237}). Therefore, one dual vector space of $S^{\delta_0}(E_0^\vee)\otimes_k\cdots\otimes_kS^{\delta_d}(E_d^\vee)$ is given by
\[\Gamma^{\delta_0}(E_0)\otimes_k\cdots\otimes_k\Gamma^{\delta_d}(E_d).\]
If $R\in S^{\delta_0}(E_0^\vee)\otimes_k\cdots\otimes_kS^{\delta_d}(E_d^\vee)$ is a multi-homogeneous polynomial of multi-degree $(\delta_0,\ldots,\delta_d)$, for any $(s_0,\ldots,s_d)\in E_0\times\cdots\times E_d$, we denote by $R(s_0,\ldots,s_d)$ the value
\[R(s_0^{\otimes \delta_0}\otimes\cdots\otimes s_d^{
\otimes \delta_d})\]
in $k$. Thus $R$ determines a function on $E_0\times\cdots\times E_d$ valued in $K$ (which we still denote by $R$ by abuse of notation). Note that, as an element of $S^{\delta_0}(E_0^\vee)\otimes_k\cdots\otimes_kS^{\delta_d}(E_d^\vee)$, $R$ is uniquely determined by the corresponding function on  $E_0\times\cdots\times E_d$ since each vector space $\Gamma^{\delta_i}(E_i)$ is spanned over $k$ by elements of the form $s_i^{\otimes \delta_i}$, $s_i\in E_i$ (see \cite[Chapitre IV, \S5, no. 5, proposition 5]{MR0274237}). This observation allows us to consider, for any $i\in\{0,\ldots,d\}$ and $s_i\in E_i$, the specification 
\[R(\cdots,\underset{\mathclap{\shortstack{\scriptsize $\uparrow$\\[-.2ex]$i$-th coordinate}}
  }{s_i},\cdots) \] of $R$ at $s_i$ as an element of
\[S^{\delta_0}(E_0^\vee)\otimes_k\cdots\otimes_k S^{\delta_{i-1}}(E_{i-1}^\vee)\otimes_k S^{\delta_{i+1}}(E_{i+1}^\vee)\otimes_k\cdots\otimes_kS^{\delta_d}(E_d^\vee)\]
or as a multi-homogeneous polynomial function on
\[E_0\times\cdots\times E_{i-1}\times E_{i+1}\times\cdots\times E_d.\]

\begin{rema}
Note that an element of $S^{\delta_0}(E_0^{\vee}) \otimes_k \cdots \otimes_k S^{\delta_d}(E_d^{\vee})$ yields
a multi-homogeneous polynomial function on $\mathbb A(E_0^{\vee}) \times_k \cdots \times_k \mathbb A(E_d^{\vee})$ and
the set of $k$-rational points of $\mathbb A(E_0^{\vee}) \times_k \cdots \times_k \mathbb A(E_d^{\vee})$ is naturally isomorphic to
$E_0 \times \cdots \times E_d$, where $\mathbb A(E_i^{\vee}) = \operatorname{Spec}(\bigoplus_{\delta=0}^{\infty} S^{\delta}(E_i^{\vee}))$
for each $i$.
\end{rema}

\section{Incidence subscheme}
Let $k$ be a field and $E$ be a finite-dimensional vector space over $k$. 
We denote $\operatorname{Proj}(\bigoplus_{\delta=0}^{\infty} S^{\delta}(E))$ by $\mathbb P(E)$.
Recall that the projective space $\mathbb P(E)$ represents the  contravariant functor from the category of $k$-schemes to that of sets, which sends a $k$-scheme $\varphi:S\rightarrow\operatorname{Spec} k$ to the set of all invertible quotient $\mathcal O_S$-modules of $\varphi^*(E)$. In particular, if we denote by $\pi_E:\mathbb P(E)\rightarrow\operatorname{Spec} k$ the structural scheme morphism, then the universal object of the representation of this functor by $\mathbb P(E)$ is a quotient $\mathcal O_{\mathbb P(E)}$-module of $\pi_E^*(E)$, which we denote by $\mathcal O_E(1)$ and which we call \emph{universal invertible sheaf} on $\mathbb P(E)$.  For any positive integer $n$, we let $\mathcal O_E(n):=\mathcal O_E(1)^{\otimes n}$ and $\mathcal O_E(-n):=(\mathcal O_E(1)^\vee)^{\otimes n}$. Note that the quotient homomorphism $\pi_E^*(E)\rightarrow\mathcal O_E(1)$ induces by passing to dual modules an injective homomorphism
\[\mathcal O_E(-1)\longrightarrow\pi_E^*(E^\vee).\]

We now consider the fibre product of projective spaces $\mathbb P(E)\times_k\mathbb P(E^\vee)$. Let \[p_1:\mathbb P(E)\times_k\mathbb P(E^\vee)\longrightarrow \mathbb P(E)\quad\text{ and }\quad p_2:\mathbb P(E)\times_k\mathbb P(E^\vee)\longrightarrow \mathbb P(E^\vee)\]
be morphisms of projection. Note that the following diagram of scheme morphisms is cartesian
\[\xymatrix{\mathbb P(E)\times_k\mathbb P(E^\vee)\ar[r]^-{p_2}\ar[d]_-{p_1}&\mathbb P(E^\vee)\ar[d]^-{\pi_{E^\vee}}\\
\mathbb P(E)\ar[r]_-{\pi_E}&\operatorname{Spec} k}\] 
The composition 
of the homomorphisms
\begin{equation}\label{Equ: incident subscheme}p_1^*(\mathcal O_E(-1))\longrightarrow p_1^*(\pi_E^*(E^\vee))\cong p_2^*(\pi_{E^\vee}^*(E^\vee))\longrightarrow p_2^*(\mathcal O_{E^\vee}(1))\end{equation}
determines a global section of the invertible sheaf \[\mathcal O_E(1)\boxtimes\mathcal O_{E^\vee}(1):= p_1^*(\mathcal O_E(1))\otimes p_2^*(\mathcal O_{E^\vee}(1)).\]

\begin{defi}
We call \emph{incidence subscheme} of $\mathbb P(E)\times_k\mathbb P(E^\vee)$ and we denote by $I_E$ the closed subscheme of $\mathbb P(E)\times_k\mathbb P(E^\vee)$ defined by the vanishing of the global section of $\mathcal O_E(1)\boxtimes\mathcal O_{E^\vee}(1)$ determined by \eqref{Equ: incident subscheme}. In particular, the cycle class of $I_E$ modulo the linear equivalence is
\[c_1(\mathcal O_E(1)\boxtimes\mathcal O_{E^\vee}(1))\cap[\mathbb P(E)\times_k\mathbb P(E^\vee)].\]
\end{defi}

The following proposition shows that the incidence subscheme can be realized as a projective bundle over $\mathbb P(E)$.

\begin{prop}\label{Pro: indcidence variety as bundle}
Let $Q_{E^\vee}$ be the quotient sheaf of $\pi_E^*(E^\vee)$ by the canonical image of $\mathcal O_E(-1)$. Then the incidence subscheme $I_E$ is isomorphic as a $\mathbb P(E)$-scheme to the projective bundle $\mathbb P(Q_{E^\vee}\otimes\mathcal O_E(1))$. Moreover, under this isomorphism, the restriction of $\mathcal O_E(1)\boxtimes\mathcal O_{E^\vee}(1)$ to $I_E$ is isomorphic to the universal invertible sheaf of the projective bundle $\mathbb P(Q_{E^\vee}\otimes\mathcal O_E(1))$.
\end{prop} 
\begin{proof}
It suffices to identify $p_1:\mathbb P(E)\times_k\mathbb P(E^\vee)\rightarrow\mathbb P(E)$ with the projective bundle \[\mathbb P(\pi_E^*(E^\vee)\otimes\mathcal O_E(1))\longrightarrow\mathbb P(E).\]
Note that the universal invertible sheaf of this projective bundle is isomorphic to $\mathcal O_E(1)\boxtimes\mathcal O_{E^\vee}(1)$. Under this identification, the vanishing locus of \eqref{Equ: incident subscheme} coincides with the projective bundle $\mathbb P(Q_{E^\vee}\otimes\mathcal O_E(1))$. 
\end{proof}

\begin{rema}\label{Rem: fibre incidence variety}
As a scheme over $\mathbb P(E)$, the incident subscheme $I_E$ also identifies with the projective bundle $\mathbb P(Q_{E^\vee})$. However, the universal invertible sheaf of this projective bundle is the restriction of $p_2^*(\mathcal O_{E^\vee}(1))$. Moreover, we can also consider the morphism of projection from the incidence subscheme to $\mathbb P(E^\vee)$. By the duality between $E$ and $E^\vee$, the incidence subscheme $I_E$ also identifies with the projective bundle of $Q_{E}:=\pi_{E^\vee}^*(E)/\mathcal O_{E^\vee}(-1)$ over $\mathbb P(E^\vee)$. In particular, if $x$ is a point of $\mathbb P(E^\vee)$, then the fibre of the incidence subscheme $I_{E}$ over $x$ identifies with \[\mathbb P((E\otimes_k\kappa(x))/x^*\mathcal O_E(-1)),\]
which is a hyperplane in $\mathbb P(E\otimes_k\kappa(x))$ defined by the vanishing locus of any non-zero element of the one-dimensional $\kappa(x)$-vector subspace of $E\otimes_k\kappa(x)$ defining the point $x$.
\end{rema}

\section{Resultants}

Let $k$ be a field and $X$ be an integral projective $k$-scheme, and $d$ be the Krull dimension of $X$. For any $i\in\{0,\ldots,d\}$, we fix a finite-dimensional vector space $E_i$ over $k$ and a closed embedding $f_i:X\rightarrow\mathbb P(E_i)$, and we denote by $L_i$ the pull-back of $\mathcal O_{E_i}(1)$ by $f_i$. For each $i\in\{0,\ldots,d\}$, we let $r_i$ be the Krull dimension of $\mathbb P(E_i)$, which identifies with $\operatorname{dim}_k(E_i)-1$. For each $i\in\{0,\ldots,d\}$, we let $\delta_i$ be the intersection number
\[\deg\big(c_1(L_0)\cdots c_1(L_{i-1})c_1(L_{i+1})\cdots c_1(L_d)\cap[X]\big)\]

Let $\mathbb P=\mathbb P(E_0)\times_k\cdots\times_k\mathbb P(E_d)$ be the product of $k$-schemes $(\mathbb P(E_i))_{i=0}^d$. The family $(f_i)_{i=0}^d$ induces a closed embedding $f:X\rightarrow \mathbb P$. Let \[\check{\mathbb P}:=\mathbb P(E_0^\vee)\times_k\cdots\times_k\mathbb P(E_d^\vee)\] 
be the product of dual projective spaces. We identify $\mathbb P\times_k\check{\mathbb P}$ with
\[(\mathbb P(E_0)\times_k\mathbb P(E_0^\vee))\times_k\cdots\times_k(\mathbb P(E_d)\times_k\mathbb P(E_d^\vee))\]
and we denote by \[I_{\mathbb P}:=I_{E_0}\times_k\cdots\times_kI_{E_d}\] the fibre product of incidence subschemes,
so that the class of $I_{\mathbb P}$ modulo the linear equivalence coincides with the intersection product
\[
c_1(r_0^*(\mathcal O_{E_0}(1)\boxtimes\mathcal O_{E_0^\vee}(1))) \cdots c_1(r_n^*(\mathcal O_{E_n}(1)\boxtimes\mathcal O_{E_n^\vee}(1)))\cap[\mathbb P\times_k\check{\mathbb P}],
\] where $r_i : \mathbb P\times_k\check{\mathbb P} \to \mathbb P(E_i)\times_k\mathbb P(E_i^\vee)$ is the $i$-th projection.
By Proposition \ref{Pro: indcidence variety as bundle} (see also Remark \ref{Rem: fibre incidence variety}), $I_{\mathbb P}$ is isomorphic to a fiber product of projective bundles
\[\mathbb P(Q_{E_0})\times_k\cdots\times_k\mathbb P(Q_{E_d}).\]

\begin{defi}
We denote by $I_X$ the fibre product 
$X\times_{\mathbb P}I_{\mathbb P}$, called the \emph{incidence subscheme} of $X\times_{k}\check{\mathbb P}$. As an $X$-scheme, it identifies with
\[\mathbb P(Q_{E_0}|_X)\times_{X} \cdots\times_{X}\mathbb P(Q_{E_d}|_X).\]
and hence is an integral closed subscheme of dimension 
\[d+(r_0-1)+\cdots+(r_d-1)=r_0+\cdots+r_d-1\]
of $\mathbb P\times_k\check{\mathbb P}$.
In particular, for any extension $K$ of $k$ and any element 
\[(x,\alpha_0,\ldots,\alpha_d)\in X(K)\times\mathbb P(E_0^\vee)(K)\times\cdots\times\mathbb P(E_d^\vee)(K),\]
if we denote by $H_i$ the hyperplane in $\mathbb P(E_{i,K})$ defined by the vanishing of $\alpha_i$, then $(x,\alpha_0,\ldots,\alpha_d)$ belongs to $I_X(K)$ if and only if $f_{i,K}(x)\in H_i$ for any $i\in\{1,\ldots,d\}$. 
In addition, 
the cycle class of $I_X$ modulo the linear equivalence is
the intersection product
\begin{equation} \label{eqn:class:incidence:subscheme}
c_1\left( p^*(L_0) \otimes q^*q_0^*(\mathcal O_{E_0^{\vee}}(1))\right)  \cdots c_1\left( p^*(L_d) \otimes q^*q_d^*(\mathcal O_{E_d^{\vee}}(1))\right)\cap[X\times_k\check{\mathbb P}], \end{equation}
where $p : X \times_k \check{\mathbb P} \to X$, $q : X \times_k \check{\mathbb P} \to \check{\mathbb P}$ and $q_i :  \check{\mathbb P} \to \mathbb P(E_i^{\vee})$ are the projections. 
\end{defi}

\if00

\begin{prop}
The direct image  by the projection $q : X \times_k\check{\mathbb P}\rightarrow\check{\mathbb P}$ of $I_X$ is a multi-homogeneous hypersurface of multi-degree $(\delta_0,\ldots,\delta_d)$.
\end{prop}
\begin{proof}
It is sufficient to see that $q_*(I_X)$ belongs to the cycle class
\[
c_1(\mathcal O_{E_0^{\vee}}(\delta_0) \boxtimes \cdots \boxtimes \mathcal O_{E_d^{\vee}}(\delta_d))\cap[\check{\mathbb P}].
\]
Note that, for any $(i_1,\ldots,i_n)\in\{0,\ldots,d\}^n$ such that $i_1,\ldots, i_n$ are distinct,
\[
q_*(c_1(p^*(L_{i_1})) \cdots c_1(p^*(L_{i_n}))\cap[X\times_k\check{\mathbb P}])\]
is equal to
\[c_1(\mathcal O_{E_{i_1}^\vee}(1)\boxtimes\cdots\boxtimes\mathcal O_{E_{i_n}^\vee}(1))\cap [\check{\mathbb P}]
\]
if $n=d$, and is equal to the zero cycle class otherwise.
Therefore, the assertion follows from \eqref{eqn:class:incidence:subscheme}.
\end{proof}

\else
\begin{prop}
The direct image  by the morphism of projection $\mathbb P\times_k\check{\mathbb P}\rightarrow\check{\mathbb P}$ of the effective algebraic cycle  in $\mathbb P\times_k\check{\mathbb P}$ associated with $I_X$ is a multi-homogeneous hypersurface of multi-degree $(\delta_0,\ldots,\delta_d)$.
\end{prop}
\begin{proof}
We reason by induction on $d$. First we assume that $d=0$. Without loss of generality, we may assume that $X=\operatorname{Spec} k'$, where $k'$ is a finite extension of $k$. Let $x$ be the image of $f_0:X\rightarrow\mathbb P(E_0)$, which corresponds to a one-dimensional $k'$-vector subspace of $E_0\otimes_kk'$
\[x^*\mathcal O_{E_0}(-1)\longrightarrow E_0^\vee\otimes_kk'. \] 
The scheme $I_X=X\times_{\mathbb P}I_{\mathbb P}$ identifies with $\mathbb P(Q_{E_0})$, where \[Q_{E_0}=(E_0^\vee\otimes_kk')/x^{*}\mathcal O_{E_0}(-1).\]
This is a hyperplane of \[\mathbb P(E_0^\vee\otimes_kk')\cong\mathbb P(E_0^\vee)\times_{\operatorname{Spec} k}\operatorname{Spec} k'.\]
Note that the field of definition of the quotient vector space $Q_{E_0}$ identifies with the residue field $\kappa(x)$ of the point $x$. Hence the projection of $\mathbb P(Q_{E_0})$ in $\mathbb P(E_0^\vee)$ is a hypersurface of degree $[\kappa(x):k]$.

We now proceed with the case $d\geqslant 1$ in supposing that the proposition is true for equidimensional projective $K$-scheme of dimension $d-1$. Without loss of generality, we may assume that $X$ is an integral projective $k$-scheme. Let $\eta$ be the generic point of $\mathbb P(E_d^\vee)$ and $k_d$ be the field of rational functions on $\mathbb P(E_d^\vee)$. Note that the point $\eta$ corresponds to a one-dimension vector subspace of \[E_d\otimes_kk_d=H^0(\mathbb P(E_d\otimes_kk_d),\mathcal O_{E_d\otimes_kk_d}(1)).\] We pick a non-zero element $s$ in this one-dimensional vector subspace and let $Z$ be the effective cycle of $X\times_{\operatorname{Spec} k}\operatorname{Spec} k_d$ defined by the vanishing of $f_{d,k_d}^*(s)$. Note that 
$I_{X}\times_{\mathbb P(E_d^\vee),\eta}\operatorname{Spec} k_d$
is isomorphic to $I_{Z_s}$ (see Remark \ref{Rem: fibre incidence variety}). 
We apply the induction hypothesis to $Z_s$ and obtain that the projection of the effective cycle associated with $I_{Z_s}$ is a multi-homogeneous hypersurface of multi-degree $(\delta_0,\ldots,\delta_{d-1})$. Therefore the projection of $I_X$ in $\check{\mathbb P}$ is also a hypersurface, which is multi-homogeneous  of multi-degree $(\delta_0,\ldots,\delta_{d-1})$ in the first $d-1$ coordinates. The same argument applied to fibres of $I_X$ over the generic points of $\mathbb P(E_0^\vee),\ldots,\mathbb P(E_{d-1}^\vee)$ shows that the hypersurface is actually multi-homogeneous in all coordinates and is of multi-degree $(\delta_0,\ldots,\delta_d)$.  
\end{proof}
\fi

\begin{defi}
Let $X$ be an integral projective $k$-scheme of dimension $d$. We call \emph{resultant of $X$ with respect to $(f_i)_{i=0}^d$} any  multi-homogeneous polynomial of multi-degree $(\delta_0,\ldots,\delta_d)$ on $E_0\times\cdots\times E_d$,   whose vanishing cycle in \[\mathbb P(E_0^\vee)\times_k\cdots\times_k\mathbb P(E_d^\vee)\] identifies with the projection of the cycle associated with the incidence subscheme $I_X$. Note that the resultant of $X$ with respect to $(f_i)_{i=0}^d$ is unique up to a factor of scalar in $k\setminus\{0\}$ as an element of $S^{\delta_0}(E_0^{\vee}) \otimes_k \cdots \otimes_k S^{\delta_d}(E_d^{\vee})$.

In general, if $X$ is a projective $k$-scheme of dimension $d$ and if 
\[\sum_{i=1}^nm_i X_i\]
is the $d$-dimensional part of the fundamental cycle of $X$, where $X_1,\ldots, X_n$ are $d$-dimensional irreducible components of $X$, and $m_i$ is the local multiplicity of $X$ at the generic point of $X_i$, we define the \emph{resultant} of $X$ with respect to $(f_i)_{i=0}^d$ as 
any multi-homogeneous polynomial of the form
\[(R_{f_0|_{X_1},\ldots,f_d|_{X_1}}^{X_1})^{m_1}\cdots (R_{f_0|_{X_n},\ldots, f_d|_{X_n}}^{X_n})^{m_n},\]
where each $R_{f_0|_{X_i},\ldots,f_d|_{X_i}}^{X_i}$ is a resultant of $X_i$ with respect to $(f_i|_{X_i})_{i=0}^d$.
\end{defi}

\begin{exem}\label{Exe: explicit construction of resultant}
We consider the particular case where $d=0$. Let $f_0:X\rightarrow\mathbb P(E_0)$ be a close embedding. We first assume that $X$ is integral. In this case $f_0$ sends $X$ to a  closed point $x$ of $\mathbb P(E_0)$. Let $\kappa(x)$ be the residue field of $x$ and $\delta_0=[\kappa(x):k]$ be the degree of $x$. Let $s_0$ be an element of $E_0$. We assume that, if we view $s_0$ as a global section of $\mathcal O_{E_0}(1)$, one has $s_0(x)\neq 0$. We construct an element $R^{X,s_0}_{f_0}\in S^{\delta_0}(E_0^\vee)$ as follows. Let \[\varphi_0:E_0\otimes_K\kappa(x)\longrightarrow\mathcal O_{E_0}(1)(x)\]
be the surjective $\kappa(x)$-linear map corresponding to the closed point $x$, and
\[\varphi_0^\vee:\mathcal O_{E_0}(-1)(x)\longrightarrow E_0^\vee\otimes_K\kappa(x)\]
be the dual $\kappa(x)$-linear map of $\varphi_0$, which is an injective linear map. Let $s_0(x)^\vee$ be the unique $\kappa(x)$-linear form on $\mathcal O_{E_0}(1)(x)$ taking the value $1$ at $s_0(x)$. We let \[R_{f_0}^{X,s_0}:=N_{\kappa(x)/K}(\varphi_0^\vee(s_0(x)^\vee))\in S^{\delta_0}(E_0^\vee),\]
which is defined as the determinant of the following homothety endomorphism of the free module $\operatorname{Sym}(E_0^{\vee})\otimes_K\kappa(x)$ of rank $\delta_0$ over the symmetric algebra $\operatorname{Sym}(E_0^\vee)$
\[\xymatrix@C+2pc{\operatorname{Sym}(E_0^{\vee})\otimes_K\kappa(x)\ar[r]^-{\varphi_0^\vee(s_0(x)^\vee)}&\operatorname{Sym}(E_0^{\vee})\otimes_K\kappa(x).}\] 
Note that \[\varphi_0^\vee(s_0(x)^\vee)(s_0\otimes 1)=s_0(x)^\vee(s_0(x))=1.\]
Therefore the following equality holds
\[R_{f_0}^{X,s_0}(s_0)=1.\]

Assume that $X$ is not irreducible.
We let $X_1,\ldots,X_n$ be irreducible components of $X$ (namely points of $X$). For each $i\in\{1,\ldots,n\}$, let $x_i=f_0(X_i)$ and $a_i$ be the local multiplicity of $X$ at $X_i$. Then\[a_1x_1+\cdots+a_nx_n\] is the decomposition of $f(X)$ as a zero-dimensional cycle in $\mathbb P(E_0)$, where $x_1,\ldots,x_n$ are closed points of $\mathbb P(E_0)$ and $a_1,\ldots,a_n$ are positive integers. If $s_0$ is a global section of $\mathcal O_{E_0}(1)$, which does not vanish on any of the points $x_1,\ldots,x_n$, we define
\[R_{f_0}^{X,s_0}:=\prod_{i=1}^n(R_{f_0|_{X_i}}^{X_i,s_0})^{a_i}.\]
Then $R_{f_0}^{X,s_0}$ is a resultant of $X$ with respect to the closed embedding $f_0$, which satisfies $R_{f_0}^{X,s_0}(s_0)=1$.
\end{exem}

\begin{exem}
Let $n$ and $m$ be positive integers, and let \[f : \mathbb P^1 \longrightarrow \mathbb P^n, \qquad (x,y) \longmapsto (x^iy^{n-i})_{i=0}^n\] and \[g : \mathbb P^1 \longrightarrow \mathbb P^m,\qquad (x,y) \longmapsto (x^jy^{m-j})_{j=0}^m\] be the Veronese embeddings of degree $n$ and $m$, respectively. Note that the resultant $R$ of $\mathbb P^1$ with respect to $f$ and $g$ is the usual resultant, that is,
{\small\[
R(a_0, \ldots, a_n, b_0, \ldots, b_m) =
\det \begin{pmatrix}
a_0 & a_{1} & \cdots  &        & a_n   &       &        & \\
    & a_0     & a_{1} & \cdots &       & a_n   &        & \\
    &         & \ddots  &        &       &       & \ddots & \\
    &         &         &  a_0   &a_{1}&\cdots &        & a_n \\
    b_0 & b_{1} & \cdots  &        & b_m   &       &        & \\
    & b_0     & b_{1} & \cdots &       & b_m   &        & \\
    &         & \ddots  &        &       &       & \ddots & \\
    &         &         &  b_0   &b_{1}&\cdots &        & b_m \\
\end{pmatrix}
\begin{array}{l}
  \\ [-10ex]
   \rdelim\}{4.7}{5ex}[$m$ rows] \\ \\ \\[5ex]  \rdelim\}{4.7}{5ex}[$n$ rows] \\ \\
\end{array}
\]}
\end{exem}

\begin{rema}\label{Rem: resultant induction formula}
Let $R_{f_0,\ldots,f_d}^X$ be a resultant of $X$ with respect to $(f_i)_{i=0}^d$. If $K/k$ is an extension and if $s$ is an element of $E_d\otimes_kK$, defining a global section of $\mathcal O_{\mathbb P(E_d\otimes_kK)}(1)$, which intersects properly with all irreducible components $X\times_{\operatorname{Spec} k}\operatorname{Spec} K$, then, viewed as a multi-homogeneous polynomial on 
\[(E_0\otimes_kK)\times\cdots\times(E_d\otimes_kK)\]
by extension of scalars, the resultant $R^X_{f_0,\ldots,f_d}$ specified on the last coordinate at $s$, is a resultant of $\operatorname{div}(s)\cap X_{K}$ with respect to $(f_{i,K})_{i=0}^{d-1}$. This observation motivates the following explicit construction of the resultant polynomial by induction.
\end{rema}

\begin{defi}\label{Def:resultant precise}
Let $(s_0,\ldots,s_d)\in E_0\times\cdots\times E_d$. We assume that, for any irreducible component $Z$ of $X$, the divisors $\operatorname{div}(s_0),\ldots,\operatorname{div}(s_d)$ intersect properly on $Z$. We denote by $R_{f_0,\ldots,f_d}^{X,s_0,\ldots,s_d}$ the unique resultant of $X$ with respect to $f_0,\ldots,f_d$ such that
\[R_{f_0,\ldots,f_d}^{X,s_0,\ldots,s_d}(s_0,\ldots,s_d)=1.\]
\end{defi}

\begin{rema}\label{Rem: resultant extension}
Let $k'/k$ be an extension of fields. For any $i\in\{0,\ldots,d\}$, the morphism $f_i:X\rightarrow\mathbb P(E_i)$ induces by base change a closed embedding $f_i'$ from $X':=X\times_{\Spec k}\Spec k'$ to $\mathbb P(E_i')$, where $E_i':=E_i\otimes_kk'$.
Note that the incidence subscheme of \[X'\times_{k'}\mathbb P(E_{0}'^\vee)\times_{k'}\cdots\times_{k'}\mathbb P(E_{d}'^\vee)\]
identifies with $I_X\times_{\Spec k}\Spec k'$. Therefore, if $R_{f_0,\ldots,f_d}^X$ is a resultant of $X$ with respect to $(f_i)_{i=0}^d$, then 
\[R_{f_0,\ldots,f_d}^{X}\otimes 1\in (S^{\delta_0}(E_0^\vee)\otimes_k\cdots\otimes_kS^{\delta_d}(E_d^\vee))\otimes_kk'\cong S^{\delta_0}(E_{0}'^\vee)\otimes_{k'}\cdots\otimes_{k'}S^{\delta_d}(E_{d}'^\vee)\]
is a resultant of $X'$ with respect to $(f_{i}')_{i=0}^d$. Similarly, if $(s_0,\ldots,s_d)$ is an element of $ E_0\times\cdots\times E_d$ such that the divisors $\operatorname{div}(s_0),\ldots,\operatorname{div}(s_d)$ intersect properly on each irreducible component of $X$, then the following equality holds
\[R_{f_{0}',\ldots,f_{d}'}^{X',s_0',\ldots,s_d'}=R_{f_0,\ldots,f_d}^{X,s_0,\ldots,s_d}\otimes 1,\]
where for each $i\in\{0,\ldots,d\}$, $s_i'$ denotes the element $s_i\otimes 1$ in $E_i'=E_i\otimes_kk'$.
\end{rema}

\section{Projection to a projective space}
Let $k$ be an infinite field, $n$ be an integer such that $n\geqslant 1$, and $V$ be a vector space of dimension $n+1$ over $k$. Let $\mathbb P(V)$ be the projective space associated with the $k$-vector space $V$ and $\mathcal O_V(1)$ be the universal invertible sheaf on $\mathbb P(V)$. Recall that for any $k$-algebra $A$, any $k$-point of $\mathbb P(V)$ valued in $A$ corresponds to a quotient invertible $A$-module of $V\otimes_kA$. In particular, if $x$ is a scheme point of $\mathbb P(V)$ and $\kappa(x)$ is the residue field of $x$, then the scheme point $x$ corresponds to a non-zero $\kappa(x)$-linear map $p_x:V\otimes_k\kappa(x)\rightarrow\kappa(x)$, which is unique up to a unique homothety $\kappa(x)\rightarrow\kappa(x)$ by an element of $\kappa(x)^{\times}$.

\begin{defi}
We call \emph{rational linear subspace} of $\mathbb P(V)$ any Zariski closed subset of $\mathbb P(V)$ defined by the vanishing of all sections in a $k$-linear subspace of $V=H^0(\mathbb P(V),\mathcal O_{V}(1))$. If $Y$ is a rational linear subspace of $\mathbb P(V)$ which is of codimension $1$, we say that $Y$ is a \emph{rational hyperplane} in $\mathbb P(V)$.
\end{defi}

\begin{exem}
\begin{enumerate}[label=\rm(\arabic*)]
\item The scheme $\mathbb P(V)$ is a rational linear subspace of $\mathbb P(V)$. It is defined by the vanishing of the zero vector in $V$.
\item Let $x$ be a rational point of $\mathbb P(V)$, which corresponds to a non-zero $k$-linear map $\pi_x:V\rightarrow k$. Then $\{x\}$ is the vanishing locus of sections in $\operatorname{Ker}(\pi_x)$ and hence is a rational linear subspace of $\mathbb P(V)$.
\item The empty subset of $\mathbb P(V)$ is a rational linear subspace, which identifies with the vanishing locus of all sections in $V$. By convention, the dimension of the empty subset of $\mathbb P(V)$ is defined as $-1$.
\end{enumerate}
\end{exem}

\begin{rema}\label{Rem: linear proj}
If $Y$ is a rational linear subspace of $\mathbb P(V)$ which is the vanishing locus of a $k$-vector  subspace $W$ of $V$, then the $k$-scheme $Y$ is isomorphic to $\mathbb P(V/W)$. We call \emph{linear projection with center $Y$} the $k$-morphisme $\pi_Y:\mathbb P(V)\setminus Y \rightarrow \mathbb P(W)$ which sends, for any commutative $k$-algebra $A$, any quotient invertible $A$-module $p_L:V\otimes_kA\rightarrow L$ in $(\mathbb P(V)\setminus Y)(A)$ to the composition \[W\otimes_kA\hooklongrightarrow V\otimes_kA\stackrel{p_L}{\longrightarrow}L,\]
which is an element of $\mathbb P(W)(A)$.

We assume that $Y=\{y\}$ is the set of one rational point of $\mathbb P(V)$, which corresponds to a non-zero $k$-linear map $p_y:V\rightarrow k$ whose kernel is $W$. Let $z$ be a scheme point of $\mathbb P(V)$, $\kappa(z)$ be the residue field of $z$, and $p_z:V\otimes_k\kappa(z)\rightarrow\kappa(z)$ be the non-zero $\kappa(z)$-linear map corresponding to the scheme point $z$. Note that $\kappa(z)$ is generated by elements of the form $p_z(f\otimes 1)/p_z(g\otimes 1)$, where $f$ and $g$ are elements of $V$ such that $p_z(g\otimes 1)\neq 0$. Assume that $y$ does not belong the Zariski closure of $\{z\}$. Then there exists at least an element $s\in V\setminus W$ such that $p_z(s\otimes 1)= 0$. Let $z'$ be the image of $z$ by the linear projection $\pi_Y$. The residue field of $z'$ identifies with the sub-extension of $\kappa(x)/k$ generated by elements of the form $p_z(f'\otimes 1)/p_z(g'\otimes 1)$, where $f'$ and $g'$ are elements of $W$ such that $p_z(g'\otimes 1)\neq 0$. As $W$ is of codimension $1$ in $V$ and $s$ is an element of $V\setminus W$ such that $p_z(s\otimes 1)=0$, we obtain that, for any $f\in V$, there exists $f'\in W$ such that $p_z(f\otimes 1)=p_z(f'\otimes 1)$. Therefore we obtain that $\kappa(z)=\kappa(z')$. In particular, if $X$ is a closed subset of $\mathbb P(V)$ which does not contain $y$, then $\pi_Y(X)$ has the same dimension as  $X$.      
\end{rema}

\begin{prop}\label{prop:projection:to:projective:space}
Let $d\in\{0,\ldots,n\}$.
Let $X$ be a Zariski closed set of $\mathbb P(V)$ such that $\dim(X) \leq d$.
Then we have the following:
\begin{enumerate}[label=\rm(\arabic*)]
\item
There is a rational linear subspace $M$ of $\mathbb P(V)$ such that $\dim(M) = n-1-d$ and $X \cap M = \emptyset$.

\item
Let $T$ be a rational linear subspace of $\mathbb P(V)$ such that $\dim(T) > n-d-1$, and that $X$ and $T$ meets properly.
Then
there is a rational linear subspace $M$ of $\mathbb P(V)$ such that $M \subseteq T$, $\dim(M) = n-1-d$ and $X \cap M = \emptyset$.

\item We assume that $X$ is irreducible and $\dim(X) = d$.
Let $M$ be a rational linear subspace of $\mathbb P(V)$ such that $\dim(M) = n-1-d$ and $M \cap X = \emptyset$, which is the vanishing locus of a vector space $W$ of $V$.
Let $\pi_M : \mathbb P(V) \setminus M \to \mathbb P(W)$ be the projection with the center $M$. Then $\pi := \rest{\pi_M}{X} : X \to \mathbb P^d_k$
is finite and surjective and $\pi^*(\mathcal O_{\mathbb P^d_k}(1)) = \rest{\mathcal O_{\mathbb P^n_k}(1)}{X}$.
\end{enumerate}
\end{prop}

\begin{proof}
(1) We prove the assertion by induction on $n-d$. 
If $n=d$, then the assertion is obvious by choosing $M$ as the empty set, so that we assume that $n > d$.
Since $X \not= \mathbb P(V)$ and $k$ is an infinite field, there is a rational point $x \in \mathbb P(V)$ which does not belong to $X$. Let $W$ be the set of sections $s\in V=H^0(\mathbb P(V),\mathcal O_V(1))$ which vanish at $x$. This is a vector subspace of $V$.
Let $\pi : \mathbb P(V) \setminus \{ x \}  \to \mathbb P(W)$
be the projection with center $\{x$\}. Since $x \not\in X$, by Remark \ref{Rem: linear proj} we obtain that $X$ and $X'$ have the same dimension. 
In particular, $\dim(X') \leq d$. As $(n-1) -d < n-d$,  by the hypothesis of induction, there is a linear
subspace $M'$ in $\mathbb P(W)$ such that $\dim (M') = n-2-d$ and $X' \cap M' = \emptyset$
Thus if we set $M = \pi^{-1}(M') \cup \{ x \}$, then one has the desired subspace.

(2) Assume that $T$ is defined by the vanishing of sections in a $k$-vector subspace $W$ of $V$. If we set $X' = X \cap T$ and $t = \dim T$, then $\dim X' \leq d - (n-t)$ and $T \simeq \mathbb P(V/W)$.
As
\[
t - (d-(n-t)) = n -d \geq 0,
\]
by (1), there is linear subspace $M$ in $T$ such that $\dim M = t - 1 - (d-(n-t))$ and $M \cap X' = \emptyset$.
Thus one has (2).

(3) Let $T$ be a linear subspace of $\mathbb P(V)$ such that $M \subseteq T$ and $\dim(T) = n -d$.
It is sufficient to show that $\dim (T \cap X) = 0$. Note that  
$M$ is a rational hyperplane in $T$,
so that if $\dim (T \cap X) \geq 1$, then $M \cap X \not= \emptyset$. Therefore $\dim (T \cap X) = 0$. 
\end{proof}


\chapter{Adelic curves and their constructions}

\section{Adelic structures}

In this section, we recall the notion of adelic curves. Let $K$ be a field.
An \emph{adelic structure} of $K$ consists of data
$((\Omega, \mathcal A, \nu), \phi)$ satisfying the following properties:
\begin{enumerate}
\renewcommand{\labelenumi}{\textup{(\arabic{enumi})}}
\item $(\Omega, \mathcal A,\nu)$ is a measure space, that is,
$\mathcal A$ is a $\sigma$-algebra of $\Omega$ and $\nu$ is a measure on $(\Omega, \mathcal A)$.

\item 
The last $\phi$ is a map
from $\Omega$ to $M_K$, where $M_K$ is the set of all absolute values of $K$. For any $\omega\in\Omega$, we denote the absolute value $\phi(\omega)$ by $|\ndot|_{\omega}$.

\item For any $\omega\in\Omega$ and any $a \in K^{\times}$, the function $(\omega \in \Omega) \mapsto \ln |a|_{\omega}$ is $\nu$-integrable.
\end{enumerate}
The field $K$ equipped with an adelic structure is called an \emph{adelic curve}. Moreover, the adelic structure $((\Omega,\mathcal A,\nu),\phi)$ is said to be \emph{proper} if
\begin{equation}\label{eqn:product:formula}
\int_{\Omega} \ln |a|_{\omega}\, \nu(\mathrm{d}\omega) = 0
\end{equation}
holds for all $a \in K^{\times}$.  If the adelic structure $((\Omega,\mathcal A,\nu),\phi)$ is proper, we also say that the adelic curve $(K,(\Omega,\mathcal A,\nu),\phi)$ is \emph{proper}. The equation \eqref{eqn:product:formula} is
called  \emph{product formula}. For details, see \cite[Chapter 3]{CMArakelovAdelic}.
We denote the set of all $\omega \in \Omega$ with $|\ndot|_{\omega}$ Archimedean
(resp. non-Archimedean) 
by $\Omega_{\infty}$ (resp. $\Omega_{\fin}$). The restriction of $\mathcal A$ to $\Omega_\infty$ (resp. $\Omega_{\operatorname{fin}}$) is denoted by $\mathcal A_\infty$ (resp. $\mathcal A_{\operatorname{fin}}$).
Note that $\Omega_{\infty}$ and $\Omega_{\fin}$ belong to $\mathcal{A}$ (see \cite[Proposition 3.1.1]{CMArakelovAdelic}).
For each $\omega \in \Omega_{\infty}$, there exist
an embedding $\iota_{\omega} : K \to \CC$ and $\kappa_{\omega} \in (0, 1]$
such that $|a|_{\omega} = |\iota_{\omega}(a)|^{\kappa_{\omega}}$ for all $a \in K$,
where $|\ndot|$ is the usual absolute value of $\CC$.
Note that the invariant $\kappa_{\omega}$ does not depend on the choice of the embedding
$\iota_{\omega} : K \to \CC$.
From now on, we always assume that $\kappa_{\omega} = 1$ for all $\omega \in \Omega_{\infty}$.

For $(a_1, \ldots, a_n) \in K^n \setminus \{ (0, \ldots, 0) \}$,
the height $h_S(a_1, \ldots, a_n)$ of $(a_1, \ldots, a_n)$ with respect to the adelic curve
$S = (K,(\Omega, \mathcal A, \nu), \phi)$ is defined to be
\begin{equation}\label{eqn:height:function:wrt:S}
h_S(a_1, \ldots, a_n) := \int_{\Omega} \ln( \max \{ |a_1|_{\omega}, \ldots, |a_n|_{\omega} \})
\nu(\mathrm{d}\omega).
\end{equation}
Note that if $S$ is proper, then $h_S(a) = 0$ for all $a \in K^{\times}$.

\begin{rema}
Many classic constructions in algebraic geometry and arithmetic geometry, such as algebraic curves, rings of algebraic integers, polarized projective varieties and arithmetic varieties, can be interpreted as adelic curves. For example, on the filed $\mathbb Q$ of rational numbers there is an adelic structure consisting of all places of $\mathbb Q$ (namely the set $\Omega_{\mathbb Q}$ of all prime numbers and $\infty$) equipped with the discrete $\sigma$-algebra and the measure $\nu$ such that $\nu(\{\omega\})=1$ for any $\omega\in\Omega_{\mathbb Q}$, where $|\ndot|_{\infty}$ is the usual absolute value on $\mathbb Q$ and $|\ndot|_p$ is the $p$-adic absolute value for any prime number $p$. The product formula for this adelic curve is just the logarithmic version of the usual product formula for rational numbers
\[\forall\,a\in\mathbb Q^{\times},\quad |a|_\infty\cdot\prod_{p}|a|_p=1.\]
We call this adelic structure \emph{the standard adelic structure on $\mathbb Q$}. We refer the readers to \cite[\S 3.2]{CMArakelovAdelic} for more examples.
\end{rema}

\begin{defi}\label{Def: covering of adelic curves}
Let $S=(K,(\Omega,\mathcal A,\nu),\phi)$ and $S'=(K',(\Omega',\mathcal A',\nu'),\phi')$ be two adelic curves. We call \emph{morphism} from $S'$ to $S$ any triplet $\alpha=(\alpha^{\#},\alpha_{\#}, I_\alpha)$, where
\begin{enumerate}[label=\rm(\arabic*)]
\item $\alpha^{\#}:K\rightarrow K'$ is a field homomorphism,
\item $\alpha_{\#}:(\Omega',\mathcal A')\rightarrow(\Omega,\mathcal A)$ is a measurable map, such that, for any $\omega'\in\Omega'$, \[|\ndot|_{\omega'}\circ\alpha^{\#}=|\ndot|_{\alpha_{\#}(\omega')},\]
and that the direct image of $\nu'$ by $\alpha_{\#}$ coincides with $\nu$, namely for any $\mathcal A$-measurable  function $f:\Omega\rightarrow\mathbb R$ which is either non-negative or integrable, one has
\[\int_{\Omega'}f(\alpha_{\#}(\omega'))\,\nu'(\mathrm{d}\omega')=\int_\Omega f(\omega)\,\nu(\mathrm{d}\omega),\]
\item \[I_\alpha:\mathscr L^1(\Omega',\mathcal A',\nu')\longrightarrow \mathscr L^1(\Omega,\mathcal A,\nu)\] is a linear map sending positive integrable functions on $(\Omega',\mathcal A',\nu')$ to positive integrable functions on $(\Omega,\mathcal A,\nu)$ such that, for any $f\in\mathscr L^1(\Omega',\mathcal A',\nu')$,
\[\int_{\Omega}I_\alpha(f)(\omega)\,\nu(\mathrm{d}\omega)=\int_{\Omega'}f(\omega')\,\nu'(\mathrm{d}\omega').\]
\end{enumerate}
If in addition for any $g\in\mathscr L^1(\Omega,\mathcal A,\nu)$, one has \[g\circ\alpha_{\#}\in\mathscr L^1(\Omega',\mathcal A',\nu')\quad\text{ and }\quad I_\alpha(g\circ\alpha_{\#})=g,\] we say that $\alpha$ is a \emph{covering} of adelic curves.
\end{defi}

\section{Algebraic coverings of adelic curves}\label{Sec: algebraic coverings}

Adelic curves are very flexible constructions. On a field there exist many adelic structures. It is also possible to construct new adelic structures from given ones.  Let $S=(K,(\Omega,\mathcal A,\nu),\phi)$ be an adelic curve. In \cite[\S 3.2]{CMArakelovAdelic} it has been explained how to construct, for any algebraic extension $L/K$, a natural adelic curve \[S\otimes_KL=(L,(\Omega_{L},\mathcal A_{L},\nu_{L}),\phi_{L})\] on $L$ such that $\Omega_{L}=\Omega\times_{M_K,\phi}M_{L}$. The projection map $\pi_{L/K}:\Omega_{L}\rightarrow\Omega$ satisfies the relation 
\[\nu=(\pi_{L/K})_*(\nu_L).\]
Moreover, for any $\omega\in\Omega$, the fibre $\pi_{L/K}^{-1}(\{\omega\})$ is equipped with a natural $\sigma$-algebra and a probability measure $\nu_{L,\omega}$, such that, for any positive $\mathcal A_{L}$-measurable function $f$ on $\Omega_{L}$, one has
\[\int_{\Omega_{L}}g(x)\,\nu_{L}(\mathrm{d}x)=\int_{\Omega}\nu(\mathrm{d}\omega)\int_{\pi_{L/K}^{-1}(\omega)}g(x)\,\nu_{L,\omega}(\mathrm{d}x).\]
In other words, the family of measures $(\nu_{L,\omega})_{\omega\in\Omega}$ form an disintegration of $\nu_{L}$ over $\nu$.
If the adelic curve $S$ is proper, then also is $S\otimes_KL$, see \cite[Proposition 3.4.10]{CMArakelovAdelic}. If we denote by $i_{K,L}:K\rightarrow L$ the inclusion map, and \[I_{L/K}:\mathscr L^1(\Omega_L,\mathcal A_L,\nu_L)\longrightarrow\mathscr L^1(\Omega,\mathcal A,\nu)\]
the linear map of fiber integrals, which sends $g\in L^1(\Omega_L,\mathcal A_L,\nu_L)$ to the function
\[(\omega\in\Omega)\longmapsto\int_{\pi_{L/K}^{-1}(\omega)}g(x)\,\nu_{L,\omega}(\mathrm{d}x),\]
then the triplet 
 $(i_{K,L},\pi_{L/K},I_{L/K})$ forms a covering of adelic curves in the sense of Definition \ref{Def: covering of adelic curves}.

\begin{lemm}\label{lemma:discrete:algebraic}
Let $K'$ be an algebraic extension of $K$ and
$S \otimes K' := (K', (\Omega', \mathcal{A}', \nu'), \phi')$.
Suppose that $K$ and $\Omega_{\operatorname{fin}}$ are countable sets. If $(\Omega_{\fin},\mathcal{A}_{\fin})$ is discrete, also is $(\Omega'_{\fin},\mathcal{A}'_{\fin})$.
\end{lemm}

\begin{proof}Since $K'/K$ is an algebraic extension, and $K$ and $\Omega_{\operatorname{fin}}$ are countable sets, we obtain that the sets $K'$ and $\Omega'_{\operatorname{fin}}$ are countable, so that it is sufficient to see that $\{ \omega' \}\in\mathcal{A}_{\fin}'$
for all $\omega' \in \Omega'_{\fin}$.

First we consider the case
where $K'$ is finite over $K$. 
Let $\omega' \in \Omega'$ and $\omega = \pi(\omega')$,
where $\pi : \Omega' \to \Omega$ is the canonical map.
Then as $\{ \omega \}$ is $\mathcal{A}_{\fin}$-measurable and
$\pi$ is measurable, $\pi^{-1}(\{\omega\})\in\mathcal{A}'_{\fin}$.
If $|\ndot|_{\omega}$ is trivial, then $\pi^{-1}(\{\omega\}) = \{ \omega' \}$,
so that the assertion is obvious. Next we assume that $|\ndot|_{\omega}$
is non-trivial. Let us see that, for any $(x, x') \in \pi^{-1}(\{\omega\})^2$
with $x \not= x'$, $|\ndot|_x$ is not equivalent to $|\ndot|_{x'}$.
Otherwise, there is $\kappa \in \RR_{>0}$ such that $|\ndot|_{x'} =
|\ndot|_x^{\kappa}$. As $|\ndot|_{\omega}$ is non-trivial, 
there is $a \in K$ such that $|a|_{\omega} < 1$. Then
\[
|a|_{\omega} = |a|_{x'} = |a|_{x}^{\kappa} = |a|_{\omega}^{\kappa},
\]
and hence $\kappa = 1$, which is a contradiction. Therefore,
there is $a' \in K'$ such that $|a'|_{\omega'} < 1$ and
$|a'|_{x} > 1$ for all $x \in \pi^{-1}(\{\omega\}) \setminus \{ \omega' \}$
(cf. \cite[the proof of Theorem~3.4]{Neukirch}).
Note that \[\Delta := \{ \chi \in \Omega' \,:\, |a'|_\chi < 1 \}\] is $\mathcal{A}'_{\fin}$-measurable,
so that $\{ \omega' \} = \pi^{-1}(\{\omega\}) \cap \Delta$ is $\mathcal{A}'_{\fin}$-measurable.

In general, for $a \in K'$, let \[(K(a),(\Omega_{K(a)},\mathcal A_{K(a)},\nu_{K(a)}),\phi_{K(a)})=S\otimes K(a)\] and let
$\pi_{K'/K(a)} : \Omega' \to \Omega_{K(a)}$ be the canonical map. 
By the previous case, $\{\pi_{K'/K(a)}(\omega')\}\in\mathcal A_{K(a)}$, so that
$\pi_{K'/K(a)}^{-1}(\{\pi_{K'/K(a)}(\omega')\})\in\mathcal A'$. 
Therefore, as $K'$ is countable,
\[
\bigcap_{a \in K'} \pi_{K'/K(a)}^{-1}(\pi_{K'/K(a)}(\omega')).
\]
belongs to $\mathcal A'$. Thus it suffices to prove
\begin{equation}\label{eqn:prop:discrete:algebraic:01}
\{\omega'\} = \bigcap_{a \in K'} \pi_{K'/K(a)}^{-1}(\pi_{K'/K(a)}(\omega')).
\end{equation}
Indeed, if $x \in \bigcap_{a \in K'} \pi_{K'/K(a)}^{-1}(\{\pi_{K'/K(a)}(\omega')\})$,
then, for any $a \in K'$, $\pi_{K'/K(a)}(x) = \pi_{K'/K(a)}(\omega')$, 
so that $|a|_x = |a|_{\omega'}$,
which means that $x = \omega'$.
\end{proof}

\section{Transcendental fibrations of adelic curves}
\label{Sec: Transcendental fibration}

The purpose of this section is to discuss the  extension of an adelic structure to a transcendental extension of the field. We fix an adelic curve $S=(K,(\Omega,\mathcal A,\nu),\phi)$. For any $\omega\in\Omega$, let $K_\omega$ be the completion of $K$ with respect to the absolute value $\phi(\omega)$. Let $B$ be a $K$-algebra. Note that $B$ is not necessarily of finite type over $K$.
We assume that $B$ is a unique factorization domain and the set $B^{\times}$ of units in $B$ coincides with $K^{\times}$. We say that two irreducible elements of $B$ are \emph{equivalent} if they differ by a unit as a factor. This defines an equivalence relation on the set of all irreducible elements of $B$. We pick a representative in each of the equivalence classes to form a subset $\mathscr P_B$ of $B$ consisting of non-equivalent irreducible elements. Let $L$ be the field of fractions of $B$.  
Recall that any non-zero element $g\in L$ can be written in a unique way as
\[c(g)\prod_{F\in\mathscr P_B}F^{\operatorname{ord}_F(g)},\]
where $c(g)$ is an element of $K^{\times}=B^{\times}$, and for each $F\in\mathscr P_B$, $\operatorname{ord}_F(g)$ is an integer. Note that $\operatorname{ord}_F(\ndot)$ is a discrete valuation on the field $L$, and $\operatorname{ord}_F(a)=0$ for any $a\in K^{\times}=B^{\times}$.

\begin{defi} \label{Def: admissible family}
For any $\omega\in\Omega$, let $S_{L,\omega}=(L,(\Omega_{L,\omega},\mathcal A_{L,\omega},\nu_{L,\omega}),\phi_{L,\omega})$ be an adelic curve such that $\nu_{L,\omega}$ is a probability measure. We say that the family $(S_{L,\omega})_{\omega\in\Omega}$ is an \emph{admissible fibration with respect to $(B,\mathscr P_B)$} over the adelic curve $S$ if the following conditions are satisfied:
\begin{enumerate}[label=
\rm(\alph*)]
\item for any $\omega\in\Omega$ and any $x\in \Omega_{L,\omega}$, the absolute value $\phi_{L,\omega}(x)$ on $L$ is an extension of $\phi(\omega)$ on $K$,
\item\label{Cond: mesurabilit} for any element $g\in B\setminus\{0\}$, any finite family $(F_j)_{j=1}^n$ of elements of $\mathscr P_B$ containing  $\{F\in\mathscr P_B\,|\,\operatorname{ord}_{F}(g)\neq 0\}$ and any $(C_j)_{j=1}^n\in\mathbb R_{\geqslant 0}^n$, the function
\[(\omega\in\Omega)\longmapsto\int_{\Omega_{L,\omega}}|g|_x\indic_{|F_1|_x\leqslant C_1,\ldots,|F_n|_x\leqslant C_n}\,\nu_{L,\omega}(\mathrm{d}x)\]
is $\mathcal A$-measurable,
\item\label{Cond: integrability} for any $\omega\in\Omega$ and any element $F$ of $\mathscr P_B$, the function
\[(\omega\in\Omega)\longmapsto\int_{\Omega_{L,\omega}}\ln|F|_x\, \nu_{L,\omega}(\mathrm{d}x)\]
is integrable with respect to $\nu$. 
\end{enumerate} 

Let $(S_{L,\omega})_{\omega\in\Omega}$ be an admissible fibration over the adelic curve $S$. We define $\Omega_L$ as the disjoint union of $(\Omega_{L,\omega})_{\omega\in\Omega}$ and let $\phi_L$ be the map from $\Omega_L$ to the set of all absolute values on $L$, whose restriction on each $\Omega_{L,\omega}$ is equal to $\phi_{L,\omega}$. Let $\pi_{L/K}:\Omega_L\rightarrow\Omega$ be the projection map, sending the elements of $\Omega_{L,\omega}$ to $\omega$. We equip $\Omega_L$ with the $\sigma$-algebra $\mathcal A_L$ generated by the projection map $\pi_{L/K}$ and all functions of the form $(x\in\Omega_{L})\mapsto |g|_x$, where $g$ runs over the set $L$.
\end{defi}

\begin{prop}\label{Pro: disintegration}
Let $f$ be a non-negative $\mathcal A_L$-measurable function on $\Omega_L$.  For any $\omega\in\Omega$, the function $f$ is $\mathcal A_{L,\omega}$-measurable on $\Omega_{L,\omega}$. Moreover, the function
\[(\omega\in\Omega)\longmapsto\int_{\Omega_{L,\omega}}f(x)\,\nu_{L,\omega}(\mathrm{d}x)\in[0,+\infty]\]
is $\mathcal A$-measurable.
\end{prop}
\begin{proof}
Let $\mathcal H$ be the set of all bounded non-negative $\mathcal A_{L}$-measurable functions $g$ on $\Omega_{L}$ which is $\mathcal A_{L,\omega}$-measurable on $\Omega_{L,\omega}$ for any $\omega\in\Omega$ and such that the function
\[(\omega\in\Omega)\longmapsto\int_{\Omega_{L,\omega}}f(x)\,\nu_{L,\omega}(\mathrm{d}x)\]
is $\mathcal A$-measurable. Note that, for any non-negative bounded $\mathcal A$-measurable function $\varphi$ on $\Omega$, one has $\varphi\circ\pi\in\mathcal H$ since it is constant on each fiber $\Omega_{L,\omega}$ and
\[\int_{\Omega_{L,\omega}}\varphi(\pi(x))\,\nu_{L,\omega}(\mathrm{d}x)=\int_{\Omega_{L,\omega}}\varphi(\omega)\,\nu_{L,\omega}(\mathrm{d}x)=\varphi(\omega).\] In particular, all non-negative constant functions belong to $\mathcal H$. Clearly, for any $(g_1,g_2)\in\mathcal H\times\mathcal H$ and any $(a_1,a_2)\in\mathbb R_{\geqslant 0}\times\mathbb R_{\geqslant 0}$, one has $a_1g_1+a_2g_2\in\mathcal H$. For any increasing sequence of functions $(g_n)_{n\in\mathbb N}$ in $\mathcal H$, the pointwise limit of $(g_n)_{n\in\mathbb N}$ belongs to $\mathcal H$.  Moreover, for functions $g_1$ and $g_2$ in $\mathcal H$ such that $g_2\geqslant g_1$, then one $g_2-g_1\in\mathcal H$.

Let $\mathcal S$ be the set of functions of the form \[(x\in\Omega_{L})\longmapsto|g|_x\indic_{|F_1|_x\leqslant C_1,\ldots,|F_n|_x\leqslant C_n}\varphi(\pi(x)),\] where $g$ is an element of $B\setminus\{0\}$, $(F_j)_{j=1}^n$ is a finite family of elements of $\mathscr P_B$ containing  $\{F\in\mathscr P_B\,|\,\operatorname{ord}_{F}(g)\neq 0\}$, $(C_j)_{j=1}^n$ is a family of positive constant and $\varphi$ is a non-negative and bounded $\mathcal A$-measurable function on $\Omega$. Clearly the set $\mathcal S$ is stable by multiplication. Note that the function sending $\omega\in\Omega$ to
\[\begin{split}&\quad\;\int_{\Omega_{L,\omega}}|g|_x\indic_{|F_1|_x\leqslant C_1,\ldots,|F_n|_x\leqslant C_n}\varphi(\pi(x))\,\nu_{L,\omega}(\mathrm{d}x)\\
&=\varphi(\omega)\int_{\Omega_{L,\omega}}|g|_x\indic_{|F_1|_x\leqslant C_1,\ldots,|F_n|_x\leqslant C_n}\,\nu_{L,\omega}(\mathrm{d}x)
\end{split}\]
takes real values and is $\mathcal A$-measurable by the condition \ref{Cond: mesurabilit} above. Therefore, $\mathcal S$ is a subset of $\mathcal H$. Since the $\sigma$-algebra $\mathcal A_L$ is generated by $\mathcal S$, by monotone class theorem (see \cite[\S2.2]{Yan}, see also \cite[\S A.1]{CMArakelovAdelic}), $\mathcal H$ contains all bounded non-negative $\mathcal A_{L}$-measurable functions. Finally, since any non-negative $\mathcal A_{L}$-measurable function $f$ can be written as the limit of an increasing sequence of bounded non-negative $\mathcal A_{L}$-measurable functions, the assertion of the proposition is true.   
\end{proof}

\begin{defi}\label{Def: adelic structure coming from a fibration}
Let $(S_{L,\omega})_{\omega\in\Omega}$ be an admissible fibration over $S$ (see Definition \ref{Def: admissible family}), where $S_{L,\omega}=(L,(\Omega_{L,\omega},\mathcal A_{L,\omega},\nu_{L,\omega}),\phi_{L,\omega})$. By Proposition \ref{Pro: disintegration}, there is a measure $\nu_{L}$ on the measurable space $(\Omega_{L},\mathcal A_{L})$ such that, for any non-negative $\mathcal A_{L}$-measurable function $f$ on $\Omega_{L}$, one has
\[\int_{\Omega_{L}}f(x)\,\nu_{L}(\mathrm{d}x)=\int_{\Omega}\nu(\mathrm{d}\omega)\int_{\Omega_{L,\omega}}f(x)\,\nu_{L,\omega}(\mathrm{d}x).\]
Therefore $S_L := (L,(\Omega_{L},\mathcal A_{L},\nu_{L}),\phi_{L})$ is an adelic curve, called the \emph{adelic curve associated with the admissible fibration $(S_{L,\omega})_{\omega\in\Omega}$.} Since $\nu_{L,\omega}$ are probability measures, if we denote by $i_{K,L}:K\rightarrow L$ the inclusion map, by $\pi_{L/K}:\Omega_L\rightarrow \Omega$ the map sending the elements of $\Omega_{L,\omega}$ to $\omega$, and by \[I_{L/K}:\mathscr L^1(\Omega_L,\mathcal A_L,\nu_L)\longrightarrow\mathscr L^1(\Omega,\mathcal A,\nu)\]
the linear map of fiber integrals, then the triplet $(i_{K,L},\pi_{L/K},I_{L/K})$ forms a covering of adelic curves in the sense of Definition \ref{Def: covering of adelic curves}.
\end{defi}

\section{Intrinsic compactification of admissible fibrations}

Let $S=(K,(\Omega,\mathcal A,\nu),\phi)$ be a proper adelic curve, $B$ be a $K$-algebra which is a unique factorization domain, and $\mathscr P_B$ be a representative family of irreducible elements as in the previous section. Let $L$ be the field of fractions of $B$ and \[\big(S_{L,\omega}=(L,(\Omega_{L,\omega},\mathcal A_{L,\omega},\nu_{L,\omega}),\phi_{L,\omega})\big)_{\omega\in\Omega}\] be an admissible fibration with respect to $(B,\mathscr P_B)$. In the previous section, we have constructed an adelic curve $S_L := (L,(\Omega_{L},\mathcal A_{L},\nu_{L}),\phi_{L})$ which fibers over $S$ and such that the measure $\nu_{L}$ disintegrates over $\nu$ by the family of measures $(\nu_{L,\omega})_{\omega\in\Omega}$ on the fibers. This construction looks similar to algebraic coverings of adelic curves. However, even in the case where the adelic structure $((\Omega,\mathcal A,\nu),\phi)$ is proper, the adelic structure $((\Omega_L,\mathcal A,\nu_L),\phi_L)$ is not necessarily proper.  In this section, we show that, under a mild condition on the admissible fibration $(S_\omega)_{\omega\in\Omega}$ over $S$, we can naturally ``compactify'' the adelic structure $((\Omega_{L},\mathcal A_{L},\nu_{L}),\phi_{L})$. For any element $F\in\mathscr P_B$, we denote by $|\ndot|_F$ the absolute value on $K(T)$ such that 
\[\forall\,g\in L^{\times},\quad |g|_F:=\operatorname{e}^{-\operatorname{ord}_F(g)}.\] 
Thus we obtain a map $\phi_{L}'$ from $\mathscr P_B$ to $M_{L}$ sending $F$ to $|\ndot|_F$.  Let $(\Omega_{L}^*,\mathcal A_{L}^*)$ be the disjoint union of the measurable spaces $(\Omega_{L},\mathcal A_{L})$ and $\mathscr P_B$ equipped with the discrete $\sigma$-algebra. Let $\phi_L^*:\Omega_L^*\rightarrow M_L$ be the map extending $\phi_L$ on $\Omega_L$ and $\phi_L'$ on $\mathscr P_B$.

\begin{prop}\label{Pro: compactification}
Let $(S_\omega)_{\omega\in\Omega}$ be an admissible fibration over $S$. We assume that, for any element $F\in\mathscr P_B$,  
\begin{equation}\label{Equ: positivity of height}h_{S_L}(F):=\int_{\Omega}\nu(\mathrm{d}\omega)\int_{\Omega_{L,\omega}}\ln|F|_x\,\nu_{L,\omega}(\mathrm{d}x)\geqslant 0.\end{equation} 
Let $\nu_{L}^*$ be the measure on $(\Omega_{L}^*,\mathcal A_{L}^*)$ which coincides with $\nu_{L}$ on $(\Omega_{L},\mathcal A_{L})$ and such that
\begin{gather*}\forall\,F\in\mathscr P_B,\quad\nu_{L}^*(\{F\})=h_{S_L}(F).\end{gather*}
Then $S^*_L:= (L,(\Omega_{L}^*,\mathcal A_{L}^*,\nu_{L}^*),\phi_L^*)$ is a proper adelic curve.
\end{prop}
\begin{proof}
For any $g\in L^{\times}$, one has
\begin{equation}\label{Equ: computation ln g}\begin{split}&\quad\;\int_{\Omega_{L}}\ln|g|_x\,\nu_{L}(\mathrm{d}x)=\int_{\Omega}\nu(\mathrm{d}\omega)\int_{\Omega_{L,\omega}}\ln|g|_x\,\nu_{L,\omega}(\mathrm{d}x)\\
&=\sum_{F\in\mathscr P_B}\operatorname{ord}_F(g)\int_{\Omega}\nu(\mathrm{d}\omega)\int_{\Omega_{L,\omega}}\ln|F|_x\,\nu_{L,\omega}(\mathrm{d}x)=\sum_{F\in\mathscr P_B}\operatorname{ord}_F(g)h_{S_L}(F).
\end{split}\end{equation} 
Thus
\[\int_{\Omega_{L}^*}\ln|g|_x\,\nu_{L}^*(\mathrm{d}x)=\int_{\Omega_L}\ln|g|_x\,\nu_L(\mathrm{d}x)+\sum_{F\in\mathscr P_B}h_{S_L}(F)\ln|g|_F=0.
\]
\end{proof}

\begin{defi}\label{def:canonical:compactification}
Under the assumption \eqref{Equ: positivity of height}, the adelic curve $S^*_L$ is called the \emph{canonical
compactification} of $S_L$. 
\end{defi}

\begin{rema}
Let $\mathcal A_B$ be the discrete $\sigma$-algebra on $\mathscr P_B$, $\nu_B$ be the  measure on $(\mathscr P_B,\mathcal A_B)$ such that 
\[\nu_B(\{F\})=h_{S_L}(F)\]
for any $F\in\mathscr P_B$, and $\phi_B:\mathscr P_B\rightarrow M_K$ be the map sending any element of $\mathscr P_B$ to the trivial absolute value on $K$. Then $S_B:=(K,(\mathscr P_B,\mathcal A_B,\nu_B),\phi_B)$ forms an adelic curve having $K$ as the underlying field. Let $S^*$ be the amalgamation of $S$ and $S_B$. Then, the inclusion map $K\rightarrow L$, the projection \[\pi_{L/K}\amalg \operatorname{Id}_{\mathscr P_B}:\Omega_L^*=\Omega_L\amalg\mathscr P_B\longrightarrow\Omega^*=\Omega\amalg\mathscr P_B\]
and the integral along fibers form a covering of adelic curves. 
\end{rema}

\section{Non-intrinsic compactification of admissible fibrations}

We keep the notation of the previous section. In this section, we assume that the family of absolute values $(|\ndot|_{F})_{F\in\mathscr P_B}$ can be included in a proper adelic structure. We will show that a weaker positivity condition than \eqref{Equ: positivity of height} would be enough to ensure the existence of (non-intrinsic) compactifications  of the adelic structure $((\Omega_{L},\mathcal A_{L},\nu_{L}),\phi_{L})$. In the rest of the subsection, we assume that there exists a \emph{proper} adelic structure $((\Omega_{L}',\mathcal A_{L}',\nu_L'),\phi_L')$ on $L$ which satisfies the following conditions:
\begin{enumerate}[label=\rm(\arabic*)]
\item $\Omega_L'$ contains $\mathscr P_B$ as a discrete measurable sub-space and $\nu_L'(\{F\})>0$ for any $F\in\mathscr P_B$,
\item for any $F\in\mathscr P_B$, one has $\phi_L'(F)=|\ndot|_F$.
\end{enumerate}
Note that the existence of such an adelic structure is is true when $K$ is of characteristic $0$ and $\Spec B$ is a smooth $K$-scheme of finite type. In this case there exists a projective $K$-scheme $X$ and an open immersion from $B$ into $X$.   Then one can construct an adelic consisting of prime divisors of $X$, by choosing a polarization on $X$. We refer the readers to \cite[\S3.2.4]{CMArakelovAdelic} for more details.

\begin{prop}\label{Pro: compactification 2}
Let $(S_\omega)_{\omega\in\Omega}$ be an admissible fibration over $S$. For any element $F\in\mathscr P_B$, let 
\begin{equation}\label{Equ: positivity of height :02}h_{S_L}(F):=\int_{\Omega}\nu(\mathrm{d}\omega)\int_{\Omega_{L,\omega}}\ln|F|_x\,\nu_{L,\omega}(\mathrm{d}x).\end{equation}
Let $\delta$ be a positive constant. We assume that
\[\forall\,F\in\mathscr P_B,\quad h_{S_L}(F)+\delta\nu_L'(\{F\})\geqslant 0.\]
Let $(\Omega_L'',\mathcal A_L'')$ be the disjoint union of $(\Omega_L,\mathcal A_L)$ and $(\Omega_L',\mathcal A_L')$, $\phi_L'':\Omega_L''\rightarrow M_L$ be the map extending $\phi_L$ and $\phi_L'$, and  $\nu_L^\delta$ be the measure on 
$(\Omega_L'',\mathcal A_L'')$ which coincides with $\nu_L$ on $(\Omega_L,\mathcal A_L)$ and coincides with
\[\delta\nu_L'+\sum_{F\in\mathscr P_B}h_{S_L}(F)\operatorname{Dirac}_F\]
on $(\Omega_L',\nu_L')$, where $\operatorname{Dirac}_F$ denotes the Dirac measure at $F$. Then $((\Omega_{L}'',\mathcal A_{L}'',\nu_{L}^\delta),\phi_L'')$ is a proper adelic structure on $L$.
\end{prop}
\begin{proof}
For any $g\in L^{\times}$, one has
\[\int_{\Omega_L^*}\ln|g|_x\,\nu_L^\delta(\mathrm{d}x)=\int_{\Omega_L}\ln|g|_x\,\nu_L(\mathrm{d}x)+\delta\int_{\Omega_L'}\ln|g|_x\,\nu_L'(\mathrm{d}x)+\sum_{F\in\mathscr P_B}h_{S_L}(F)\ln|g|_F.
\]
By \eqref{Equ: computation ln g}, one has
\[\int_{\Omega_L}\ln|g|_x\,\nu_L(\mathrm{d}x)+\sum_{F\in\mathscr P_B}h_{S_L}(F)\ln|g|_F=0.\]
Moreover, since $((\Omega_L',\mathcal A_L',\nu_L'),\phi_L')$ is a proper adelic structure, one has
\[\int_{\Omega_L'}\ln|g|_x\,\nu_L'(\mathrm{d}x)=0.\]
Therefore we obtain 
\[\int_{\Omega_L''}\ln|g|_x\,\nu_L^\delta(\mathrm{d}x)=0.\]
\end{proof}

\section{Purely transcendental fibration of adelic curves}
\label{subsec:adelic:structure:of:purely:transcendental:fibration}
In this section, we apply the results obtained in previous sections to the study of adelic structures on a purely transcendental extension of the underlying field of an adelic curve. Let $S=(K,(\Omega,\mathcal A,\nu),\phi)$ be an adelic curve and $I$ be a non-empty set. We consider the polynomial ring $K[\boldsymbol{T}_I]$ spanned by $I$, where $\boldsymbol{T}_I=(T_i)_{i\in I}$ denotes the variables. Let $\mathbb N^{\oplus I}$ be the set of vectors $\boldsymbol{d}=(d_i)_{i\in I}\in\mathbb N^I$ such that $d_i=0$ for all but a finite number of $i\in I$. For any vector $\boldsymbol{d}=(d_i)_{i\in I}\in\mathbb N^{\oplus I}$,  we denote by $\boldsymbol{T}^{
\boldsymbol{d}}$ the monomial
\[\prod_{i\in I,\, d_i>0}T_i^{d_i}.\]
If $g$ is an element of $K[\boldsymbol{T}_I]$, for any $\boldsymbol{d}\in\mathbb N^{\oplus I}$ we denote by $a_{\boldsymbol{d}}(g)$ the coefficient of $\boldsymbol{T}^{\boldsymbol{d}}$ in the writing of $g$ as a $K$-linear combination of monomials.
For convenience, $K[\boldsymbol{T}_I]$ means $K$ in the case where $I = \emptyset$.

\begin{lemm}\label{lem:UFD:K:Lambda}
\begin{enumerate}[label=\rm(\arabic*)]
\item
Let $J$ be a subset of $I$. If $f$ and $g$ are two elements of $K[\boldsymbol{T}_I]$ such that $fg$ belongs to $K[\boldsymbol{T}_J]$, then both polynomials $f$ and $g$ belong to $K[\boldsymbol{T}_J]$.

\item The ring $K[\boldsymbol{T}_I]$ is a unique factorization domain and $K[\boldsymbol{T}_I]^{\times}=K^{\times}$. 
\end{enumerate}
\end{lemm}

\begin{proof}
(1) For $i \in I$ and $\varphi\in K[\boldsymbol{T}_I]$, the degree of $\varphi$ with
respect to $x_i$ is denoted by $\deg_i(\varphi)$. Note that the function $\deg_i(\ndot)$ satisfies the equality $\deg_i(fg) = \deg_i(f) + \deg_i(g)$, so that $\deg_i(f) = \deg_i(g) = 0$ once $i \in I \setminus J$,
which means that $g$ and $h$  belong to $K[\boldsymbol{T}_J]$.

\medskip
(2)
For any finite subset $J$ of $I$, it is well-known that $K[\boldsymbol{T}_J]$ is a unique factorization domain.
Moreover, for $f \in K[\boldsymbol{T}_I] \setminus \{ 0 \}$,
there is a finite subset $J$ of $I$ such that $f \in K[\boldsymbol{T}_J]$.
Thus the first assertion follows from (1). The second assertion is a direct consequence of (1) in the particular case where $J=\varnothing$.
\end{proof}

 Let $L=K(\boldsymbol{T}_I)$ be the field of fractions of $K[\boldsymbol{T}_I]$. As in \S\ref{Sec: Transcendental fibration}, we pick in each equivalent class of irreducible polynomials in $K[\boldsymbol{T}_I]$, a representative to form a subset $\mathscr P_{K[\boldsymbol{T}_I]}$. For each element $F\in\mathscr P_{K[\boldsymbol{T}_I]}$, we let $\mathrm{ord}_F(\ndot)$ be the discrete valuation on $L$ defined by $F$ and let $|\ndot|_F:=\mathrm{e}^{-\operatorname{ord}_F(\ndot)}$ be the corresponding absolute value. Let $\deg(\ndot)$ be the degree function on $K[\boldsymbol{T}_I]$. Note that for any $(f,g)\in K[\boldsymbol{T}_I]^2$ one has
\[\deg(f+g)\leqslant\max\{\deg(f),\deg(g)\},\quad
\deg(fg)=\deg(f)+\deg(g).\]
Therefore the function $-\deg(\ndot)$ extends to a discrete valuation on $L$. Denote by $|\ndot|_\infty$ the corresponding absolute value, defined as
\[|\ndot|_{\infty}=\mathrm{e}^{\deg(\ndot)}.\]
Note that the following product formula holds
\[\forall\, g\in L\setminus\{0\},\quad \ln|g|_{\infty}+\sum_{F\in\mathscr P_{K[\boldsymbol{T}_I]}}\deg(F)\ln|g|_F=0.\]
In other words, if we equip $\Omega_L':=\mathscr P_{K[\boldsymbol{T}_I]}\amalg\{\infty\}$ with the discrete $\sigma$-algebra $\mathcal A_L'$ and the measure $\nu_L'$ such that 
\[\nu_L'(\{\infty\})=1\text{ and }\nu_L'(\{F\})=\deg(F)\] for any $F\in\mathscr P_{K[\boldsymbol{T}_I]}$, then $(L,(\Omega_L',\mathcal A_L',\nu_L'),\phi_L')$ forms a proper adelic curve, where \[\phi_L':\mathscr P_{K[\boldsymbol{T}_I]}\amalg\{\infty\}\rightarrow M_L\] sends $x$ to $|\ndot|_x$.\medskip

\begin{rema}
Let $\boldsymbol{X}_{I \cup \{ \infty \}} = \{ X_i \}_{i \in I} \cup \{ X_{\infty} \}$ be the variables indexed by $I \cup \{ \infty \}$. Let $\varphi : K[\boldsymbol{X}_{I \cup \{ \infty \}}] \to K[\boldsymbol{T}_I]$ be the homomorphism given by 
$\varphi(f) = f((T_i)_{i \in I}, 1)$. If $f$ is an irreducible homogeneous polynomial in $K[\boldsymbol{X}_{I \cup \{ \infty \}}]$ and $f \not= X_{\infty}$, then $\varphi(f)$ is an irreducible polynomial in $K[\boldsymbol{T}_I]$. Moreover, for any irreducible polynomial $g$ in $K[\boldsymbol{T}_I]$,
there is an irreducible homogeneous polynomial $f$ in $K[\boldsymbol{X}_{I \cup \{ \infty \}}]$ such that $\varphi(f) = g$.
Note that the above $|\ndot|_{\infty}$ comes from the irreducible polynomial $X_{\infty}$, so that
the corresponding element is $1 = \varphi(X_{\infty})$.
\end{rema}

\begin{lemm}[Gauss's Lemma]
\label{lem:Gauss:lemma}
Let $|\ndot|$ be a non-Archimedean absolute value on $K$. 
We fix $\pmb{e} = (e_i)_{i \in I} \in \RR_{>0}^{I}$. For $\boldsymbol{d} = (d_i)_{i\in I} \in\mathbb N^{\oplus I}$, we set $\pmb{e}^{\boldsymbol{d}} := \prod_{i \in I} e_i^{d_i}$.
We denote by $|\ndot|_{\pmb{e},L}$ the function on $K[\boldsymbol{T}_I]$ sending $f\in K[\boldsymbol{T}_I]$ to \[\max_{\boldsymbol{d}\in\mathbb N^{\oplus I}}|a_{\boldsymbol{d}}(f)|\pmb{e}^{\boldsymbol{d}}.\]
Then, for any $(f,g)\in K[\boldsymbol{T}_I]^2$ one has
\[|fg|_{\pmb{e},L}=|f|_{\pmb{e},L}\cdot|g|_{\pmb{e},L}\quad\text{and}\quad|f+g|_{\pmb{e},L}\leqslant\max\{|f|_{\pmb{e},L},|g|_{\pmb{e},L}\}. \] In particular,
$|\ndot|_{\pmb{e},L}$ extends to an absolute value on $L=K(\boldsymbol{T}_I)$.
\end{lemm}

\begin{proof}
If we set $f = \sum_{\boldsymbol{d}'\in\mathbb N^{\oplus I}} a_{\boldsymbol{d}'} \boldsymbol{T}^{\boldsymbol{d}'}$ and
$g = \sum_{\boldsymbol{d}''\in\mathbb N^{\oplus I}} b_{\boldsymbol{d}''} \boldsymbol{T}^{\boldsymbol{d}''}$, then
\[
fg = \sum_{\boldsymbol{d}\in\mathbb N^{\oplus I}} \left( \sum_{\substack{\boldsymbol{d}', \boldsymbol{d}'' \in\mathbb N^{\oplus I}, \\ \boldsymbol{d}' + \boldsymbol{d}'' = \boldsymbol{d}}}
a_{\boldsymbol{d}'} b_{\boldsymbol{d}''} \right) \boldsymbol{T}^{\boldsymbol{d}}
\quad\text{and}\quad
f + g = \sum_{\boldsymbol{d}\in\mathbb N^{\oplus I}} (a_{\boldsymbol{d}} + b_{\boldsymbol{d}}) \boldsymbol{T}^{\boldsymbol{d}}
\]
Thus it is easy to see
\begin{equation}\label{eqn:lem:Gauss:lemma:01}
\begin{cases}
|fg|_{\pmb{e},L}\leqslant|f|_{\pmb{e},L}\cdot|g|_{\pmb{e},L}, \\
|f+g|_{\pmb{e},L}\leqslant\max\{|f|_{\pmb{e},L},|g|_{\pmb{e},L}\}.
\end{cases}
\end{equation}

Let $\Sigma_f = \{ \boldsymbol{d}'\in\mathbb N^{\oplus I} \mid |a_{\boldsymbol{d}'}|\pmb{e}^{\boldsymbol{d}'} =
|f|_{\pmb{e},L} \}$ and
$\Sigma_g = \{ \boldsymbol{d}''\in\mathbb N^{\oplus I} \mid |b_{\boldsymbol{d}''}|\pmb{e}^{\boldsymbol{d}''} =
|g|_{\pmb{e},L} \}$. Let $\leqslant_{\mathrm{lex}}$ be the lexicographic order on $\mathbb N^{\oplus I}$.
We choose $\pmb{\delta}(f) \in \Sigma_f$ and $\pmb{\delta}(g) \in \Sigma_g$ such that
$\boldsymbol{d}' \leqslant_{\mathrm{lex}} \pmb{\delta}(f)$ and
$\boldsymbol{d}'' \leqslant_{\mathrm{lex}} \pmb{\delta}(g)$
for all $\boldsymbol{d}' \in \Sigma_f$ and $\boldsymbol{d}'' \in \Sigma_g$.

\begin{enonce}{Claim}\label{claim:lem:Gauss:lemma:01}
One has $|a_{\boldsymbol{d}'}||b_{\boldsymbol{d}''}| \leqslant |a_{\pmb{\delta}(f)}||b_{\pmb{\delta}(g)}|$
for all $\boldsymbol{d}', \boldsymbol{d}'' \in \mathbb N^{\oplus I}$ with 
$\boldsymbol{d}' + \boldsymbol{d}'' = \pmb{\delta}(f) + \pmb{\delta}(g)$.
Moreover, the equality holds if and only if $\boldsymbol{d}' = \pmb{\delta}(f)$
and $\boldsymbol{d}'' = \pmb{\delta}(g)$.
\end{enonce}

\begin{proof}
As $|a_{\boldsymbol{d}'}| \pmb{e}^{\boldsymbol{d}'} \leqslant |f|_{\pmb{e}, L}$ and
$|b_{\boldsymbol{d}''}| \pmb{e}^{\boldsymbol{d}''} \leqslant |g|_{\pmb{e}, L}$, one has
\[
|a_{\boldsymbol{d}'}||b_{\boldsymbol{d}''}| \leqslant \frac{|f|_{\pmb{e}, L} |g|_{\pmb{e}, L}}{\pmb{e}^{\boldsymbol{d}'+ \boldsymbol{d}''}} = \frac{|a_{\pmb{\delta}(f)}|\pmb{e}^{\pmb{\delta}(f)} |b_{\pmb{\delta}(g)}|\pmb{e}^{\pmb{\delta}(g)}}{\pmb{e}^{\boldsymbol{d}'+ \boldsymbol{d}''}} = |a_{\pmb{\delta}(f)}||b_{\pmb{\delta}(g)}|.
\]
We assume that the equality holds. Then $\boldsymbol{d}' \in \Sigma_f$ and $\boldsymbol{d}'' \in \Sigma_g$, so that
$\boldsymbol{d}' \leqslant_{\mathrm{lex}} \pmb{\delta}(f)$ and $\boldsymbol{d}'' \leqslant_{\mathrm{lex}} \pmb{\delta}(g)$.
Therefore, one has the assertion because $\boldsymbol{d}' + \boldsymbol{d}'' = \pmb{\delta}(f) + \pmb{\delta}(g)$.
\end{proof}

The above claim implies that
\[
\left| \sum_{\substack{\boldsymbol{d}', \boldsymbol{d}'' \in\mathbb N^{\oplus I}, \\ \boldsymbol{d}' + \boldsymbol{d}'' = \pmb{\delta}(f) + \pmb{\delta}(g)}}
a_{\boldsymbol{d}'} b_{\boldsymbol{d}''} \right| \pmb{e}^{\pmb{\delta}(f) + \pmb{\delta}(g)}= |a_{\pmb{\delta}(f)}|\pmb{e}^{\pmb{\delta}(f)} |b_{\pmb{\delta}(g)}|\pmb{e}^{\pmb{\delta}(g)} = |f|_{\pmb{e}, L} |g|_{\pmb{e}, L},
\]
which means that $|fg|_{\pmb{e}, L} \geqslant |f|_{\pmb{e}, L}|g|_{\pmb{e}, L}$, as required.
\end{proof}

For any $\omega\in\Omega\setminus \Omega_{\infty}$, let $|\ndot|_{\omega,L}$ be the absolute value on $L$ such that 
\[\forall\,g=\sum_{\boldsymbol{d}\in\mathbb N^{\oplus I}} a_{\boldsymbol{d}}(g)\boldsymbol{T}_I^{\boldsymbol{d}}\in K[\boldsymbol{T}_I],\quad
|g|_{\omega,L}:=\sup_{\boldsymbol{d}\in\mathbb N^{\oplus I}}|a_{\boldsymbol{d}}(g)|_\omega.\]
By Lemma~\ref{lem:Gauss:lemma}, this absolute value is an extension of $|\ndot|_{\omega}$ on $K$. Let \[((\Omega_{L,\omega},\mathcal A_{L,\omega},\nu_{L,\omega}),\phi_{L,\omega})\] be the adelic structure on $L$ which consists of a single copy of the absolute value $|\ndot|_{\omega,L}$, equipped with the unique probability measure. We denote by $S_{L,\omega}$ the adelic curve $(L,(\Omega_{L,\omega},\mathcal A_{L,\omega},\nu_{L,\omega}),\phi_{L,\omega})$.

\begin{prop}\label{Pro: S is an admissible family}
If $\Omega_{\infty} = \emptyset$, then  
family $(S_{L,\omega})_{\omega\in\Omega}$ is an admissible fibration over $S$. 
\end{prop}
\begin{proof}
Let 
$g$ be a non-zero element of $K[\boldsymbol{T}_I]$, $(F_j)_{j=1}^n$ be a finite family of elements of $\mathscr P_{K[\boldsymbol{T}_I]}$ containing $\{F\in \mathscr P_{K[\boldsymbol{T}_I]}\,|\,\operatorname{ord}_F(g)\neq 0\}$ and $(C_j)_{j=1}^n$ be a family of non-negative constants. One has
\[\int_{\Omega_{L,\omega}}|g|_x\indic_{|F_1|_x\leqslant C_1,\ldots,|F_n|_x\leqslant C_n}\,\nu_{L,\omega}(\mathrm{d}x)=\max_{\boldsymbol{d}\in\mathbb N^{\oplus I}}|a_{\boldsymbol{d}}(g)|_\omega\cdot\prod_{j=1}^n\prod_{\boldsymbol{d}\in\mathbb N^{\oplus I}}\indic_{|a_{\boldsymbol{d}}(F_j)|_\omega\leqslant C_j}.\]
Therefore the function
\[(\omega\in\Omega)\longmapsto\int_{\Omega_{L,\omega}}|g|_x\indic_{|F_1|_x\leqslant C_1,\ldots,|F_n|_x\leqslant C_n}\,\nu_{L,\omega}(\mathrm{d}x).\]
is $\mathcal A$-measurable. Moreover, for any element $F$ of $\mathscr P_{K[\boldsymbol{T}_I]}$, one has
\begin{equation}\label{Equ: local integral for finite places}\int_{\Omega_{L,\omega}}\ln|F|_x\,\nu_{L,\omega}(\mathrm{d}x)=\max_{\boldsymbol{d}\in\mathbb N^{\oplus I},\,a_{\boldsymbol{d}}(F)\neq 0}\ln|a_{\boldsymbol{d}}(F)|_\omega.\end{equation}
Therefore the function 
\[(\omega\in\Omega)\longmapsto\int_{\Omega_{L,\omega}}\ln|F|_x\,\nu_{L,\omega}(\mathrm{d}x)\]
is $\nu$-integrable. 
\end{proof}

\begin{rema}
In the case where $\Omega_{\infty} = \emptyset$ and
the adelic curve $S$ is proper, for any $\boldsymbol{d}$ such that $a_{\boldsymbol{d}}(F)\neq 0$, one has 
\[\int_{\omega\in\Omega}\ln|a_{\boldsymbol{d}}(F)|_\omega\,\nu(\mathrm{d}\omega)=0,\]
and hence 
\[h_{S_L}(F)=\int_{\Omega}\nu(\mathrm{d}\omega)\int_{\Omega_{L,\omega}}\ln|F|_x\,\nu_{L,\omega}(\mathrm{d}x)\geqslant 0.\]
\end{rema}

\section{Arithmetic adelic structure}
\label{sec:Arithmetic adelic structure}

In this section, we provides a ``standard'' construction of an adelic structure for a countable field of characteristic zero. More precisely, for any countable field $E$ of characteristic zero, we will construct an adelic curve $S_E=(E,(\Omega_E,\mathcal A_E,\nu_E),\phi_E)$, which satisfies the following properties:
\begin{enumerate}[label=\rm(\arabic*)]
\item\label{Item: propre} $S_E$ is proper.

\item For any $\omega \in \Omega_E$, the absolute value $\phi_E(\omega)$ is not trivial.

\item The set $\Omega_{E,\fin}$ of $\omega\in\Omega$ such that $\phi_E(\omega)$ is non-Archimedean is infinite but countable.

\item\label{Item: Northcott} Let $E^{\mathrm{ac}}$ be an algebraic closure of $E$.
If $E_0$ is a subfield of $E^{\mathrm{ac}}$ such that $E_0$ is finitely generated over $\QQ$,
then 
\[\{ a \in E^{\mathrm{ac}} \mid \text{$h_{S_E \otimes_E E^{\mathrm{ac}}} (1,a) \leqslant C$ and $[E_0(a) : E_0] \leqslant \delta$} \}\] is finite for all $C \in \RR_{\geqslant 0}$ and
$\delta \in \ZZ_{\geqslant 1}$.
\end{enumerate}

\begin{defi}\label{Def: arithmetic adelic structure}
Let $K$ be a countable field of characteristic $0$. An adelic structure of $K$ which satisfies the above conditions \ref{Item: propre}--\ref{Item: Northcott} is said to be \emph{arithmetic}.
\end{defi}

\begin{rema}
Note that the condition \ref{Item: Northcott} is analogous to Northcott's property  in Diophantine geometry. In Arakelov geometry of adelic curve, we say that an adelic curve $S=(K,(\Omega,\mathcal A,\nu),\phi)$ \emph{has Northcott property} if the set 
\[\{a\in K\,|\,h_S(1,a)\leqslant C\}\]
is finite for any $C\geqslant 0$ (see \cite[Definition 3.5.2]{CMArakelovAdelic}). In the case where the adelic curve $S$ is proper and has Northcott property, an analogue of Northcott's theorem holds (see \cite[Definition 3.5.3]{CMArakelovAdelic})
\end{rema}

In the remaining of the section, we fix a countable field $K$ of characteristic $0$ and a countable non-empty set $I$.  We equip $K$ with an adelic structure $((\Omega,\mathcal A,\nu),\phi)$ to form an adelic curve, which we denote by $S$.
We also fix a family $(\iota_\omega)_{\omega\in\Omega_\infty}$ of embeddings from $K$ to $\mathbb C$ such that $|\ndot|_\omega=|\iota_\omega(\ndot)|$ for any $\omega\in\Omega_\infty$ and that the map $(\omega\in\Omega_\infty)\mapsto\iota_\omega(a)$ is measurable for each $a\in K$ (see Theorem \ref{thm:measurable:family:embeddings}). For any element $f\in K[\boldsymbol{T}_I]$, we denote by $\iota_\omega(f)$ the polynomial in $\mathbb C[\boldsymbol{T}_I]$ defined as
\[\iota_\omega(f):=\sum_{\boldsymbol{d}\in\mathbb N^{\oplus I}}\iota_\omega(a_{\boldsymbol{d}}(f))\boldsymbol{T}_I^{\boldsymbol{d}}.\]
This defines a ring homomorphism from $K[\boldsymbol{T}_I]$ to $\mathbb C[\boldsymbol{T}_I]$, which extends to a homomorphism of fields form $K(\boldsymbol{T}_I)$ to $\mathbb C(\boldsymbol{T}_I)$, which we still denote by $\iota_\omega(\ndot)$.

\begin{enonce}[remark]{Notation}
For convenience, for any $f\in K[\boldsymbol{T}_I]$, the complex polynomial $\iota_{\omega}(f)\in\mathbb C[\boldsymbol{T}_I]$ is often denoted by $f_{\omega}$.
\end{enonce}

For any $t\in[0,1]$, we denote by $e(t)$ the complex number $\mathrm{e}^{2\pi t\sqrt{-1}}$.  For any $\omega\in\Omega_\infty$, we denote by $\Omega_{L,\omega}$ the set
\[
 \Omega_{L,\omega}:=\left\{ (t_i)_{i\in I} \in  [0,1]^I  \left|
\begin{array}{l}
\text{\small $( e(t_i))_{i\in I}$ is algebraically} \\
\text{\small independent over $\iota_{\omega}(K)$} 
\end{array} 
\right\}\right.. 
\]
Note that by definition one has
\begin{equation}\label{Equ: transcendental points}[0,1]^I\setminus\Omega_{L,\omega}=\bigcup_{f\in K[\boldsymbol{T}_I]\setminus\{0\}}\{(t_i)_{i\in I}\in[0,1]^I\,:\,f_\omega\big((e(t_i))_{i\in I}\big)=0\}.\end{equation}
We equip $[0,1]^I$ with the product $\sigma$-algebra (namely the smallest $\sigma$-algebra making measurable the projection maps to the coordinates) and the product of the uniform probability measure on $[0,1]$, denoted by $\eta_I$ (see \cite[\S4.2]{MR0651017} for the product of an arbitrary family of probability spaces).

\begin{lemm}\label{Lem: outside Omega negligible}
For any $\omega\in\Omega_\infty$, the subset $\Omega_{L,\omega}$ of $[0,1]^I$ is measurable, and $[0,1]^I\setminus\Omega_{L,\omega}$ is $\eta_I$-negligible.
\end{lemm}
\begin{proof} The measurability of $\Omega_{L,\omega}$ follows from \eqref{Equ: transcendental points}.

For any non-zero element of $K[\boldsymbol{T}_I]$, let
\[V_I(f)=\{(t_i)_{i\in I}\in[0,1]^I\,:\,f_\omega\big((e(t_i))_{i\in I}\big)=0\}.\]
Since $K$ and $I$ are countable,  $K[\boldsymbol{T}_I]$ is a countable set. Therefore, by \eqref{Equ: transcendental points}, to prove the second statement it suffices to show that $\eta_I(V_I(f))=0$. We first treat the case where $I$ is a finite set. Without loss of generality, we assume that $I=\{1,\ldots,n\}$, where $n\in\mathbb N$. The case where $n=0$ (namely $I=\varnothing$) is trivial since in this case $V_I(f)$ is empty. Assume that $n\geqslant 1$. For $t\in[0,1]$, let $f_t$ be the polynomial 
\[\iota_{\omega}(f)(T_1,\ldots,T_{n-1},e(t))\in \iota_\omega(K)(e(t))[T_1,\ldots,T_{n-1}].\]
Then by Fubini's theorem, one has 
\[\eta_{\{1,\ldots,n\}}(V_{\{1,\ldots,n\}}(f))=\int_{[0,1]}\eta_{\{1,\ldots,n-1\}}(V_{\{1,\ldots,n-1\}}(f_t))\,\mathrm{d}t=0,\]
where the second equality comes from the induction hypothesis.

We now consider the general case. Let $J$ be a finite subset of $I$ such that $f\in K[(T_i)_{i\in J}]$. By the definition of the product measure, one has
\[\eta_I(V_I(f))=\eta_J(V_J(f))=0.\]
\end{proof}

For any $\omega\in\Omega_\infty$, we equip $\Omega_{L,\omega}$ with the restriction of the product $\sigma$-algebra on $[0,1]^I$ and the restriction of the product probability measure $\eta_I$ to obtain a probability space denoted by $(\Omega_{L,\omega},\mathcal A_{L,\omega},\nu_{L,\omega})$.
Let $\phi_{L,\omega}:\Omega_{L,\omega}\rightarrow M_L$ be the map sending $x=(t_i)_{i\in I}\in\Omega_{L,\omega}$ to the absolute value
\[(f\in L)\longmapsto |f|_x:=\Big|f_\omega\big((e(t_i))_{i\in I}\big)\Big|.\]
Thus we obtain an adelic curve $S_{L,\omega}:=(L,(\Omega_{L,\omega},\mathcal A_{L,\omega},\nu_{L,\omega}),\phi_{L,\omega})$. 

We recall Jensen's formula for Mahler measure of polynomials (see \cite{MR1554908} for a proof). 
\begin{lemm}[Jensen's formula]
Let \[P(T)=a_d(T-\alpha_1)\cdots(T-\alpha_d)\in\mathbb C[T]\] be a complex polynomial of one variable $T$, with $a_d\in\mathbb C\setminus\{0\}$ and $(\alpha_1,\ldots,\alpha_d)\in\mathbb C^d$. One has
\[\int_0^1\ln|P(e(t))|\,\mathrm{d}t=\ln|a_d|+\sum_{j=1}^d\ln(\max\{1,|\alpha_j|\})\geqslant\ln|a_d|.\]
\end{lemm}

\begin{prop}\label{prop:admissible:fibration:K:T}
The family of adelic curves $(S_{L,\omega})_{\omega\in\Omega}$ is an admissible fibration over the adelic curve $S$. Moreover, in the case where the adelic curve $S$ is proper, for any $F\in\mathscr P_{K[T]}$, one has 
\[h_{S_L}(F):=\int_{\Omega}\nu(\mathrm{d}\omega)\int_{\Omega_{L,\omega}}\ln|F|_x\,\nu_{L,\omega}(\mathrm{d}x)\geqslant 0.\]
\end{prop}
\begin{proof}
{\bf Step 1.} By construction, for any $\omega\in\Omega$ and any $x\in\Omega_{L,\omega}$, the absolute value $\phi_{L,\omega}(x)$ on $L$ extends the absolute value $\phi(\omega)$ on $K$. 

{\bf Step 2.} Let $g$ be a non-zero element of $K[\boldsymbol{T}_I]$, $(F_j)_{j=1}^n$ be elements of $\mathscr P_{K[\boldsymbol{T}_I]}$ containing \[\{F\in\mathscr P_{K[\boldsymbol{T}_I]}\,:\,\operatorname{ord}_F(g)\neq 0\},\] and $(C_j)_{j=1}^n\in\mathbb R_{\geqslant 0}^n$. We show that the function
\begin{equation}\label{Equ: interal sur Omega}(\omega\in\Omega)\longmapsto\int_{\Omega_{L,\omega}}|g|_x\indic_{|F_1|_x\leqslant C_1,\ldots,|F_n|_x\leqslant C_n}\,\nu_{L,\omega}(\mathrm{d}x)\end{equation}
is $\mathcal A$-measurable. We choose a finite subset $J$ of $I$ such that $g$, $F_1,\ldots,F_n$ belong to $K[(T_i)_{i\in J}]$. By Lemma \ref{Lem: outside Omega negligible}, one has
\[\begin{split}
&\quad\int_{\Omega_{L,\omega}}|g|_x\indic_{|F_1|_x\leqslant C_1,\ldots,|F_n|_x\leqslant C_n}\,\nu_{L,\omega}(\mathrm{d}x)\\&=\int_{[0,1]^I}\Big|g_\omega\big((e(t_i))_{i\in I}\big)\Big|\prod_{j=1}^n\indic_{|F_{j,\omega}((e(t_i))_{i\in I})|\leqslant C_j}\,\eta_I(\mathrm{d}(t_i)_{i\in I})\\
&=\int_{[0,1]^J}\Big|g_\omega\big((e(t_i))_{i\in J}\big)\Big|\prod_{j=1}^n\indic_{|F_{j,\omega}((e(t_i))_{i\in I})|\leqslant C_j}\,\eta_J(\mathrm{d}(t_i)_{i\in J})
\end{split}\]
Note that $[0,1]^J$ is a separable 
compact metric space. By the criterion of measurability for functions on product measurable space proved in \cite[Lemma 9.2]{MR44536} and the measurability of integrals with parameter (see \cite[Lemma 1.26]{MR1464694}), we obtain the measurability of the function \eqref{Equ: interal sur Omega} on $\Omega_\infty$. The measurability of this function on $\Omega\setminus\Omega_\infty$ follows from Proposition \ref{Pro: S is an admissible family}.

{\bf Step 3.}  It remains to show that the function
\begin{equation}\label{Equ: int ln F}(\omega\in\Omega)\longmapsto \int_{\Omega_{L,\omega}}\ln|F|_x\,\nu_{L,\omega}(\mathrm{d}x)\end{equation}
is well defined and is integrable for any $F\in\mathscr P_{K[T]}$. By Proposition \ref{Pro: S is an admissible family} again, it suffices to show its integrability on $\Omega_\infty$.   Let 
\[\Theta:=\{\boldsymbol{d}\in\mathbb N^{\oplus I}\,:\,a_{\boldsymbol{d}}(F)\neq 0\}.\]
One has
\[\ln|F|_x\leqslant\max_{\boldsymbol{d}\in\Theta}\ln|a_{\boldsymbol{d}}(F)|_\omega+\ln(\operatorname{card}(\Theta)).\]
Therefore, for $\omega\in\Omega_\infty$, the integral
\[\int_{\Omega_{L,\omega}}\ln|F|_x\,\nu_{L,\omega}(\mathrm{d}x)\]
is well defined and the following inequality holds:
\begin{equation}\label{Equ: upper bound of the integral}\int_{\Omega_{L,\omega}}\ln|F|_x\,\nu_{L,\omega}(\mathrm{d}x)\leqslant\max_{\boldsymbol{d}\in\Theta}\ln|a_{\boldsymbol{d}}(F)|_\omega+\ln(\operatorname{card}(\Theta)).\end{equation} Moreover, by an argument similar to that in Step 2, it can be shown that the function 
\[(\omega\in\Omega_\infty)\longrightarrow \int_{\Omega_{L,\omega}}\ln|F|_x\,\nu_{L,\omega}(\mathrm{d}x)\]
is measurable. Finally, by writing 
\[\int_{\Omega_{L,\omega}}\ln|F|_x\,\nu_{L,\omega}(\mathrm{d}x)\]
as successive integrals, and then by applying Jensen's formula in a recursive way, we obtain that there exists $\boldsymbol{d}_0\in\Theta$ such that 
\begin{equation}\label{Equ: lower bound for int ln}\forall\,\omega\in\Omega_\infty,\quad \int_{\Omega_{L,\omega}}\ln|F|_x\,\nu_{L,\omega}(\mathrm{d}x)\geqslant\ln|a_{\boldsymbol{d}_0}(F)|_\omega.\end{equation} 
Combining this inequality with \eqref{Equ: upper bound of the integral} and the fact that $\nu(\Omega_\infty)<+\infty$ (see \cite[Proposition 3.1.2]{CMArakelovAdelic}), we obtain the integrability of the function \eqref{Equ: int ln F} on $\Omega_\infty$. Finally, applying \eqref{Equ: local integral for finite places} to $\omega\in\Omega\setminus\Omega_\infty$, the inequality \eqref{Equ: lower bound for int ln} leads to
\[h_{S_L}(F)\geqslant\int_{\omega\in\Omega}\ln|a_{\boldsymbol{d_0}}(F)|_\omega
\,\nu(\mathrm{d}\omega)=0\]
provided that the adelic curve $S$ is proper. The proposition is thus proved.  
\end{proof}

\begin{rema}
Note that, for $f \in L$,
\[
h_{S_L}(f) = \int_{\Omega_{\infty}} \nu(\mathrm{d}\omega) \int_{\Omega_{L,\omega}} \ln \Big|f_\omega\big(( e(t_i) )_{i \in I}\big)\Big|\, \eta_I (\mathrm{d} (t_i)_{i \in I}) + \int_{\Omega_{\rm{fin}}} 
\ln |f|_{\omega} \,\nu(\mathrm{d}\omega).
\]
Thus $h_{S_L}(1) = 0$ and $h_{S_L}(T_i) = 0$ for all $i \in I$.
\end{rema}

\begin{defi}\label{def:S:L:lambda}
As a corollary, to the admissible fibration $(S_{L,\omega})_{\omega\in\Omega}$ one can associate an adelic structure $((\Omega_L,\mathcal A_L,\nu_L),\phi_L)$ on $L$ as in Definition \ref{Def: adelic structure coming from a fibration}. We fix $\lambda \in \RR_{\geqslant 0}$.
Let $S_L^\lambda:=(L,(\Omega_L^\lambda,\mathcal A_L^\lambda,\nu_L^\lambda),\phi_L^\lambda)$ be an adelic curve with underlying field $L$ such that
\begin{enumerate}[label=\rm(\arabic*)]
\item $(\Omega_L^\lambda,\mathcal A_L^\lambda,\nu_L^\lambda)$ is the disjoint union of $(\Omega_L,\mathcal A_L,\nu_L)$ and $\mathscr P_{K[\boldsymbol{T}_I]} \cup \{ \infty \}$ equipped with the discrete $\sigma$-algebra and the measure satisfying
\[
\nu_L^\lambda(\{ F \}) = h_{S_L}(F) + \lambda \deg(F) \quad\text{and}\quad
\nu_L^\lambda(\{ \infty \}) = \lambda
\]
for any $F\in\mathscr P_{K[\boldsymbol{T}_I]}$.
\item the map $\phi_L^\lambda:\Omega_L^\lambda\rightarrow M_L$ extends $\phi_L$ and the map \[(x\in\mathscr P_{K[\boldsymbol{T}_I]} \cup \{\infty\} )\longmapsto|\ndot|_x.\] 
\end{enumerate}
The adelic curve $S^\lambda_L$ is called the \emph{$\lambda$-twisted compactification} of $S_L$.
\end{defi}

\begin{rema}\label{Rem:S:L:lambda}
Note that if $\lambda = 0$, then $S_L^{\lambda} = S_L^*$.
Moreover, if $K$ and $\Omega_{\rm{fin}}$ are countable and
$\mathcal A_{\Omega_{\rm{fin}}}$ is discrete, then
$L$ and $\Omega_{L,\rm{fin}}^*$ are countable and
$\mathcal A_{\Omega_{L,\rm{fin}}^*}$ is discrete.
\end{rema}

\begin{prop}\label{prop:proper:L:lambda}
The adelic curve $S_L^\lambda = (L,(\Omega_L^\lambda,\mathcal A_L^\lambda,\nu_L^\lambda),\phi_L^\lambda)$ is proper.
\end{prop}

\begin{proof}
If $\lambda = 0$, then the assertion follows from Proposition~\ref{Pro: compactification} and Proposition~\ref{prop:admissible:fibration:K:T}.
Note that
\begin{equation}\label{Equ: deg equality}
\deg(g) = \sum_{F \in \mathscr P_{K[\boldsymbol{T}_I]}} \deg(F) \ord_F(g)
\end{equation}
for $g \in L^{\times}$,
so that
\[
\sum_{F \in \mathscr P_{K[\boldsymbol{T}_I]}} (h_{S_L}(F) + \lambda \deg(F)) (-\ord_F(g)) + \lambda \deg(g)
= \sum_{F \in \mathscr P_{K[\boldsymbol{T}_I]}} h_{S_L}(F) (-\ord_F(g)),
\]
as required.
\end{proof}

\begin{rema}
The above result can be considered as a particular case of Proposition \ref{Pro: compactification 2}. In fact, if we equip $\mathscr P_{K[\boldsymbol{T}_I]}\cup\{\infty\}$ with the discrete $\sigma$-algebra $\mathcal A'$ and the measure $\nu'$ such that $\nu'(\{\infty\})=1$ and $\nu'(\{F\})=\deg(F)$, then \[(L,(\mathscr P_{K[\boldsymbol{T}_I]}\cup\{\infty\},\mathcal A',\nu'),\phi')\] 
forms an adelic curve, where $\phi'$ sends $x\in\mathscr P_{K[\boldsymbol{T}_I]}\cup\{\infty\}$ to the absolute value $|\ndot|_x$. Then the equality \eqref{Equ: deg equality} shows that this adelic curve is proper. Note that the restriction of $\nu^\lambda_L$ on $\mathscr P_{K[\boldsymbol{T}_I]}\cup\{\infty\}$ coincides with 
\[\lambda\nu_L'+\sum_{F\in\mathscr P_{K[\boldsymbol{T}_I]}}h_{S_L}(F)\operatorname{Dirac}_F.\]
Therefore the statement of Proposition \ref{prop:proper:L:lambda} follows from  Proposition \ref{Pro: compactification 2}. 
\end{rema}

\begin{lemm}\label{lem:deg:mu:formula}
\begin{enumerate}
\renewcommand{\labelenumi}{\textup{(\arabic{enumi})}}
\item
If $F_0, \ldots, F_r \in K[\boldsymbol{T}_I]$ with $(F_0, \ldots, F_r) \not= (0, \ldots, 0)$, then
{\allowdisplaybreaks
\begin{multline*}
\kern3em h_{S_L^\lambda}(F_0, \ldots, F_r) \leqslant
\int_{\Omega_{L,\infty}} \ln \max \{ |F_0|_{x}, \ldots, |F_r|_{x} \}
\nu_{L,\infty}(\mathrm{d} x) \\
\kern10em + \int_{\Omega_{\fin}} \ln \max \{ |F_0|_{\omega}, \ldots, |F_r|_{\omega} \}
\nu_{\fin}(\mathrm{d} \omega) \\
+ \lambda \max \{ \deg(F_0), \ldots, \deg(F_r) \}.
\end{multline*}}
Moreover, if $\mathrm{G.C.D}(F_0, \ldots, F_r) = 1$, then the equality holds.

\item
Fix $n \in I$ and let $I' = I \setminus \{ n \}$ and $L' = K(\boldsymbol{T}_{I'})$.
For $F \in K[\boldsymbol{T}_I] \setminus \{ 0 \}$, if we set 
$F = a_0 T_n^d + a_1 T_n^{d-1} + \cdots + a_d$ such that $a_0, a_1, \ldots, a_d \in K[\boldsymbol{T}_{I'}]$ and $a_0 \not= 0$,
then \[h_{S_{L'}^\lambda}(a_0, \ldots, a_d) \leqslant h_{S_L}(F) + \deg(F) (\lambda + \ln(2) 
\nu(\Omega_{\infty})).\] 
\end{enumerate}
\end{lemm}

\begin{proof}
(1) Note that 
\[
\max \{ |F_0|_{\xi}, \ldots, |F_r|_{\xi} \} \begin{cases}
\leqslant 1 & \text{in general}, \\
= 1 & \text{if $\mathrm{G.C.D}(F_0, \ldots, F_r) = 1$},
\end{cases}
\]
for $\xi \in\mathscr P_{K[\boldsymbol{T}_I]}$,
so that the assertion follows.

\medskip
(2)
Note that $d \leqslant \deg(F)$. We set $f = F/a_0$. 
For $y \in \Omega_{L', \infty}$,
let \[f_{y} = T_n^d + \iota_{y}(a_1/a_0)T_n^{d-1} + \cdots + \iota_{y}(a_d/a_0) = (T_n - \alpha_1) \cdots (T_n - \alpha_d)\] be the irreducible decomposition in $\CC[T_n]$.
Then,
\[
\iota_{y}(a_k/a_0) = (-1)^k \sum_{1 \leqslant i_1 < \cdots < i_k \leqslant d} \alpha_{i_1} \cdots \alpha_{i_k},
\]
so that {\allowdisplaybreaks
\begin{align*}
|a_k/a_0|_{y} & \leqslant \sum_{1 \leqslant i_1 < \cdots < i_k \leqslant d}  |\alpha_{i_1}| \cdots |\alpha_{i_k}| \leqslant \sum_{1 \leqslant i_1 < \cdots < i_k \leqslant d} \max \{1, |\alpha_{i_1}|\} \cdots\max \{1,  |\alpha_{i_k}|\} \\
& \leqslant \sum_{1 \leqslant i_1 < \cdots < i_k \leqslant d}  \max \{1, |\alpha_{1}|\} \cdots\max \{1,  |\alpha_{d}|\} \\
& \leqslant 2^{\deg(F)} \max \{1, |\alpha_{1}|\} \cdots\max \{ 1, |\alpha_{d}|\}
\end{align*}}
because $ \binom{d}{k} \leqslant 2^d \leqslant 2^{\deg(F)}$,
and hence one has
\[
\max \{ 1, |a_k/a_0|_{y}\} \leqslant 2^{\deg(F)} \max \{1, |\alpha_{1}|\} \cdots\max \{1,  |\alpha_{d}|\}.
\]
On the other hand, by Jensen's formula,
\[
\int_{0}^1 \ln |f_{y}(e(t_n))| \,\mathrm{d}t_n = \sum_{i=1}^d \ln \max \{ 1, |\alpha_i| \}.
\]
Therefore, one obtains
\[
\ln \max\{ 1, |a_k/a_0|_{y} \} \leqslant \int_{0}^1 \ln |f_{y}(e(t_n))| \,\mathrm{d}t_n +\deg(F)\ln(2)
\]
for all $k=1, \ldots, d$,
so that
{\allowdisplaybreaks
\begin{align*}
\ln \max\{ |a_0|_{y}, |a_1|_{y}, \ldots, |a_d|_{y} \} \\
& \kern-10em =
\ln |a_0|_{y} +
\ln \max\{ 1, |a_1/a_0|_{y}, \ldots, |a_d/a_0|_{y} \} \\
& \kern-10em \leqslant \ln |a_0|_{y} + \int_{0}^1 \ln |f_{y}(e(t_n))|\, \mathrm{d}t_n + \deg(F)\ln(2) \\
& \kern-10em = \int_{0}^1 \ln |F_{y}(e(t_n))|\, \mathrm{d}t_n + \deg(F)\ln(2).
\end{align*}}
Thus, by 
Fubini's theorem, 
{\allowdisplaybreaks
\begin{align*}
\int_{\Omega_{L,\infty}} \ln |F|_{x} \nu_{L, \infty}(\mathrm{d} x) &  =
\int_{\Omega_{L', \infty} \times [0,1]} \ln |F_{y}(e(t_n))|\, 
\nu_{L'}(\mathrm{d}y)\, \mathrm{d} t_n \\
&  \kern-8em= \int_{\Omega_{L', \infty}} \left(
\int_{0}^1 \ln |F_{y}(e(t_n))|\, \mathrm{d}t_n\right)\, \nu_{L'}(\mathrm{d} y) \\
&  \kern-8em\geqslant \int_{\Omega_{L', \infty}} \left( \ln \max\{|a_0|_{y},  \ldots, |a_d|_{y} \} - \deg(F) \ln(2) \right) \, \nu_{L'}(\mathrm{d} y) \\
& \kern-8em = \int_{\Omega_{L', \infty}} \ln \max\{ |a_0|_{y}, \ldots, |a_d|_{y}\}\,  \nu_{L'}(\mathrm{d} y) - \deg(F) \ln(2)\nu_{L'}(\Omega_{L', \infty}) \\
& \kern-8em = \int_{\Omega_{L', \infty}} \ln \max\{ |a_0|_{y}, \ldots, |a_d|_{y}\} \, \nu_{L'}(\mathrm{d} y) - \deg(F) \ln(2)\nu(\Omega_{\infty}).
\end{align*}}
On the other hand, note that \[
|F|_{\omega} = \max \{ |a_0|_{\omega}, \ldots, |a_d|_{\omega}\}
\]
for $\omega \in \Omega_{\rm{fin}}$, so that
{\allowdisplaybreaks
\begin{multline*}
\int_{\Omega_{L', \infty}} \ln \max\{ |a_0|_{y}, \ldots, |a_d|_{y}\} \, \nu_{L'}(\mathrm{d} y) \\
+
\int_{\Omega_{\fin}} \ln \max\{ |a_0|_{\omega}, \ldots, |a_d|_{\omega}\}\,  \nu(\mathrm{d} \omega) \\
\leqslant h_{S_L}(F) + \deg(F) \ln(2)\,\nu(\Omega_{\infty}).
\end{multline*}}{\allowdisplaybreaks
\begin{align*}
h_{S_{L'}^{\lambda}}(a_0,\ldots, a_d) & \leqslant  \int_{\Omega_{L', \infty}} \ln \max\{ |a_0|_{y}, \ldots, |a_d|_{y}\}  \,\nu_{L'}(\mathrm{d} y) \\
& \kern5em + \int_{\Omega_{\mathrm{fin}}}
\ln \max \{ |a_0|_{\omega}, \ldots, |a_d|_{\omega} \}\,\nu (\mathrm{d} \omega) \\
& \kern10em + \lambda\max \{\deg(a_0), \ldots, \deg(a_d) \} \\
& \leqslant h_{S_L}(F) + \deg(F) (\lambda + \ln(2)\nu(\Omega_{\infty})),
\end{align*}}
as required.
\end{proof}

Fix $n \in I$ and let $I' = I \setminus \{ n \}$ and $L' = K(\boldsymbol{T}_{I'})$.
For $F \in K[\boldsymbol{T}_I] \setminus \{ 0 \}$, we set
$F = a_0 T_n^{d} + \cdots + a_d$ such that $a_0, \ldots, a_d \in K[\boldsymbol{T}_{I'}]$ and $a_0 \not= 0$. We define $\bbnu(F)$ to be
\[\bbnu(F) := F/a_0 = T_n^{d} + (a_1/a_0)T_n^{d-1} + \cdots + (a_d/a_0).\] Note that $\bbnu(F)$ is a monic polynomial over $L'$.

\begin{prop}\label{prop:finiteness:deg:mu:bounded}
If $S_{L'}^{\lambda}$ has Northcott's property, then,
for $C \in \RR$ and $\delta \in \ZZ_{\geqslant 1}$, then the set 
\[\{ \bbnu(F) \mid \text{$F \in K[\boldsymbol{T}_I] \setminus \{ 0 \}$, $h_{S_L}(F) \leqslant C$ and
$\deg(F) \leqslant \delta$} \}\] is  finite.
\end{prop}

\begin{proof}
Let $\Theta := \{ F \in K[\boldsymbol{T}_I] \setminus \{ 0 \} \mid \text{$h_{S_L}(F) \leqslant C$ and
$\deg(F) \leqslant \delta$} \}$ and $\vartheta : \Theta \to \PP^{\delta}(L')$
be a map given by the following way:
for \[F = a_0 T_n^{d} + \cdots + a_d \in \Theta\quad (\text{$a_0, \ldots, a_d \in K[\boldsymbol{T}_I]$ and $a_0 \not= 0$}),\] \[\vartheta(F) := 
\overbrace{(a_0 : \cdots : a_d : 0 : \cdots : 0)}^{\delta+1}
\in \PP^{\delta}(L').\]
By Lemma~\ref{lem:deg:mu:formula},
{\allowdisplaybreaks
\begin{multline*}
h_{S_{L'}^{\lambda}}(\vartheta(F)) \leqslant h_{S_L}(F) + \deg(F) (\lambda + \ln(2) \nu(\Omega_{\infty})) \leqslant C + \delta(\lambda + \ln(2) \nu(\Omega_{\infty})).
\end{multline*}}
Thus the assertion of the proposition is a consequence of Northcott's property of $S_{L'}^{\lambda}$.
\end{proof}

\begin{prop}\label{prop:Northcott:property}
If $S$ has Northcott's property, $\operatorname{card}(I) < \infty$ and $\lambda > 0$, then $S_L^{\lambda}$ has also Northcott's property.
\end{prop}

\begin{proof}
We prove it by induction on $\operatorname{card}(I)$. If $\operatorname{card}(I)=0$, then the assertion is obvious because
$S_{L}^{\lambda} = S$.
Fix $n \in I$ and let $I' = I \setminus \{ n \}$ and $L' = K(\boldsymbol{T}_{I'})$.
It is sufficient to see that $\{ f \in L^{\times} \mid h_{S_L^{\lambda}} (f, 1 ) \leqslant C \}$ is finite for any $C$.
For $f \in L^{\times}$, let us choose $F_1, F_2 \in K[\boldsymbol{T}_{I}] \setminus \{ 0 \}$
such that $f = F_1/F_2$, and $F_1$ and $F_2$ are relatively prime.
We set
\[
\begin{cases}
F_1 = a_{10}T_n^{d_1} + \cdots + a_{1d_1}\quad(\text{$a_{10}, \ldots, a_{1d_1} \in K[\boldsymbol{T}_{I'}]$ and $a_{10} \not= 0$}),\\
F_2 = a_{20}T_n^{d_2} + \cdots + a_{2d_2}\quad(\text{$a_{20}, \ldots, a_{2d_2} \in K[\boldsymbol{T}_{I'}]$ and $a_{20} \not= 0$}). 
\end{cases}
\]

\begin{enonce}{Claim}\label{claim:thm:Northcott:property:01}
If $h_{S_L^{\lambda}}(f , 1) \leqslant C$, then one has the following:
\begin{enumerate}
\renewcommand{\labelenumi}{\textup{(\arabic{enumi})}}
\item $\max \{ \deg(F_1), \deg(F_2) \} \leqslant C/\lambda$ and $\max \{ h_{S_L}(F_1), h_{S_L}(F_2) \} \leqslant C$.
\item $h_{S_{L'}^{\lambda}}(a_{10} , a_{20}) \leqslant  C$.
\end{enumerate}
\end{enonce}

\begin{proof}
(1) As $C \geqslant h_{S_L^{\lambda}}(f, 1) = h_{S_L^{\lambda}}(F_1, F_2)$ and
$F_1$ and $F_2$ are relatively prime, by (1) in Lemma~\ref{lem:deg:mu:formula}, one has
{\allowdisplaybreaks
\begin{multline}\label{eqn:thm:Northcott:property:01}
C \geqslant \lambda \max \{ \deg(F_1), \deg(F_2) \} + \int_{\Omega_{L, \infty}} \ln \max \{ |F_1|_{x}, |F_2|_{x} \}\,
\nu_{L}(\mathrm{d} x) \\
+ \int_{\Omega_{\mathrm{fin}}} \ln \max \{ |F_1|_{\omega}, |F_2|_{\omega} \}\, \nu(\mathrm{d} \omega).
\end{multline}}
Thus, 
{\allowdisplaybreaks
\begin{align*}
C & \geqslant \lambda \max \{ \deg(F_1), \deg(F_2) \} + \max \{ h_{S_L}(F_1), h_{S_L}(F_2) \} 
\end{align*}}
Therefore, (1) follows because $h_{S_L}(F_1), h_{S_L}(F_2) \geqslant 0$. 

\medskip
(2) By (1) in Lemma~\ref{lem:deg:mu:formula},
{\allowdisplaybreaks
\begin{align*}
h_{S_{L'}^{\lambda}}(a_{10}, a_{20}) & \leqslant \lambda \max \{ \deg(a_{10}), \deg(a_{20}) \} \\
& \kern5em + \int_{\Omega_{L', \infty}} \ln \max \{ |a_{10}|_{y}, |a_{20}|_{y} \}\,\nu_{L'}(\mathrm{d} y) \\
& \kern8em + \int_{\Omega_{\mathrm{fin}}} \ln \max \{ |a_{10}|_{\omega}, |a_{20}|_{\omega} \} \,\nu(\mathrm{d} \omega).
\end{align*}}
Therefore, by \eqref{eqn:thm:Northcott:property:01},
it is sufficient to see the following:
{\allowdisplaybreaks
\begin{equation}\label{eqn:thm:Northcott:property:02}
\int_{\Omega_{L, \infty}} \ln \max \{ |F_1|_{x}, |F_2|_{x} \}
\nu_{L,\infty}(\mathrm{d} x)
\geqslant
\int_{\Omega_{L', \infty}} \ln \max \left\{ |a_{10}|_{y}, |a_{20}|_{y}
\right\} \,
\nu_{L'}(\mathrm{d} y)
\end{equation}}
and
{\allowdisplaybreaks
\begin{equation}\label{eqn:thm:Northcott:property:03}
\int_{\Omega_{\mathrm{fin}}} \ln \max \{ |F_1|_{\omega}, |F_2|_{\omega} \} \,\nu(\mathrm{d} \omega)  \geqslant
\int_{\Omega_{\mathrm{fin}}} \ln \max \{ |a_{10}|_{\omega}, |a_{20}|_{\omega} \} \,\nu(\mathrm{d} \omega).
\end{equation}}
Indeed, by Jensen's formula together with Fubini's theorem, 
{\allowdisplaybreaks
\begin{align*}
\int_{\Omega_{L, \infty}} \ln \max \{ |F_1|_{x}, |F_2|_{x} \}\,
\nu_{L}(\mathrm{d} x) \\
&\kern-8em= \int_{\Omega_{L', \infty} \times [0,1]} \ln \max_{i=1,2} \{ |F_{iy}(e(t_n))| \}
\,\nu_{L'}(\mathrm{d}y) \mathrm{d}t_n \\
&\kern-8em= \int_{\Omega_{L', \infty}} \left(\int_{0}^1 \ln \max_{i=1,2} \{ |F_{iy}(e(t_n))| \}\mathrm{d}t_{n} \right) \,\nu_{L'}(\mathrm{d}y) \\
& \kern-8em\geqslant \int_{\Omega_{L', \infty}} \max_{i=1,2} \left\{ \int_{0}^1 \ln |F_{iy}(e(t_n))|\mathrm{d}t_{n} \right\}\,\nu_{L'}(\mathrm{d}y) \\
& \kern-8em\geqslant \int_{\Omega_{L', \infty}} \max_{i=1,2} \left\{ \ln |a_{i0y}| \right\} \,\nu_{L'}(\mathrm{d}y) \\
&\kern-8em= \int_{\Omega_{L', \infty}} \ln \max \left\{ |a_{10}|_{y}, |a_{20}|_{y}
\right\} 
\,\nu_{L'}(\mathrm{d} y),
\end{align*}}
as required for \eqref{eqn:thm:Northcott:property:02}. Further, since $|F_1|_{\omega} \geqslant |a_{10}|_{\omega}$ and $|F_2|_{\omega} \geqslant |a_{20}|_{\omega}$, one has \eqref{eqn:thm:Northcott:property:03}.
\end{proof}

If we set 
\[
\begin{cases}
\Delta= \{ \bbnu(F) \mid \text{$F \in K[\boldsymbol{T}_I] \setminus \{ 0 \}$, $h_{S_L}(F) \leqslant C$ and $\deg(F) \leqslant C/\lambda$ }\}, \\
\Delta' = \{ a \in K(\boldsymbol{T}_{I'}) \mid \text{$h_{S_{L'}^{\lambda}}(a, 1) \leqslant C$} \},
\end{cases}
\]
then, by Proposition~\ref{prop:finiteness:deg:mu:bounded} together with the hypothesis of
induction, $\Delta$ and $\Delta'$ are finite.
Moreover, by Claim~\ref{claim:thm:Northcott:property:01},
if $h_{S_L^{\lambda}}(f, 1) \leqslant C$, then
\[\bbnu(F_1), \bbnu(F_2) \in \Delta\quad\text{and}\quad 
a_{10}/a_{20} \in \Delta'.\]
Thus the assertion follows because $f = (a_{10}/a_{20})(\bbnu(F_1)/\bbnu(F_2))$.
\end{proof}

\begin{rema}
\begin{enumerate}
\renewcommand{\labelenumi}{\textup{(\arabic{enumi})}}
\item Note that $h_{S_L^{\lambda}}(1, T_n) = \lambda$ for all $n \in I$,
so that Northcott's property does not hold for $S_L^{\lambda}$ if $\Lambda$ is
infinite.

\item Let $S_{\QQ}$ be the standard adelic structure of $\QQ$. Then, it is easy to see that
\[
h_{(S_{\QQ})_{\QQ(T)}^*}(1, T^n-1) = \int_{0}^1 \ln \max \{1, |e(nt) - 1|\}\, \mathrm{d}t \leqslant \ln 2
\]
for all $n \geqslant 0$, so that the Northcott's property does not hold for $S_{\QQ(T)}^*$.
\end{enumerate}
\end{rema}

\begin{theo}\label{thm:Northcott:S:Lambda:L}
We use the same notation as in Section~\ref{subsec:adelic:structure:of:purely:transcendental:fibration}.
We assume that $S$ has Northcott's propery and $\lambda > 0$.
Let $E$ be an algebraic closure  of $L = K(\boldsymbol{T}_{I})$.
If $E_0$ is a subfield of $E$ such that $E_0$ is finitely generated over $K$,
then $S_L^{\lambda} \otimes_L E$ has Northcott's property over $E_0$, that is, \[\Big\{ a \in E \mid \text{$h_{S_L^{\lambda} \otimes E}(1 , a) \leqslant C$ and $[E_0(a) : E_0] \leqslant \delta$} \Big\}\] is finite
for any $C \in \RR_{\geqslant 0}$ and $\delta\in\ZZ_{\geqslant 1}$.
\end{theo}

\begin{proof}
Since $E_0$ is finitely generated over $K$ and $E$ is algebraic over $L$,
we can choose a finite subset $I'$ of $I$ such that
$E_0(\boldsymbol{T}_{I'})$ is finite over $K(\boldsymbol{T}_{I'})$.
It is sufficient to see that the set 
\begin{equation}\label{eqn:lemma:Northcott:S:Lambda:L:01}
\Big\{ \alpha \in E \mid
\text{$h_{S_L^{\lambda} \otimes E}(1, \alpha) \leqslant C$ and $[ K(\boldsymbol{T}_{I'})(\alpha):K(\boldsymbol{T}_{I'})] \leqslant \delta$}\Big\}
\end{equation} 
is finite for any $C \in \RR_{\geqslant 0}$ and $\delta \in \ZZ_{\geqslant 1}$. Indeed, note that 
\begin{align*}
[K(\boldsymbol{T}_{I'})(\alpha):K(\boldsymbol{T}_{I'})] & \leqslant [E_0(\boldsymbol{T}_{I'})(\alpha):K(\boldsymbol{T}_{I'})] \\
& = [E_0(\boldsymbol{T}_{I'})(\alpha):E_0(\boldsymbol{T}_{I'})][E_0(\boldsymbol{T}_{I'}):K(\boldsymbol{T}_{I'})] \\
& \leqslant [E_0(\alpha) : E_0][E_0(\boldsymbol{T}_{I'}):K(\boldsymbol{T}_{I'})],
\end{align*}
so that
\begin{multline*}
\Big\{ a \in E \mid
\text{$h_{S_L^{\lambda} \otimes E}(1, a) \leqslant C$ and $[ E_0(a):E_0] \leqslant \delta$}\Big\} \\
\kern-5em \subseteq
\Big\{ \alpha \in E \mid
\text{$h_{S_L^{\lambda} \otimes E}(1, \alpha) \leqslant C$ and }\\
\kern10em \text{$[K(\boldsymbol{T}_{I'})(\alpha):K(\boldsymbol{T}_{I'})] \leqslant \delta[ E_0(\boldsymbol{T}_{I'}):K(\boldsymbol{T}_{I'})]$}\Big\}.
\end{multline*}
Let $\alpha$ be an element of the set \eqref{eqn:lemma:Northcott:S:Lambda:L:01}.
Let $f(t)$ be the minimal polynomial of $\alpha$ over $K(\boldsymbol{T}_{I'})$.
As $K(\boldsymbol{T}_{I})$ is a regular extension over $K(\boldsymbol{T}_{I'})$,
\[
K(\boldsymbol{T}_{I})[t]/f(t)K(\boldsymbol{T}_{I})[t] \simeq \left(K(\boldsymbol{T}_{I'})[t]/f(t)K(\boldsymbol{T}_{I'})[t]\right) \otimes_{K(\boldsymbol{T}_{I'})} K(\boldsymbol{T}_{I})
\]
is an integral domain, so that $f(t)$ is irreducible over $K(\boldsymbol{T}_{I})$, and hence $f(t)$ is also the minimal polynomial of $\alpha$ over $K(\boldsymbol{T}_{I})$.
We set
\[f = t^d + a_1t^{d-1} + \cdots + a_d\quad(a_1, \ldots, a_d \in K(\boldsymbol{T}_{I'})).\]
Then, in the same arguments as \cite[Theorem~3.5.3]{CMArakelovAdelic},
one has
\[
h_{S_L^{\lambda}}(1, a_1, \ldots, a_d) \leqslant \delta C + (\delta - 1) \ln(2) \nu(\Omega_{\infty}),
\]
so that
$h_{S_{K(\boldsymbol{T}_{I'})}^{\lambda}}(1, a_1, \ldots, a_d) \leqslant \delta C + (\delta - 1) \ln(2) \nu(\Omega_{\infty})$.
Therefore, the assertion is a consequence of Proposition~\ref{prop:Northcott:property}.
\end{proof}

\begin{theo}\label{thm:adelic:structure:countable}
If $E$ is a countable field of characteristic zero, then 
$E$ has an arithmetic adelic structure (see Definition \ref{Def: arithmetic adelic structure}). 
\end{theo}

\begin{proof}
We denote by $S$ the standard adelic curve with $\mathbb Q$ as underlying field. Recall that the measure space of $S$ is given by the set of all places of $\mathbb Q$ equipped with the discrete $\sigma$-algebra and the counting measure.
Let $\{ x_n \}_{n=1}^N$ be a transcendental basis of $E$ over $\QQ$. Note that $N$ might be $+\infty$. Moreover, $E$ is algebraic over $L:=\QQ((x_n)_{n=1}^N)$. Let $\lambda$ be a positive number.
Starting from the adelic curve $S$, by the way in Subsetion~\ref{subsec:adelic:structure:of:purely:transcendental:fibration},
let $S_{L}^{\lambda}$ be the $\lambda$-twisted compactification of $S_L$.
We claim that the adelic curve $S_L^\lambda\otimes_L E$ satisfies the properties (1) -- (4) characterizing an arithmetic adelic curve.
The property (1) follows from Proposition~\ref{prop:proper:L:lambda} and \cite[Proposition~3.4.10]{CMArakelovAdelic}.
The property (2) is obvious. 
For (3), see Lemma~\ref{lemma:discrete:algebraic} and Remark~\ref{Rem:S:L:lambda}. 
Finally the property (4) follows from Theorem~\ref{thm:Northcott:S:Lambda:L}.
\end{proof}

\subsection{Density of Fermat property over arithmetic function fields}
In this subsection, let us consider a simple application of Theorem~\ref{thm:adelic:structure:countable}
together with Faltings' theorem \cite{MR1307396}.
Let $K$ be a field. We denote by $\mu(K)$ the subgroup of $K^{\times}$
consisting of roots of unity in $K$, that is, \[\mu(K) := \{ a \in K \mid \text{$a^n = 1$ for some $n \in \ZZ_{>0}$} \}.\]
Let $N$ be a positive integer and let $F_N := \Spec(\ZZ[X, Y]/(X^N + Y^N - 1))$.
We say that \emph{$F_N$ has Fermat's property over $K$} if $x, y \in \mu(K) \cup \{ 0 \}$ for all $(x, y) \in F_N(K)$.
Then one has the following theorem.

\begin{theo}\label{thm:density:Fermat:conj}
If $K$ is an arithmetic function field, then
\[
\lim_{m\to\infty} \frac{\#\{ N \in \ZZ \mid \text{$1 \leq N \leq m$ and $F_N$ has Fermat's property over $K$} \}}{m}
= 1.
\]
\end{theo}

\begin{proof}
Let $S$ be a proper adelic structure of $K$ with Northcott's property (cf. Theorem~\ref{thm:adelic:structure:countable}).
Let us begin with the following claim:

\begin{enonce}{Claim}
\begin{enumerate}
\renewcommand{\labelenumi}{\textup{(\arabic{enumi})}}
\item
For $x, y \in K$, $h_S(x,y,1) = 0$ if and only if $x, y \in \mu(K) \cup \{ 0 \}$.

\item If $N \geq 4$, then there is a positive integer $m_0$ such that
$F_{Nm}$ has Fermat's property of every integer $m \geq m_0$.
\end{enumerate}
\end{enonce}

\begin{proof}
(1) We assume that $h_S(x,y,1) = 0$ for $x, y \in K$. Then $h_S(x^n, y^n, 1) = n h_S(x, y, 1) = 0$ for all $n \in \ZZ_{>0}$, so that,
by Northcott's property, \[\{ (x^n, y^n) \mid n \in \ZZ_{>0} \}\] is finite. Therefore, there are $n, n' \in \ZZ_{>0}$ such that
$n > n'$ and $(x^n, y^n) = (x^{n'}, y^{n'})$, and hence $x, y \in \mu(K) \cup \{ 0 \}$.
The converse is obvious.

\smallskip
(2) First of all, 
note that $F_N(K)$ is finite by Faltings' theorem \cite{MR1307396}.
We set
\[
\begin{cases}
H := \max \{h_S(x, y, 1) \mid (x, y) \in F_N(K) \}, \\
a := \inf \{ h_S(x, y, 1) \mid \text{$x, y \in K$ and $h_S(x, y, 1) > 0$} \}.
\end{cases}
\]
Note that $a > 0$ by Northcott's property. For a positive integer $m$ with $m \geq \exp(H/a)$, we assume that $h_S(x, y, 1) > 0$ for some $(x, y) \in F_{Nm}(K)$.
Then, as $(x^m, y^m) \in F_N(K)$, 
\[
H \geq h_S(x^m, y^m, 1) = m h_S(x, y, 1) \geq m a,
\]
so that $\exp(H/a) \geq \exp(m)$, and hence $m \geq \exp(m)$. This is a contradiction.
Therefore, $h_S(x, y, 1) = 0$ for all $(x, y) \in F_{Nm}(K)$. Thus, by (1), $F_{Nm}$ has Fermat's property.
\end{proof}

By (2) together with \cite[Lemma~5.16]{IKMFaltings},
one can conclude the assertion of the theorem.
\end{proof}

In the case where $K = \QQ$, it was proved by \cite{MR734070, MR777765, MR766439} (cf. \cite{MR1719329}). A general number field case is treated in \cite{IKMFaltings}. The above theorem gives an evidence of the following conjecture:

\begin{enonce}{Conjecture}[Fermat's conjecture over an arithmetic function field] Let K be an arithmetic function field. 
Then is there a positive integer $N_0$ depending on $K$ such that $F_N$ has Fermat's property over K for all $N \geq N_0$?
\end{enonce}

\section{Polarized adelic structure}\label{sec:Polarized adelic structure}
In this subsection, we recall an adelic structure induced by a polarization of a field.
Let $K$ be a finitely generated field over $\mathbb Q$ and $n$ be the transcendental degree of $K$ over $\QQ$.
Let $\mathscr B\rightarrow\Spec\mathbb Z$ be a normal projective arithmetic variety such that the function field of $\mathscr B$ is $K$.
Note that $\dim \mathscr B = n+1$.
Let \[\big(\mathscr B; \overline{\mathscr H}_1 = (\mathscr H_1, h_1), \ldots, \overline{\mathscr H}_n = (\mathscr H_n, h_n)\big)\] be data with the following properties:
\begin{enumerate}
\renewcommand{\labelenumi}{\textup{(\arabic{enumi})}}
\item $\mathscr H_1, \ldots, \mathscr H_n$ are invertible $\OO_{\mathscr B}$-modules that
are nef along all fibers of $\mathscr B \to \Spec(\ZZ)$.

\item The second entries $h_1, \ldots, h_n$ are semipositive metrics of $\mathscr H_1, \ldots, \mathscr H_n$ on
$\mathscr B(\CC)$, respectively.

\item For each $i=1,\ldots,n$, the associated height function with $\overline{\mathscr H}_i$ is non-negative
\end{enumerate}
According to \cite{MR1779799}, the data $(\mathscr B; \overline{\mathscr H}_1, \ldots, \overline{\mathscr H}_n)$
is called a \emph{polarization of $K$}.

Let $x$ be a $\CC$-valued point of $\mathscr B$, that is, there are a unique scheme point $p_x \in \mathscr B$ and
a unique homomorphism $\phi_x :\OO_{\mathscr B, p_x} \to \CC$ such that
$x$ is given by $\phi_x$.
We say $x$ is \emph{generic} if $p_x$ is the generic point of $\mathscr B$.
We denote the set of all generic $\CC$-valued points by $\mathscr B(\CC)_{\mathrm{gen}}$.
Note that the measure of $\mathscr B(\CC) \setminus \mathscr B(\CC)_{\mathrm{gen}}$ is zero.

The polarization $(\mathscr B; \overline{\mathscr H}_1, \ldots, \overline{\mathscr H}_n)$ yields a proper
adelic structure of $K$ in the following way. 
First of all, we set
\[
\begin{cases}
\Omega_{\infty} := \mathscr B(\CC)_{\mathrm{gen}},\\
\Omega \setminus \Omega_{\infty} :=
\text{the set of all prime divisors on $\mathscr B$}.
\end{cases}
\]
For each element of $\omega \in \Omega$, $|\ndot|_{\omega}$ is give by 
\[
\begin{cases}
| f |_x := |\phi_x(f)| & \text{if $x \in \Omega_{\infty}$},\\
|f|_{\Gamma} := \exp\left(-\ord_{\Gamma}(f) 
\right) &
\text{if $\Gamma \in \Omega \setminus \Omega_{\infty}$}
\end{cases}
\]
for $f \in K$. Note that $\Omega_{\infty}$ is a measurable subset of a projective space, so that
one can give the standard measurable space structure and
a measure on $\Omega_{\infty}$ is given by $c_1(\overline{\mathscr H}_1)\wedge \cdots \wedge c_1(\overline{\mathscr H}_n)$.
The measurable space structure on $\Omega \setminus \Omega_{\infty}$ is discrete and a measure $\nu$ on $\Omega \setminus \Omega_{\infty}$ is given by $\nu(\{ \Gamma \}) = (\overline{\mathscr H}_1 \cdots \overline{\mathscr H}_n \cdot (\Gamma, 0))$.
This adelic structure is called the \emph{polarized adelic structure by the polarization 
$(\mathscr B; \overline{\mathscr H}_1, \ldots, \overline{\mathscr H}_n)$}.

\begin{enonce}[remark]{Example}
Let $h$ be the metric of $\OO_{\PP^1_{\CC}}(1)$ on $\PP^1_{\CC} = \Proj(\CC[T_0, T_1])$ given by
\[
|aT_0 + bT_1|_{h}(\zeta_0, \zeta_1) := \frac{|a\zeta_0 + b \zeta_1|}{\max \{ |\zeta_0|, |\zeta_1| \}}.
\]
Then $(\OO_{\PP^1_{\ZZ}}(1), h)$ gives rise to a semipositive metrized invertible $\OO_{\PP^1_{\ZZ}}$-module, so
that
\[
\Big((\PP^1_{\ZZ})^n; p_1^*(\OO_{\PP^1_{\ZZ}}(1), h), \ldots, p_n^*(\OO_{\PP^1_{\ZZ}}(1), h)\Big)
\]
yields to an adelic structure of the purely transcendental extension $\QQ(x_1, \ldots, x_n)$
over $\QQ$, where $p_i : (\PP^1_{\ZZ})^n \to \PP^1_{\ZZ}$ is the projection to the $i$-th factor. Note that it is nothing more than the adelic structure described in Section~\ref{subsec:adelic:structure:of:purely:transcendental:fibration} and
Section~\ref{sec:Arithmetic adelic structure}.
\end{enonce}


\chapter{Local intersection number and local height}
In this chapter, we fix a field $k$
equipped with an absolute value $|\ndot|$, such that $k$ is complete under the topology induced by the absolute value $|\ndot|$. 
In the case where $|\ndot|$ is Archimedean, $k$ is equal to $\mathbb R$ or $\mathbb  C$. In this case we always assume that $|\ndot|$ is the usual absolute value on $\mathbb R$ or $\mathbb C$. Note that the absolute value $|\ndot|$ extends in a unique way to any algebraic extension of $k$ (see \cite{Neukirch} Chapter II, Theorem 6.2). In particular, we fix an algebraic closure $k^{\operatorname{ac}}$, on which the absolute value $|\ndot|$ extends in a unique way.
Throughout this chapter, we denote the pair $(k, |\ndot|)$ by $v$. In the case where $|\ndot|$ is non-Archimedean, we denote by $\mathfrak o_v$ the valuation ring of $v=(k,|\ndot|)$, and by $\mathfrak m_v$ the maximal ideal of $\mathfrak o_v$.

\section{Reminder on completion of an algebraic closure}We denote by $\mathbb C_k$ the completion of an algebraic closure $k^{\operatorname{ac}}$ of $k$, on which the absolute value $|\ndot|$ extends by continuity. Recall that $\mathbb C_k$ is algebraically closed. A proof for the case where $k=\mathbb Q_p$ can for example be found in \cite[(10.3.2)]{MR2392026}, by using Krasner's lemma. The positive characteristic case is quite similar, but a supplementary argument is needed to show that there is no inseparable algebraic extension of $\mathbb C_k$. For the convenience of the readers, we include the proof here (see also \cite[Theorem~17.1]{MR2444734} for another proof). 

\begin{lemm}\label{Lem: completion is perfect}
Let $K$ be a field equipped with an absolute value $|\ndot|$ and $\widehat{K}$ be the completion of $K$. If the field $K$ is perfect, then also is $\widehat{K}$. 
\end{lemm}
\begin{proof}
Clearly it suffices to treat the case where the characteristic of $K$ is $p>0$. To prove that the completed field $\widehat{K}$ is perfect, we need to show that any element $a$ of $\widehat{K}$ has a $p$-th root in $\widehat{K}$. We choose a sequence $(a_n)_{n\in\mathbb N}$ of elements of $K$ which converges to $a$. Since $K$ is supposed to be perfect, for each $n\in\mathbb N$ we can choose $b_n\in K$ such that $b_n^p=a_n$. For any $(n,m)\in\mathbb N^2$ one has
\[|b_n-b_m|^p=|(b_n-b_m)^p|=|b_n^p-b_m^p|=|a_n-a_m|.\]
Hence $(b_n)_{n\in\mathbb N}$ is a Cauchy sequence in $K$, which converges to an element $b\in\widehat{K}$. Therefore
\[b^p=\lim_{n\rightarrow+\infty}b_n^p=\lim_{n\rightarrow+\infty}a_n=a,\]
as required. 
\end{proof}

\begin{prop}
The field $\mathbb C_k$ is algebraically closed. 
\end{prop}
\begin{proof}
It suffices to treat the case where the absolute value $|\ndot|$ is non-Archimedean.
We begin with proving that the field $\mathbb C_k$ is separably closed. Let $\mathbb C_k^{s}$ be a separable closure of $\mathbb C_k$, on which $|\ndot|$ extends in a unique way. Let $\alpha$ be a non-zero element of $\mathbb C_k^s$ and 
\[f(T)=T^r+a_1T^{r-1}+\cdots+a_r\in\mathbb C_k[T]\]
be the minimal polynomial of $\alpha$. Assume that $r\geqslant 2$. Let $\alpha_2,\ldots,\alpha_r$ be conjugates of $\alpha$ in $\mathbb C_k^s$ which are different from $\alpha$, and let
\[\varepsilon=\min_{j\in\{2,\ldots,r\}}|\alpha-\alpha_j|.\] Since $k^{\operatorname{ac}}$ is dense in $\mathbb C_k$, there exists a polynomial
\[g(T)=T^r+b_1T^{r-1}+\cdots+b_r\in k^{\operatorname{ac}}[T]\]
such that 
\[\max_{i\in\{1,\ldots,r\}}|\alpha|^{r-i}|b_i-a_i|<\varepsilon^r.\]
Since $k^{\operatorname{ac}}$ is algebraically closed, there exist elements $\beta_1,\ldots,\beta_r$ such that
\[g(T)=(T-\beta_1)\cdots(T-\beta_r).\]
One has
\[\prod_{i=1}^r|\alpha-\beta_i|=|g(\alpha)|=|g(\alpha)-f(\alpha)|\leqslant\max_{i\in\{1,\ldots,r\}}|\alpha|^{r-i}|b_i-a_i|<\varepsilon^n.\]
Hence there exists $\beta\in\{\beta_1,\ldots,\beta_r\}$ such that $|\alpha-\beta|<\varepsilon$. However, for any $\sigma\in\operatorname{Gal}(\mathbb C_k^s/\mathbb C_k)$, one has 
\[|\alpha-\beta|=|\sigma(\alpha-\beta)|=|\sigma(\alpha)-\beta|.\]
This implies $|\alpha-\sigma(\alpha)|<\varepsilon$, which leads to a contradiction. Therefore one has $r=1$, or equivalently, $\alpha\in\mathbb C_k$.

To show that $\mathbb C_k$ is algebraic closed, it suffices to check that $\mathbb C_k$ does not admit any algebraic inseparable extension, or equivalently, $\mathbb C_k$ is a perfect field. Note that any algebraic closed field is perfect (see \cite[Chapitre V, \S1, no.5, Proposition 5]{MR643362}). Hence the result follows from Lemma \ref{Lem: completion is perfect}.
\end{proof}

\section{Continuous metrics}
If $X$ is a  projective $k$-scheme, we denote by $X^{\mathrm{an}}$ the analytification of $X$. If  $k=\mathbb C$ and $|\ndot|$ is the usual absolute value, then $X^{\mathrm{an}}$ is a complex analytic space; if $|\ndot|$ is non-Archimedean, then the analytification $X^{\mathrm{an}}$ is defined in the sense of Berkovich (see \cite[\S4.3]{MR1070709}). Recall that any element $x$ of $X^{\mathrm{an}}$ consists of a scheme point of $X$ and an absolute value $|\ndot|_x$ on the residue field of the scheme point, which extends the absolute value $|\ndot|$ on $k$. We denote by $\widehat{\kappa}(x)$ the completion of the residue field of the scheme point with respect to the absolute value $|\ndot|_x$, on which the absolute value extends by continuity. In the remaining of the section, we fix a projective $k$-scheme $X$.
\begin{defi}\label{Defi: continuous metric}
Let $L$ be an invertible $\mathcal O_X$-module. We call \emph{continuous metric} on $L$ any family $\varphi=(|\ndot|_{\varphi}(x))_{x\in X^{\mathrm{an}}}$, where for each $x\in X^{\mathrm{an}}$, $|\ndot|_\varphi(x)$ is a norm on $L(x):=L\otimes_{\mathcal O_X}\widehat{\kappa}(x)$, such that, for any section $s$ of $L$ on a Zariski open subset $U$ of $X$, the map $|s|_\varphi$ from $U^{\mathrm{an}}$ to $\mathbb R_{\geqslant 0}$ sending $(x\in U^{\mathrm{an}})$ to $|s(x)|_{\varphi}(x)$ is a continuous function on $U^{\mathrm{an}}$. If $\varphi$ and $\psi$ are continuous metrics on $L$, we define
\[d(\varphi,\psi):=\sup_{x\in X^{\mathrm{an}}}\bigg|\ln\frac{|\ndot|_{\varphi}(x)}{|\ndot|_\psi(x)}\bigg|,\]
where
\[\frac{|\ndot|_{\varphi}(x)}{|\ndot|_\psi(x)}:=\frac{|\ell|_{\varphi}(x)}{|\ell|_\psi(x)}
\quad\text{for any $\ell\in L(x)\setminus\{0\}$.}\]

\end{defi}

\begin{exem}\phantomsection\label{Exa: Fubini-Study}
\begin{enumerate}[label=\rm(\arabic*)]
\item\label{Item: quotient metrics} Let $L$ be an invertible $\mathcal O_X$-module and $n$ be a positive integer. Let $(E,\|\ndot\|)$ be a finite-dimensional normed vector space over $k$. We assume that $p:E\otimes_{k}\mathcal O_X\rightarrow L^{\otimes n}$ is a surjective homomorphism of $\mathcal O_X$-modules, which induces a $k$-morphism $f:X\rightarrow\mathbb P(E)$ such that $L^{\otimes n}$ is isomorphic to $f^*(\mathcal O_E(1))$, where $\mathcal O_E(1)$ denotes the universal invertible sheaf on the projective space $\mathbb P(E)$ (see \cite[II.(4.2.3)]{EGA}).  For each point $x\in X^{\mathrm{an}}$ the norm $\mathopen{\|}{\ndot}\mathclose{\|}$ induces a quotient norm $|\ndot|(x)$ on $L(x)$ such that, for any $\ell\in L(x)\setminus\{0\}$,
\[|\ell|(x)=\inf_{\begin{subarray}{c}s\in E,\,\lambda\in\widehat{\kappa}(x)^{\times}\\
p(s)(x)=\lambda\ell^{\otimes n}
\end{subarray}}\big(|\lambda|_x^{-1}\mathopen{\|}s\mathclose{\|}\big)^{1/n}.\] 
The quotient norms $(|\ndot|(x))_{x\in X^{\mathrm{an}}}$ define a continuous metric on $L$, called the \emph{quotient metric induced by $\mathopen{\|}{\ndot}\mathclose{\|}$}. By definition, if $\|\ndot\|_1$ and $\|\ndot\|_2$ are two norms on $E$, and if $\varphi_1$ and $\varphi_2$ are quotient metrics induced by $\|\ndot\|_1$ and $\|\ndot\|_2$, respectively, then one has
\begin{equation}\label{Equ: distance of quotient metrics}d(\varphi_1,\varphi_2)\leqslant d(\|\ndot\|_1,\|\ndot\|_2):=\sup_{s\in E\setminus\{0\}}\Big|\ln\|s\|_1-\ln\|s\|_2\Big|.\end{equation}
\item Let $L$ be an invertible $\mathcal O_X$-module and $\varphi=(|\ndot|_{\varphi}(x))_{x\in X^{\mathrm{an}}}$ be a continuous metric on $L$. The dual norms of $|\ndot|_{\varphi}(x)$ on $L(x)^\vee$ form a continuous metric on $L^\vee$, which we denote by $\varphi^{\vee}$. Recall that for any $\ell\in L(x)\setminus\{0\}$, one has
\[|\ell^\vee|_{\varphi^\vee}=|\ell|_{\varphi}^{-1},\]
where $\ell^\vee$ denotes the linear form on $L(x)$ such that $\ell^\vee(\lambda\ell)=\lambda$ for any $\lambda\in\widehat{\kappa}(x)$.
\item Let $L_1$ and $L_2$ be invertible $\mathcal O_X$-modules, and $\varphi_1$ and $\varphi_2$ be continuous metrics on $L_1$ and $L_2$ respectively. Then the tensor product norms of $|\ndot|_{\varphi_1}(x)$ and $|\ndot|_{\varphi_2}(x)$ form a continuous metric on  $L_1\otimes L_2$, which we denote by $\varphi_1\otimes\varphi_2$. Note that, for any $\ell_1\in L_1(x)$ and $\ell_2\in L_2(x)$, one has
\[|\ell_1\otimes\ell_2|_{\varphi_1\otimes\varphi_2}(x)=|\ell_1|_{\varphi_1}(x)\cdot|\ell_2|_{\varphi_2}(x).\]
\item Let $f:Y\rightarrow X$ be a $k$-morphism of projective $k$-schemes. Let $L$ be an invertible $\mathcal O_X$-module, equipped with a continuous metric $\varphi$. Then the metric $\varphi$ induces by pull-back a continuous metric $f^*(\varphi)$ on $f^*(L)$ such that, for any $y\in Y^{\operatorname{an}}$ and any $\ell\in L(f^{\operatorname{an}}(y))$, one has
\[|f^*(\ell)|_{f^*(\varphi)}(y)=|\ell|_{\varphi}(f^{\operatorname{an}}(y)).\]
The metric $f^*(\varphi)$ is called the \emph{pull-back} of $\varphi$ by $f$.

\item\label{Item: extension of scalars metric} Let $k'/k$ be an extension of fields. We assume that the absolute value $|\ndot|$ extends to $k'$ and that the field $k'$ is complete with respect to the topology induced by the extended absolute value.   Let $X_{k'}$ be the fiber product $X\times_{\Spec k}\Spec k'$. We denote by $\pi:X_{k'}\rightarrow X$ the morphism of projection. Then the map \begin{equation}\label{Equ: pi natural}\pi^{\natural}:X_{k'}^{\operatorname{an}}\longrightarrow X^{\mathrm{an}},\end{equation}
sending any point $x'=(j(x'),|\ndot|_{x'})\in X_{k'}^{\mathrm{an}}$ to the pair consisting of the scheme point $\pi(j(x'))$ of $X$ and the restriction of $|\ndot|_{x'}$ to the residue field of $\pi(j(x'))$, is continuous (see \cite[Proposition 2.1.17]{CMArakelovAdelic}), where $j:X_{k'}^{\mathrm{an}}\rightarrow X_{k'}$ denotes the map sending a point in the analytic space to its underlying scheme point.

 Let $L$ be an invertible $\mathcal O_X$-module, equipped with a continuous metric  $\varphi$. Let $L_{k'}$ be the pull-back of $L$ by the morphism of projection $\pi$. The continuous metric $\varphi$ induces  a continuous metric $\varphi_{k'}$ on $L_{k'}$ such that, for any $x'\in X_{k'}^{\operatorname{an}}$ and any $\ell\in L(\pi^{\natural}(x'))$, one has
\[\forall\,a\in\widehat{\kappa}(x'),\quad |a\otimes\ell|_{\varphi_{k'}}(x')=|a|_{x'}\cdot|\ell|_{\varphi}(\pi^{\natural}(x')). \]
In particular, if $\psi$ is another continuous metric on $L$, then one has
\begin{equation}\label{Equ: distance extension of scalars}
d(\varphi_{k'},\psi_{k'})\leqslant d(\varphi,\psi).
\end{equation} 
\end{enumerate}  
\end{exem}

\begin{defi}\label{Def: extension of sclars}
Let $(E,\|\ndot\|)$ be a finite-dimensional normed vector space over $k$. We assume that the norm $\|\ndot\|$ is either ultrametric or induced by an inner product. Let $k'/k$ be an extension of fields, on which the absolute value $|\ndot|$ extends. We assume that field $k'$ is complete with respect to the extended absolute value. We denote by $\|\ndot\|_{k'}$ the following norm on $E_{k'}:=E\otimes_kk'$.
\begin{enumerate}[label=\rm(\arabic*)]
\item In the case where the absolute value $|\ndot|$ is non-Archimedean and the norm $\|\ndot\|$ is ultrametric, $\|\ndot\|_{k'}$ is the $\varepsilon$-extension of scalars of the norm $\|\ndot\|$. Namely, for any $t=s_1\otimes \lambda_1+\cdots+s_m\otimes \lambda_m\in E\otimes_kk'$
\[ \|t\|_{k'}:=\sup_{f\in E^\vee\setminus\{0\}}\frac{|\lambda_1f(s_1)+\cdots+\lambda_mf(s_m)|}{\|f\|_*},\]
where $\|\ndot\|_*$ denotes the dual norm of $\|\ndot\|$, which is defined as
\[\|f\|_*:=\sup_{x\in E\setminus\{0\}}\frac{|f(s)|}{\|s\|}.\]
This is an ultrametric norm on $E_{k'}$ such that $\|s\otimes a\|_{k'}=\|s\|\cdot|a|$ (see \cite[Proposition 1.3.1]{CMArakelovAdelic}). Moreover, if $(e_i)_{i=1}^r$ is an orthonormal basis of $(E,\|\ndot\|)$, then $(e_i\otimes 1)_{i=1}^r$ is an orthonormal basis of $(E_{k'},\|\ndot\|_{k'})$ (see \cite[Proposition 1.3.13]{CMArakelovAdelic}).  
\item In the case where the absolute value $|\ndot|$ is Archimedean, $k=\mathbb R$, $k'=\mathbb C$, and $\|\ndot\|$ is induced by an inner product $\emptyinnprod$, $\|\ndot\|_{\mathbb C}$ is the orthogonal extension of scalars of $\|\ndot\|$. Namely, for any $(s,t)\in E\times E$,
\[\|s\otimes 1+t\otimes\sqrt{-1}\|_{\mathbb C}:=(\|s\|^2+\|t\|^2)^{1/2}.\]
Clearly, for any $s\in E$ one has $\|s\otimes 1\|_{\mathbb C}=\|s\|$.
Note that the norm $\|\ndot\|_{\mathbb C}$ is induced by an inner product $\emptyinnprod_{\mathbb C}$ on $E_{\mathbb C}$ such that, for any $u=s\otimes 1+t\otimes\sqrt{-1}$ and $u'=s'\otimes 1+t'\otimes\sqrt{-1}$ in $E_{\mathbb C}$,
\[\mathopen{\langle}u,u'\mathclose{\rangle}=\langle s,s'\rangle+\langle t,t'\rangle+\sqrt{-1}(\langle s,t'\rangle-\langle t,s'\rangle).\]
Moreover, if $(e_i)_{i=1}^r$ is an orthonormal basis of $(E,\|\ndot\|)$, then $(e_i\otimes 1)_{i=1}^r$ is an orthonormal basis of $(E_{\mathbb C},\|\ndot\|_{\mathbb C})$.
\end{enumerate}
\end{defi}

\begin{rema}
\label{Rem: extension of normss}
Let $n$ be a positive integer. Assume that $p:E\otimes_{k}\mathcal O_X\rightarrow L^{\otimes n}$ is a surjective homomorphism of $\mathcal O_X$-modules, which induces a $k$-morphism $f:X\rightarrow\mathbb P(E)$ such that $L^{\otimes n}=f^*(\mathcal O_E(1))$. We equip $L$ with the quotient metric $\varphi$ induced by $\|\ndot\|$. In the case where the absolute value $|\ndot|$ is non-Archimedean, for any point $x\in X^{\mathrm{an}}$, the norm $|\ndot|_{\varphi^{\otimes n}}(x)$ on $L^{\otimes n}(x)$ coincides with the quotient norm on $L^{\otimes n}(x)$ induced by the norm $\|\ndot\|_{\widehat{\kappa}(x)}$ on $E\otimes_{k}{\widehat{\kappa}(x)}$ and the quotient map $p_x:E\otimes_{k}{\widehat{\kappa}(x)}\rightarrow L^{\otimes n}$. We refer the readers to \cite[Proposition 1.3.26 (i)]{CMArakelovAdelic} for a proof. As for the Archimedean case with $k=\mathbb R$ and $\widehat{\kappa}(x)=\mathbb C$, note that, if $s$ and $t$ are elements of $E$ and $a$ and $b$ are complex numbers such that
\[p_x(s)=a\ell^{\otimes n},\quad p_x(t)=b\ell^{\otimes n},\]
where $\ell$ is a fixed non-zero element of $L(x)$.   Then one has
\[p_x(s\otimes 1+t\otimes\sqrt{-1})=(a+b\sqrt{-1})\ell^{\otimes n}\]
and hence
\[\frac{(\|s\|^2+\|t\|^2)^{\frac 12}}{|a+b\sqrt{-1}|}\geqslant\frac{(\|s\|^2+\|t\|^2)^{\frac 12}}{|a|+|b|}\geqslant \frac{1}{\sqrt{2}}\frac{\|s\|+\|t\|}{|a|+|b|}\geqslant \frac{1}{\sqrt{2}}|\ell|_{\varphi^{\otimes n}}(x).\]
Therefore, the quotient norm on $L^{\otimes n}$ induced by $\|\ndot\|_{\widehat{\kappa}(x)}$ and the quotient map \[p_x:E\otimes_k\widehat{\kappa}(x)\longrightarrow L^{\otimes n}(x),\]
which is bounded from above by $|\ndot|_{\varphi^{\otimes n}}(x)$ by definition, is actually bounded from below by $\frac{1}{\sqrt{2}}|\ndot|_{\varphi^{\otimes n}}(x)$.

Let $k'/k$ be a valued extension of $(k,|\ndot|)$ which is complete. By extension of scalars, we obtain a surjective homomorphism of $\mathcal O_{X_{k'}}$-modules \[p_{k'}:E_{k'}\otimes_{k'}\mathcal O_{X_{k'}}\longrightarrow L_{k'}^{\otimes n},\] which corresponds to the $k'$-morphism $f_{k'}:X_{k'}\rightarrow\mathbb P(E_{k'})$. Let $\varphi$ be the quotient metric on $L$ induced by $\|\ndot\|$.  In the case where $|\ndot|$ is non-Archimedean, it turns out that the quotient metric on $L_{k'}$ induced by $\|\ndot\|_{k'}$ coincides with $\varphi_{k'}$. This fact follows from \cite[Proposition 1.3.15 (i)]{CMArakelovAdelic} and the above identification of the quotient metric to a family of quotient norms. In the Archimedean case with $k=\mathbb R$ and $k'=\mathbb C$, by the above estimate, in general the quotient metric $\varphi'$ on $L_{\mathbb C}$ induced by $\|\ndot\|_{\mathbb C}$ is different from  $\varphi_{\mathbb C}$. The above estimate actually shows that, for any $x\in X_{\mathbb C}^{\operatorname{an}}$ one has
\[2^{-\frac{1}{2n}}|\ndot|_{\varphi_{\mathbb C}}(x)\leqslant |\ndot|_{\varphi'}(x)\leqslant|\ndot|_{\varphi_{\mathbb C}}(x).\]
Note that the metric $\varphi_{\mathbb C}$ is still a quotient metric. In fact, if we consider the $\pi$-extension of scalars $\|\ndot\|_{\mathbb C,\pi}$ on $E_{\mathbb C}$ defined as
\[\forall\,t\in E_{\mathbb C},\quad\|t\|_{\mathbb C,\pi}:=\inf_{t=s_1\otimes\lambda_1+\cdots+s_m\otimes\lambda_m}\sum_{i=1}^m|\lambda_i|\cdot\|s_i\|.\]
Then the metric $\varphi_{\mathbb C}$ identifies with the quotient metric induced by $\|\ndot\|_{\mathbb C,\pi}$.
\end{rema}

\begin{defi}\label{Def: Fubini-Study metric}

Let $L$ be an invertible $\mathcal O_X$-module and $n$ be a positive integer. 
Let $(E,\|\ndot\|)$ be a finite-dimensional normed vector space over $k$. We assume that the norm $\|\ndot\|$ is either ultrametric or induced by an inner product. Let $p:E\otimes_k\mathcal O_X\rightarrow L^{\otimes n}$ be a surjective homomorphism of $\mathcal O_X$-module, which induces a $k$-morphism $f:X\rightarrow\mathbb P(E)$ such that $L^{\otimes n}$ is isomorphic to $f^*(\mathcal O_E(1))$. For each point $x\in X^{\mathrm{an}}$, the norm $\|\ndot\|_{\widehat{\kappa}(x)}$ on $E\otimes_k\widehat{\kappa}(x)$ induces by quotient a norm $|\ndot|(x)$ on $L^{\otimes n}(x)$. There then exists a unique continuous metric $\varphi$ on $L$ such that $|\ndot|_{\varphi^{\otimes n}}(x)=|\ndot|(x)$ for any $x\in X^{\mathrm{an}}$. The metric $\varphi$ is called the \emph{orthogonal quotient metric induced by $\|\ndot\|$}.
Note that, in the case where $|\ndot|$ is non-Archimedean or $(k,|\ndot|)$ is $\mathbb C$ equipped with the usual absolute value, the orthogonal quotient metric identifies with the quotient metric induced by $\|\ndot\|$ introduced in Example \ref{Exa: Fubini-Study} \ref{Item: quotient metrics}. Moreover, for any complete valued extension $k'/k$, the metric $\varphi_{k'}$ identifies with the orthogonal quotient metric induced by $\|\ndot\|_{k'}$ (see Remark \ref{Rem: extension of normss} above). 
\end{defi}

\begin{defi}
Let $L$ be a semi-ample invertible $\mathcal O_X$-module and $\varphi$ be a continuous metric on $L$. If there exists a sequence of quotient metrics $\varphi_n$ on $L$ such that \[\lim_{n\rightarrow+\infty}d(\varphi_n,\varphi)=0,\] we say that the metric $\varphi$ is \emph{semi-positive} (see \cite[\S2.2]{MR3783789}). In the case where $|\ndot|$ is Archimedean and $k=\mathbb C$, this definition is equivalent to the plurisubharmonicity of the metric $\varphi$ (see for example \cite[Theorem 3.5]{MR1254133}). 
\end{defi}

\begin{rema}\label{Rem: positivity extension scalar}
Let $L$ be an invertible $\mathcal O_X$-module. Let $k'/k$ be a complete valued extension of $k$, $X_{k'}$ be the fiber product $X\times_{\Spec k}\Spec k'$ and $\pi^{\natural}:X_{k'}^{\mathrm{an}}\rightarrow X^{\mathrm{an}}$ be the map defined in \eqref{Equ: pi natural}. If $\varphi$ and $\psi$ are two continuous metrics on $L$, then the metrics $\varphi_{k'}$ and $\psi_{k'}$ satisfy the relation (see \eqref{Equ: distance extension of scalars})
\[d(\varphi_{k'},\psi_{k'})\leqslant d(\varphi,\psi).\]
Therefore, if $\varphi$ is a semi-positive metric on $L$, then $\varphi_{k'}$ is also a semi-positive metric.
\end{rema}

\begin{defi}\label{Def: Fubini-Study}
Let $L$ be a very ample invertible $\mathcal O_X$-module and $\varphi$ be a continuous metric on $L$. For any positive integer $m$, the continuous metric $\varphi$ induces a seminorm $\|\ndot\|_{\varphi^{\otimes m}}$ on $H^0(X,L^{\otimes m})$ as follows:
\[\forall\,s\in H^0(X,L^{\otimes m}),\quad \|s\|_{\varphi^{\otimes m}}=\sup_{x\in X^{\operatorname{an}}}|s|_{\varphi^{\otimes m}}(x).\]
This seminorm is a norm notably when the scheme $X$ is reduced.
For each point $x\in X^{\mathrm{an}}$, the seminorm $\|\ndot\|_{\varphi^{\otimes m}}$ induces a quotient seminorm $|\ndot|_{\varphi^{(m)}}(x)$ on $L(x)$ such that, for any $\ell\in L(x)\setminus\{0\}$
\[|\ell|_{\varphi^{(m)}}(x)=\inf_{\begin{subarray}{c}
s\in H^0(X,L^{\otimes m}),\,\lambda\in\widehat{\kappa}(x)^{\times}\\
s(x)=\lambda\ell^{\otimes m}
\end{subarray}}(|\lambda|^{-1}\|s\|_{\varphi^{\otimes m}})^{1/m}.\]
This seminorm is actually a norm and is bounded from below by $|\ndot|_\varphi(x)$. The norms $(|\ndot|_{\varphi^{(m)}}(x))_{x\in X^{\operatorname{an}}}$ form a continuous metric on $L$, which we denote by $\varphi^{(m)}$.  
\end{defi}

\begin{prop}\label{Pro: distance varphi 1mvarphi2m}
Let $L$ be a very ample invertible $\mathcal O_X$-module. If $\varphi_1$ and $\varphi_2$ are two continuous metrics on $L$, then the following inequalities hold:
\[\forall\,m\in\mathbb N_{\geqslant 1},\quad d(\varphi_1^{(m)},\varphi_2^{(m)})\leqslant d(\varphi_1,\varphi_2).\]
\end{prop}
\begin{proof}
By definition, one has
\[\sup_{\begin{subarray}{c}
s\in H^0(X,L^{\otimes m})\\
\|s\|_{\varphi_1^{\otimes m}}\neq 0
\end{subarray}}\bigg|\ln\frac{\|s\|_{\varphi_1^{\otimes m}}}{\|s\|_{\varphi_2^{\otimes m}}}\bigg|\leqslant d(\varphi_1^{\otimes m},\varphi_2^{\otimes m})=m\,d(\varphi_1,\varphi_2).\]
Therefore,
\[d(\varphi_1^{(m)},\varphi_2^{(m)})\leqslant \frac 1md(\varphi_1^{\otimes m},\varphi_2^{\otimes m})\leqslant d(\varphi_1,\varphi_2).\]
\end{proof}

\begin{rema}\label{Rem: Fubini-Study of quotient metric}Let $(E,\|\ndot\|)$ be a finite-dimensional vector space over $k$, $m$ be a positive integer and $p:\pi^*(E)\rightarrow L^{\otimes m}$ be a surjective homomorphism of $\mathcal O_X$-modules, where $\pi:X\rightarrow\Spec k$ denotes the structural morphism of schemes. Let $\varphi$ be the quotient metric induced by $\|\ndot\|$. Note that $p$ induces by adjunction between $\pi^*$ and $\pi_*$ a $k$-linear map $\alpha:E\rightarrow H^0(X,L^{\otimes m})$. Let $s$ be an element of $H^0(X,L^{\otimes m})$. For any $x\in X^{\mathrm{an}}$, one has
\[|s|_{\varphi^{\otimes m}}(x)=\inf_{\begin{subarray}{c}t\in E,\,\lambda\in\widehat{\kappa}(x)^{\times}\\
\alpha(t)(x)=\lambda s(x)
\end{subarray}}\frac{\|t\|}{|\lambda|_x}.\]
In particular, for any $s$ in the image of the linear map $\alpha$, one has
\[\|s\|_{\varphi^{\otimes m}}\leqslant\inf_{t\in E,\,\alpha(t)=s}\|t\|.\]
Therefore, for $x\in X^{\mathrm{an}}$ and $\ell\in L(x)\setminus\{0\}$, one has 
\[|\ell|_{\varphi^{(m)}}(x)=\inf_{\begin{subarray}{c}s\in H^0(X,L^{\otimes m}),\,\lambda\in\widehat{\kappa}(x)^{\times}\\
s(x)=\lambda\ell^{\otimes m}
\end{subarray}}\Big(\frac{\|s\|_{\varphi^{\otimes m}}}{|\lambda|_x}\Big)^{1/m}\leqslant\inf_{\begin{subarray}{c}t\in E,\,\lambda\in\widehat{\kappa}(x)^{\times}\\
\alpha(t)(x)=b\ell^{\otimes m}
\end{subarray}}\Big(\frac{\|t\|}{|\lambda|_x}\Big)^{1/m}=|\ell|_{\varphi}(x).\]
Combining with the inequality $|\ell|_{\varphi^{(m)}}(x)\geqslant|\ell|_{\varphi}(x)$, we obtain the equality $\varphi^{(m)}=\varphi$.
\end{rema}

\begin{prop}\label{Pro: varphi 1 as a continuous metric}
Let $L$ be a very ample invertible $\mathcal O_X$-module, equipped with a continuous metric $\varphi$. Let $\|\ndot\|$ be a norm on the vector space $H^0(X,L^{\otimes n})$. For any $a>0$, let $\|\ndot\|_{a}$ be the norm on $H^0(X,L^{\otimes n})$ defined by
\[\forall\,s\in H^0(X,L^{\otimes n}),\quad\|s\|_{a}=\max\{\|s\|_{\varphi},a\|s\|\}=\max\Big\{\sup_{x\in X^{\operatorname{an}}}|s|_{\varphi}(x),a\|s\|\Big\},\]
and let $\varphi_{a}$ be the quotient metric on $L$ induced by $\|\ndot\|_a$. Then, for any $x\in X^{\mathrm{an}}$
\begin{equation}\label{Equ: metric varphi1 bounded by varphi eps}|\ndot|_{\varphi^{(1)}}(x)\leqslant|\ndot|_{\varphi_{a}}(x),\end{equation}
and there exists $a_0>0$ such that $\varphi_{a}=\varphi^{(1)}$ when $0<a\leqslant a_0$.
\end{prop}
\begin{proof}
By definition, one has $\|\ndot\|_{a}\geqslant\|\ndot\|_{\varphi}$. Hence the inequality \eqref{Equ: metric varphi1 bounded by varphi eps} holds.

Let $N_{\|\ndot\|_{\varphi}}$ be the null space of the norm $\|\ndot\|_{\varphi}$, which is defined as
\[N_{\|\ndot\|_{\varphi}}=\{s\in H^0(X,L)\,|\,\|s\|_{\varphi}=0\}.\] Let $E$ be the quotient vector space $H^0(X,L)/N_{\|\ndot\|_{\varphi}}$ and $\pi:H^0(X,L)\rightarrow E$ be the projection map. We denote by $\|\ndot\|_E$ the quotient norm of $\|\ndot\|$ on $E$ and $\|\ndot\|_{\varphi,E}$ be the quotient seminorm of $\|\ndot\|_{\varphi}$ on $E$, which is actually a norm satisfying the relation
\begin{equation}\label{Equ: quotient seminorm}\forall\,s\in H^0(X,L),\quad \|\pi(s)\|_{\varphi,E}=\|s\|_{\varphi}.\end{equation}
Since all norms on $E$ are equivalent, there exists $C>0$ such that $\|\ndot\|_{E}\leqslant C\|\ndot\|_{\varphi,E}$. Therefore, for any $x\in X^{\mathrm{an}}$, and any $\ell\in L(x)\setminus\{0\}$ one has
\[\begin{split}|\ell|_{\varphi_a}(x)&=\inf_{\begin{subarray}{c}
s\in H^0(X,L^{\otimes n}),\,\lambda\in\widehat{\kappa}(x)^{\times}\\
s(x)=\lambda\ell^{\otimes n}
\end{subarray}}\Big(\frac{\max\{\|s\|_{\varphi},a\|s\|\}}{|\lambda|}\Big)^{\frac 1n}\\
&=\inf_{\begin{subarray}{c}
s\in H^0(X,L^{\otimes n}),\,\lambda\in\widehat{\kappa}(x)^{\times}\\
s(x)=\lambda\ell^{\otimes n}
\end{subarray}}\Big(\frac{\max\{\|\pi(s)\|_{\varphi,E},a\|\pi(s)\|_E\}}{|\lambda|}\Big)=|\ell|_{\varphi^{(1)}}(x)
\end{split}\]
once $a<C^{-1}$, where the second equality comes from the fact that $s(x)=0$ when $s\in N_{\|\ndot\|_{\varphi}}$.
\end{proof}

\begin{prop}\label{Pro: convergence of varphin}
Let $L$ be a very ample invertible $\mathcal O_X$-module, equipped with a semi-positive continuous metric $\varphi$. Then one has
\[\lim_{m\rightarrow+\infty}d(\varphi^{(m)},\varphi)=0.\]
\end{prop}
\begin{proof}First of all, for positive integers $m$ and $m'$, one has 
\[\forall\,x\in X^{\mathrm{an}},\;\forall\,\ell\in L(x)\setminus\{0\},\quad |\ell|_{\varphi^{(m+m')}}^{m+m'}(x)\leqslant |\ell|_{\varphi^{(m)}}^m\cdot|\ell|_{\varphi^{(m')}}^{m'}. \] 
Therefore \[(m+m')d(\varphi^{(m+m')},\varphi)\leqslant md(\varphi^{(m)},\varphi)+m'd(\varphi^{(m')}).\]
By Fekete's lemma we obtain that the sequence
\[d(\varphi^{(m)},\varphi),\quad m\in\mathbb N,\;m\geqslant 1\]
converges to a non-negative real number, which is also equal to 
\[\inf_{m\in\mathbb N,\;m\geqslant 1}d(\varphi^{(m)},\varphi).\] 
Moreover,
since the metric $\varphi$ is semi-positive, there exist a sequence of positive integers $(m_n)_{n\in\mathbb N}$, a sequence of finite-dimensional normed vector spaces $((E_n,\|\ndot\|_n))_{n\in\mathbb N}$ and surjective homomorphisms of $\mathcal O_X$-modules $p_n:E_n\otimes_k\mathcal O_X\rightarrow L^{\otimes m_n}$ such that, if we denote by $\varphi_n$ the quotient metric on $L$ induced by $\|\ndot\|_n$, then one has
\[\lim_{n\rightarrow+\infty}d(\varphi_n,\varphi)=0.\]
By Remark \ref{Rem: Fubini-Study of quotient metric}, one has $\varphi_n^{(m_n)}=\varphi_n$ and hence 
\[d(\varphi^{(m_n)},\varphi)\leqslant d(\varphi^{(m_n)},\varphi_n)+d(\varphi_n,\varphi)=d(\varphi^{(m_n)},\varphi_n^{(m_n)})+d(\varphi_n,\varphi)\leqslant 2d(\varphi_n,\varphi),\]
where the last inequality comes from Proposition \ref{Pro: distance varphi 1mvarphi2m}. By taking the limite when $n\rightarrow +\infty$, we obtain that
\[\inf_{m\in\mathbb N,\,m\geqslant 1}d(\varphi^{(m)},\varphi)=0.\]
\end{proof}

\begin{defi}
Let $(L,\varphi)$ be a metrized invertible $\mathcal O_X$-module. We say that $(L,\varphi)$ is \emph{integrable} if there exist ample invertible $\mathcal O_X$-modules $L_1$ and $L_2$ equipped with semi-positive metrics $\varphi_1$ and $\varphi_2$ respectively, such that $L=L_1\otimes L_2^\vee$ and $\varphi=\varphi_1\otimes\varphi_2^\vee$. 
\end{defi}

\begin{defi}
We assume that $v$ is non-Archimedean.
Let $(L,\varphi)$ be a metrized invertible $\mathcal O_X$-module.
We say $\varphi$ is a \emph{model metric} if there are a positive integer $n$ and
a model $(\mathscr X, \mathscr L)$ of $(X, L^{\otimes n})$ 
such that $\varphi^{\otimes n}$ coincides with the metric arising from 
the model $(\mathscr X, \mathscr L)$ (cf. \cite[Subsection~2.3.2]{CMArakelovAdelic}).
In the above definition, we may assume that $\mathscr X$ is flat over $\mathfrak o_v$ (for details, see \cite[Subsection~2.3.2]{CMArakelovAdelic}).
In the case where $L$ is nef, if $\mathscr L$ is nef along the special fiber of
$\mathscr X \to \Spec(\mathfrak o_v)$, then 
the model $(\mathscr X, \mathscr L)$ is said to be \emph{nef} and
$\varphi$ is called a \emph{nef model metric}.
\end{defi}

\begin{rema}\label{rem:model:nef:metric}
Let $(\mathscr X, \mathscr L)$ be a model of $(X, L)$, $\mathscr X_{\mathrm{red}}$ be the reduced scheme associated with $\mathscr X$ and $\mathscr L_{\mathrm{red}} := \rest{\mathscr L}{X_{\mathrm{red}}}$.
For $x \in X^{\an}$, the morphism $\Spec (\mathfrak o_x) \to \mathscr X$ factors through
$\Spec (\mathfrak o_x) \to \mathscr X_{\mathrm{red}} \to \mathscr X$, and hence 
$\varphi_{\mathscr L}$ coincides with $\varphi_{\mathscr L_{\mathrm{red}}}$.
Moreover, $\mathscr L$ is nef with respect to $\mathscr{X} \to \Spec(\mathfrak o_v)$ if and only if
$\mathscr L_{\mathrm{red}}$ is nef with respect to $\mathscr{X}_{\mathrm{red}} \to \Spec(\mathfrak o_v)$.
\end{rema}

\begin{defi}
Let $(L,\varphi)$ be a metrized invertible $\mathcal O_X$-module.
We say $(L,\varphi)$ is \emph{smooth} if one of the following conditions is satisfied:
\begin{enumerate}[label=\rm(\roman*)]
\item if $v$ is Archimedean, $\varphi$ is a $C^{\infty}$-metric;

\item if $v$ is  non-Archimedean, $\varphi$ is a model metric.
\end{enumerate}
If $L$ is nef
and $v$ is non-Archimedean,
then $\varphi$ is said to be \emph{$M$-semi-positive}
if there is a sequence $\{ \varphi_m\}_{m=1}^{\infty}$ of nef model metrics of $L$ 
such that $\lim\limits_{m\to\infty} d(\varphi, \varphi_m) = 0$.
\end{defi}

\begin{lemm}\label{lemma:finitely:presented:model}
We assume that $v$ is non-Archimedean.
Let $L$ be an invertible $\OO_X$-module and $(\mathscr X, \mathscr L)$ be a model of $(X, L)$.
Then there is a model $(\mathscr X', \mathscr L')$ of $(X, L)$ with the following properties:
\begin{enumerate}[label=\rm(\arabic*)]
\item $\mathscr X' \to \Spec(\mathfrak o_v)$ is finitely presented, that is, $(\mathscr X', \mathscr L')$ is a coherent model of $(X,L)$  (cf. \cite[Subsection~2.3.2]{CMArakelovAdelic}).
\item $\mathscr X$ is a closed subscheme of $\mathscr X'$.
\item The special fiber of $\mathscr X' \to \Spec(\mathfrak o_v)$ coincides with
the special fiber of $\mathscr X \to \Spec(\mathfrak o_v)$.
\item $\rest{\mathscr L'}{\mathscr X} = \mathscr L$.
\end{enumerate}
\end{lemm}

\begin{proof}
By \cite[Corollary~5.16 in Chapter~II]{Hart77}, there are a polynomial ring $A := \mathfrak o_v[T_0, \ldots, T_N]$ over $\mathfrak o_v$ and
a homogeneous ideal $I$ of $A$ such that $\mathscr X = \Proj (A/I)$. We set $R := A/I$. Let $p : A \to R$ and $\pi : A \to A \otimes_{\mathfrak o_v} \mathfrak o_v/\mathfrak m_v = (\mathfrak o_v/\mathfrak m_v)[T_0, \ldots, T_N]$ be the natural homomorphisms.
There are homogeneous elements $h_1, \ldots, h_e$ of $R$ and $g_{ij} \in R_{(h_ih_j)}$ ($i, j \in \{ 1, \ldots, e \}$) such that $\mathscr X = \bigcup_{i=1}^e D_+(h_i)$ and
$\{ g_{ij} \}_{i, j \in \{ 1, \ldots, e \}}$ gives transition functions of $\mathscr L$,
where $R_{(h)}$ (for a homogenous element $h$) is the homogeneous localization with respect to $h$.
We choose a homogeneous element $H_i$ of $A$ such that $p(H_i) = h_i$.
Since \[\emptyset = \bigcap_{i=1}^e V_+(h_i) = V_+(h_1 R + \cdots + h_e R),\] we have
$R_+ \subseteq \operatorname{rad}(h_1 R + \cdots + h_e R)$
by \cite[Lemma~3.35 in Section~2.3]{MR1917232}, that is,
there is a positive integer $a$ such that $p(T_0)^a, \ldots, p(T_N)^a \in h_1 R + \cdots + h_e R$, so that
\begin{equation}\label{eqn:thm:equiv:semi-positive:M:semi-positive:01}
T_0^a, \ldots, T_N^a \in H_1 A + \cdots + H_e A + I.
\end{equation}
We also choose $G_{ij} \in A_{(H_iH_j)}$ such that $p(G_{ij}) = g_{ij}$ and $G_{ii} = 1$. As $g_{ij}g_{jl} = g_{il}$ on $R_{(h_ih_jh_l)}$, one can see
\begin{equation}\label{eqn:thm:equiv:semi-positive:M:semi-positive:02}
G_{ij}G_{jl} - G_{il} \in I_{(H_iH_jH_l)}
\end{equation}
for all $i,j, l\in \{ 1, \ldots, e \}$. Let $S = \mathfrak o_v \setminus \{ 0 \}$.
Since $I_S$ and $\pi(I)$ are homogeneous ideals of $k[T_0, \ldots, T_N]$ and $(\mathfrak o_v/\mathfrak m_v)[T_0, \ldots, T_N]$, respectively, $I_S$ and $\pi(I)$ are finitely generated ideals.
Therefore, by using \eqref{eqn:thm:equiv:semi-positive:M:semi-positive:01} and \eqref{eqn:thm:equiv:semi-positive:M:semi-positive:02}, one can find a finitely generated 
homogeneous ideal $I'$ of $A$ such that
\[
\begin{cases}
I' \subseteq I,\quad I'_S = I_S,\quad \pi(I') = \pi(I), \\
T_0^a, \ldots, T_N^a \in H_1 A + \cdots + H_e A + I', \\
G_{ij}G_{jl} - G_{il} \in I'_{(H_iH_jH_l)}\quad(\forall\ i, j, l \in \{ 1, \ldots, e \}).
\end{cases}
\]
Let $R' := A/I'$, $\mathscr X' := \Proj(R')$ and $p' : A \to R'$ be the natural homomorphism. 
Obviously $\mathscr X$ is a closed subscheme of $\mathscr X'$.
We set $h'_i = p'(H_i)$ and $g'_{ij} = p'(G_{ij})$.
Then $p'(T_0)^a, \ldots, p'(T_N)^a \in h'_1 R' + \cdots + h'_eR'$, which means that
$\mathscr X' = \bigcup_{i=1}^e D_+(h'_i)$ by \cite[Lemma~3.35 in Section~2.3]{MR1917232}.
Moreover, $g'_{ij}g'_{il} = g'_{il}$.
In particular, $g'_{ij}g'_{ji} = g'_{ii} = 1$, so that $g'_{ij} \in {R'}_{(h'_ih'_j)}^{\times}$.
This means that $\{ g'_{ij} \}_{i,j \in \{ 1, \ldots, e \}}$ gives rise to an invertible $\OO_{\mathscr X'}$-module
$\mathscr L'$ such that $\rest{\mathscr L'}{\mathscr X} = \mathscr L$.
Moreover, $(\mathscr X', \mathscr L')$ is a model of $(X, L)$ and
the special fiber of $\mathscr X' \to \Spec(\mathfrak o_v)$ is same as
the special fiber of $\mathscr X' \to \Spec(\mathfrak o_v)$, as required.
\end{proof}

\begin{prop}\label{prop:ample:every:fiber:ample}
Let $\mathscr X \to \Spec(\mathfrak o_v)$ be a model of $X$ and $\mathscr L$ be an invertible $\OO_{\mathscr X}$-module.
If $\mathscr L$ is ample on every fiber of $\mathscr X \to \Spec(\mathfrak o_v)$, then $\mathscr L$ is ample.
\end{prop}

\begin{proof}
By Lemma~\ref{lemma:finitely:presented:model}, there are a coherent model of $\mathscr X'$ of $X$ and an invertible $\OO_{\mathscr X'}$-module $\mathscr L'$
such that 
$\mathscr X$ is a closed subscheme of $\mathscr X'$, $\rest{\mathscr L'}{\mathscr X} = \mathscr L$ and
the special fiber of $\mathscr X' \to \Spec(\mathfrak o_v)$  coincides with the special fiber of $\mathscr X \to \Spec(\mathfrak o_v)$.
Note that $\mathscr L'$ is ample on every fiber of $\mathscr X' \to \Spec(\mathfrak o_v)$, and hence
$\mathscr L'$ is ample by \cite[IV-3, Corollaire~(9.6.4)]{EGA} because
$\mathscr X' \to \Spec(\mathfrak o_v)$ is finitely presented. Therefore $\mathscr L$ is ample.
\end{proof}

\begin{theo}
\label{thm:equiv:semi-positive:M:semi-positive}
We assume that $v$ is non-Archimedean and $|\ndot|$ is not trivial.
Let $L$ be a semi-ample invertible $\OO_X$-module and $\varphi$ be a continuous metric of $L$.
Then $\varphi$ is semi-positive if and only if $\varphi$ is $M$-semi-positive. 
\end{theo}

\begin{proof}
First we assume that $\varphi$ is semi-positive.
By Remark~\ref{rem:model:nef:metric}, we may assume that $X$ is reduced.
As $L$ is semi-positive, there is a positive integer $n_0$ such that $L^{\otimes n_0}$ is generated by global sections,
so we may assume that $L$ is generated by global sections, and hence
$L^{\otimes n}$ is generated by global sections for all $n \geqslant 1$.
Fix $\lambda \in \mathopen{]}0,1\mathclose{[}$ such that $\lambda < \sup \{ |a| : a \in k^{\times}, |a| < 1 \}$.
By \cite[Proposition~1.2.22]{CMArakelovAdelic},
there is a finitely generated lattice $\mathcal E_n$ of $H^0(X, L^{\otimes n})$ such that
$d(\|\ndot\|_{\mathcal E_n}, \|\ndot\|_{n\varphi}) \leqslant \log(\lambda^{-1})$.
Note that there is a morphism $f_n : X \to \PP(H^0(X, L^{\otimes n}))$ with $f_n^*(\OO_{\PP(H^0(X, L^{\otimes n}))}(1)) = L^{\otimes n}$, so we can find a morphism 
$\mathcal F_n : \mathscr X_n \to \PP(\mathcal E_n)$ over $\mathfrak o_v$
such that $\mathscr X_n$ is flat and projective over $\mathfrak o_v$ and
$\mathcal F_n$ is an extension of $f_n$ over $\mathfrak o_v$. If we set $\mathscr L_n = \mathcal F_n^*(\OO_{\PP(\mathcal E_n)}(1))$,
then $(\mathscr X_n, \mathscr L_n)$ is a flat model of $(X, L^{\otimes n})$.
As $\mathcal E_n \otimes_{\mathfrak o_v} \OO_{\PP(\mathcal E_n)} \to \OO_{\PP(\mathcal E_n)}(1)$ is surjective,
one also has the sujectivity of $\mathcal E_n \otimes_{\mathfrak o_v} \OO_{\mathscr X_n} \to \mathscr L_n$.
Therefore, by \cite[Proposition~2.3.12]{CMArakelovAdelic}, the model metric $\varphi_{\mathscr L_n}$ coincides
with the quotient metric induced by $\|\ndot\|_{\mathcal E_n}$.
Therefore, if we denote by $\varphi_n$ the quotient metic induced by $\|\ndot\|_{n\varphi}$, then, by \cite[Proposition~2.2.20]{CMArakelovAdelic}, 
\[
d(\varphi_{\mathscr L_n}, \varphi_n) \leqslant d(\|\ndot\|_{\mathcal E_n}, \|\ndot\|_{n\varphi})
\leqslant \log(\lambda^{-1}),
\]
which implies
\begin{multline*}
d((1/n)\varphi_{\mathscr L_n}, \varphi) \leqslant d((1/n)\varphi_{\mathscr L_n}, (1/n)\varphi_n) + d( (1/n)\varphi_n, \varphi) \\
\leqslant (1/n) \log(\lambda^{-1}) + d( (1/n)\varphi_n, \varphi),
\end{multline*}
and hence $\lim\limits_{n\to\infty} d((1/n)\varphi_{\mathscr L_n}, \varphi) = 0$.
Thus $\varphi$ is $M$-semi-positive because $\mathscr L_n$ is nef.

\medskip
Let us see the converse.
Let $(\mathscr X, \mathscr L)$ be a model of $(X, L)$ such that $\mathscr L$ is nef along
the special fiber of $\mathscr X \to \Spec(\mathfrak o_v)$.
Let $\varphi_{\mathscr L}$ be the metric arising from the model $(\mathscr X, \mathscr L)$.
It is sufficient to see that $\varphi_{\mathscr L}$ is semi-positive.
Let $\mathscr A$ be an ample invertible $\OO_{\mathscr X}$-module.
Then, for $n \geqslant 1$, $\mathscr A \otimes \mathscr L^{\otimes n}$ is ample
on every fiber of $\mathscr X \to \Spec(\mathfrak o_v)$, and hence,
by Proposition~\ref{prop:ample:every:fiber:ample}, $\mathscr A \otimes \mathscr L^{\otimes n}$ is ample
on $\mathscr X$ for all $n \geqslant 1$.
Therefore, by \cite[Proposition~2.3.17]{CMArakelovAdelic}, $\varphi_{\mathscr L}$ is semi-positive.
\end{proof}

\section{Green functions}

In this section, we fix a projective $k$-scheme $X$.

\begin{defi}
Let $D$ be a Cartier divisor on $X$. We call \emph{Green function} of $D$ any real-valued continuous function on $(X\setminus\Supp(D))^{\mathrm{an}}$ such that, for any regular meromorphic function  $f\in\Gamma(U,\mathscr M_X^{\times})$ which defines the Cartier divisor locally on a Zariski open subset $U$, the function $g+\log|f|$ on $(U\setminus\Supp(D))^{\mathrm{an}}$ extends to a continuous function on $U^{\mathrm{an}}$. A pair $(D,g)$ consisting of a Cartier divisor $D$ on $X$ and a Green function $g$ of $D$ is called a \emph{metrized Cartier divisor}. We denote by $\widehat{\operatorname{Div}}(X)$ the set of all metrized Cartier divisors on $X$.
Further $g$ is said to be \emph{smooth} if $(\OO_X(D), |\ndot|_g)$ is smooth.
A smooth Green function of $D=0$ is called a \emph{smooth function} on $X^{\an}$.      
\end{defi}

\begin{exem}\label{Rem: Green function of several points}
In the case where $D$ is the zero Cartier divisor, Green functions of $D$ are continuous functions on $X^{\mathrm{an}}$. In particular, if the Krull dimension of $X$ is zero, then $X^{\mathrm{an}}$ consists of isolated points. In this case any Cartier divisor $D$ on $X$ is trivial (see Remark \ref{remark:Cartier:div:on:0:dim}) and hence Green functions identify with elements in the real vector space spanned by $X^{\operatorname{an}}$.

In the case where $D$ is a principal Cartier divisor, namely a Cartier divisor of the form $\operatorname{div}(f)$, where $f$ is a regular meromorphic function, then by definition $-\ln|f|$ is a Green function of $\operatorname{div}(f)$. We denote by $\widehat{\operatorname{div}}(f)$ the pair $(\operatorname{div}(f),-\ln|f|)$. Such a metrized Cartier divisor is said to be \emph{principal}.

\end{exem}

\begin{rema}\label{Rem: metric induced by a green function}
Metrized Cartier divisors are closely related to metrized invertible sheafs. Let $D$ be a Cartier divisor on $X$. We denote by $\mathcal O_X(D)$ the sub-$\mathcal O_X$-module of $\mathscr M_X$ generated by $-D$. Let $(U_i)_{i\in I}$ be an open covering of $X$ such that, on each $U_i$ the Cartier divisor is defined by a regular meromorphic function $s_i$. Then the restriction of $\mathcal O_X(D)$ at $U_i$ is given by $\mathcal O_{U_i}s_i^{-1}$. If $g$ is a Green function of $D$, then it induces a continuous metric $\varphi_g=(|\ndot|_g(x))_{x\in X^{\mathrm{an}}}$ on $\mathcal O_X(D)$ such that \[|s_i^{-1}|_{g}:=\exp(-g-\ln|s_i|)\text{ on $U_i^{\mathrm{an}}$}.\]
Note that the metric of the canonical regular meromorphic section (see Definition \ref{Def: divisor of a section}) is given by 
\[|s_D|_g=|s_i\otimes s_i^{-1}|_g=\exp(-g)\text{ on $U_i$}.\]

Conversely, given an invertible $\mathcal O_X$-module $L$, any non-zero rational section $s$ of $L$ defines a Cartier divisor $\operatorname{div}(L;s)$. Moreover, if $\varphi$ is a continuous metric on $L$, then $-\ln|s|_{\varphi}$ is a Green function of $\operatorname{div}(L;s)$. 
We denote by $\operatorname{\widehat{div}}(\overline L;s)$ (or by $\operatorname{\widehat{div}}(s)$ for simplicity) the metrized Cartier divisor $(\operatorname{div}(L;s),-\ln|s|_\varphi)$.  

The above relation between metrized Cartier divisors and metrized invertible sheaves is important to define the following composition law on the set of metrized Cartier divisors. Let $(D_1,g_1)$ and $(D_2,g_2)$ be metrized Cartier divisors. Note that $\mathcal O_X(D_1+D_2)$ is canonically isomorphic to $\mathcal O_X(D_1)\otimes_{\mathcal O_X}\mathcal O_{X}(D_2)$. Moreover, under the canonical isomorphism 
\[\mathcal O_{X}(D_1+D_2)\stackrel{\sim}{\longrightarrow}\mathcal O_X(D_1)\otimes_{\mathcal O_X}\mathcal O_{X}(D_2),\]
the regular meromorphic section $s_{D_1+D_2}$ corresponds to $s_{D_1}\otimes s_{D_2}$. We equip the invertible sheaf $\mathcal O_{X}(D_1)$ and $\mathcal O_X(D_2)$ with the metrics $\varphi_{g_1}=(|\ndot|_{g_1}(x))_{x\in X^{\operatorname{an}}}$ and $\varphi_{g_2}=(|\ndot|_{g_2}(x))_{x\in X^{\operatorname{an}}}$ respectively, and $\mathcal O_X(D_1+D_2)$ with the tensor product metric $\varphi_{g_1}\otimes\varphi_{g_2}$. We then denote by $g_1+g_2$ the Green function in the metrized Cartier divisor $\widehat{\operatorname{div}}(s_{D_1+D_2})$. Clearly, for any $x\in \big(X\setminus(\operatorname{Supp}(D_1)\cup\operatorname{Supp}(D_2))\big)^{\operatorname{an}}$, one has
\[(g_1+g_2)(x)=g_1(x)+g_2(x).\]
Note that the set $\widehat{\operatorname{Div}}(X)$ of metrized Cartier divisors equipped with this composition law forms a commutative group.
\end{rema}

\begin{defi}
Let $(A,g)$ be a metrized Cartier divisor such that $\mathcal O_X(A)$ is an ample invertible $\mathcal O_X$-module (namely the Cartier divisor $A$ is ample). We say that the Green function $g$ is \emph{plurisubharmonic} if the metric $|\ndot|_g$ on $\mathcal O_X(A)$ is semi-positive. We refer to \cite[\S6.8]{ChamDucros} and \cite[\S6]{MR3975640} for a local version of positivity conditions.

We say that a metrized Cartier divisor  $(D, g)$ is \emph{integrable}
 if there are ample Cartier divisors $A_1$ and $A_2$ together with
plurisubharmonic Green functions $g_1$ and $g_2$ of $A_1$ and $A_2$, respectively, such that
$(D, g) = (A_1, g_1) - (A_2, g_2)$.  We denote by $\widehat{\operatorname{Int}}(X)$ the set of all integrable metrized Cartier divisors. This is a subgroup of the group $\widehat{\operatorname{Div}}(X)$ of metrized Cartier divisors.
\end{defi}

\begin{rema}\label{Rem: extension scalar Green function}
Let $k'/k$ be a valued extension which is complete. Let $X_{k'}$ be the fiber product $X\times_{\Spec k}\Spec k'$, and $\pi:X_{k'}\rightarrow X$ be the morphism of projection. Let $(D,g)$ be a metrized Cartier divisor on $X$. Then the pull-back $D_{k'}$ of $D$ by the morphisme $\pi$ is well defined (see Definition \ref{Def: pull back} and Remark \ref{Rem: extension of scalars}). Note that $\mathcal O_{X_{k'}}(D_{k'})$ is isomorphic with the pull-back of $\mathcal O_X(D)$ by $\pi$, and the canonical meromorphic section $s_{D_{k'}}$ of $D_{k'}$ identifies with the pull-back of $s_D$ by $\pi$. Let $\varphi_g$ be the continuous metric on $\mathcal O_{X}(D)$ induced by the Green function $g$. We denote by $g_{k'}$ the Green function of $D_{k'}$ defined as
\[g_{k'}=-\ln|s_{D_{k'}}|_{\varphi_{g,k'}},\]
where $\varphi_{g,k'}$ is the continuous metric on $\pi^*(\mathcal O_X(D))\cong\mathcal O_{X_{k'}}(D_{k'})$ induced by $\varphi_g$ (see Example \ref{Exa: Fubini-Study} \ref{Item: extension of scalars metric}). Note that, for any element $x'\in X_{k'}^{\mathrm{an}}$ such that $\pi^{\natural}(x')\in (X\setminus\operatorname{Supp}(D))^{\mathrm{an}}$, one has 
\[g_{k'}(x')=g(\pi^\natural(x')).\]

Moreover, the composition of $g$ with the restriction of $\pi^{\natural}$ to $(X_{k'}\setminus\operatorname{Supp}(D_{k'}))^{\mathrm{an}}$ forms a Green function of $D_{k'}$. We denote by $g_{k'}$ this Green function. By Remark \ref{Rem: positivity extension scalar}, if $\mathcal O_X(D)$ is semi-ample and $g$ is plurisubharmonic, then $g_{k'}$ is also plurisubharmonic. If $(D,g)$ is integrable, then $(D_{k'},g_{k'})$ is also integrable. Therefore the correspondance
$(D,g)\mapsto (D_{k'},g_{k'})$
defines a group homomorphism from $\widehat{\operatorname{Div}}(X)\rightarrow\widehat{\operatorname{Div}}(X_{k'})$, whose restriction to $\widehat{\operatorname{Int}}(X)$ defines a group homomorphism $\widehat{\operatorname{Int}}(X)\rightarrow\widehat{\operatorname{Int}}(X_{k'})$.  
\end{rema}

\begin{theo}
\label{thm:approximation:smooth:funcs}
Let $X$ be a $d$-dimensional projective and integral scheme over $k$. Let
$D$ be a nef and effective Cartier divisor and $g$ be a Green function of $D$ such that either
\begin{enumerate}[label=\rm(\alph*)]
\item if $v$ is Archimedean, the metric of $|\ndot|_g$ of $\OO_X(D)$ is $C^{\infty}$ and semi-positive, or
\item if $v$ is non-Archimedean, the metric of $|\ndot|_g$ of $\OO_X(D)$ is a nef model metric.
\end{enumerate}
Then there is a sequence $( \psi_n )_{n\in\mathbb N}$ of smooth functions on $X^{\an}$ with the following properties:
\begin{enumerate}[label=\rm(\arabic*)]
\item for all $n \in \NN$, $\psi_n \leqslant g$, $\psi_n \leqslant \psi_{n+1}$.
\item for each point $x \in X^{\an}$, $\sup \{ \psi_n(x) ; n \in \NN \} = g(x)$.

\item for all $n \in \NN$, $g - \psi_n$ is a Green function of $D$ such that either
\begin{enumerate}[label=\rm(\arabic{enumi}.\alph*)]
\item if $v$ is Archimedean, the metric of $|\ndot|_{g-\psi_n}$ of $\OO_X(D)$ is $C^{\infty}$ and semi-positive, or
\item if $v$ is non-Archimedean, the metric of $|\ndot|_{g-\psi_n}$ of $\OO_X(D)$ is a nef model metric.
\end{enumerate}
\end{enumerate}
\end{theo}

\begin{proof}
This theorem is nothing more than \cite[Th\'eor\`eme~3.1]{MR2244803}. In the case where $v$ is non-Archimedean, it is proved under
the additional assumption that $v$ is discrete. However, their proof works well by slight modifications.
For reader's convenience, we reprove it here. 

\medskip
We may assume that $v$ is non-Archimedean.
If the theorem holds for $(mD, mg)$ for some positive number $m$, then
it also holds for $(D, g)$, so that we may assume that 
there is a flat model $(\mathscr X, \mathscr L)$ of $(X, \OO_X(D))$ such that
$|\ndot|_g = |\ndot|_{\varphi_{\mathscr L}}$ and $\mathscr L$ is nef along the special fiber of $\mathscr X \to \Spec(\oo_v)$.
By Lemma~\ref{lem:eq:regular:mero:generic:total}, there is a Cartier divisor $\mathscr D$ on $\mathscr X$
such that $\OO_{\mathscr X}(\mathscr D) = \mathscr L$, $\rest{\mathscr D}{X} = D$ and
$g$ is the Green function arising from $(\mathscr X, \mathscr D)$.
Let $\mathscr X = \bigcup_{i=1}^N \Spec (\mathscr A_i)$ be an affine open covering of $\mathscr X$ such that
$\mathscr D$ is given by a local equation $f_i$ on $\Spec (\mathscr A_i)$.
Since $D$ is effective, one has $f_i \in (\mathscr A_i)_S$, that is, $s_i f_i \in \mathscr A_i$ for some $s_i \in S$,
where $S := \mathfrak o_v \setminus \{ 0 \}$, so that
if we set $s = s_1 \cdots s_N$, then $s f_i \in \mathscr A_i$ for all $i = 1, \ldots, N$.
Let \[g' := g - \log |s|,\quad \mathscr L' := \mathscr L \otimes \OO_{\mathscr X} s^{-1}\quad\text{and}\quad
\mathscr D' := \mathscr D + \operatorname{div}(s).\]
Then $\mathscr D'$ is effective, $\OO_{\mathscr X}(\mathscr D') = \mathscr L'$ and $|\ndot|_{g'} = |\ndot|_{\mathscr L'}$. Thus, if the theorem holds for $g'$, then one has the assertion for $g$, and hence we may further assume that $\mathscr D$ is
effective.

Fix $a \in S$ such that $|a| < 1$, and set
\[ \psi_n = \min \{ g, -n \log |a|\}\quad (\forall\ n \in \NN).\]
The properties (1) and (2) are obvious, so
we need to see (3).
Let $\mathscr I_n$ be the ideal sheaf of $\OO_{\mathscr X}$ generated by a local equation of $\mathscr D$ and $a^n$.
Let $p_n : \mathscr Y_n \to \mathscr X$ be the blowing-up in terms of the ideal sheaf $\mathscr I_n$.
Note that $\mathscr I_n \OO_{\mathscr Y}$ is a locally principal ideal sheaf of $\OO_{\mathscr Y_n}$
whose support is contained in the special fiber of $\mathscr Y_n \to \Spec(\mathfrak o_v)$, that is,
there is an effective Cartier divisor $\mathscr E_n$ on $\mathscr Y_n$ such that
$\OO_{\mathscr Y_n}(-\mathscr E_n) = \mathscr I_n \OO_{\mathscr Y_n}$ and $\rest{\mathscr E_n}{X} = 0$.
Obviously $\psi_n$ is a smooth function arising from the model $(\mathscr Y_n, \mathscr E_n)$.
Therefore, it is sufficient to show that $p_n^*(\mathscr D) - \mathscr E_n$ is nef along the special fiber
$\mathscr Y_n \to \Spec(\mathfrak o_v)$.
Let $\mathscr X = \bigcup_{i=1}^N \Spec (\mathscr A_i)$ be an affine open covering of $\mathscr X$ as before.
Note that
$\mathscr D$ is given by $f_i \in \mathscr A_i$ on $\Spec \mathscr A_i$ for each $i$.
Then \[p_n^{-1}(\Spec \mathscr A_i) = \Proj (\mathscr A_i[T_0, T_1]/(f_i T_0 - a^n T_1)).\]
If we set $p_n^{-1}(\Spec \mathscr A_i)_\alpha = \{ T_\alpha \not = 0 \}$ for $\alpha \in \{ 0, 1 \}$, 
then $f_i = a^n (T_1/T_0)$ on $p_n^{-1}(\Spec \mathscr A_i)_0$ and
$a^n = f_i (T_0/T_1)$ on $p_n^{-1}(\Spec \mathscr A_i)_1$, so that
\begin{equation}
\label{eqn:thm:approximation:smooth:funcs:01}
\begin{cases}
\rest{\OO_{\mathscr Y_n}(-\mathscr E_n)}{p_n^{-1}(\Spec \mathscr A_i)_0} = a^n \OO_{p_n^{-1}(\Spec \mathscr A_i)_0}, \\
\rest{\OO_{\mathscr Y_n}(-\mathscr E_n)}{p_n^{-1}(\Spec \mathscr A_i)_1} = f_i \OO_{p_n^{-1}(\Spec \mathscr A_i)_1}.
\end{cases}
\end{equation}
Therefore, one can see that
$p_n^*(\mathscr D) - \mathscr E_n$ and $\mathrm{div}(a^n) - \mathscr E_n$ are effective.
Let us see $(p_n^*(\mathscr D) - \mathscr E_n \cdot C) \geqslant 0$
for any  irreducible curve $C$ on the special fiber of $\mathscr Y_n \to \Spec(\mathfrak o_v)$.
Let $\xi$ be the generic point of $C$. We choose $i$ such that $\xi \in p_n^{-1}(\Spec \mathscr A_i)$.
If $\xi \not\in \Supp(p_n^*(\mathscr D) - \mathscr E_n)$, then the assertion is obvious because $p_n^*(\mathscr D) - \mathscr E_n$ is effective.
Otherwise, by \eqref{eqn:thm:approximation:smooth:funcs:01}, $\xi \in p_n^{-1}(\Spec \mathscr A_i)_0$.
Then, by \eqref{eqn:thm:approximation:smooth:funcs:01} again, $\xi \not\in \Supp(\mathrm{div}(a^n) - \mathscr E_n)$,
so that $((\mathrm{div}(a^n) - \mathscr E_n) \cdot C) \geqslant 0$ by the reason of
the effectivity of $\mathrm{div}(a^n) - \mathscr E_n$.
Note that $p_n^*(\mathscr D) - \mathscr E_n$ is linearly equivalent to $p_n^*(\mathscr D) + (\mathrm{div}(a^n) - \mathscr E_n)$. Thus it is sufficient to show that $(p_n^*(\mathscr D) \cdot C) \geqslant 0$,
which is obvious because of the projection formula and the nefness of $\mathscr D$.
\end{proof}

\section{Local measures}\label{sec:local:measures}

In this section, we assume that $k$ is \emph{algebraically closed}. Let $X$ be a 
projective $k$-scheme and let $d$ be the dimension of $X$. Assume given a family $(L_i)_{i=1}^d$ of semi-ample invertible $\mathcal O_X$-modules. For any $i\in\{1,\ldots,d\}$, let $\varphi_i$ be a semi-positive continuous metric on $L_i$. 
First we assume that $X$ is integral.
In the case where $|\ndot|$ is Archimedean (and hence $k=\mathbb C$), by Bedford-Taylor theory \cite{MR445006} one can construct a Borel measure
\[c_1(L_1,\varphi_1)\cdots c_1(L_{d},\varphi_d)\]
having
\[\deg(c_1(L_1)\cdots c_1(L_d)\cap[X])\]
as its total mass. In the non-Archimedean case, an analoguous measure has been proposed by Chambert-Loir \cite{MR2244803}, assuming that the field $k$ admits a dense countable subfield (see also \cite[\S5]{ChamDucros} for a general non-Archimedean analogue of Bedford-Taylor theory). 
In any case, the measure $c_1(L_1,\varphi_1)\cdots c_1(L_{d},\varphi_d)$ is often denoted by
$\mu_{(L_1,\varphi_1)\cdots (L_{d},\varphi_d)}$.
Note that the measure $\mu_{(L_1,\varphi_1)\cdots (L_d,\varphi)}$ is additive with respect to each $(L_i,\varphi_i)$. More precisely, if $i\in\{1,\ldots,d\}$ and if $(M_i,\psi_i)$ is another semi-positively metrized invertible $\mathcal O_X$-module, then the measure
\[\mu_{(L_1,\varphi_1)\cdots(L_{i-1},\varphi_{i-1})(L_i\otimes M_i,\varphi_i\otimes\psi_i)(L_{i+1},\varphi_{i+1})\cdots(L_d,\varphi_d)}\]
is equal to
\[\mu_{(L_1,\varphi_1)\cdots(L_d,\varphi_d)}+\mu_{(L_1,\varphi_1)\cdots(L_{i-1},\varphi_{i-1})(M_i,\psi_i)(L_{i+1},\varphi_{i+1})\cdots(L_d,\varphi_d)}\]
Moreover, for any permutation $\sigma:\{1,\ldots,d\}\rightarrow\{1,\ldots,d\}$, one has
\[\mu_{(L_{\sigma(1)},\varphi_{\sigma(1)})\cdots(L_{\sigma(d)},\varphi_{\sigma(d)})}=\mu_{(L_1,\varphi_1)\cdots(L_d,\varphi_d)}.\]

In general, let $X_1, \ldots, X_n$ be irreducible components of $X$ which are of dimension $d$, and $\eta_1, \ldots, \eta_n$ the generic points of $X_1, \ldots, X_n$, respectively.
Let $\xi_i : X_i \hookrightarrow X$ be the canonical closed embedding for each $i$.
Then a measure $\mu_{(L_1,\varphi_1)\cdots (L_{d},\varphi_d)}$ on $X^{\mathrm{an}}$ is defined to be
\begin{multline}\label{def:measure:assocated:invertible:sheaves}
\mu_{(L_1,\varphi_1)\cdots (L_{d},\varphi_d)} := \\
\sum_{j=1}^n \length_{\OO_{X,\eta_j}}(\OO_{X,\eta_j})(\xi^{\mathrm{an}}_j)_*\Big(c_1\big(\xi_j^*(L_1,\varphi_1)\big)\cdots c_1\big(\xi_j^*(L_{d},\varphi_d)\big)\Big).
\end{multline}

\begin{defi}
Let $(L_1,\varphi_1),\ldots,(L_d,\varphi_d)$ be a family of integrable metrized invertible $\mathcal O_X$-modules. For each $i\in\{1,\ldots,d\}$, we let $(L_i',\varphi_i')$ and $(L_i'',\varphi_i'')$ be ample invertible $\mathcal O_X$-modules equipped with semi-positive metrics, such that $L_i=L_i'\otimes (L_i'')^\vee$ and $\varphi_i=\varphi_i'\otimes(\varphi_i'')^\vee$. We define a signed Radon measure $\mu_{(L_1,\varphi_1)\cdots(L_d,\varphi_d)}$ on $X^{\mathrm{an}}$ as follows:
\[\mu_{(L_1,\varphi_1)\cdots(L_d,\varphi_d)}:=\sum_{I\subseteq\{1,\ldots,d\}}(-1)^{\operatorname{card}(I)}\mu_{(L_{1,I},\varphi_{1,I})\cdots(L_{d,I},\varphi_{d,I})},\]
where $(L_{j,I},\varphi_{j,I})=(L_j'',\varphi_j'')$ if $j\in I$, and $(L_{j,I},\varphi_{j,I})=(L_j',\varphi_j')$ if $j\in\{1,\ldots,d\}\setminus I$. 
\end{defi}

\begin{exem}\label{Exe:chambert-loir measure}We recall the explicit construction of Chambert-Loir's measure in a particular case as explained in \cite[\S2.3]{MR2244803}. 
Assume that the absolute value $|\ndot|$ is non-Archimedean and that  the $k$-scheme $X$ is integral and normal. Let $k^\circ$ be the valuation ring of $(k,|\ndot|)$ and $\mathfrak m$ be the maximal ideal of $k^\circ$. Suppose given an integral model of $X$, namely, a flat and normal projective $k^\circ$-scheme $\mathscr X$ such that \[\mathscr X\times_{\Spec k^\circ}\Spec k\cong X.\]
Let $\mathscr X_{\mathfrak m}$ be the fibre of $\mathscr X$ over the closed point of $\Spec k^\circ$. It turns out that the reduction map from $X^{\mathrm{an}}$ to $\mathscr X_{\mathfrak m}$ is surjective. Let $Z_1,\ldots,Z_n$ be irreducible components of $\mathscr X_{\mathfrak m}$. For any $j\in\{1,\ldots,n\}$, there exists a unique point $z_j\in X^{\mathrm{an}}$ whose reduction identifies with the generic point of $Z_j$.

Assume that each metric $\varphi_j$ is induced by an integral model $\mathscr L_i$, which is an invertible sheaf on $\mathscr X$ such that $\mathscr L_i|_X\cong L_i$. Then the measure
\[c_1(L_1,\varphi_1)\cdots c_1(L_{d},\varphi_d)\]
is given by
\[\sum_{j=1}^d \operatorname{mult}_{Z_j}(\mathscr X_{\mathfrak m})\deg(c_1(\mathscr L_1|_{\mathscr X_{\mathfrak m}})\cdots c_1(\mathscr L_{d}|_{\mathscr X_{\mathfrak m}})\cap[Z_j])\operatorname{Dirac}_{z_j},\]  
where $\operatorname{mult}_{Z_j}(\mathscr X_{\mathfrak m})$ is the multiplicity of $Z_j$ in $\mathscr X_{\mathfrak m}$, and $\operatorname{Dirac}_{z_j}$ denotes the Dirac measure at $z_j$.

\end{exem}

\begin{rema}\label{rem:approximation:CL:measure}
We assume that $X$ is integral. Let 
$\{ \varphi_{1, n} \}_{n=1}^{\infty}, \ldots, \{ \varphi_{d, n} \}_{n=1}^{\infty}$ be sequences of semi-positive
metrics of $L_1, \ldots, L_d$, respectively such that \[\lim\limits_{n\to\infty} d(\varphi_{i, n}, \varphi_i) = 0\]
for all $i=1, \ldots, d$. Then, by using \cite[Corollary~(3.6)]{Demagbook} and \cite[Corollaire~(5.6.5)]{ChamDucros}, one can see
\[
\lim_{n\to\infty} \int_{X^{\an}} f \mu_{(L_1, \varphi_{1, n}) \cdots (L_d, \varphi_{d, n})} =
\int_{X^{\an}} f \mu_{(L_1, \varphi_{1}) \cdots (L_d, \varphi_{d})}
\]
for any smooth function $f$ on $X^{\an}$.
\end{rema}

\begin{defi}
Let $\overline D_1=(D_1,g_1),\ldots,\overline D_d=(D_d,g_d)$ be a family of  integrable metrized Cartier divisors on $X$. For any $i\in\{1,\ldots,d\}$, we write $(D_i,g_i)$ as the difference of two metrized Cartier divisors $(D_i',g_i')-(D_i'',g_i'')$, where $D_i'$ and $D_i''$ are ample, and $g_i'$ and $g_i''$ are plurisubharmonic. We define a signed Radon measure $\mu_{\overline D_1\cdots \overline D_d}$ on $X^{\mathrm{an}}$ 
to be
\begin{equation*}\mu_{\overline D_1\cdots\overline D_d}
:=\sum_{I\subseteq\{1,\ldots,d\}}(-1)^{\operatorname{card}(I)}\mu_{\overline{D}_{1,I}\cdots\overline D_{d,I}},
\end{equation*}
where $\overline D_{j,I}=(D_i'',g_{i}'')$ if $j\in I$, and $\overline D_{j,I}=(D_i',g_i')$ if $j\in\{1,\ldots,d\}\setminus I$.
Note that this signed measure does not depend on the choice of the decompositions.

Let $X_1, \ldots, X_n$ be irreducible components of $X$ and $\eta_1, \ldots, \eta_n$ be the generic points
of $X_1, \ldots, X_n$, respectively.
Let $\xi_j : X_j \hookrightarrow X$ be the canonical closed embedding.
Then it is easy to see
\begin{equation}\label{eqn:measure:non:reduced:scheme}
\mu_{(D_1, g_1) \cdots (D_d, g_d)} = \sum_{j=1}^n \length_{\OO_{X, \eta_j}}(\OO_{X, \eta_j})
(\xi_j^{\mathrm{an}})_*\Big(\mu_{\xi_j^*(D_1, g_1) \cdots \xi_j^*(D_d, g_d)}\Big).
\end{equation}
\end{defi}

\begin{prop}\label{prop:CH:measure:projection:formula}
Let $\pi : Y \to X$ be a surjective morphism between integral projective schemes over $k$.
We set $e = \dim X$ and $d = \dim Y$. 
Let $(L_1, \varphi_1), \ldots, (L_d, \varphi_d)$ be integrable metrized invertible $\OO_X$-modules.
Then one has the following:
\begin{enumerate}[label=\rm(\arabic*)]
\item If $d > e$, then $\pi_*(\mu_{\pi^*(L_1, \varphi_1) \cdots \pi^*(L_d, \varphi_d)}) = 0$.
\item If $d = e$, then 
$\pi_*(\mu_{\pi^*(L_1, \varphi_1) \cdots \pi^*(L_d, \varphi_d)}) = (\deg \pi)\mu_{(L_0, \varphi_0) \cdots (L_d, \varphi_d)}$.
\end{enumerate}
\end{prop}

\begin{proof}
We may assume that $L_1, \ldots, L_d$ are ample and $\varphi_1, \ldots, \varphi_d$ are semi-positive.
If $\varphi_1, \ldots, \varphi_d$ are smooth, then the assertion is well-known (cf. \cite[Proposition~10.4]{GKTropical}).
Let $\{ \varphi_{1,n} \}_{n=1}^{\infty}, \ldots, \{ \varphi_{d,n} \}_{n=1}^{\infty}$ be regularizations of
$\varphi_1, \ldots, \varphi_d$, that is, $\varphi_{1,n}, \ldots, \varphi_{d,n}$ are smooth and
semi-positive for $i=1, \ldots, d$ and $n \geqslant 1$, and
$\lim\limits_{n\to\infty} d(\varphi_{i}, \varphi_{i, n}) = 0$ for $i=1, \ldots, d$
(for example, see \cite{MR3028754} for the Archimedean case and Theorem~\ref{thm:equiv:semi-positive:M:semi-positive} for the non-Archimedean case).
Let $f$ be a smooth function on $X^{\an}$. Then, by using \cite[Corollary~(3.6)]{Demagbook} and \cite[Corollaire~(5.6.5)]{ChamDucros}, one can see that
\begin{align*}
& \lim_{n\to\infty} \int_{X^{\an}} \pi^*(f) \mu_{\pi^*(L_{1}, \varphi_{1,n}) \cdots \pi^*(L_d, \varphi_{d,n})}
= \int_{X^{\an}} \pi^*(f) \mu_{\pi^*(L_{1}, \varphi_{1}) \cdots \pi^*(L_d, \varphi_{d})} \\
\intertext{and if $d = e$, then}
& \lim_{n\to\infty} \int_{Y^{\an}} f \mu_{(L_{1}, \varphi_{1,n}) \cdots (L_d, \varphi_{d,n})}
= \int_{X^{\an}} f \mu_{(L_{1}, \varphi_{1}) \cdots (L_d, \varphi_{d})}.
\end{align*}
Thus the assertions follow.
\end{proof}

\begin{rema}Let $X$ and $Y$ be two projective schemes over $\Spec k$, of Krull dimension $d$ and $n$, respectively. Let $\overline L_1,\ldots,\overline L_d$ be integrable metrized invertible $\mathcal O_X$-modules, $\overline M_1,\ldots,\overline M_d$ be integrable $\mathcal O_Y$-modules. We consider the fiber product $X\times_kY$ and let $\pi_1:X\times_kY\rightarrow X$ and $\pi_2:X\times_kY\rightarrow Y$ be the two morphisms of projection.   In the case where $k$ is Archimedean, the analytic space $(X\times_k Y)^{\mathrm{an}}$ is homeomorphic to $X^{\mathrm{an}}\times Y^{\mathrm{an}}$ and the measure \[\mu_{\pi_1^*(\overline L_1)\cdots\pi_1^*(\overline L_d)\pi_2^*(\overline M_1)\cdots\pi_2^*(\overline M_n)}\]
on $(X\times_kY)^{\operatorname{an}}$ identifies with
\[\mu_{\overline L_1\cdots\overline L_d}\otimes\mu_{\overline M_1\cdots\overline M_n}.\]
In the case where $|\ndot|$ is non-Archimedean, in general the topological space $(X\times_kY)^{\mathrm{an}}$ is not homeomorphic to $X^{\operatorname{an}}\times Y^{\operatorname{an}}$. However, there is a natural continuous map \[\alpha:(X\times_kY)^{\operatorname{an}}\longrightarrow X^{\mathrm{an}}\times Y^{\mathrm{an}}.\] Then the following equality holds (see \cite[\S2.8]{MR2244803})
\begin{equation*}
\alpha_*\Big(\mu_{\pi_1^*(\overline L_1)\cdots\pi_1^*(\overline L_d)\pi_2^*(\overline M_1)\cdots\pi_2^*(\overline M_n)}\Big)=\mu_{\overline L_1\cdots\overline L_d}\otimes\mu_{\overline M_1\cdots\overline M_n}.
\end{equation*} 
In particular, if $g$ is a measurable function on $Y^{\mathrm{an}}$ which is integrable with respect to $\mu_{\overline M_1\cdots\overline M_n}$, one has
\begin{equation}\label{Equ:Fubini}\int_{(X\times_kY)^{\mathrm{an}}}(g\circ\pi_2^{\operatorname{an}})\,\mathrm{d}\mu_{\pi_1^*(\overline L_1)\cdots\pi_1^*(\overline L_d)\pi_2^*(\overline M_1)\cdots\pi_2^*(\overline M_n)}=\int_{Y^{\mathrm{an}}}g\,\mathrm{d}\mu_{\overline M_1\cdots\overline M_n}.\end{equation}
\end{rema}

\begin{defi}
Let $E$ be a finite-dimensional vector space over $k$. We say that a norm $\|\ndot\|$ on $E$ is \emph{orthonormally decomposable} if 
\begin{enumerate}[label=\rm(\arabic*)]
\item in the case where $|\ndot|$ is non-Archimedean, the norm $\|\ndot\|$ is ultrametric, and $(E,\|\ndot\|)$ admits an orthonormal basis $(e_j)_{j=0}^{r}$, namely, \[\forall\,(\lambda_j)_{j=0}^{r}\in k^{r+1},\quad \|\lambda_0e_0+\cdots+\lambda_{r}e_{r}\|=\max_{j\in\{0,\ldots,r\}}|\lambda_j|;\]
\item in the case where $|\ndot|$ is Archimedean, the norm $\|\ndot\|$ is induced by an inner product $\mathopen{\langle}\cdot,\cdot\mathclose{\rangle}$. 
\end{enumerate}
Note that for each valued extension $(k',|\ndot|')$ of $(k,|\ndot|)$, there is a unique norm $\|\ndot\|_{k'}$ on $E\otimes_kk'$, which is either ultrametric or induced by an inner product, such that any orthonormal basis of $(E,\|\ndot\|)$ is also an orthonormal basis of the extended normed vector space $(E\otimes_kk',\|\ndot\|_{k'})$ (see Definition \ref{Def: extension of sclars}). 
\end{defi}

\begin{rema}\label{Rem: approximation by good norms}
Let $E$ be a finite-dimensional vector space over $k$, and $\|\ndot\|$ be an orthonormally decomposable norm on $E$. For any $s\in E$, the real number $\|s\|$ belongs to the image of the absolute value $|\ndot|$. In particular, if $s$ is non-zero, then there exists $\lambda\in k$ such that $\|\lambda s\|=1$. 

In the case where the absolute value $|\ndot|$ is non-Archimedean, it is \emph{not} true that any ultrametrically normed vector space admits an orthonormal basis (see \cite[Example 2.3.26]{MR2598517}). However, if $(E,\|\ndot\|)$ is a finite-dimensional ultrametrically normed vector space over $k$, for any $\alpha\in\mathbb R$ such that $0<\alpha<1$, there exists an $\alpha$-orthogonal basis of $E$ (cf. \cite[\S 2.3]{MR2598517}, see also \cite[\S1.2.6]{CMArakelovAdelic} for details), namely a basis $(e_i)_{i=1}^r$ such that, for any $(\lambda_i)_{i=1}^r\in k^r$, 
\[ \alpha\max_{i\in\{1,\ldots,r\}}|\lambda_i|\cdot\|e_i\|\leqslant\|\lambda_1e_1+\cdots+\lambda_re_r\|\leqslant\max_{i\in\{1,\ldots,r\}}|a_i|\cdot\|e_i\|.\]
Moreover, since $k$ is assumed to be algebraically closed, in the case where absolute value $|\ndot|$ is non-trivial, the image of $|\ndot|$ is dense in $\mathbb R$. In fact, if $a$ is an element of $k$ such that $|a|\neq 1$, for any non-zero rational number $p/q$ with $p\in\mathbb Z$ and $q\in\mathbb Z_{>0}$, any element $x\in k$ satisfying the polynomial equation
\[x^q=a^p\]
has $|a|^{p/q}$ as absolute value. Therefore, by possibly delating the vectors $(e_i)_{i=1}^r$ we may assume that 
\[\alpha\leqslant\|e_i\|\leqslant 1\]
for any $i\in\{1,\ldots,r\}$. Therefore, if we denote by $\|\ndot\|_\alpha$ the norm on $E$ under which $(e_i)_{i=1}^r$ is an orthonormal basis of $E$, then for any $x=\lambda_1e_1+\cdots+\lambda_re_r$ in $E$, one has
\[\|x\|_\alpha=\max_{i\in\{1,\ldots,r\}}|\lambda_i|\leqslant\alpha^{-1}\max_{i\in\{1,\ldots,r\}}|\lambda_i|\cdot\|e_i\|\leqslant\alpha^{-2}\|x\|,\]
and
\[\|x\|\leqslant\max_{i\in\{1,\ldots,r\}}|\lambda|\cdot\|e_i\|\leqslant\max_{i\in\{1,\ldots,r\}}|\lambda_i|=\|x\|_\alpha.\]
Therefore, one has
\[d(\|\ndot\|_\alpha,\|\ndot\|):=\sup_{x\in E\setminus\{0\}}\Big|\ln\|x\|_\alpha-\ln\|x\|\Big|\leqslant -2\ln(\alpha).\]
Thus we can approximate the  ultrametric norm $\|\ndot\|$ by a sequence of ultrametric norms which are orthonormally decomposable.
\end{rema}

\begin{prop}\label{Pro: decomposable quotient norm}Let $(E,\|\ndot\|)$ be a finite-dimensional vector space over $k$, equipped with an orthonormally decomposable norm. Then any element $s_0\in E$ such that $\|s_0\|=1$ belongs to an orthonormal basis. Moreover, for any quotient vector space $G$ of $E$, the quotient norm on $F$ is orthonormally decomposable.
\end{prop}
\begin{proof}
The statement is classic when $|\ndot|$ is Archimedean, which follows from the Gram-Schmidt process. In the following, we assume that $|\ndot|$ is non-Archimedean. Let $k^\circ$ be the valuation ring of $(k,|\ndot|)$.

Let $(e_j)_{j=0}^r$ be an orthonormal basis of $(E,\|\ndot\|)$. Without loss of generality, we may assume that $s_0=\lambda_0e_0+\cdots+\lambda_re_r$ with $(\lambda_0,\ldots,\lambda_r)\in (k^\circ)^{r+1}$ and $|\lambda_0|=1$. We then construct an upper triangular matrix $A$ of size $(r+1)\times(r+1)$, such that the first row $A$ is $(\lambda_0,\ldots,\lambda_r)$ and the diagonal coordinates of $A$ are elements of absolute value $1$ in $k$. Then the matrix $A$ belongs to $\operatorname{GL}_{r+1}(k^\circ)$. Let $(s_j)_{j=0}^r$ be the basis of $E$ such that
\[(s_0,\ldots,s_r)^{T}=A(e_0,\ldots,e_r)^T.\]
For any $j\in\{0,\ldots,r\}$, one has $\|s_j\|=1$. Moreover, for any $(b_0,\ldots,b_r)\in k^r$, one has
\[b_0s_0+\cdots+b_rs_r=(b_0,\ldots, b_r)A(e_0,\ldots,e_r)^T.\]
Let $(c_0,\ldots,c_r)=(b_0,\ldots,b_r)A$. Since $(e_0,\ldots,e_r)$ is an orthonormal basis, one has 
\[\|b_0s_0+\cdots+b_rs_r\|=\max_{j\in\{0,\ldots,r\}}|c_j|.\]
Note that $(b_0,\ldots,b_r)=(c_0,\ldots,c_r)A^{-1}$. Since $A^{-1}$ belongs to $\operatorname{GL}_{r+1}(k^\circ)$, one has
\[\forall\,i\in\{0,\ldots,r\},\quad |b_i|\leqslant\max_{j\in\{0,\ldots,r\}}|c_j|.\]
Therefore one obtains 
\[\|b_0s_0+\cdots+b_rs_r\|\geqslant\max_{i\in\{0,\ldots,r\}}|b_i|.\]
Combined with the strong triangle inequality, we obtain
\[\|b_0s_0+\cdots+b_rs_r\|=\max_{i\in\{0,\ldots,r\}}|b_i|.\]
Therefore $(s_j)_{j=0}^r$ is an orthonormal basis of $(E,\|\ndot\|)$. In particular, the image of $(s_1,\ldots,s_{r})$ in $E/ks_0$ forms an orthonormal basis of $E/ks_0$ with respect to $\|\ndot\|$. Therefore the quotient norm on $E/ks_0$ is orthonormally decomposable. By induction  we can show that all quotient norms of $\|\ndot\|$ are orthonormally decomposable.
\end{proof}

In the remaining of this section, we fix a finite-dimensional vector space $E$ equipped with an orthonormally decomposable norm $\|\ndot\|$. We also choose an orthonormal basis $(e_j)_{j=1}^r$ of $(E,\|\ndot\|)$. Let $\mathbb P(E)$ be the projective space of $E$ and $\mathcal O_E(1)$ be the universal invertible sheaf on $\mathbb P(E)$. We equip $\mathcal O_E(1)$ with the orthogonal quotient metric $(|\ndot|(x))_{x\in\mathbb P(E)^{\mathrm{an}}}$ (see Definition \ref{Def: Fubini-Study metric}) and denote by $\overline{\mathcal O_E(1)}$ the corresponding metrized invertible sheaf. Recall that each point $x\in\mathbb P(E)^{\mathrm{an}}$ corresponds to a one-dimensional quotient vector space 
\[E\otimes_K\widehat{\kappa}(x)\longrightarrow\mathcal O_E(1)(x),\]
where $\widehat{\kappa}(x)$ denotes the completed residue field of $x$. Then the norm $|\ndot|(x)$ on $\mathcal O_E(1)(x)$ is by definition the quotient norm of $\|\ndot\|_{\widehat{\kappa}(x)}$. 

\begin{defi}Assume that $|\ndot|$ is non-Archimedean. We denote by $\xi$ the point in $\mathbb P(E)^{\mathrm{an}}$ which is the generic point of $\mathbb P(E)^{\mathrm{an}}$ equipped with the absolute value \[\textstyle|\ndot|_\xi:k\big(\frac{e_0}{e_r},\ldots,\frac{e_{r-1}}{e_r}\big)\longrightarrow\mathbb R_{\geqslant 0}\]
such that, for any  
\[ P=\sum_{\boldsymbol{a}=(a_0,\ldots,a_{r-1})\in\mathbb N^d}\lambda_{\boldsymbol{a}}\Big(\frac{e_0}{e_r}\Big)^{a_0}\cdots\Big(\frac{e_{r-1}}{e_r}\Big)^{a_{r-1}}\in \textstyle{k\big[\frac{e_0}{e_r},\ldots,\frac{e_{r-1}}{e_r}\big]},\]
one has 
\[|P|_\xi=\max_{\boldsymbol{a}\in\mathbb N^d}|\lambda_{\boldsymbol{a}}|.\]
Note that the point $\xi$ does not depend on the choice of the orthonormal basis $(e_j)_{j=0}^r$. In fact, the norm $\|\ndot\|$ induces a symmetric algebra norm on $k[E]$ (which is often called a \emph{Gauss norm}) and hence defines an absolute value on the fraction field of $k[E]$. The restriction of this absolute value to the field of rational functions on $\mathbb P(E)$ identifies with $|\ndot|_\xi$. Hence $\xi$ is called the \emph{Gauss point} of $\mathbb P(E)^{\mathrm{an}}$.
\end{defi}  

\begin{prop}\label{Pro: non Archimedean projective space} Assume that the absolute value $|\ndot|$ is non-Archimedean. 
The following equality holds
\[c_1(\overline{\mathcal O_E(1)})^r=\operatorname{Dirac}_\xi,\]
where $\operatorname{Dirac}_\xi$ denotes the Dirac measure at $\xi$.
\end{prop}
\begin{proof} Let $k^\circ$ be the valuation ring of $(k,|\ndot|)$, $\mathfrak m$ be the maximal ideal of $k^\circ$, and $\kappa=k^\circ/\mathfrak m$ be the residue field of $k^\circ$. 
Let $\mathcal E$ be the free $k^\circ$-module generated by $\{e_0,\ldots,e_r\}$. Then $\mathbb P(\mathcal E)$ is a projective flat $k^\circ$-scheme such that
\[\mathbb P(\mathcal E)\times_{\Spec k^\circ}\Spec k\cong\mathbb P(E).\]
Note that the fibre product
\[\mathbb P(\mathcal E)\times_{\Spec k^\circ}\Spec\kappa\]
is isomorphic to $\mathbb P(\mathcal E\otimes_{k^\circ}\kappa)$, which is an integral $\kappa$-scheme. Therefore, one has (see Example \ref{Exe:chambert-loir measure})
\[c_1(\overline{\mathcal O_E(1)})^r=\deg(c_1(\mathcal O_{\mathcal E_{\kappa}}(1))^r\cap\mathbb P(\mathcal E_{\kappa}))\operatorname{Dirac}_\xi=\operatorname{Dirac}_\xi.\]
\end{proof}

\begin{rema} \label{Rem: mesure complexe} Assume that $k=\mathbb C$ and $|\ndot|$ is the usual absolute value. Let $(E,\|\ndot\|)$ be a Hermitian space and \[\mathbb{S}(E^\vee,\|\ndot\|_*)=\{\alpha\in E^\vee\,|\,\|\alpha\|_*=1\}\] be the unit sphere in $E$, where $\|\ndot\|_*$ denotes the dual norm of $\|\ndot\|$, which is also a Hermitian norm. Note that $\mathbb P(E)^{\mathrm{an}}$ identifies with the quotient of $\mathbb S(E^\vee,\|\ndot\|_*)$ by the action of the unit sphere $\mathbb S(\mathbb C)=\{z\in\mathbb C\,|\,|z|=1\}$ in $\mathbb C$. We equip the universal invertible sheaf $\mathcal O_E(1)$ with the orthogonal quotient metric induced by $\|\ndot\|$ and equip $\mathbb S(E^\vee,\|\ndot\|_*)$ with the unique $U(E^\vee,\|\ndot\|_*)$-invariant Borel probability measure
$\eta_{\mathbb S(E^\vee,\|\ndot\|_*)}$ which is locally equivalent to Lebesgue measure. Then the measure \[c_1(\overline{\mathcal O_E(1)})^{\dim_{\mathbb C}(E)-1}\] identifies with the direct image of $\eta_{\mathbb S(E^\vee,\|\ndot\|_*)}$  by the projection map from $\mathbb S(E^\vee,\|\ndot\|_*)$ to $\mathbb P(E)^{\mathrm{an}}$ (see for example \cite[(1.4.7)]{MR1260106} for more details).
\end{rema}

\begin{theo}\label{thm:integrability:Green:func}
Let $\overline{L} = (L, \varphi),\overline{L}_1 = (L_1, \varphi_1), \ldots, \overline{L}_d = (L_d, \varphi_d)$
be integrable metrized invertible $\OO_X$-modules.
Let $s$ be a regular meromorphic section of $L$.
Then $g =-\log |s|_\varphi$ is integrable with respect to $\mu_{\overline{L}_1 \cdots \overline{L}_d}$.
\end{theo}

\begin{proof}
The proof of this theorem is same as \cite[Th\'eor\`eme~4.1]{MR2244803}.
We prove it without using the local intersection numbers.

Clearly we may assume that $X$ is integral, $L, L_1, \ldots, L_d$ are ample and
$\overline{L}, \overline{L}_1, \ldots, \overline{L}_d$ are semi-positive.
Let $\mathcal I$ be the ideal sheaf of $\OO_X$ given by
\[\mathcal I_x = \{ a \in \OO_{X,x} \mid a s_x \in L_x \}.\]
Choose a positive number $m$ and a non-zero section $t_1 \in H^0(X, \mathcal{I} L^{\otimes m}) \setminus \{ 0 \}$.
If we set $t_2 = t_1 \otimes s$, 
then $s = t_2 \otimes t_1^{-1}$ and $t_2 \in H^0(X, L^{\otimes m+1}) \setminus \{ 0 \}$ and
$g = -\log |t_2|_{(m+1)\varphi} + \log |t_1|_{m \varphi}$, so that we may assume that $s \in H^0(X, L)\setminus \{ 0 \}$.
Let $\varphi'$ be a metric of $L$ such that either
(a) if $v$ is Archimedean, $\varphi'$ is $C^{\infty}$ and semi-positive, or
(b) if $v$ is non-Archimedean, $\varphi'$ is a nef model metric.
Then $-\log |s|_{\varphi} + \log |s|_{\varphi'}$ is a continuous function, so that
we may assume that $\varphi = \varphi'$. 
By Theorem~\ref{thm:approximation:smooth:funcs}, 
there is a sequence $\{ \psi_n \}_{n\in\mathbb N}$ of smooth functions on $X^{\an}$ with the following properties:
\begin{enumerate}[label=\rm(\arabic*)]
\item for all $n \in \NN$, $\psi_n \leqslant g$, $\psi_n \leqslant \psi_{n+1}$.
\item for each point $x \in X^{\an}$, $\sup \{ \psi_n(x) ; n \in \NN \} = g(x)$.

\item for all $n \in \NN$, $g - \psi_n$ is a Green function of $D$ such that either
\begin{enumerate}[label=\rm(\arabic{enumi}.\alph*)]
\item if $v$ is Archimedean, the metric of $|\ndot|_{g-\psi_n}$ of $L$ is $C^{\infty}$ and semi-positive, or
\item if $v$ is non-Archimedean, the metric of $|\ndot|_{g-\psi_n}$ of $L$ is a nef model metric.
\end{enumerate}
\end{enumerate}
We prove the assertion by induction on the number
\[
e := \operatorname{Card}\{ i \in \{ 1, \ldots, d \} \mid \text{$\varphi_i$ is not smooth} \}.
\]
If $e = 0$, that is, $\varphi_i$ is smooth for all $i$, then the assertion is obvious.
We assume that $e > 0$. Obviously we may assume that $\varphi_1$ is not smooth.
Let $\varphi'_1$ be a semi-positive and smooth metric of $L_1$.
If we choose a continuous function $\vartheta$ such that $|\ndot|_{\varphi_1} = \exp(-\vartheta)|\ndot|_{\varphi'_1}$, then $c_1(\overline{L}_1) = c_1(\overline{L}'_1) + dd^c(\vartheta)$, where
$\overline{L}'_1 = (L_1,\varphi'_1)$.

Let us consider the following integral:
\[
I_n := \int_{X^{\an}} \psi_n c_1(\overline{L}_1) \cdots c_1(\overline{L}_d).
\]
Note that $\psi_n$ and $\vartheta$ are locally written by differences of plurisubharmonic functions, so that,
by \cite[Proposition~2.3]{MR2244803},
\begin{align*}
I_n & = \int_{X^{\an}} \psi_n c_1(\overline{L}'_1) \cdots c_1(\overline{L}_d) + \int_{X^{\an}} \psi_n dd^c(\vartheta) c_1(\overline{L}_2)\cdots c_1(\overline{L}_d) \\
 & = \int_{X^{\an}} \psi_n c_1(\overline{L}'_1) \cdots c_1(\overline{L}_d) + \int_{X^{\an}} \vartheta dd^c(\psi_n) c_1(\overline{L}_2) \cdots c_1(\overline{L}_d).
\end{align*}
By the hypothesis of induction, 
\[
\lim_{n\to\infty} \int_{X^{\an}} \psi_n c_1(\overline{L}'_1) \cdots c_1(\overline{L}_d)
\]
exists. Moreover, by the same arguments as the last part of \cite[Th\'eor\`em~4.1]{MR2244803},
one can see
\begin{multline*}
\lim_{n\to\infty}\int_{X^{\an}} \vartheta dd^c(\psi_n) c_1(\overline{L}_2) \cdots c_1(\overline{L}_d) \\
= \int_{X^{\an}} \vartheta c_1(\overline{L}) c_1(\overline{L}_2) \cdots c_1(\overline{L}_d) -
\int_{\operatorname{div}(s)^{\an}} \vartheta c_1(\overline{L}_2) \cdots c_1(\overline{L}_d).
\end{multline*}
Therefore $\lim_{n\to\infty} I_n$ exists, as required.
\end{proof}

\section{Local intersection number over an algebraically closed field} 

Let $k$ be an algebraically closed field equipped with a non-trivial absolute value $|\ndot|$ such that $k$ is complete with respect to the topology defined by $|\ndot|$.
The pair $(k, |\ndot|)$ is denoted by $v$.
Let $X$ be a projective scheme over $k$ and $d$ be its dimension. 
Recall that any element $x$ of $X^{\mathrm{an}}$ consists of a scheme point of $X$ and an absolute value $|\ndot|_x$ of the residue field of the scheme point. We denote by $\widehat{\kappa}(x)$ the completion of the residue field of the scheme point with respect to the absolute value $|\ndot|_x$, on which the absolute value extends by continuity.

\begin{defi}\label{def:local:intersection}

Let $(D_0, g_0), \ldots, (D_d, g_d)$ be integrable metrized Cartier divisors on $X$.
We assume that $D_0, \ldots, D_d$ intersect properly, that is, $(D_0, \ldots, D_d) \in \IP_X$ (see Definition~\ref{def:meet:properly}).
According to \cite{MR2244803}, we define \emph{the local intersection number $\big((D_0, g_0) \cdots (D_d, g_d)\big)_v$ at $v$} as follows.

In the case where $d = 0$, one has $X=\Spec(A)$ for some $k$-algebra with $\dim_k(A) < \infty$.
By Remark~\ref{remark:Cartier:div:on:0:dim} and Example~\ref{Rem: Green function of several points},
\[A = \bigoplus_{x \in \Spec(A)} A_x\;\text{ and }\;(D_0, g_0) = \sum_{x \in \Spec(A)} (0, a_x),\]
where $a_x \in \mathbb{R}$ for all $x \in \Spec(A)$.
Then 
\begin{equation}\label{eqn:def:local:intersection:dim:0}
\big((D_0, g_0)\big)_v := \sum_{x \in \Spec(A)} \length_{A_x}(A_x)\,a_x.
\end{equation}
Note that $\length_{A_x}(A_x) = \dim_k (A_x)$ because $k$ is algebraically closed.

If $d > 0$ and $\sum_{i=1}^n a_i Z_i$ is the cycle associated with $D_d$ (cf. Remark~\ref{remark:expansion:Cartier:div:as:cycle}), then the local intersection number $\big((D_0, g_0) \cdots (D_d, g_d)\big)_v$ is defined in a recursive way with respect to $d=\dim(X)$ as
\begin{equation}\label{eqn:def:local:intersection}
\sum_{i=1}^n  a_i \Big(\rest{(D_0, g_0)}{Z_i} \cdots \rest{(D_{d-1}, g_{d-1})}{Z_i}\Big)_v + \int_{X^{\an}} g_{d}(x)\, \mu_{(D_0, g_0) \cdots (D_{d-1}, g_{d-1})}(\mathrm{d}x).
\end{equation} 
For the integrability of $g_d$ with respect to the measure $\mu_{(D_0, g_0) \cdots (D_{d-1}, g_{d-1})}$, see Theorem~\ref{thm:integrability:Green:func}.

\end{defi}

\begin{prop}\label{prop:formula:local:intersection:irreducible:components}
Let $X_1, \ldots, X_\ell$ be irreducible components of $X$ and $\eta_1, \ldots, \eta_\ell$ be the generic points of $X_1, \ldots, X_\ell$, respectively.
Then
\[
\big((D_0, g_0) \cdots (D_d, g_d)\big)_v = \sum_{j=1}^{\ell} \operatorname{length}_{\mathcal O_{X, \eta_j}}(\mathcal O_{X, \eta_j}) \big(\rest{(D_0, g_0)}{X_j} \cdots \rest{(D_d, g_d)}{X_j}\big)_v.
\]
\end{prop}

\begin{proof}
In the case where $d = 0$, the assertion is obvious.
We assume that $d > 0$. By the definition of $\mu_{(D_0, g_0) \cdots (D_{d-1}, g_{d-1})}$ (cf. Section~\ref{sec:local:measures}), if we set
\[b_j = \operatorname{length}_{\mathcal O_{X, \eta_j}}(\mathcal O_{X, \eta_j}),\] then one has
\[
\int_{X^{\mathrm{an}}} g_d(x)\, \mu_{(D_0, g_0) \cdots (D_{d-1}, g_{d-1})}(\mathrm{d}x) \\
= \sum_{j=1}^{\ell}  b_j \int_{X_j^{\mathrm{an}}}
g_d(x)\, \mu_{\rest{(D_0, g_0)}{X_j} \cdots \rest{(D_{d-1}, g_{d-1})}{X_j}}(\mathrm{d}x).
\]
If $\sum_{i=1}^n a_i Z_i$ and $\sum_{i=1}^n a_{ji} Z_i$ are the cycles associated with $D_d$ and $\rest{D_d}{X_j}$, respectively,
then,  by \eqref{eqn:remark:expansion:Cartier:div:as:cycle:01}, 
$a_i = \sum_{j=1}^{\ell} b_j a_{ji}$, so that
\begin{multline*}
\sum_{i=1}^n a_i \Big(\rest{(D_0, g_0)}{Z_i} \cdots \rest{(D_{d-1}, g_{d-1})}{Z_i}\Big)_v \\
= \sum_{i=1}^n  \sum_{j=1}^{\ell} b_j a_{ji} \Big(\rest{(D_0, g_0)}{Z_i} \cdots \rest{(D_{d-1}, g_{d-1})}{Z_i}\Big)_v \\
= \sum_{j=1}^{\ell} b_j \sum_{i=1}^n  a_{ji}\Big(\rest{(D_0, g_0)}{Z_i} \cdots \rest{(D_{d-1}, g_{d-1})}{Z_i}\Big)_v.
\end{multline*}
Therefore, since
\begin{multline*}
\big(\rest{(D_0, g_0)}{X_j} \cdots \rest{(D_d, g_d)}{X_j}\big)_v = \sum_{i=1}^n  a_{ji}\Big(\rest{(D_0, g_0)}{Z_i} \cdots \rest{(D_{d-1}, g_{d-1})}{Z_i}\Big)_v \\
+ \int_{X_j^{\mathrm{an}}}
g_d(x)\,\mu_{\rest{(D_0, g_0)}{X_j} \cdots \rest{(D_{d-1}, g_{d-1})}{X_j}}(\mathrm{d}x),
\end{multline*}
one has the desired formula.
\end{proof}

\begin{prop}\label{prop:multilinear:symmetric:semiample:case}
Let $(D_0, g_0) \ldots, (D_i, g_i), (D'_i, g'_i), \ldots, (D_d, g_d)$ be integrable metrized Cartier divisors on $X$ such that
$(D_0, \ldots, D_i, \ldots, D_d)$ and $(D_0, \ldots, D'_i, \ldots, D_d)$ belong to $\IP_X$.
Then one has the following:
\begin{enumerate}[label=\rm(\arabic*)]
\item
The local intersection pairing is multi-linear, that is, 
\[
\kern2em\begin{cases}
\big( (D_0, g_0) \cdots (D_i + D'_i, g_i + g'_i) \cdots (D_d, g_d) \big)_v \\
\kern0.7em = \big( (D_0, g_0) \cdots (D_i, g_i) \cdots (D_d, g_d) \big)_v +
\big( (D_0, g_0) \cdots (D'_i, g'_i) \cdots (D_d, g_d) \big)_v. \\[2ex]
\big( (D_0, g_0) \cdots (-D_i, -g_i) \cdots (D_d, g_d) \big)_v = - \big( (D_0, g_0) \cdots (D_i, g_i) \cdots (D_d, g_d) \big)_v.
\end{cases}
\]

\item 

We assume that $D_0, \ldots, D_d$ are ample and $g_0, \ldots, g_d$ are plurisubharmonic.
For each $i$, let $\{ g_{i,n} \}_{n=1}^{\infty}$ be a sequence of plurisubharmonic Green functions of $D_i$
such that $\lim_{n\to\infty} \| g_i - g_{i,n}\|_{\sup} = 0$. Then
\[
\lim_{n\to\infty} \big( (D_0, g_{0,n}) \cdots (D_{d}, g_{d,n}) \big)_v = \big( (D_0, g_0) \cdots (D_d, g_d) \big)_v
\]

\item 
The local intersection pairing is symmetric, that is, for any bijection $\sigma:\{0,\ldots,d\}\rightarrow\{0,\ldots,d\}$ one has
\[
 \big( (D_{\sigma(0)}, g_{\sigma(0)}) \cdots \cdots (D_{\sigma(d)}, g_{\sigma(d)}) \big)_v 
= \big( (D_0, g_0) \cdots (D_d, g_d) \big)_v.\]\end{enumerate}
\end{prop}

\begin{proof}
Clearly we may assume that $X$ is integral.
We prove (1), (2) and (3) by induction on $d$. In the case $d=0$, the assertion is obvious, so that we assume $d > 0$.

\medskip
(1) If $0 \leqslant i < d$, the assertions follow from the hypothesis of induction and
the multi-linearity of the measure $\mu_{(D_0, g_0) \cdots (D_{d-1}, g_{d-1})}$ with respect to $(D_0, g_0)$, $\ldots$, $(D_{d-1}, g_{d-1})$, so that we may assume that $i=d$.
Let $D_d = a_1 Z_1 + \cdots + a_n Z_n$ and $D'_d = a'_1 Z_1 + \cdots + a'_n Z_n$ be
the decompositions of $D_d$ and $D'_{d}$ as cycles.
Then $D_d + D'_d = (a_1 + a'_1) Z_1 + \cdots + (a_n + a'_n) Z_n$ and
$-D_d = (-a_1)Z_1 + \cdots + (-a_n)Z_n$, so that the assertions are obvious.

\medskip

(2) By \eqref{eqn:def:local:intersection} and the hypothesis of induction,
it is sufficient to see 
\[
\lim_{n\to\infty} \int_{X^{\an}} g_{d,n} \mu_{(D_0, g_{0,n}) \cdots (D_{d-1}, g_{d-1, n})} =
\int_{X^{\an}} g_{d} \mu_{(D_0, g_{0}) \cdots (D_{d-1}, g_{d-1})},
\]
which follows from \cite[Corollary~(3.6)]{Demagbook} and \cite[Corollaire~(5.6.5)]{ChamDucros}.

\medskip
(3) We may assume that $D_0, \ldots, D_d$ are ample and $g_0, \ldots, g_d$ are plurisubharmonic.
By (2) together with regularizations of metrics, we may further assume that metrics $|\ndot|_{g_0}, \ldots, |\ndot|_{g_d}$ are smooth. It suffices to prove the assertion in the particular case where $\sigma$ is a transposition exchanging two indices $i$ and $j$ with $i<j$. If $j < d$, then the assertion follows from the hypothesis of induction,
so that we may assume that $j = d$. If $i < d-1$, then
\begin{multline*}
\big( (D_0, g_0) \cdots (D_i, g_i) \cdots (D_{d-1}, g_{d-1}) \cdot (D_d, g_d) \big)_v \\
= \big( (D_0, g_0) \cdots (D_{d-1}, g_{d-1}) \cdots (D_i, g_i) \cdot (D_d, g_d) \big)_v.\end{multline*}
by the hypothesis of induction. Therefore we may assume that $i = d-1$.
Let $D_d = a_1 Z_1 + \cdots + a_n Z_n$ and $\rest{D_{d-1}}{Z_i} = a_{i1} Z_{i1} + \cdots + a_{in} Z_{in}$
be the decomposition as cycles. Then
\begin{multline*}
\big( (D_0, g_0) \cdots (D_{d-1}, g_{d-1}) \cdot (D_d, g_d) \big)_v \\
= \sum_{i,j} a_i a_{ij} \big(\kern-.3em\rest{(D_0, g_0)}{Z_{ij}} \cdots \rest{(D_{d-2}, g_{d-2})}{Z_{ij}}\big)_v \\
\kern7em+ \sum_{i} a_i \int_{Z_i^{\an}} g_{d-1}(x)\,\mu_{\rest{(D_0, g_0)}{Z_i} \cdots \rest{(D_{d-2}, g_{d-2})}{Z_i}}(\mathrm{d}x) \\
+ \int_{X^{\an}} g_d(x)\, \mu_{(D_0, g_0) \cdots (D_{d-1}, g_{d-1})}(\mathrm{d}x).
\end{multline*}
In the same way, if $D_{d-1} = a'_1 Z'_1 + \cdots + a'_n Z'_n$ and $\rest{D_{d}}{Z'_i} = a'_{i1} Z'_{i1} + \cdots + a'_{in} Z'_{in}$
be the decomposition as cycles, then
\begin{multline*}
\big( (D_0, g_0) \cdots (D_{d}, g_{d}) \cdot (D_{d-1}, g_{d-1}) \big)_v \\
= \sum_{i,j} a'_i a'_{ij} \big(\kern-.3em\rest{(D_0, g_0)}{Z'_{ij}} \cdots \rest{(D_{d-2}, g_{d-2})}{Z'_{ij}} \big)_v \\
\kern7em + \sum_{i} a'_i \int_{(Z'_i)^{\an}} g_{d}(x)\, \mu_{\rest{(D_0, g_0)}{Z'_i} \cdots \rest{(D_{d-2}, g_{d-2})}{Z'_i}} (\mathrm{d}x)\\
+ \int_{X^{\an}} g_{d-1}(x)\, \mu_{(D_0, g_0) \cdots (D_{d-2}, g_{d-2}) \cdot (D_{d}, g_{d})}(\mathrm{d}x).
\end{multline*}
By \cite[Proposition~5.2 (2)]{MArakelov}, one has $\sum_{ij} a_i a_{ij}Z_{ij} = \sum_{ij} a'_i a'_{ij}Z'_{ij}$ as cycles, so that it is sufficient to show that
\begin{multline*}
\sum_{i} a_i \int_{Z_i^{\an}} g_{d-1}(x)\,\mu_{\rest{(D_0, g_0)}{Z_i} \cdots \rest{(D_{d-2}, g_{d-2})}{Z_i}}(\mathrm{d}x)
+ \int_{X^{\an}} g_d(x)\, \mu_{(D_0, g_0) \cdots (D_{d-1}, g_{d-1})}(\mathrm{d}x) \\
= \sum_{i} a'_i \int_{(Z'_i)^{\an}} g_{d}(x)\,\mu_{\rest{(D_0, g_0)}{Z'_i} \cdots \rest{(D_{d-2}, g_{d-2})}{Z'_i}}(\mathrm{d}x) \\
+ \int_{X^{\an}} g_{d-1}(x)\, \mu_{(D_0, g_0) \cdots (D_{d-2}, g_{d-2}) \cdot (D_{d}, g_{d})}(\mathrm{d}x),
\end{multline*}
which is nothing more than \cite[Theorem~5.6]{MArakelov} for the Archimedean case and \cite[Proposition~11.5]{GKTropical} for the non-Archimedean case.
\end{proof}

\begin{prop}\label{prop:intersection:finite:morphism}
Let $\pi : Y \to X$ be a surjective morphism of integral projective schemes over $k$.
We set $e = \dim X$ and $d = \dim Y$. 
Let $(D_0, g_0), \ldots, (D_d, g_d)$ be integrable metrized Cartier divisors on $X$ such that
$(\pi^*(D_0), \ldots, \pi^*(D_d)) \in \IP_Y$.
Then one has the following:
\begin{enumerate}[label=\rm(\arabic*)]
\item If $d > e$, then $(\pi^*(D_0, g_0) \cdots \pi^*(D_d, g_d))_v = 0$.
\item If $d = e$ and $(D_0,\ldots,D_d) \in \IP_X$, then 
\[(\pi^*(D_0, g_0) \cdots \pi^*(D_d, g_d))_v = (\deg \pi)((D_0, g_0) \cdots (D_d, g_d))_v.\]
\end{enumerate}
\end{prop}

\begin{proof}
We prove (1) and (2) by induction on $e$. 
If $e = 0$, then (2) is obvious. For (1), as $\pi^*(D_0, g_0) = (0, a_0), \ldots, \pi^*(D_d, g_d) = (0, a_d)$
for some $a_0, \ldots, a_d \in \RR$, then
\[
(\pi^*(D_0, g_0) \cdots \pi^*(D_d, g_d))_v = \int_{X^{\an}} a_d \mu_{(0, a_0)\cdots(0, a_d)} = 0,
\]
as desired.

We assume $e > 0$.
Let $D_d = a_1  Z_1 + \cdots + a_n Z_n$ and $\pi^*(D_d) = b_{1} Z'_{1} + \cdots + b_{N}Z'_{N}$ 
be the decompositions as cycles.
By \eqref{eqn:def:local:intersection},
\begin{multline*}
\big(\pi^*(D_0, g_0) \cdots \pi^*(D_d, g_d)\big)_v = \sum_{j=1}^N b_j \big(\kern-.3em\rest{\pi^*(D_0, g_0)}{Z'_j} \cdots \rest{\pi^*(D_{d-1}, g_{d-1})}{Z'_j}\big)_v \\
+ \int_{Y^{\an}} g_{d}(\pi^{\mathrm{an}}(y)) \,\mu_{\pi^*(D_0, g_0) \cdots \pi^*(D_d, g_d)}(\mathrm{d}y),
\end{multline*}
Note that if $e < d$, then $\dim \pi(Z'_j) < \dim Z'_j$
and $\pi_*(\mu_{\pi^*(D_0, g_0) \cdots \pi^*(D_d, g_d)}) = 0$ by Proposition~\ref{prop:CH:measure:projection:formula}, so that one has (1).

Next we assume that $e = d$.
For each $i$, we set $J_i = \{ j \in \{ 1, \ldots, N \} \mid \pi(Z'_j) = Z_i \}$.
We set $J_0 = \{ 1, \ldots, N \} \setminus (J_1 \cup \cdots \cup J_n)$. By the hypothesis of induction for (1),
$\big(\kern-.3em\rest{\pi^*(D_0, g_0)}{Z'_j} \cdots \rest{\pi^*(D_{d-1}, g_{d-1})}{Z'_j}\big)_v = 0$ for all
$j \in J_0$, so that, by the hypothesis of induction for (2) and Proposition~\ref{prop:CH:measure:projection:formula}, the above equation implies 
\begin{align*}
\big(\pi^*(D_0, g_0) \cdots \pi^*(D_d, g_d)\big)_v \\
& \kern-5em = \sum_{i=1}^n \sum_{j \in J_i} b_j \big(\kern-.3em\rest{\pi^*(D_0, g_0)}{Z'_j} \cdots \rest{\pi^*(D_{d-1}, g_{d-1})}{Z'_j}\big)_v \\
& + \int_{Y^{\an}} g_{d}(\pi^{\mathrm{an}}(y)) \,\mu_{\pi^*(D_0, g_0) \cdots \pi^*(D_d, g_d)}(\mathrm{d}y) \\
& \kern-5em = \sum_{i=1}^n \big(\kern-.3em\rest{(D_0, g_0)}{Z_i} \cdots \rest{(D_{d-1}, g_{d-1})}{Z_i}\big)_v \sum_{j \in J_i} b_j \deg(\rest{\pi}{Z'_j})  \\
& + \deg(\pi) \int_{X^{\an}} g_{d}(x)\, \mu_{(D_0, g_0) \cdots (D_d, g_d)}(\mathrm{d}x).
\end{align*}
Therefore, the assertion follows because $\sum_{j \in J_i} b_j \deg(\rest{\pi}{Z'_j}) = \deg(\pi)a_i$ (cf. \cite[Lemma~1.12]{MArakelov}).
\end{proof}

\begin{prop}\label{prop:intersection:principal:div}
Let $f$ be a regular meromorphic function on $X$ and $(D_1, g_1), \ldots, (D_d, g_d)$
be integrable metrized Cartier divisors on $X$ such that $(\operatorname{div}(f), D_1, \ldots, D_d) \in \IP_X$.
If we set $D_1 \cdots D_{d} = \sum_{x \in X_{(0)}} a_x x$ as cycle, then
\begin{equation} \label{Equ: intersection principal divisor}\big(\widehat{\operatorname{div}}(f) \cdot (D_1, g_1) \cdots (D_{d}, g_{d})\big)_v = \sum_{x \in X_{(0)}} a_x \big(-\log |f|(x^{\an})\big),\end{equation}
where $X_{(0)}$ is the set of all closed point of $X$ and $x^{\an}$ is the associated absolute value at $x$.
Note that in the case where $\dim(X) = 0$, the above formula means that 
\[\big(\widehat{\operatorname{div}}(f) \big)_v = 0.\]
\end{prop}

\begin{proof}
Let $X = a_1 X_1 + \cdots + a_n X_n$ be the decomposition as cycles.
Then
\begin{equation*} \big(\widehat{\operatorname{div}}(f) \cdot (D_1, g_1) \cdots (D_{d}, g_{d})\big)_v 
= 
\sum_{i=1}^n a_i \big({\widehat{\operatorname{div}}(f)}|_{X_i} \cdot \rest{(D_1, g_1)}{X_i} \cdots \rest{(D_{d}, g_{d})}{X_i}\big)_v
\end{equation*}
and
\[
D_1 \cdots D_{d} = \sum_{i=1}^n a_i \big(\rest{D_1}{X_i} \cdots\rest{D_d}{X_i}\big),
\]
so that we may assume that $X$ is integral.

We prove the equality \eqref{Equ: intersection principal divisor} by induction on $d=\dim(X)$. In the case where $\dim(X) = 0$, the assertion is obvious
because $f$ is a unit.
We assume that $\dim(X) \geqslant 1$. Let $D_d = a_1 Z_1 + \cdots + a_n Z_n$ 
be the decomposition as cycles. Let $\sum_{x \in X_{(0)}} b_{ix} x$ be the decomposition of
$\rest{D_1}{Z_i} \cdots \rest{D_{d-1}}{Z_i} = D_1 \cdots D_{d-1} \cdot Z_i$ as cycles.
Then
\[
\sum_{i=1}^n a_i \sum_{x \in X_{(0)}} b_{ix} x = \sum_{x \in X_{(0)}} a_x x,
\]
so that $a_x = \sum_{i=1}^n a_i b_{ix}$. On the other hand, by hypothesis of induction,
\[
\big(\widehat{\operatorname{div}}(f)|_{Z_i} \cdot \rest{(D_1, g_1)}{Z_i} \cdots \rest{(D_{d-1}, g_{d-1})}{Z_i}\big)_v = \sum_{x \in X_{(0)}} b_{ix} (-\log |f|(x^{\an})).
\]
Therefore,
\begin{multline*}
\sum_{x \in X_{(0)}} a_x (-\log |f|(x^{\an})) \\
= \sum_{x \in X_{(0)}} \Big( \sum_{i=1}^n a_i b_{ix} \Big) (-\log |f|(x^{\an})) 
=\sum_{i=1}^n a_i \sum_{x \in X} b_{ix} (-\log |f|(x^{\an})) \\
= \sum_{i=1}^n a_i \big(\widehat{\operatorname{div}}(f)|_{Z_i} \cdot \rest{(D_1, g_1)}{Z_i} \cdots \rest{(D_{d-1}, g_{d-1})}{Z_i}\big)_v.
\end{multline*}
Note that $\mu_{(\widehat{\operatorname{div}}(f) \cdot (D_1, g_1) \cdots (D_d, g_d))} = 0$, and hence the assertion follows by \eqref{eqn:def:local:intersection}.
\end{proof}

\begin{prop}\label{prop:intersection:trivial:divisor}
Let $(D_0, g_0), \ldots, (D_{d-1}, g_{d-1}), (0,g)$ be integrable metrized Cartier divisors on $X$ with $(D_0, \ldots, D_{d-1}, 0) \in \IP_X$.
We assume that $D_0, \ldots, D_{d-1}$ are semiample and $g_0, \ldots, g_{d-1}$ are plurisubharmonic.
Then
\[
\big((D_0, g_0) \cdots (D_{d-1}, g_{d-1}) \cdot (0, g)\big)_v = \int_{X^{\an}} g(x) \mu_{(D_0, g_0) \cdots (D_{d-1}, g_{d-1})}(\mathrm{d}x).
\]
In particular,
\begin{multline*}
\min \{ g(x) \mid x \in X^{\an} \} (D_0 \cdots D_{d-1}) \\\
\leqslant
\big((D_0, g_0) \cdots (D_{d-1}, g_{d-1}) \cdot (0, g)\big)_v \\
\leqslant \max \{ g(x) \mid x \in X^{\an} \} (D_0 \cdots D_{d-1}).
\end{multline*}
\end{prop}

\begin{proof}
This is trivial by the definition.
\end{proof}

\begin{coro}\label{cor:intersection:trivial:divisor}
Let $(D_0, g_0), \ldots, (D_d, g_d)$ be integrable arithmetic Cartier divisors on $X$ with $(D_0, \ldots, D_d) \in \IP_X$.
We assume that $D_0, \ldots, D_d$ are semiample and $g_0, \ldots, g_d$ are plurisubmarmonic.
Let $g'_0, \ldots, g'_d$ be another plurisubharmonic Green functions of $D_0, \ldots, D_d$, respectively.
Then one has
\begin{multline*}
\left| \big((D_0, g'_0) \cdots (D_d, g'_d)\big)_{v} - \big((D_0, g_0) \cdots (D_d, g_d)\big)_{v} \right| \\
\leqslant
\sum_{i=0}^d \max \{|g'_i - g_i|(x) \mid x \in X^{\an} \} (D_0 \cdots D_{i-1} \cdot D_{i+1} \cdots D_d).
\end{multline*}
\end{coro}

\begin{proof}
By using Proposition~\ref{prop:multilinear:symmetric:semiample:case},
\begin{multline*}
\big((D_0, g'_0) \cdots (D_d, g'_d)\big) - \big((D_0, g_0) \cdots (D_d, g_d)\big) \\
=
\sum_{i=0}^d ((D_0, g_0) \cdots (D_{i-1}, g_{i-1}) \cdot (0, g'_i - g_i) \cdot (D_{i+1}, g'_{i+1}) \cdots (D_d, g'_d) \big),
\end{multline*}
so that the assertion follows from Proposition~\ref{prop:intersection:trivial:divisor}.
\end{proof}

\begin{prop}\label{prop:local:intersection:PP}
We assume that $X = \PP^d_k$ and $L = \OO_{\PP^d}(1)$. Let $\{ T_0, \ldots, T_d \}$ be a basis of $H^0(\PP^d_k, \OO_{\PP^d}(1))$ over $k$. We view $(T_0 : \cdots : T_d)$ as a homogeneous coordinate of
$\PP^d_k$.
Let $\|\ndot\|$ be a norm of $H^0(\PP^d_k, \OO_{\PP^d}(1))$ given by
\[
\| a_0 T_0 + \cdots + a_d T_d \| = \begin{cases}
\sqrt{|a_0|^2 + \cdots + |a_d|^2} & \text{if $v$ is Archimedean},\\[2ex]
\max \{ |a_0|, \ldots, |a_d| \} & \text{if $v$ is non-Archimedean}.
\end{cases}
\]
Let $\varphi$
be the orthogonal quotient metric of $\OO_{\PP^d}(1)$ given by the surjective homomorphism
$H^0(\PP^d_k, \OO_{\PP^d}(1)) \otimes \OO_{\PP^d} \to \OO_{\PP^d}(1)$ and the above norm $\|\ndot\|$.
We set $H_i = \{ T_i = 0 \}$ and $h_i = -\log |T_i|_{\varphi}$. Then
\[
((H_0, h_0) \cdots (H_d, h_d))_v = \begin{cases}
\adeg\Big(\widehat{c}_1(\OO_{\PP^d_{\ZZ}}(1), \varphi)^{d+1}\Big) & \text{if $v$ is Archimedean},\\[2ex]
0 & \text{if $v$ is non-Archimedean},
\end{cases}
\]
where $\adeg\Big(\widehat{c}_1(\OO_{\PP^d_{\ZZ}}(1), \varphi)^{d+1}\Big)$ is the self-intersection number of
the arithmetic first Chern class $\widehat{c}_1(\OO_{\PP^d_{\ZZ}}(1), \varphi)$ on 
the $d$-dimensional projective space $\PP^d_{\ZZ}$ over $\ZZ$.
\end{prop}

\begin{proof}
If we set
\[
a_m := \int_{\PP^m_k} -\log |T_m|_{\varphi}(x)\,\mu_{(\OO_{\PP^m}(1), \varphi_{\rm{FS}})^m}(\mathrm{d}x)
\]
for a positive integer $m$, then
\[
((H_0, h_0) \cdots (H_d, h_d))_v = \sum_{m=1}^d a_m.
\]
In the following, we set $x_i = T_i/T_0$.

\medskip
{\bf $\bullet$ Archimedean case} :\quad
The algorithms of the calculation is exactly same as one on $\PP^d_{\ZZ}$, so that we have the assertion.

\bigskip
{\bf $\bullet$ non-Archimedean case} :\quad
If we set
$|f|_* = \max_{i_1, \ldots, i_m} \{ |c_{i_1, \ldots, i_m}| \}$ for \[f = \sum_{i_1,\ldots,i_m} c_{i_1, \ldots, i_m} x_1^{i_1} \cdots x_m^{i_m} \in k[x_1, \ldots, x_m],\] then $|\ndot|_*$ extends to an absolute value of $k(x_1, \ldots, x_m)$ (cf. Lemma~\ref{lem:Gauss:lemma}). We set $U = \{ T_m  \not= 0\}$. Note that if $\xi \in U^{\an}$, then
\[
|T_m|_{\varphi}(\xi) = \frac{|x_m|_\xi}{\max\{1, |x_1|_\xi, \ldots, |x_m|_\xi\}}.
\]
Let $\oo_v$ be the valuation ring of $v$. Note that $\varphi$ coincides with
the metric of the model $(\PP^d_{\oo_v}, \OO_{\PP^d_{\oo_v}}(1))$ by \cite[Proposition~2.3.12]{CMArakelovAdelic},
so that $\mu_{(\OO_{\PP^m}(1), \varphi)^m} = \delta_{|\ndot|_*}$. Thus
\[
a_m = -\log \frac{|x_m|_*}{\max\{1, |x_1|_*, \ldots, |x_m|_*\}} = 0,
\]
and hence the assertion follows.
\end{proof}

\section{Local intersection number over a general field}\label{Sec: Local intersection number over a general field} 
 
In this section, we consider the local intersection product and local height formula in the non-necessarily algebraically closed case. We fix in this section a complete valued field $v=(k,|\ndot|)$ such that $|\ndot|$ is not trivial. Let $\mathbb C_k$ be the completion of an algebraic closure of $k$. Note that the absolute value $|\ndot|$ extends naturally to $\mathbb C_k$ and the valued field $(\mathbb C_k,|\ndot|)$ is both algebraically closed and complete. We denote by $v^{\operatorname{ac}}$ the couple $(\mathbb C_k,|\ndot|)$. We also fix a projective morphism $\pi:X\rightarrow\Spec k$ and we denote by $X_{\mathbb C_k}$ the fiber product $X\times_{\Spec k}\Spec \mathbb C_k$. Let $d$ be the Krull dimension of $X$, which is also equal to the Krull dimension of $X_{\mathbb C_k}$.

\begin{defi}\label{Def: intersection pairing general case}
Let $(D_0,g_0),\ldots,(D_d,g_d)$ be a family of metrized Cartier divisor on $X$ such that $D_0,\ldots,D_d$ intersect properly and that $g_0,\ldots,g_d$ are integrable Green functions.  By Remark \ref{Rem: extension of scalars}, the Cartier divisors $D_{0,\mathbb C_k},\ldots,D_{d,\mathbb C_k}$ intersect properly. Moreover, by Remark \ref{Rem: extension scalar Green function}, the Green functions $g_{0,\mathbb C_k},\ldots,g_{d,\mathbb C_k}$ are integrable. We then define the local intersection number of $(D_0,g_0),\ldots,(D_d,g_d)$ as
\[\big( (D_0, g_0) \cdots (D_d, g_d) \big)_v:=\big( (D_{0,\mathbb C_k}, g_{0,\mathbb C_k}) \cdots (D_{d,\mathbb C_k}, g_{d,\mathbb C_k}) \big)_{v^{\operatorname{ac}}}\]  
\end{defi}

Several properties of the local intersection number follow directly from the results of the previous section. We gather them below.

\begin{rema}\label{rem:properties:local:intersection:general:field}
Recall that $\operatorname{\widehat{Int}}(X)$ denotes the group of integrable metrized Cartier divisors on $X$. Let $\widehat{\IP}_X$ be the subset of $\operatorname{\widehat{Int}}(X)^{d+1}$ consisting of elements
\[\big((D_0,g_0),\ldots,(D_d,g_d)\big)\]
such that the Cartier divisors $D_0,\ldots,D_d$ intersect properly. 
\begin{enumerate}[label=\rm(\arabic*)]
\item The set $\widehat{\IP}_X$ forms a symmetric multi-linear subset of the group $\operatorname{\widehat{Int}}(X)^{d+1}$. Moreover, the function of local intersection number
\[\big( (D_0, g_0) \cdots (D_d, g_d) \big)\longmapsto\big( (D_0, g_0) \cdots (D_d, g_d) \big)_v\]
form a symmetric multi-linear map from $\widehat{\IP}_X$ to $\mathbb R$. These statements follow from Proposition \ref{prop:multilinear:symmetric:semiample:case}.

\item
Let $\pi : Y \to X$ be a surjective morphism of geometrically integral projective schemes over $k$.
We set $e = \dim X$ and $d = \dim Y$. 
Let $(D_0, g_0), \ldots, (D_d, g_d)$ be integrable metrized Cartier divisors on $X$ such that
$(\pi^*(D_0), \ldots, \pi^*(D_d)) \in \IP_Y$.
Then one has the following:
\begin{enumerate}[label=\rm(\roman*)]
\item If $d > e$, then $(\pi^*(D_0, g_0) \cdots \pi^*(D_d, g_d))_v = 0$.
\item If $d = e$ and $(D_0,\ldots,D_d) \in \IP_X$, then 
\[(\pi^*(D_0, g_0) \cdots \pi^*(D_d, g_d))_v = (\deg \pi)((D_0, g_0) \cdots (D_d, g_d))_v.\]
\end{enumerate}
We refer to Proposition \ref{prop:intersection:finite:morphism} for a proof.

\item Let $f$ be a regular meromorphic function on $X$ and $(D_1,g_1),\ldots,(D_d,g_d)$ be integrable metrized Cartier divisors on $X$ such that $(\operatorname{div}(f),D_1,\ldots,D_d)\in\IP_X$. Suppose that 
\[D_1\cdots D_d=\sum_{x\in X_{(0)}}a_xx\]
as a cycle, then
\[\big(\widehat{\operatorname{div}}(f) \cdot (D_1, g_1) \cdots (D_{d}, g_{d})\big)_v = \sum_{x \in X_{(0)}} a_x[\kappa(x):k]_s \big(-\log |f|(x^{\an})\big),\]
where $[\kappa(x):k]_s$ denotes the separable degree of the residue field $\kappa(x)$ over $k$. We refer to Proposition \ref{prop:intersection:principal:div} for more details.
\item Let $\big((D_0,g_0),\ldots,(D_d,g_d)\big)$ be an element of $\widehat{\IP}_X$. We assume that $D_0,\ldots,D_{d-1}$ are semi-ample, $g_0,\ldots,g_{d-1}$ are plurisubharmonic, and $D_{\textcolor{mred}{d}}=0$. Then one has
\[\delta\min_{x\in X^{\mathrm{an}}}g_d(x)\leqslant \big((D_0,g_0),\ldots,(D_d,g_d)\big)_v\leqslant\delta\max_{x\in X^{\mathrm{an}}}g_d(x),\]
where $\delta=(D_0,\ldots,D_{d-1})$. See Proposition \ref{prop:intersection:trivial:divisor} for more details.

\item Let $\big((D_0,g_0),\ldots,(D_d,g_d)\big)$ and $\big((D_0,g_0'),\ldots,(D_d,g_d')\big)$ be two elements of $\widehat{\IP}_X$ having the same family of underlying Cartier divisors. One has
\[\begin{split}\Big|\big((D_0,g_0),\ldots,(D_d,g_d)\big)_v&-\big((D_0,g_0'),\ldots,(D_d,g_d')\big)_v\Big|\\&\quad\leqslant\sum_{i=0}^d\max_{x\in X^{\operatorname{an}}}|g_i'-g_i|(x)(D_0\cdots D_{i-1}\cdot D_{i+1}\cdots D_d).\end{split}\] See Corollary \ref{cor:intersection:trivial:divisor} for more details.
\end{enumerate}
\end{rema}

\section{Local height}

In this section, we fix a complete valued field $v=(k,|\ndot|)$ and a projective scheme $X$ over $\Spec k$. Let $d$ be the dimension of $X$.
\begin{defi}\label{Def: intersection number}
Let $\overline L_i=(L_i,\varphi_i)$, $i\in\{0,\ldots,d\}$ be a family of metrized invertible $\mathcal O_X$-modules, where each $L_i$ is an invertible $\mathcal O_X$-module, and $\varphi_i$ is a continuous and integrable metric on $L_i$. For any $i\in\{0,\ldots,d\}$, we let $s_i$ be a regular meromorphic section of $L_i$ on $X$. Assume that the Cartier divisors $\operatorname{div}(s_0),\ldots,\operatorname{div}(s_d)$ intersect properly. We define the \emph{local height} of $X$ with respect to the family  of metrized invertible $\mathcal O_X$-modules $(\overline L_i)_{i=0}^d$ and the family of regular meromorphic sections $(s_i)_{i=0}^d$ as the local intersection number (see Definition \ref{Def: intersection pairing general case})
\[h_{\overline L_0,\ldots,\overline L_d}^{s_0,\ldots,s_d}(X):=\big(\operatorname{\widehat{div}}(s_0)\cdots\operatorname{\widehat{div}}(s_d)\big)_v.\] 
\end{defi}

\begin{enonce}[remark]{Notation}\label{Not: simplified notation}
We often encounter the situation where each $\overline L_i$ is the pull-back by a projective morphism $f_i:X\rightarrow Y_i$ of a metrized invertible $\mathcal O_{Y_i}$-module $\overline{M_i}$ and $s_i$ is the pull-back of a regular meromorphic section $t_i$. In such a situation, for simplicity of notation, 
we often use the expressions $h_{\overline M_0,\ldots,\overline M_{d}}^{t_0,\ldots,t_{d}}(X)$ or $h_{\overline L_0,\ldots,\overline L_{d}}^{t_0,\ldots,t_{d}}(X)$  to denote $h_{\overline L_0,\ldots,\overline L_d}^{s_0,\ldots,s_{d}}(X)$.
\end{enonce}

\begin{rema} We keep the notation of Definition \ref{Def: intersection number} in assuming that the field $k$ is algebraically closed.
Let $X_1,\ldots,X_n$ be irreducible components of $X$, considered as reduced closed subscheme of $X$. For any $j\in\{1,\ldots,n\}$, let $\operatorname{mult}_{X_j}(X)$ be the multiplicity of the component $X_j$, which is by definition the length of the Artinian local ring of $\mathcal O_X$ at the generic point of $X_j$. Then, for any $j\in\{1,\ldots,n\}$, the divisors on $X_j$ associated with the restricted sections $(s_i|_{X_j})_{i=0}^d$ intersect properly on $X_j$. 

Assume firstly that $d=0$. In this case, each $X_j$ consists of a closed point $x_j$ of $X$, which is actually a rational point since $k$ is supposed to be algebraically closed. Hence $X_j^{\operatorname{an}}$  only contains one point, which we denote by $x_j^{\operatorname{an}}$. Note that $s_0$ does not vanish at any of the closed points $X_j$. By definition, $h_{\overline L_0}^{s_0}(X)$ is equal to 
\begin{equation}\label{Equ: height of one point}-\sum_{j=1}^n\operatorname{mult}_{X_j}(X)\ln|s_0|_{\varphi_0}(x_j^{\operatorname{an}}).\end{equation}

In the case where $d\geqslant 1$, the induction formula in Definition \ref{def:local:intersection} 
for local intersection number leads to the following formula for the local height. 
\begin{equation}\label{Equ: height recursive formula}\begin{split}h_{\overline L_0,\ldots,\overline L_d}^{s_0,\ldots,s_d}(X)=
\sum_{i=1}^n & a_i h_{\overline L_0,\ldots,\overline L_{d-1}}^{s_0,\ldots,s_{d-1}}(Z_i)\\&-\int_{X^{\mathrm{an}}}\ln|s_d|_{\varphi_d}(x) \,\mu_{(L_0,\varphi_0)\cdots (L_{d-1},\varphi_{d-1})}(\mathrm{d}x),\end{split}\end{equation}
where $\sum_{i=1}^n a_i Z_i$ is the cycle associated with $\operatorname{div}(L_d; s_d)$.
\end{rema}

\begin{defi}
Let $(E,\|\ndot\|)$ be a finite-dimensional normed vector space over $k$, and $r$ be the rank of $E$. We denote by $\|\ndot\|_{\det}$ the norm on the one-dimensional vector space $\det(E):=\Lambda^r(E)$ such that,
\[\forall\,\eta\in\det(E),\quad
\|\eta\|_{\det}:=\inf_{\eta=t_1\wedge\cdots\wedge t_r}\|t_1\|\cdots\|t_r\|.\]
Note that, if the norm $\|\ndot\|$ is ultrametric or induced by an inner product, for any complete valued extension $k'$ of $k$, one has (see Definition \ref{Def: extension of sclars})
\begin{equation}\label{Equ: determinant norm}\|\ndot\|_{k',\det}=\|\ndot\|_{\det,k'},\end{equation}
if we identify $\det(E)\otimes_kk'$ with $\det(E\otimes_kk')$.
We refer the readers to \cite[Proposition 1.3.19]{CMArakelovAdelic} for a proof.
\end{defi}

\begin{prop}\label{Pro: height projective space}
Let $E$ be a finite-dimensional vector space over $k$, equipped with a norm $\|\ndot\|$ which is either ultrametric or induced by an inner product, $r=\operatorname{dim}_k(E)$, and $L=\mathcal O_E(1)$ be the universal invertible sheaf on $\mathbb P(E)$. We equip $L$ with the orthogonal quotient metric $\varphi$ induced by $\|\ndot\|$. Let $(s_j)_{j=0}^r$ be a basis of $E$ over $k$. If $|\ndot|$ is non-Archimedean, then
\[h_{\overline{L},\ldots,\overline{L}}^{s_0,\ldots,s_r}(\mathbb P(E))=-\ln\|s_0\wedge\cdots\wedge s_r\|_{\det};\] if $|\ndot|$ is Archimedean, then
\[h_{\overline{L},\ldots,\overline{L}}^{s_0,\ldots,s_r}(\mathbb P(E))=-\ln\|s_0\wedge\cdots\wedge s_r\|_{\det}+\sigma_r,\]
where
\[\sigma_r=\frac 12\sum_{m=1}^r\sum_{\ell=1}^m\frac{1}{\ell}\]
is the $r$-th Stoll number.
\end{prop}

\begin{proof}
First, the metric $\varphi_{\mathbb C_k}$ identifies with the orthogonal quotient metric induced by $\|\ndot\|_{\mathbb C_k}$. Therefore, by \eqref{Equ: determinant norm} we may assume without loss of generality that $k$ is algebraically closed.

By Remark \ref{Rem: approximation by good norms}, one can find a sequence $(\|\ndot\|_n)_{n\in\mathbb N}$ of orthonormally decomposable norms such that 
\[\lim_{n\rightarrow+\infty}d(\|\ndot\|_n,\|\ndot\|)=0.\]
By \eqref{Equ: distance of quotient metrics}, if we denote by $\varphi_n$ the orthogonal quotient metric on $L$ induced by $\|\ndot\|_n$, then one has 
\[\lim_{n\rightarrow+\infty}d(\varphi_n,\varphi)=0.\]
By Corollary \ref{cor:intersection:trivial:divisor}, one has
\[\lim_{n\rightarrow+\infty}h_{(L,\varphi_n),\ldots,(L,\varphi_n)}^{s_0,\ldots,s_r}(\mathbb P(E))=h_{\overline{L},\ldots,\overline{L}}^{s_0,\ldots,s_r}(\mathbb P(E)).\]
Moreover, by \cite[Proposition 1.1.64]{CMArakelovAdelic} one has
\[0\leqslant d(\|\ndot\|_{n,\det},\|\ndot\|_{\det})\leqslant r d(\|\ndot\|_n,\|\ndot\|) \]
and hence
\[\lim_{n\rightarrow+\infty}d(\|\ndot\|_{n,\det},\|\ndot\|_{\det})=0.\]
Therefore, without loss of generality, we may assume that the norm $\|\ndot\|$ is orthonormally decomposable.

We reason by induction on $r$. In the case where $r=0$, the vector space $E$ is one-dimensional, and $s_0$ is a non-zero element of $E$. One has 
\[h_{\overline L}^{s_0}(\mathbb P(E))=-\ln\|s_0\|.\]
We now assume that $r\geqslant 1$. Let $G$ be the quotient vector space of $E$ by $ks_r$. Note that the quotient norm $\|\ndot\|_{\operatorname{quot}}$ on $G$ is orthonormally decomposable (see Proposition \ref{Pro: decomposable quotient norm}). For $j\in\{0,\ldots,r-1\}$, let $\overline{s}_j$ be the class of $s_j$ in $G$. We can also view $\overline{s}_j$ as the restriction of $s_j$ to the closed subscheme $\mathbb P(G)$ of $\mathbb P(E)$. We apply the induction hypothesis to $(G,\|\ndot\|_{\operatorname{quot}})$ and obtain (see Notation \ref{Not: simplified notation})
\[h_{\overline L,\ldots,\overline L}^{s_0,\ldots,s_{r-1}}(\mathbb P(G))=-\ln\|\overline{s}_0\wedge\cdots\wedge\overline{s}_{r-1}\|_{\operatorname{quot},\det}\]
when $|\ndot|$ is non-Archimedean and
\[h_{\overline L,\ldots,\overline L}^{s_0,\ldots,s_{r-1}}(\mathbb P(G))=-\ln\|\overline{s}_0\wedge\cdots\wedge\overline{s}_{r-1}\|_{\operatorname{quot},\det}+\sigma_{r-1}\]
We now compute the integral
\[-\int_{\mathbb P(E)^{\mathrm{an}}}\ln|s_r|_{\varphi}\,\mathrm{d}\mu_{\overline{L}^r}.\]
We first consider the case where $|\ndot|$ is non-Archimedean. By Proposition \ref{Pro: non Archimedean projective space} one has
\[\int_{\mathbb P(E)^{\mathrm{an}}}\ln|s_r|_{\varphi}\,\mathrm{d}\mu_{\overline{L}^r}=-\ln|s_r|_{\varphi}(\xi)=-\ln\|s_r\|,\]
where $\xi$ denotes the Gauss point of $\mathbb P(E)^{\mathrm{an}}$. 
Therefore, by \cite[Proposition 1.2.51]{CMArakelovAdelic} we obtain
\[h_{\overline{L},\ldots,\overline{L}}^{s_0,\ldots,s_r}(\mathbb P(E))=-\ln\|\overline s_0\wedge\cdots\wedge\overline s_{r-1}\|_{\operatorname{quot},\det}-\ln\|s_r\|=-\ln\|s_0\wedge\cdots\wedge s_r\|_{\det}.\]
In the case where $|\ndot|$ is Archimedean, by \cite[\S1.4.3]{MR1260106} Remark (iii), one has 
\[-\int_{\mathbb P(E)^{\mathrm{an}}}\ln|s_r|_{\varphi_r}\,\mathrm{d}\mu_{\overline{L}^r}=-\ln\|s\|+\frac 12\sum_{\ell=1}^r\frac{1}{\ell}.\]
Therefore \[\begin{split}h_{\overline{L},\ldots,\overline{L}}^{s_0,\ldots,s_r}(\mathbb P(E))&=-\ln\|\overline{s}_0\wedge\cdots\wedge\overline{s}_{r-1}\|_{\operatorname{quot},\det}-\ln\|s_r\|+\frac 12\sum_{m=1}^r\sum_{\ell=1}^m\frac{1}{\ell}\\
&=-\ln\|s_0\wedge\cdots\wedge s_r\|_{\det}+\frac 12\sum_{m=1}^r\sum_{\ell=1}^m\frac{1}{\ell}.
\end{split}\]
\end{proof}

In the remaining of the section, we consider a family \[(E_i,\|\ndot\|_i),\quad i\in\{0,\ldots,d\}\] of finite-dimensional vector spaces over $k$ equipped with norms which are either ultrametric or induced by inner products. For each $i\in\{0,\ldots,d\}$, we let  $(E_i^\vee,\|\ndot\|_{i,*})$ be the dual normed vector space of $(E_i,\|\ndot\|_i)$, $r_i:=\operatorname{dim}_k(E_i)-1$, $(s_{i,j})_{j=0}^{r_i}$ be a basis of $E_i$ over $k$, and $(\alpha_{i,j})_{j=0}^{r_i}$ be the dual basis of $(s_{i,j})_{j=0}^{r_i}$, namely
\[\alpha_{i,j}(s_{i,j})=1\quad \text{and}\quad \alpha_{i,j}(s_{i,\ell})=0\quad \text{if $j\neq\ell$}.\] 

Let $\check{\mathbb P}$ be the product projective space
\[\mathbb P(E_0^\vee)\times_k\cdots\times_k\mathbb P(E_d^\vee).\]
For any $i\in\{0,\ldots,d\}$, let $\pi_i:\check{\mathbb P}\rightarrow\mathbb P(E_i^\vee)$ be the morphism of projection to the $i^{\text{th}}$ coordinate, and $L_i=\pi_i^*(\mathcal O_{E_i^\vee}(1))$. We equip $L_i$ with the orthogonal quotient metric induced by $\|\ndot\|_{i,*}$, which we denote by $\varphi_i$. Let $(\delta_0,\ldots,\delta_d)$ be an element of $\mathbb N^{d+1}$,
\[L=\pi_0^*(\mathcal O_{E_0^\vee}(\delta_0))\otimes\cdots\otimes\pi_d^*(\mathcal O_{E_d^\vee}(\delta_d))=L_0^{\otimes\delta_0}\otimes\cdots\otimes L_d^{\otimes\delta_d}.\]
We equip $L$ with the metric
\[\varphi:=\varphi_0^{\otimes\delta_0}\otimes\cdots\otimes\varphi_d^{\otimes\delta_d}.\]
Let $R$ be a non-zero element of 
\[S^{\delta_0}(E_0^\vee)\otimes_k\cdots\otimes_k S^{\delta_d}(E_d^\vee),\]
which is considered as a global section of $L$, and also as a multi-homogenous polynomial of multi-degree $(\delta_0,\ldots,\delta_d)$ on $E_0\times\cdots\times E_d$. For any $i\in\{0,\ldots,d\}$, let 
\[{\overline{\boldsymbol{L}}_i}=(\underbrace{\overline L_i,\ldots,\overline L_i}_{r_i\text{ copies}}),\quad\boldsymbol{\alpha}_i:=(\alpha_{i,j})_{j=1}^{r_i}.\]
The purpose of this section is to compute the  local height $h_{\overline L,\overline{\boldsymbol{L}}_0,\ldots,\overline{\boldsymbol{L}}_d}^{R,\boldsymbol{\alpha}_0,\ldots,\boldsymbol{\alpha}_d}(\check{\mathbb P})$.

\begin{prop}\label{Pro: hight of multi-projective} Assume that the sections $R$ and 
\[\alpha_{i,j},\quad i\in\{0,\ldots,d\},\quad j\in\{1,\ldots,r_i\}\]
intersect property on $\check{\mathbb P}$. If the absolute value $|\ndot|$ is non-Archimedean, then
\[h_{\overline L,\overline{\boldsymbol{L}}_0,\ldots,\overline{\boldsymbol{L}}_d}^{R,\boldsymbol{\alpha}_0,\ldots,\boldsymbol{\alpha}_d}(\check{\mathbb P})=-\ln|R(s_{0,0},\ldots,s_{d,0})|-\sum_{i=0}^d\delta_i\ln\|\alpha_{i,0}\wedge\cdots\wedge\alpha_{i,r_i}\|_{i,*,\det};\]
if the absolute value $|\ndot|$ is Archimedean, then
\[h_{\overline L,\overline{\boldsymbol{L}}_0,\ldots,\overline{\boldsymbol{L}}_d}^{R,\boldsymbol{\alpha}_0,\ldots,\boldsymbol{\alpha}_d}(\check{\mathbb P})=-\ln|R(s_{0,0},\ldots,s_{d,0})|-\sum_{i=0}^d\delta_i\big(\ln\|\alpha_{i,0}\wedge\cdots\wedge\alpha_{i,r_i}\|_{i,*,\det}-\sigma_{r_i}\big).\]
\end{prop}
\begin{proof} 
By the same argument as in the beginning of the proof of Proposition \ref{Pro: height projective space}, we may assume without loss of generality that $k$ is algebraically closed and that all norms $\|\ndot\|_i$ are orthonormally decomposable.

We reason by induction on $r_0+\cdots+r_d$. Consider first the case where $r_0=\cdots=r_d=0$. One has
\[h_{\overline L}^R(\check{\mathbb P})=-\ln|R(s_{0,0},\ldots,s_{d,0})|.\]
In the following, we assume that $r_0+\cdots+r_d>0$. Let $i$ be an element of $\{0,\ldots,d\}$ such that $r_i>0$. We consider the quotient vector space $G_i^\vee=E_i^\vee/k\alpha_{i,r_i}$. For $j\in\{0,\ldots,r_i-1\}$, let $\overline \alpha_{i,j}$ be the class of $\alpha_{i,j}$ in $G_i$. Let $\overline{\boldsymbol{\alpha}}_i:=(\overline \alpha_{i,j})_{j=1}^{r_i-1}$ and
\[\check{\mathbb P}'=\mathbb P(E_0)\times_k\cdots\times\mathbb P(E_{i-1})\times_k\mathbb P(G_i)\times_k\mathbb P(E_{i+1})\times_k\cdots\times_k\mathbb P(E_d).\]
By the same argument as in Proposition \ref{Pro: non Archimedean projective space}, we obtain that, in the case where the absolute value $|\ndot|$ is non-Archimedean, one has 
\[\mu_{\overline L\,\overline L_0^{r_0}\cdots\overline L_{i-1}^{r_{i-1}}\overline L_i^{r_i-1}\overline{L}_{i+1}^{r_{i+1}}\cdots\overline L_d^{r_d}}=\delta_i\operatorname{Dirac}_\xi,\]
where $\operatorname{Dirac}_\xi$ denotes the Dirac measure at the Gauss point $\xi$ of $\check{\mathbb P}^{\operatorname{an}}$. Hence,
 by \eqref{Equ: height recursive formula}, one has
\[\begin{split}h_{\overline L,\overline{\boldsymbol{L}}_0,\ldots,\overline{\boldsymbol{L}}_d}^{R,\boldsymbol{\alpha}_0,\ldots,\boldsymbol{\alpha}_d}(\check{\mathbb P})&=h_{\overline L,\overline{\boldsymbol{L}}_0,\ldots,\overline{\boldsymbol{L}}_{i-1},\overline{\boldsymbol{L}}_i',\overline{\boldsymbol{L}}_{i+1},\ldots,\overline{\boldsymbol{L}}_d}^{R,\boldsymbol{\alpha}_0,\ldots,\boldsymbol{\alpha}_{i-1},\overline{\boldsymbol{\alpha}}_i,\boldsymbol{\alpha}_{i+1},\ldots,\boldsymbol{\alpha}_d}(\check{\mathbb P}')-\delta_i\ln|\alpha_{i,r_i}|_{\varphi_i}(\xi)\\&=h_{\overline L,\overline{\boldsymbol{L}}_0,\ldots,\overline{\boldsymbol{L}}_{i-1},\overline{\boldsymbol{L}}_i',\overline{\boldsymbol{L}}_{i+1},\ldots,\overline{\boldsymbol{L}}_d}^{R,\boldsymbol{\alpha}_0,\ldots,\boldsymbol{\alpha}_{i-1},\overline{\boldsymbol{\alpha}}_i,\boldsymbol{\alpha}_{i+1},\ldots,\boldsymbol{\alpha}_d}(\check{\mathbb P}')-\delta_i\ln\|\alpha_{i,r_i}\|_{i,*},
\end{split}\]
where
\[\overline{\boldsymbol{L}}_i':=(\underbrace{\overline L_i,\ldots,\overline L_i}_{r_i-1\text{ copies}}).\]
By the induction hypothesis, we obtain 
\[\begin{split}h_{\overline L,\overline{\boldsymbol{L}}_0,\ldots,\overline{\boldsymbol{L}}_d}^{R,\boldsymbol{\alpha}_0,\ldots,\boldsymbol{\alpha}_d}(\check{\mathbb P})&=-\ln|R(s_{0,0},\ldots,s_{d,0})|-\sum_{j\in\{0,\ldots,d\}\setminus\{i\}}\delta_j\ln\|\alpha_{j,0}\wedge\cdots\wedge\alpha_{j,r_j}\|_{j,*,\det}\\
&\qquad\quad -\delta_i\ln\|\overline{\alpha}_{i,0}\wedge\cdots\wedge\overline{\alpha}_{i,r_i-1}\|_{i,*,\operatorname{quot},\det}-\delta_i\ln\|\alpha_{i,r_i}\|_{i,*}\\
&=-\ln|R(s_{0,0},\ldots,s_{d,0})|-\sum_{j=0}^d\delta_j\ln\|\alpha_{j,0}\wedge\cdots\wedge{\alpha_{j,r_j}}\|_{j,*,\det},
\end{split}\]
where the last equality comes from \cite[Proposition 1.2.51]{CMArakelovAdelic}.

In the case where $|\ndot|$ is Archimedean, by \cite[\S1.4.3]{MR1260106} Remark (iii) one has 
\[h_{\overline L,\overline{\boldsymbol{L}}_0,\ldots,\overline{\boldsymbol{L}}_d}^{R,\boldsymbol{\alpha}_0,\ldots,\boldsymbol{\alpha}_d}(\check{\mathbb P})=h_{\overline L,\overline{\boldsymbol{L}}_0,\ldots,\overline{\boldsymbol{L}}_{i-1},\overline{\boldsymbol{L}}_i',\overline{\boldsymbol{L}}_{i+1},\ldots,\overline{\boldsymbol{L}}_d}^{R,\boldsymbol{\alpha}_0,\ldots,\boldsymbol{\alpha}_{i-1},\overline{\boldsymbol{\alpha}}_i,\boldsymbol{\alpha}_{i+1},\ldots,\boldsymbol{\alpha}_d}(\check{\mathbb P}')-\delta_i\Big(\ln\|\alpha_{i,r_i}\|_{i,*}-\frac 12\sum_{\ell=1}^{r_i}\frac{1}{\ell}\Big).\]
Thus the induction hypothesis leads to 
\[\begin{split}h_{\overline L,\overline{\boldsymbol{L}}_0,\ldots,\overline{\boldsymbol{L}}_d}^{R,\boldsymbol{\alpha}_0,\ldots,\boldsymbol{\alpha}_d}(\check{\mathbb P})&=-\ln|R(s_{0,0},\ldots,s_{d,0})|-\hspace{-1.5ex}
\sum_{\begin{subarray}{c}j\in\{0,\ldots,d\}\\j\neq i\end{subarray}}\hspace{-1.5ex}
\delta_j\big(\ln\|\alpha_{j,0}\wedge\cdots\wedge\alpha_{j,r_j}\|_{j,*,\det}-\sigma_{r_j}\big)\\
&\hspace{-8ex} -\delta_i\big(\ln\|\overline{\alpha}_{i,0}\wedge\cdots\wedge\overline{\alpha}_{i,r_i-1}\|_{i,*,\operatorname{quot},\det}-\sigma_{r_i-1}\big) -
\delta_i\Big(\ln\|\alpha_{i,r_i}\|_{i,*}-\frac 12\sum_{\ell=1}^{r_i}\frac{1}{\ell}\Big)\\&=-\ln|R(s_{0,0},\ldots,s_{d,0})|-\sum_{j=0}^d\delta_j\big(\ln\|\alpha_{j,0}\wedge\cdots\wedge\alpha_{j,r_j}\|-\sigma_{r_j}\big),
\end{split}\]
as required.
\end{proof}

\section{Local height of the resultant}\label{Sec: local height of resultant}

The purpose of this subsection is to relate local heights of a projective variety and its resultant. As in the previous section, $v=(k,|\ndot|)$ denotes a complete valued field such that $|\ndot|$ is not trivial. We fix a projective $k$-scheme $X$ and we let $d$ be the Krull dimension of $X$. Let $(E_i)_{i=0}^d$ be a family of finite-dimensional vector spaces over $k$. For each $i\in\{0,\ldots,d\}$, we denote by $r_i:=\dim_k(E_i)-1$ and let $\|\ndot\|_{i}$ be a norm on $E_i$, which is supposed to be either ultrametric or induced by an inner product. Let $f_i:X\rightarrow\mathbb P(E_i)$ be a closed immersion. We pick elements $s
_0,\ldots,s_d$ of $E_0,\ldots,E_d$ respectively, such that   \[\operatorname{div}(s_0|_X),\ldots,\operatorname{div}(s_d|_X)\] intersect properly on $X$. For simplicity of notation, we denote by 
\[\boldsymbol{s}:=(s_0,\ldots,s_d).\]
Let $\check{\mathbb P}:=\mathbb P(E_0^\vee)\times_k\cdots\times_k\mathbb P(E_d^\vee)$, and let \[p:X\times_k\check{\mathbb P}\longrightarrow X\quad\text{and}\quad q:X\times_k\check{\mathbb P}\longrightarrow\check{\mathbb P}\]
be morphisme of projections.
For any $i\in\{0,\ldots,d\}$, let $\pi_i:\check{\mathbb P}\rightarrow\mathbb P(E_i^\vee)$ be the projection to the $i$-th coordinate and let $q_i=\pi_i\circ q$.

For $i\in\{0,\ldots,d\}$, let $\overline L_i$ be $L_i=\pi_i^*({\mathcal O_{E_i^\vee}(1)})$ equipped with the pull-back of the orthogonal quotient metric on $\mathcal O_{E_i^\vee}(1)$ associated with $\|\ndot\|_{i,*}$, and let \[\overline{\boldsymbol{L}}_i:=(\underbrace{\overline{L}_i,\ldots,\overline{L}_i}_{r_i\text{ copies}})\]
and \[(\alpha_{i,0},\boldsymbol{\alpha}_i=(\alpha_{i,1},\ldots,\alpha_{i,r_i}))\] be a basis of $E_i^\vee$ such that \[\alpha_{i,0}(s_i)=1\quad\text{and}\quad \alpha_{i,j}(s_i)=0\text{ for $j\in\{1,\ldots,r_i\}$.}\]  For simplicity, we denote by $R$ the resultant
 \[R_{f_0,\ldots,f_d}^{X,s_0,\ldots,s_d}\]
as in Definition \ref{Def:resultant precise}, considered as a global section of 
\[L=\pi_0^*(L_0)^{\otimes\delta_0}\otimes\cdots\otimes\pi_d^*(L_d)^{\otimes\delta_d},\]
where 
\[\delta_i:=\deg\big(c_1(L_0)\cdots c_1(L_{i-1})c_1(L_{i+1})\cdots c_1(L_d)\cap[X]\big).\] Note that one has $R(s_0,\ldots,s_d)=1$. Moreover, the Cartier divisors \[\operatorname{div}(R), \operatorname{div}(\pi_0^*(\alpha_{0,1})),\ldots, \operatorname{div}(\pi_0^*(\alpha_{0,r_0})),\ldots,\operatorname{div}(\pi_d^*(\alpha_{d,1})),\ldots, \operatorname{div}(\pi_d^*(\alpha_{d,r_d}))\] intersect properly.

\begin{lemm}\label{Lem: h div R h IX}
Assume that the field $k$ is algebraically closed and $X$ is integral. One has
\[h_{\overline{\boldsymbol{L}}_0,\ldots,\overline{\boldsymbol{L}}_d}^{\pi_0^*(\boldsymbol{\alpha}_0),\ldots,\pi_d^*(\boldsymbol{\alpha}_d)}\big(\operatorname{div}(R)\big)=h_{q^*(\boldsymbol{\overline L}_0),\ldots,q^*(\boldsymbol{\overline L}_d)}^{q_0^*(\boldsymbol{\alpha}_0),\ldots,q_d^*(\boldsymbol{\alpha}_d)}(I_X).\]
\end{lemm}
\begin{proof}
The projection $q:I_X\rightarrow\operatorname{div}(R)$ is a birational morphism (see the proof of \cite[Proposition 3.1]{MR1264417}). Hence the equality follows from the induction formula \eqref{Equ: height recursive formula} and \cite[Proposition 2.4.11 (4)]{MR3498148}.
\end{proof}

\begin{defi}\label{Def: epsilon tensor product}
Assume that the absolute value $|\ndot|$ is non-Archimedean. We equip each symmetric power $S^{\delta_i}(E_i^\vee)$ with the $\varepsilon$-symmetric power norm of $\|\ndot\|_{i,*}$, namely the quotient norm of the $\varepsilon$-tensor power of $\|\ndot\|_{i,*}$ (see Remark \ref{Rem: extension of normss} for the definition of the dual norm $\|\ndot\|_{i,*}$).
Recall that the $\varepsilon$-tensor power of the norm $\|\ndot\|_{i,*}$ is  the norm $\|\ndot\|_{i,*,\varepsilon}$ on $(E_i^{\vee})^{\otimes_k\delta_i}$ defined as (see \cite[Definition 1.1.52]{CMArakelovAdelic})
\[\|T\|_{i,*,\varepsilon}=\sup_{\begin{subarray}{c}(t_1,\ldots,t_{\delta_i})\in E_i^{\delta_i}\\
\forall\,j\in\{1,\ldots,\delta_i\},\;t_j\neq 0
\end{subarray}}\frac{|T(t_1,\ldots,t_{\delta_i})|}{\|t_1\|_i\cdots\|t_{\delta_i}\|_i}.\]
We then equip the vector space $S^{\delta_0}(E_0^\vee)\otimes_k\cdots\otimes_k S^{\delta_d}(E_d^\vee)$ with the $\varepsilon$-tensor product of the $\varepsilon$-symmetric power norms, which we denote simply by $\|\ndot\|$. 
\end{defi}

\begin{rema}\label{Rem: extensions scalars}
Note that, by \cite[Definition 1.1.58]{CMArakelovAdelic}, the norm $\|\ndot\|$ also identifies with the quotient norm by the canonical quotient map
\[(E_0^{\vee})^{\otimes_k\delta_0}\otimes_k\cdots\otimes_k(E_d^{\vee})^{\otimes_k\delta_d}\longrightarrow S^{\delta_0}(E_0^\vee)\otimes_k\cdots\otimes_k S^{\delta_d}(E_d^\vee)\]
of the $\varepsilon$-tensor product of $\delta_i$ copies of $\|\ndot\|_{i,*}$, $i\in\{0,\ldots,d\}$. By Propositions 1.3.20 and 1.3.21 of \cite{CMArakelovAdelic}, we obtain that, for any complete valued extension $k'$ of $k$, the norm $\|\ndot\|_{k'}$ on
\[\big(S^{\delta_0}(E_0^\vee)\otimes_k\cdots\otimes_k S^{\delta_d}(E_d^\vee)\big)\otimes_{k}k'\cong S^{\delta_0}(E_{0,k'}^\vee)\otimes_{k'}\cdots\otimes_{k'}S^{\delta_d}(E_{d,k'}^\vee)\]
identifies with the $\varepsilon$-tensor product of $\delta_i$ copies of $\|\ndot\|_{i,k',*}$, $i\in\{0,\ldots,d\}$.
\end{rema}

\begin{lemm}\label{Lem: h div r et h p}
In the case where $|\ndot|$ is non-Archimedean and $k$ is algebraically closed, one has
\begin{equation}\label{Equ: passage a hauteur de P}h_{\overline{\boldsymbol{L}}_0,\ldots,\overline{\boldsymbol{L}}_d}^{\pi_0^*(\boldsymbol{\alpha}_0),\ldots,\pi_d^*(\boldsymbol{\alpha}_d)}(\operatorname{div}(R))=h_{\overline L,\overline{\boldsymbol{L}}_0,\ldots,\overline{\boldsymbol{L}}_d}^{R,\pi_0^*(\boldsymbol{\alpha}_0),\ldots,\pi_d^*(\boldsymbol{\alpha}_d)}(\check{\mathbb P})+\ln\|R\|.\end{equation}
\end{lemm}
\begin{proof}
Let $\xi$ be the Gauss point of $\check{\mathbb P}^{\mathrm{an}}$. It suffices to observe that
\[|R|_{\varphi}(\xi)=\|R\|,\]
where $\varphi$ is tensor product of orthogonal quotient metrics. In fact, if we consider the Veronese-Segre embedding
\[\check{\mathbb P}\longrightarrow\mathbb P(S^{\delta_0}(E_0^\vee))\times_k\cdots\times_k\mathbb P(S^{\delta_d}(E_d^\vee))\longrightarrow\mathbb P(S^{\delta_0}(E_0^\vee)\otimes_k\cdots\otimes_kS^{\delta_d}(E_d^\vee)),\]
then the metric $\varphi$ identifies with the quotient metric induced by $\|\ndot\|$ (see \cite[Proposition 1.1.58]{CMArakelovAdelic}). Moreover, one has
\[\mu_{\overline L_0^{r_0}\cdots \overline L_d^{r_d}}=\operatorname{Dirac}_\xi.\]
Therefore the equality \eqref{Equ: passage a hauteur de P} follows from the induction formula \eqref{Equ: height recursive formula}.
\end{proof}

\begin{lemm}\label{Lem: h div R etc}
In the case where $|\ndot|$ is Archimedean and $k=\mathbb C$, one has
\[\begin{split}h_{\overline{\boldsymbol{L}}_0,\ldots,\overline{\boldsymbol{L}}_d}^{\pi_0^*(\boldsymbol{\alpha}_0),\ldots,\pi_d^*(\boldsymbol{\alpha}_d)}&(\operatorname{div}(R))=h_{\overline L,\overline{\boldsymbol{L}}_0,\ldots,\overline{\boldsymbol{L}}_d}^{R,\pi_0^*(\boldsymbol{\alpha}_0),\ldots,\pi_d^*(\boldsymbol{\alpha}_d)}(\check{\mathbb P})\\&+\int_{\mathbb S_0\times\cdots\times \mathbb S_d}\ln|R(z_0,\ldots,z_d)|\,\eta_{\mathbb S_0}(\mathrm{d}z_0)\otimes\cdots\otimes\eta_{\mathbb S_d}(\mathrm{d}z_d),
\end{split}\]
where $\mathbb S_i$ is the unit sphere of $(E_{i,\mathbb C},\|\ndot\|_{i,\mathbb C})$, and $\mu_{\mathbb S_i}$ is the $U(E_{i,\mathbb C},\|\ndot\|_{i,\mathbb C})$-invariant Borel probability measure on $\mathbb S_i$.
\end{lemm}
\begin{proof}
This is a direct consequence of the induction formula \eqref{Equ: height recursive formula} and Remark \ref{Rem: mesure complexe}.
\end{proof}

\begin{lemm}\label{Lem: h s in h IX}
Assume that the field $k$ is algebraically closed and $X$ is integral. For any $i\in \{0,\ldots,d\}$, we equip $\mathcal O_{E_i}(1)$ with the orthogonal quotient metric induced by $\|\ndot\|_i$, and denote by $M_i'$ the restriction of $\mathcal O_{E_i}(1)$ to $X$ and equip it with the restricted metric. If $|\ndot|$ is non-Archimedean, then one has
\[h_{\overline{M_0'},\ldots,\overline{M_d'}}^{s_0,\ldots,s_d}(X)=h_{q^*(\boldsymbol{\overline L}_0),\ldots,q^*(\boldsymbol{\overline L}_d)}^{q_0^*(\boldsymbol{\alpha}_0),\ldots,q_d^*(\boldsymbol{\alpha}_d)}(I_X) +\sum_{i=0}^d\delta_i\ln\|\alpha_{i,0}\wedge\cdots\wedge\alpha_{i,r_i}\|_{i,*,\det};\]
if $|\ndot|$ is Archimedean, then one has
\begin{multline*}h_{\overline{M_0'},\ldots,\overline{M_d'}}^{s_0,\ldots,s_d}(X)=h_{q^*(\boldsymbol{\overline L}_0),\ldots,q^*(\boldsymbol{\overline L}_d)}^{q_0^*(\boldsymbol{\alpha}_0),\ldots,q_d^*(\boldsymbol{\alpha}_d)}(I_X)+\sum_{i=0}^d\delta_i\big(\ln\|\alpha_{i,0}\wedge\cdots\wedge\alpha_{i,r_i}\|_{i,*,\det}-\sigma_{r_i-1}\big),
\end{multline*}
where 
\[\sigma_{r_i-1}=\frac 12\sum_{m=1}^{r_i-1}\sum_{\ell=1}^m\frac{1}{\ell}.\]
\end{lemm}
\begin{proof}
For $i\in\{0,\ldots,d\}$, let $t_i$ be the global section of $\mathcal O_{E_i}(1)\boxtimes\mathcal O_{E_i^\vee}(1)$ on $\mathbb P(E_i)\times_k\mathbb P(E_i^\vee)$ defining the incidence subscheme. Then $t_i$ corresponds to the restriction of the trace element of $E_i\otimes_kE_i^\vee$ via the Segre embedding  \[\mathbb P(E_i)\times_k\mathbb P(E_i^\vee)\longrightarrow\mathbb P(E_i\otimes_kE_i^\vee).\]
Let $\boldsymbol{t}=(t_0,\ldots,t_{d})$. For any $i\in\{0,\ldots,d\}$, let
\[(s_i,s_{i,1},\ldots,s_{i,r_i})\]
be the dual basis of $(\alpha_{i,j})_{j=0}^{r_i}$. By definition one has
\[t_i=s_i\otimes\alpha_{i,0}+s_{i,1}\otimes\alpha_{i,1}+\cdots+s_{i,r_i}\otimes\alpha_{i,r_i}.\] 
For $i\in\{0,\ldots,d\}$, let $L_i:=q_i^*(\mathcal O_{E_i^\vee}(1))$,  $M_i=p^*(\mathcal O_{E_i}(1)|_X)$ and $N_i= L_i\otimes M_i$. We use two methods to compute the following local height of $X\times\check{\mathbb P}$ (see Notation \ref{Not: simplified notation})
\[h_{\overline{\boldsymbol{N}},\overline{\boldsymbol{L}}_0,\ldots,\overline{\boldsymbol{L}}_d}^{\boldsymbol{t},\boldsymbol{\alpha}_0,\ldots,\boldsymbol{\alpha}_d}(X\times_k\check{\mathbb P}),\]
where $\overline{\boldsymbol{N}}=(\overline N_0,\ldots,\overline{N}_d)$.
 We will show by induction that 
\begin{equation}\label{Equ: height of product}h_{\overline{\boldsymbol{N}},\overline{\boldsymbol{L}}_0,\ldots,\overline{\boldsymbol{L}}_d}^{\boldsymbol{t},\boldsymbol{\alpha}_0,\ldots,\boldsymbol{\alpha}_d}(X\times_k\check{\mathbb P})=h_{\overline{M_0'},\ldots,\overline{M_d'}}^{s_0,\ldots,s_d}(X)-\sum_{i=0}^d\delta_i\ln\|\alpha_{i,0}\wedge\cdots\wedge\alpha_{i,r_i}\|_{i,*,\det}\end{equation}
if $|\ndot|$ is non-Archimedean, and
\begin{equation}\label{Equ: height of product Archimedean}h_{\overline{\boldsymbol{N}},\overline{\boldsymbol{L}}_0,\ldots,\overline{\boldsymbol{L}}_d}^{\boldsymbol{t},\boldsymbol{\alpha}_0,\ldots,\boldsymbol{\alpha}_d}(X\times_k\check{\mathbb P})=h_{\overline{M_0'},\ldots,\overline{M_d'}}^{s_0,\ldots,s_d}(X)-\sum_{i=0}^d\delta_i\big(\ln\|\alpha_{i,0}\wedge\cdots\wedge\alpha_{i,r_i}\|_{i,*,\det}-\sigma_{r_i}\big)\end{equation}
if $|\ndot|$ is Archimedean.
Let $i\in\{0,\ldots,d\}$ be such that $r_i>0$. Let $G_i^\vee=E_i^\vee/k\alpha_{i,r_i}$, $\overline{\boldsymbol{\alpha}}_i=(\overline \alpha_{i,j})_{j=1}^{r_i-1}$, and
\[\check{\mathbb P}'=\mathbb P(E_0)\times_k\cdots\times\mathbb P(E_{i-1})\times_k\mathbb P(G_i)\times_k\mathbb P(E_{i+1})\times_k\cdots\times_k\mathbb P(E_d).\]
Then, with the notation \[\overline{\boldsymbol{L}}_i':=(\underbrace{\overline L_i,\ldots,\overline L_i}_{r_i-1\text{ copies}}),\]
by \eqref{Equ: height recursive formula} one can write $h_{\overline{\boldsymbol{N}},\overline{\boldsymbol{L}}_0,\ldots,\overline{\boldsymbol{L}}_d}^{\boldsymbol{t},\boldsymbol{\alpha}_0,\ldots,\boldsymbol{\alpha}_d}(X\times_k\check{\mathbb P})$ as
\[h_{\overline{\boldsymbol{N}},\overline {\boldsymbol L}_0,\ldots,\overline {\boldsymbol L}_{i-1},\overline {\boldsymbol L}_i',\overline {\boldsymbol L}_{i+1},\ldots,\overline {\boldsymbol L}_{d}}^{\boldsymbol{t},\boldsymbol{\alpha}_0,\ldots,\boldsymbol{\alpha}_{i-1},\overline{\boldsymbol{\alpha}}_i,\boldsymbol{\alpha}_{i+1},\ldots,\boldsymbol{\alpha}_d}(X\times\check{\mathbb P})-\int_{(X\times\check{\mathbb P})^{\mathrm{an}}}\hspace{-2mm}\ln|\alpha_{i,r_i}|\,\mathrm{d}\mu_{\overline N_0\cdots\overline N_d\overline L_0^{r_0}\cdots\overline L_{i-1}^{r_{i-1}}\overline L_i^{r_i-1}\overline L_{i+1}^{r_{i+1}}\cdots\overline L_d^{r_d}},\]
which is equal to 
\[h_{\overline{\boldsymbol{N}},\overline {\boldsymbol L}_0,\ldots,\overline {\boldsymbol L}_{i-1},\overline {\boldsymbol L}_i',\overline {\boldsymbol L}_{i+1},\ldots,\overline {\boldsymbol L}_{d}}^{\boldsymbol{t},\boldsymbol{\alpha}_0,\ldots,\boldsymbol{\alpha}_{i-1},\overline{\boldsymbol{\alpha}}_i,\boldsymbol{\alpha}_{i+1},\ldots,\boldsymbol{\alpha}_d}(X\times\check{\mathbb P})-\int_{(X\times\check{\mathbb P})^{\mathrm{an}}}\ln|\alpha_{i,r_i}|\,\mathrm{d}\mu_{\overline M_0\cdots\overline{M}_{i-1}\overline M_{i+1}\cdots\overline M_d\overline L_0^{r_0}\cdots\overline L_d^{r_d}}.\]
If $|\ndot|$ is non-Archimedean, it identifies with 
\[h_{\overline{\boldsymbol{N}},\overline {\boldsymbol L}_0,\ldots,\overline {\boldsymbol L}_{i-1},\overline {\boldsymbol L}_i',\overline {\boldsymbol L}_{i+1},\ldots,\overline {\boldsymbol L}_{d}}^{\boldsymbol{t},\boldsymbol{\alpha}_0,\ldots,\boldsymbol{\alpha}_{i-1},\overline{\boldsymbol{\alpha}}_i,\boldsymbol{\alpha}_{i+1},\ldots,\boldsymbol{\alpha}_d}(X\times\check{\mathbb P})-\delta_i\ln\|\alpha_{i,r_i}\|_{i,*}\]
In the case where $|\ndot|$ is Archimedean, it equals
\[h_{\overline N,\overline {\boldsymbol L}_0,\ldots,\overline {\boldsymbol L}_{i-1},\overline {\boldsymbol L}_i',\overline {\boldsymbol L}_{i+1},\ldots,\overline {\boldsymbol L}_{d}}^{\boldsymbol{t},\boldsymbol{\alpha}_0,\ldots,\boldsymbol{\alpha}_{i-1},\overline{\boldsymbol{\alpha}}_i,\boldsymbol{\alpha}_{i+1},\ldots,\boldsymbol{\alpha}_d}(X\times\check{\mathbb P})-\delta_i\Big(\ln\|\alpha_{i,r_i}\|_{i,*}-\frac 12\sum_{\ell=1}^{r_i}\frac{1}{\ell}\Big).\]
Hence by induction we obtain \eqref{Equ: height of product} and \eqref{Equ: height of product Archimedean} according to the nature of $|\ndot|$.

Now let $\boldsymbol{t}'=(t_0,\ldots,t_{d-1})$ and $\overline{\boldsymbol{N}}'=(\overline N_0,\ldots,\overline N_{d-1})$, still by \eqref{Equ: height recursive formula} one can write $h_{\overline{\boldsymbol{N}},\overline{\boldsymbol{L}}_0,\ldots,\overline{\boldsymbol{L}}_d}^{\boldsymbol{t},\boldsymbol{\alpha}_0,\ldots,\boldsymbol{\alpha}_d}(X\times_k\check{\mathbb P})$ as
\[\begin{split}&\quad\; h_{\overline{\boldsymbol{N}},\overline{\boldsymbol{L}}_0,\ldots,\overline{\boldsymbol{L}}_d}^{\boldsymbol{t}',\boldsymbol{\alpha}_0,\ldots,\boldsymbol{\alpha}_d}(\operatorname{div}(t_d)) -\int_{(X\times_k\check{\mathbb P})^{\mathrm{an}}}\ln|t_d|\,\mathrm{d}\mu_{\overline N_0\cdots\overline{N}_{d-1}\overline L_0^{r_0}\cdots\overline L_d^{r_d}}\\
&=h_{\overline{\boldsymbol{N}},\overline{\boldsymbol{L}}_0,\ldots,\overline{\boldsymbol{L}}_d}^{\boldsymbol{t}',\boldsymbol{\alpha}_0,\ldots,\boldsymbol{\alpha}_d}(\operatorname{div}(t_d)) -\int_{(X\times_k\check{\mathbb P})^{\mathrm{an}}}\ln|t_d|\,\mathrm{d}\mu_{\overline M_0\cdots\overline M_{d-1}\overline L_0^{r_0}\cdots\overline L_d^{r_d}}
\end{split}\] 
Note that for any element $z\in(X\times_k\check{\mathbb P})^{\mathrm{an}}$ represented by
\[(\beta,x_0,\ldots,x_d)\in E_{d,\widehat{\kappa}(z)}^\vee\times E_{0,\widehat{\kappa}(z)}\cdots\times E_{d,\widehat{\kappa}(z)}\]
one has
\begin{equation}\label{Equ: ln td}\ln|t_d|(z)=\ln\frac{|\beta(x_d)|_{z}}{\|\beta\|_{d,\widehat{\kappa}(z)}\cdot\|x_d\|_{d,\widehat{\kappa}(z)}}.\end{equation}
In the case where $|\ndot|$ is non-Archimedean, this leads to \[\int_{(X\times_k\check{\mathbb P})^{\mathrm{an}}}\ln|t_d|\,\mathrm{d}\mu_{\overline M_0\cdots\overline M_{d-1}\overline L_0^{r_0}\cdots\overline L_d^{r_d}}=0\]
by using \eqref{Equ:Fubini} and 
\[\int_{\mathbb P(E_{d,\widehat{\kappa}(z)}^\vee)^{\mathrm{an}}}\ln|\beta|\,\mathrm{d}\mu_{\overline{\mathcal O_{E_d}(d)}^{r_d}}=\ln\|\beta\|_{d,*,\widehat{\kappa}(x)}.\]
In the case where $|\ndot|$ is Archimedean, by \cite[\S1.4.3]{MR1260106} Remark (iii), \eqref{Equ: ln td} leads to 
\[-\int_{(X\times_k\check{\mathbb P})^{\mathrm{an}}}\ln|t_d|\,\mathrm{d}\mu_{\overline M_0\cdots\overline{M}_{d-1}\overline L_0^{r_0}\cdots\overline L_d^{r_d}}=\frac {\delta_d}2\sum_{\ell=1}^{r_i}\frac{1}{\ell}.\]
Then by induction we obtain 
\begin{equation}\label{Equ: h X times P}h_{\overline{\boldsymbol{N}},\overline{\boldsymbol{L}}_0,\ldots,\overline{\boldsymbol{L}}_d}^{\boldsymbol{t},\boldsymbol{\alpha}_0,\ldots,\boldsymbol{\alpha}_d}(X\times_k\check{\mathbb P})=h_{\overline{\boldsymbol{L}}_0,\ldots,\overline{\boldsymbol{L}}_d}^{\boldsymbol{\alpha}_0,\ldots,\boldsymbol{\alpha}_d}(I_X)\end{equation}
when $|\ndot|$ is non-Archimedean and
\begin{equation}\label{Equ: h X times P2}h_{\overline{\boldsymbol{N}},\overline{\boldsymbol{L}}_0,\ldots,\overline{\boldsymbol{L}}_d}^{\boldsymbol{t},\boldsymbol{\alpha}_0,\ldots,\boldsymbol{\alpha}_d}(X\times_k\check{\mathbb P})=h_{\overline{\boldsymbol{L}}_0,\ldots,\overline{\boldsymbol{L}}_d}^{\boldsymbol{\alpha}_0,\ldots,\boldsymbol{\alpha}_d}(I_X)+\frac 12\delta_i\sum_{i=0}^d\sum_{\ell=1}^{r_i}\frac{1}{\ell}\end{equation}
when $|\ndot|$ is Archimedean. Combining \eqref{Equ: h X times P} with \eqref{Equ: height of product}, and \eqref{Equ: h X times P2} with \eqref{Equ: height of product Archimedean}, we obtain the result.
\end{proof}

\begin{theo}\label{Thm: equality of local height} For any $i\in \{0,\ldots,d\}$, we equip $\mathcal O_{E_i}(1)$ with the orthogonal quotient metric induced by $\|\ndot\|_i$, and denote by $M_i'$ the restriction of $\mathcal O_{E_i}(1)$ to $X$ and equip it with the restricted metric. In the case where $|\ndot|$ is non-Archimedean, one has
\[h_{\overline{M_0'},\ldots,\overline{M_d'}}^{s_0,\ldots,s_d}(X)=\ln\|R\|,\]
where the norm $\|\ndot\|$ was introduced in Definition \ref{Def: epsilon tensor product}. In the case where $|\ndot|$ is Archimedean, one has
\begin{multline*}h_{\overline{M_0'},\ldots,\overline{M_d'}}^{s_0,\ldots,s_d}(X)=\int_{\mathbb S_0\times\cdots\times \mathbb S_d}\hspace{-1mm}\ln|R(z_0,\ldots,z_d)|\,\eta_{\mathbb S_0}(\mathrm{d}z_0)\otimes\cdots\otimes\eta_{\mathbb S_d}(\mathrm{d}z_d)+\frac 12\sum_{i=0}^d\delta_i\sum_{\ell=1}^{r_i}\frac{1}{\ell},\end{multline*}
where $\mathbb S_i$ is the unit sphere of $(E_{i,\mathbb C},\|\ndot\|_{i,\mathbb C})$, and $\eta_{\mathbb S_i}$ is the $U(E_{i,\mathbb C},\|\ndot\|_{i,\mathbb C})$-invariant Borel probability measure on $\mathbb S_{i,\sigma}$.
\end{theo}
\begin{proof} By Remark \ref{Rem: resultant extension}, \[R\otimes 1\in (S^{\delta_0}(E_0^\vee)\otimes_k\cdots\otimes_kS^{\delta_d}(E_d^\vee))\otimes_k\mathbb C_k\]
is the resultant of $X_{\mathbb C_k}$ with respect to $f_{0,\mathbb C_k},\ldots,f_{d,\mathbb C_k}$, which takes value $1$ at $(s_0,\ldots,s_d)$. Therefore, by extension of scalars, we may assume without loss of generality that $k$ is algebraically closed and $X$ is integral.

We treat firstly the non-Archimedean case.
By Lemma \ref{Lem: h s in h IX}, one has
\[h_{\overline{M_0'},\ldots,\overline{M_d'}}^{s_0,\ldots,s_d}(X)=h_{q^*(\boldsymbol{\overline L}_0),\ldots,q^*(\boldsymbol{\overline L}_d)}^{q_0^*(\boldsymbol{\alpha}_0),\ldots,q_d^*(\boldsymbol{\alpha}_d)}(I_X)+\sum_{i=0}^d\delta_i\ln\|\alpha_{i,0}\wedge\cdots\wedge\alpha_{i,r_i}\|_{i,*,\det}.\]
By Lemma \ref{Lem: h div R h IX}, this is also equal to
\[h_{\overline{\boldsymbol{L}}_0,\ldots,\overline{\boldsymbol{L}}_d}^{\pi_0^*(\boldsymbol{\alpha}_0),\ldots,\pi_d^*(\boldsymbol{\alpha}_d)}(\operatorname{div}(R))+\sum_{i=0}^d\delta_i\ln\|\alpha_{i,0}\wedge\cdots\wedge\alpha_{i,r_i}\|_{i,*,\det}.\]
By Lemma \ref{Lem: h div r et h p}, it is equal to 
\[h_{\overline L,\overline{\boldsymbol{L}}_0,\ldots,\overline{\boldsymbol{L}}_d}^{R,\pi_0^*(\boldsymbol{\alpha}_0),\ldots,\pi_d^*(\boldsymbol{\alpha}_d)}(\check{\mathbb P})+\ln\|R\|+\sum_{i=0}^d\delta_i\ln\|\alpha_{i,0}\wedge\cdots\wedge\alpha_{i,r_i}\|_{i,*,\det}.\]
By Proposition \ref{Pro: hight of multi-projective} and the relation (see Definition \ref{Def:resultant precise})
\[R(s_0,\ldots,s_d)=1,\]
we obtain
\[h_{\overline{M_0'},\ldots,\overline{M_d'}}^{s_0,\ldots,s_d}(X)=\ln\|R\|.\]

The case where $|\ndot|$ is Archimedean is quite similar. We have
\[\begin{split}
&\quad\; h_{\overline{M_0'},\ldots,\overline{M_d'}}^{s_0,\ldots,s_d}(X)=h_{q^*(\boldsymbol{\overline L}_0),\ldots,q^*(\boldsymbol{\overline L}_d)}^{q_0^*(\boldsymbol{\alpha}_0),\ldots,q_d^*(\boldsymbol{\alpha}_d)}(I_X)+\sum_{i=0}^d\delta_i\big(\ln\|\alpha_{i,0}\wedge\cdots\wedge\alpha_{i,r_i}\|_{i,*,\det}-\sigma_{r_i-1}\big)
\\&=h_{\overline{\boldsymbol{L}}_0,\ldots,\overline{\boldsymbol{L}}_d}^{\pi_0^*(\boldsymbol{\alpha}_0),\ldots,\pi_d^*(\boldsymbol{\alpha}_d)}(\operatorname{div}(R))+\sum_{i=0}^d\delta_i\big(\ln\|\alpha_{i,0}\wedge\cdots\wedge\alpha_{i,r_i}\|_{i,*,\det}-\sigma_{r_i-1}\big)\\
&=h_{\overline L,\overline{\boldsymbol{L}}_0,\ldots,\overline{\boldsymbol{L}}_d}^{R,\pi_0^*(\boldsymbol{\alpha}_0),\ldots,\pi_d^*(\boldsymbol{\alpha}_d)}(\check{\mathbb P})\\
& \qquad\quad+\int_{\mathbb S_0\times\cdots\times \mathbb S_d}\ln|R(z_0,\ldots,z_d)|\,\eta_{\mathbb S_0}(\mathrm{d}z_0)\otimes\cdots\otimes\eta_{\mathbb S_d}(\mathrm{d}z_d)\\
&\qquad\quad+\sum_{i=0}^d\delta_i\big(\ln\|\alpha_{i,0}\wedge\cdots\wedge\alpha_{i,r_i}\|_{i,*,\det}-\sigma_{r_i}\big)+\frac 12\sum_{i=0}^d\delta_i\sum_{\ell=1}^{r_i}\frac{1}{\ell}\\
&=\int_{\mathbb S_0\times\cdots\times \mathbb S_d}\ln|R(z_0,\ldots,z_d)|\,\eta_{\mathbb S_0}(\mathrm{d}z_0)\otimes\cdots\otimes\eta_{\mathbb S_d}(\mathrm{d}z_d)+\frac 12\sum_{i=0}^d\delta_i\sum_{\ell=1}^{r_i}\frac{1}{\ell},
\end{split}\]
where the first equality comes from Lemma \ref{Lem: h s in h IX}, the second one  from Lemma \ref{Lem: h div R h IX}, the third one  from Lemma \ref{Lem: h div R etc}, and the last one from Proposition \ref{Pro: hight of multi-projective}.
\end{proof}

\begin{rema} Note that the result of Theorem \ref{Thm: equality of local height} does not depend on the choice of the vectors $\boldsymbol{\alpha}_0,\ldots,\boldsymbol{\alpha}_d$. 
If we are only interested in the equalities in the theorem, we could choose $\boldsymbol{\alpha}_0,\ldots,\boldsymbol{\alpha}_d$ carefully to make the computation simpler. However, the formulae in the lemmas proving the theorem are of their proper interest, especially in the computations of height of homogeneous hypersurfaces in multi-projective spaces, and hence are worth to be detailed.
\end{rema}

\begin{prop}\label{prop:invariance:local:intersection:field:extension}
Assume that the absolute value $|\ndot|$ is non-Archimedean. Let $K$ be an extension of $k$, on which the absolute value extends. We assume that $K$ is complete with respect to the extended absolute value. Let $X$ be a projective scheme over $\Spec k$, $d$ be the dimensional of $X$, and $\overline D_i=(D_i,g_i)$ be a family of integrable metrised Cartier divisors, where $i\in\{0,\ldots,d\}$, such that $D_0,\ldots,D_d$ intersect properly. For each $i\in\{0,\ldots,d\}$, let $\overline D_{i,K}:=(D_{i,k},g_{i,k})$. Then the following equality holds:
\begin{equation}\label{Equ: intersection invariant by extenison of scalars}
(\overline D_0\cdots\overline D_d)_{(k,|\ndot|)}=(\overline D_{0,K}\cdots\overline D_{d,K})_{(K,|\ndot|)}.
\end{equation}
\end{prop}
\begin{proof}
{\bf Step 1}: In this step, we assume that $D_0,\ldots,D_d$ are very ample, and, for each $i\in\{0,\ldots,d\}$, there exist a positive integer $m_i$ and an ultrametric norm $\|\ndot\|_i$ on $E_i=H^0(X,\mathcal O_X(m_iD_i))$, such that $\varphi_{g_i}$ identifies with the quotient metric induced by $\|\ndot\|_i$.

For each $i\in\{0,\ldots,d\}$, let $f_i:X\rightarrow\mathbb P(E_i)$ be the canonical closed embedding. Note that $\mathcal O_X(m_iD_i)\cong f_i^*(\mathcal O_{E_i}(1))$. In order to simplify the notation, we let $L_i$ be the line bundle $\mathcal O_X(m_iD_i)$ and $s_i$ be the canonical regular meromorphic section of $L_i$. Let $R$ be the resultant
\[R_{f_0,\ldots,f_d}^{X,s_0,\ldots,s_d},\]
which is considered as an element of
\[S^{\delta_0}(E_0^\vee)\otimes_k\cdots\otimes_kS^{\delta_d}(E_d^\vee),\]
and \[\delta_i=(D_0\cdots D_{i-1}D_{i+1}\cdots D_d).\]
Then, by Theorem \ref{Thm: equality of local height}, the equality 
\[(\overline D_0\cdots\overline D_d)_{(k,|\ndot|)}=\ln\|R\|\]
holds, where $\|{\ndot}\|$ denotes the $\varepsilon$-tensor product of $\varepsilon$-tensor powers of $\|\ndot\|_{i,*}$.
Similarly, by Remarks  \ref{Rem: resultant extension} and \ref{Rem: extensions scalars}, one has
\[(\overline{D}_{0,K},\ldots,\overline D_{d,K})_{(K,|\ndot|)}=\ln\|R\otimes 1\|_{K}.\]
By \cite[Proposition 1.3.1 (1)]{CMArakelovAdelic}, one has $\|R\otimes 1\|_K=\|R\|$. Hence the equality \eqref{Equ: intersection invariant by extenison of scalars} follows.

{\bf Step 2}: In this step, we still assume that $D_0,\ldots,D_d$ are very ample. However, the Green functions $g_0,\ldots,g_d$ are only supposed to be plurisubharmonic.  

For any $i\in\{0,\ldots,d\}$ and any positive integer $m$, let $g_i^{(m)}$ be the Green function associated with the quotient metric $\varphi_{g_i}^{(m)}$ as in Definition \ref{Def: Fubini-Study}, and let $\overline D_i^{(m)}=(D_i,g_i^{(m)})$. By Proposition \ref{Pro: convergence of varphin}, we obtain that, for any $i\in\{0,\ldots,d\}$,
\begin{equation}\label{Equ: convergence of Green functions}\lim_{m\rightarrow+\infty}\sup_{x\in X^{\mathrm{an}}}|g_i^{(m)}-g_i|(x)=0,\end{equation}
Therefore, by Corollary \ref{cor:intersection:trivial:divisor} (see also \S\ref{Sec: Local intersection number over a general field}), we obtain 
\begin{equation}\label{Equ: intersection m}\lim_{m\rightarrow+\infty}(\overline D_0^{(m)}\cdots\overline D_d^{(m)})_{(k,|\ndot|)}=(\overline D_0\cdots\overline D_d)_{(k,|\ndot|)}.\end{equation}
Moreover, \eqref{Equ: convergence of Green functions} leads to 
\[\lim_{m\rightarrow+\infty}\sup_{x\in X_K^{\mathrm{an}}}|g_{i,K}^{(m)}-g_{i,K}|(x)=0.\]
Hence, similarly to \eqref{Equ: intersection m}, we have
\[\lim_{m\rightarrow+\infty}(\overline D_{0,K}^{(m)}\cdots\overline D_{d,K}^{(m)})_{(K,|\ndot|)}=(\overline D_{0,K}\cdots\overline D_{d,K})_{(K,|\ndot|)}\]
Note that, by \cite[Proposition 1.3.16]{CMArakelovAdelic}, $g_{i,K}^{(m)}$ is also the Green function associated with a quotient metric. Therefore, by the result in Step 1, we obtain that
\[(\overline D_0^{(m)}\cdots\overline D_d^{(m)})_{(k,|\ndot|)}=(\overline D_{0,K}^{(m)}\cdots\overline D_{d,K}^{(m)})_{(K,|\ndot|)}\]
for any $m$, so that, by passing to limit when $m\rightarrow+\infty$, we obtain \eqref{Equ: intersection invariant by extenison of scalars}.

{\bf Step 3:} We now treat the general case. For each $i\in\{-1,0,\ldots,d\}$, we consider the following condition $(C_r)$:
\begin{quote}
\it For any $i\in\{0,\ldots,d\}$ such that $1\leqslant i\leqslant r$, the Cartier divisor $D_i$ is very ample and the Green function $g_i$ is plurisubharmonic. 
\end{quote}
We will show by inverted induction on $r$ that, under the condition $(C_r)$, the equality \eqref{Equ: intersection invariant by extenison of scalars} holds. Note that the initial case where $r=d$ is proved in Step 2. We suppose that the equality \eqref{Equ: intersection invariant by extenison of scalars} is true under the condition $(C_r)$ and will prove it under the condition $(C_{r-1})$. Since $\overline D_r$ is integrable, there exists very ample Cartier divisors $A_r'$ and $A_r''$, and plurisubharmonic Green functions $h_r'$ and $h_r''$ of $A_r'$ and $A_r''$, respectively, such that 
\[(D_r,g_r)=(A_r',h_r')-(A_r'',h_r'').\]
By Claim \ref{Claim: decomposition} (see also Remark \ref{Rem: further explanation to Claim}), there exists a very ample Cartier divisor $B_r$ such that \[(D_0,\ldots,D_{r-1},B_r+A_r',D_{r+1},\ldots,D_d) \in\mathcal {IP}_X^{(d)}.\] Since $\mathcal {IP}_X^{(d)}$ is a multilinear subset of $\operatorname{Div}(X)^{n+1}$, we obtain that \[(D_0,\ldots,D_{r-1},B_r+A_r'',D_{r+1},\ldots,D_d) \in\mathcal {IP}_X^{(d)}.\]  We pick arbitrarily a plurisubharmonic Green function $l_r$ on $B_r$. Let 
\[\overline D_r'=(B_r+A_r',l_r+h_r'),\quad \overline D_r''=(B_r+A_r'',l_r+h_r'')\]Then the induction hypothesis shows that 
\begin{gather*}(\overline D_0\cdots\overline D_{r-1}\overline D_r'\overline D_{r+1}\cdots\overline D_d)_{(k,|\ndot|)}=(\overline D_{0,K}\cdots\overline D_{r-1,K}\overline D_{r,K}'\overline D_{r+1,K}\cdots\overline D_{d,K})_{(K,|\ndot|)},\\
(\overline D_0\cdots\overline D_{r-1}\overline D_r''\overline D_{r+1}\cdots\overline D_d)_{(k,|\ndot|)}=(\overline D_{0,K}\cdots\overline D_{r-1,K}\overline D_{r,K}''\overline D_{r+1,K}\cdots\overline D_{d,K})_{(K,|\ndot|)}.\end{gather*}
Taking the difference, we obtain \eqref{Equ: intersection invariant by extenison of scalars}
\end{proof}

\begin{rema}\label{Rem: extension of scalars archimedean}
If $K$ is a subfield of $\mathbb C_k$, the assertion of Proposition~\ref{prop:invariance:local:intersection:field:extension}
is obvious by its definition (cf. Definition~\ref{Def: intersection pairing general case}). In particular, the statement of Proposition \ref{prop:invariance:local:intersection:field:extension} is also true when $|\ndot|$ is Archimedean.
Proposition~\ref{prop:invariance:local:intersection:field:extension} guarantees the invariance of intersection number under any field extension.
\end{rema}

\section{Trivial valuation case}
In this section, we fix a field $k$ and equip it with the trivial absolute value $|\ndot|$, namely $|a|=1$ for any $a\in k^{\times}$. Let $K=k(T)$ be the field of rational functions over $k$, and $u$ be a positive constant such that $u\not=1$. 
By Lemma \ref{lem:Gauss:lemma}, there exists a non-Archimedean absolute value $|\ndot|_u$ on $K$ which extends the above absolute value $|\ndot|$ on $k$, such that, 
\[\forall\,f=a_0+a_1T+\cdots+a_nT^n\in k[T],\quad |f|_u=\max_{i\in\{0,\ldots,n\}}|a_i|u^i.\]  
Note that $|\ndot|_u$ is not trivial.

\begin{defi}\label{Def: intersection in trivial valuation case}
Let $X$ be a projective scheme of dimension $d$ over $\Spec k$. If $\overline D_i=(D_i,g_i)$, $i\in\{0,\ldots,d\}$, is a family of integrable metrized Cartier divisors, such that $D_0,\ldots,D_d$ intersect properly. We denote by $(\overline D_0\cdots\overline D_d)_{(k,|\ndot|)}$ the intersection number
\[((D_{0,K},g_{0,K})\cdots(D_{d,K},g_{d,K}))_{(K,|\ndot|_u)}\]
\end{defi}

\begin{enonce}{Notation and assumptions}\label{Notation and assumptions}\rm
Let $((E_i,\|\ndot\|_i))_{i=0}^d$ be a family of finite-dimensional ultrametrically normed vector space over $k$. For any $i\in\{0,\ldots,d\}$, let $r_i=\dim_{k}(E_i)-1$, $f_i:X\rightarrow\mathbb P(E_i)$ be a closed immersion, and $s_i$ an element of $E_i$, viewed as a global section of $\mathcal O_{E_i}(1)$. We assume that the restriction of $s_i$ to $X$ defines a regular meromorphic section of $L_i:=\mathcal O_{E_i}(1)|_X$ and that the Cartier divisors
\[D_i=\operatorname{div}(s_i|_X), \quad i\in\{0,\ldots, d\}\]
intersect properly. We equip each $D_i$ with the Green function associated with the quotient metric induced by $\|\ndot\|_i$. Let $R$ be the resultant
\[R=R_{f_0,\ldots,f_d}^{X,s_0,\ldots,s_d}\in S^{\delta_0}(E_0^\vee)\otimes_k\cdots\otimes_kS^{\delta_d}(E_d^\vee),\]
where 
\[\delta_i=(D_0\cdots D_{i-1}D_{i+1}\cdots D_d).\]
\end{enonce}

\begin{prop} Under Notation and assumptions \ref{Notation and assumptions},
the following equality holds
\begin{equation}\label{Equ: local intersection number trivial valuation}(\overline D_0\cdots\overline D_d)_{(k,|\ndot|)}=\ln\|R\|,\end{equation}
where $\|\ndot\|$ denotes the $\varepsilon$-tensor product of $\varepsilon$-tensor power norms of $\|\ndot\|_{i,*}$. 
\end{prop}
\begin{proof}
Under the isomorphism of $K$-vector spaces
\[(S^{\delta_0}(E_0^\vee)\otimes_k\cdots\otimes_kS^{\delta_d}(E_d^\vee))\otimes_kK\cong S^{\delta_0}(E_{0,K}^\vee)\otimes_K\cdots\otimes_KS^{\delta_d}(E_{d,K}^\vee),\]
the element $R\otimes 1$ coincides with the resultant (see Remark  \ref{Rem: resultant extension})
\[R_{f_{0,K},\ldots,f_{d,K}}^{X_K,s_{0}\otimes 1,\ldots,s_{d}\otimes 1}.\]
By Theorem \ref{Thm: equality of local height} and Remark \ref{Rem: extensions scalars}, one has 
\[(\overline D_0\cdots\overline D_d)_{(k,|\ndot|)}=\ln\|R\otimes 1\|_K.\]
By \cite[Proposition 1.3.1 (1)]{CMArakelovAdelic}, one has $\|R\otimes 1\|_K=\|R\|$. Hence we obtain the equality \eqref{Equ: local intersection number trivial valuation}.
\end{proof}

\begin{coro}
Let $X$ be a projective scheme of dimension $d$ over $\Spec k$. If $\overline D_i=(D_i,g_i)$, $i\in\{0,\ldots,d\}$, is a family of integrable metrized Cartier divisors, such that $D_0,\ldots,D_d$ intersect properly. Then the intersection number $(\overline D_0\cdots\overline D_d)_{(k,|\ndot|)}$ does not depend on
the choice of $u$.
\end{coro}
\begin{proof}
By the multi-linearity of the intersection number, it suffices to treat the case where all Cartier divisors $D_i$ are very ample and all $g_i$ are plurisubharmonic. Moreover, by Proposition \ref{Pro: convergence of varphin} and Corollary \ref{cor:intersection:trivial:divisor} we can further reduce the problem to the case of Notation and assumptions \ref{Notation and assumptions}. In that case the assertion follows from \eqref{Equ: local intersection number trivial valuation}. 
\end{proof}

\begin{rema}
By using Remark~\ref{rem:properties:local:intersection:general:field}, one has the following properties.

\begin{enumerate}[label=\rm(\arabic*)]
\item The set $\widehat{\IP}_X$ forms a symmetric multi-linear subset of the group $\operatorname{\widehat{Int}}(X)^{d+1}$. Moreover, the function of local intersection number
\[\big( (D_0, g_0) \cdots (D_d, g_d) \big)\longmapsto\big( (D_0, g_0) \cdots (D_d, g_d) \big)_v\]
form a symmetric multi-linear map from $\widehat{\IP}_X$ to $\mathbb R$. 
\item 
Let $\pi : Y \to X$ be a surjective morphism of geometrically integral projective schemes over $k$.
We set $e = \dim X$ and $d = \dim Y$. 
Let $(D_0, g_0), \ldots, (D_d, g_d)$ be integrable metrized Cartier divisors on $X$ such that
$(\pi^*(D_0), \ldots, \pi^*(D_d)) \in \IP_Y$.
Then one has the following:
\begin{enumerate}[label=\rm(\roman*)]
\item If $d > e$, then $(\pi^*(D_0, g_0) \cdots \pi^*(D_d, g_d))_v = 0$.
\item If $d = e$ and $(D_0,\ldots,D_d) \in \IP_X$, then 
\[(\pi^*(D_0, g_0) \cdots \pi^*(D_d, g_d))_v = (\deg \pi)((D_0, g_0) \cdots (D_d, g_d))_v.\]
\end{enumerate}

\item Let $f$ be a regular meromorphic function on $X$ and $(D_1,g_1),\ldots,(D_d,g_d)$ be integrable metrized Cartier divisors on $X$ such that $(\operatorname{div}(f),D_1,\ldots,D_d)\in\IP_X$. Then
\[\big(\widehat{\operatorname{div}}(f) \cdot (D_1, g_1) \cdots (D_{d}, g_{d})\big)_v = 0. \]
Note that $-\log|f|(x^{\mathrm{an}}) = 0$ for any $x \in X_{(0)}$ in Remark~\ref{rem:properties:local:intersection:general:field} 
because $|\ndot|$ is trivial.
\end{enumerate}

Let $(L_0, \varphi_0), \ldots, (L_d, \varphi_d)$ be a family of integrable metrized invertible $\OO_X$-modules.
By the property (3), the local intersection number $\big( (L_0, \varphi_0) \cdots (L_d, \varphi_d) \big)_v$ is well-defined.
\end{rema}

\begin{rema}
In \cite{Chen_Moriwaki2020}, an intersection product of metrized divisors has been introduced in the setting of curves over a trivially valued field $(k,|\ndot|)$. Let $X$ be a regular projective curve over $\operatorname{Spec} k$. Recall that the Berkovich space $X^{\operatorname{an}}$ is an infinite tree
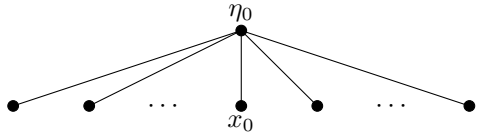
\vspace{3mm}
\begin{center}
\begin{tikzpicture}
\filldraw(0,1) circle (2pt) node[align=center, above]{$\eta_0$};
\filldraw(-3,0) circle (2pt) ;
\draw (-1,0) node{$\cdots$};
\filldraw(-2,0) circle (2pt) ;
\filldraw(-0,0) circle (2pt) node[align=center, below]{$x_0$} ;
\filldraw(1,0) circle (2pt) ;
\draw (2,0) node{$\cdots$};
\filldraw(3,0) circle (2pt) ;
\draw (0,1) -- (0,0);
\draw (0,1) -- (-3,0);
\draw (0,1) -- (1,0);
\draw (0,1) -- (-2,0);
\draw (0,1) -- (3,0);
\end{tikzpicture}
\end{center}
\vspace{3mm}
where the root point $\eta_0$ corresponds to the generic point of $X$ together with the trivial absolute value on $\kappa(\eta)$, and each leaf $x_0$ corresponds to the closed point $x$ together with the trivial absolute value on $\kappa(x)$. Moreover, each  branch $\mathopen{]}\eta_0,x_0\mathclose{[}$ is parametrized by  $\mathrm{]}0,+\infty\mathrm{[}$, where $t\in\mathrm{]}0,+\infty\mathrm{[}$ corresponds to the generic point $\eta$ together with the absolute value \[|\ndot|_{x,t}=\exp(-t\operatorname{ord}_x(\ndot)).\]
We denote by $t(\ndot):X^{\mathrm{an}}\rightarrow [0,+\infty]$ the parametrization map, where $t(\eta_0)=0$ and $t(x_0)=+\infty$.
Let $D$ be a Cartier divisor on $X$. Recall that a Green function $g$ of $D$ is of the form 
\[g=g_D+\varphi_g,\] 
where $g_D$ is the canonical Green function of $D$, which is defined as 
\[g_D(\xi)=\operatorname{ord}_x(D)t(\xi),\]
and $\varphi_g$ is a continuous real-valued function on $X^{\mathrm{an}}$ (which is hence bounded since $X^{\mathrm{an}}$ is compact). Then, the intersection number of two integrable metrized Cartier divisor $\overline D_0=(D_0,g_0)$ and $\overline D_1=(D_1,g_1)$ has been defined as
\begin{equation}\label{Equ: alternative intersection product}g_1(\eta_0)\deg(D_0)+g_0(\eta_0)\deg(D_1)-\sum_{x\in X^{(1)}}[\kappa(x):k]\int_0^{+\infty}\varphi'_{g_0\circ\xi_x}(t)\varphi'_{g_1\circ\xi_x}(t)\,\mathrm{d}t,
\end{equation}
where $X^{(1)}$ is the set of closed points of $X$, $\xi_x:[0,+\infty]\rightarrow[\eta_0,x_0]$ is the map sending $t\in[0,+\infty]$ to the point in $[\eta_0,x_0]$ of parameter $t$, and the function $\varphi_{g_1\circ\xi_x}'(\ndot)$ should be considered as right-continuous version of the Radon-Nikodym density of the function $\varphi_{g_1\circ\xi_x}(\ndot)$ with respect to the Lebesgue measure.

Let $(L, \varphi_0)$ and $(L_1, \varphi_1)$ be integrable metrized invertible $\OO_X$-modules.
By \cite[Remark~7.3]{Chen_Moriwaki2020}, the above local intersection number with respect to $(L, \varphi_0)$ and $(L_1, \varphi_1)$ is well-defined.
To destinguish this intersection number with the intersection number defined in Definition~\ref{Def: intersection in trivial valuation case},
it is denoted by $((L_0, \varphi_0) \cdot (L_1, \varphi_2))_v'$.

 Let $(E,\|\ndot\|)$ be a ultrametrically normed vector space and $f:X\rightarrow\mathbb P(E)$ be a closed embedding. Let $L$ be the pull-back of the universal invertible sheaf $\mathcal O_E(1)$ by $f$. We equip $L$ with the quotient metric induced by the norm $\|\ndot\|$. Let $E\otimes_{k}\kappa(\eta)\rightarrow L(\eta)$ be the universal quotient $\kappa(\eta)$-linear map and 
\[L(\eta)^\vee\longrightarrow E_{\kappa(\eta)}^\vee\]
be the dual linear map. By this injective linear map we identify $L(\eta)^\vee$ with a one-dimensional vector subspace of $E_{\kappa(\eta)}^\vee$. Consider an element $s\in E\setminus\{0\}$ which is viewed as a global section of $\mathcal O_E(1)$. We suppose that $s$ does not vanish at $f(\eta)$. Let $x$ be a closed point of $X$, $t\in\mathopen{[}0,+\infty\mathclose{[}$, and $\xi$ the element of $\mathopen{[}\eta_0,x_0\mathclose{[}$ having $t$ as its parameter, one has 
\[|s|(\xi)=\frac{|\beta(s)|_{x,t}}{\|\beta\|_{x,t,*}},\]
where $\beta$ is an arbitrary non-zero element of $L(\eta)^\vee$, and $\|\ndot\|_{x,t}$ is the norm on $E_{\kappa(\eta)}=E\otimes_k\kappa(\eta)$ constructed from $\|\ndot\|$ by extension of scalars to $(\kappa(\eta),|\ndot|_{x,t})$. In particular, if we pick an orthogonal basis of $(E,\|\ndot\|)$ of the form $\{s_0,s_1,\ldots,s_r\}$ with $s_0=s$, and pick a vector $\beta\in L(\eta)^\vee$ of the form  
\[\beta=\lambda_0s_0^\vee+\lambda_1s_1^\vee+\cdots+\lambda_rs_r^\vee,\quad (\lambda_0,\ldots,\lambda_r)\in\kappa(\eta)^r,\]
then one has
\[|s|(\xi)=\frac{|\lambda_0|_{x,t}}{
\displaystyle\max_{i\in\{0,\ldots,r\}}|\lambda_i|_{x,t}\cdot\|s_i\|^{-1}}=
\min_{i\in\{0,\ldots,r\}}
\|s_i\|\exp(\operatorname{ord}_x(\lambda_i)-\operatorname{ord}_x(\lambda_0))^t.\]
Therefore, the Green function $g$ of $\operatorname{div}(s)$ corresponding to the quotient metric is given by 
\[g(\xi)=-\min_{i\in\{0,\ldots,r\}}\Big(t\big(\operatorname{ord}_x(\lambda_i)-\operatorname{ord}_x(\lambda_0)\big)+\ln\|s_i\|\Big)\text{ on $\mathopen{]}\eta_0,x_0\mathclose{[}$}.\] 
Moreover, one has
\[g(\eta_0)=-\min_{\begin{subarray}{c}i\in\{0,\ldots,r\}\\
\lambda_i\neq 0
\end{subarray}}\ln\|s_i\|\]
and
\begin{equation}\label{Equ: computation of varphi g}\varphi_g(\xi)=-\min_{i\in\{0,\ldots,r\}}
\bigg(t\Big(\operatorname{ord}_x(\lambda_i)-\min_{j\in\{0,\ldots,r\}}\operatorname{ord}_x(\lambda_j)\Big)+\ln\|s_i\|\bigg).\end{equation}

\end{rema}

We now illustrate the comparison of \eqref{Equ: alternative intersection product} and the local intersection product introduced in Definition \ref{Def: intersection in trivial valuation case} in the particular case. 

\begin{prop}
Let $(E,\|\ndot\|)$ be a finite-dimensional 
ultrametrically normed vector space over $k$. 
Let $X$ be a regular projective curve over $k$ and $L$ be an invertible $\OO_X$-module.
We assume that there is a surjective homomorphism $E \otimes_k \OO_X \to L$.
Let $\varphi$ be the Fubini-Study metric of $L$ induced by the above homomorphism and $(E,\|\ndot\|)$.
If either (1) $X = \PP^1_k$ and $L=\OO_{\PP^1}(1)$, or (2) $\dim E = 2$,  then
\[
((L, \varphi) \cdot (L, \varphi))_v = ((L, \varphi) \cdot (L, \varphi))_v'.
\]
\end{prop}
\begin{proof}
(1) (the case where $X = \PP^1$ and $L = \OO_{\PP^1}(1)$)
Let $E \to H^0(X, L)$ be the natural homomorphism, which is surjective because
$E \otimes_k \OO_X \to L$ is surjective. 
Let $\|\ndot\|'$ be the quotient norm of $H^0(X, L)$ by $\|\ndot\|$ of $E$.
Then the Fubini-Study metric of $L$ induced by $(H^0(X, L), \|\ndot\|')$ coincides with $\varphi$,
so that we may assume that $E = H^0(X, L)$ and $X = \PP(E)$.

Let $\{s_0,s_1\}$ be an orthogonal basis of $E$.
Let $\overline D_0=(D_0,g_0)$ and $\overline D_1=(D_1,g_1)$ be the metrized Cartier divisors $\operatorname{\widehat{\mathrm{div}}}(s_0)$ and $\operatorname{\widehat{\mathrm{div}}}(s_1)$, respectively.
We need to prove
\begin{equation}\label{Equ: equality of two definitions}\begin{split}(\overline D_0\cdot\overline D_1)_{(k,|\ndot|)}=g_1(\eta_0)&\deg(D_0)+g_0(\eta_0)\deg(D_1)\\
&-\sum_{x\in\mathbb P(E)^{(1)}}[\kappa(x):k]\int_0^{+\infty}\varphi_{g_0\circ\xi_x}'(t)\varphi_{g_1\circ\xi_x}'(t)\,\mathrm{d}t.
\end{split}\end{equation}
Note that $\{s_0,s_1\}$ forms an orthogonal basis of $E\otimes_kK$ (see \cite[Proposition 1.3.13]{CMArakelovAdelic}). Therefore, Proposition \ref{Pro: height projective space} leads to
\[(\overline D_0\cdot\overline D_1)_{(k,|\ndot|)}=-\ln\|s_0\wedge s_1\|_{K,\det}=-\ln\|s_0\|_K-\ln\|s_1\|_K=-\ln\|s_0\|-\ln\|s_1\|,\]
where the last equality comes from \cite[Proposition 1.3.1]{CMArakelovAdelic}. 
Note that $\varphi = \varphi_{g_0}=\varphi_{g_1}$, which we denote by $\varphi_g$.

Note that the field of rational functions of $\mathbb P(E)$ is given by the subfield of $k(E)$ generated by $\tau=s_1/s_0$. Moreover, the universal one-dimensional $k(\tau)$-linear subspace
\[L^\vee\otimes_kk(\tau)\longrightarrow E^\vee\otimes_kk(\tau)\]
is spanned by the vector 
\[\beta=s_0^\vee+\tau s_1^\vee,\]
where $\{s_0^\vee,s_1^\vee\}$ is the dual basis of $\{s_0,s_1\}$. Let $x_0$ and $x_1$ be the vanishing point of $s_0$ and $s_1$, respectively. These are rational points of $\mathbb P(E)$. By \eqref{Equ: computation of varphi g}, we obtain that, for $x\in\mathbb P(E)^{(1)}\setminus\{x_0,x_1\}$, one has
\begin{equation}\label{Equ: varphi g xi x}\varphi_{g\circ\xi_x}(t)=-\min\{\ln\|s_0\|,\ln\|s_1\|\},\end{equation}
and
\begin{gather*}
\varphi_{g\circ\xi_{x_0}}(t)=-\min\{t+\ln\|s_0\|,\ln\|s_1\|\},\\
\varphi_{g\circ\xi_{x_1}}(t)=-\min\{\ln\|s_0\|,t+\ln\|s_1\|\}.
\end{gather*}
Therefore, one obtains that   
\begin{gather*}
\varphi_{g\circ\xi_{x_0}}'=-\indic_{[0,\max\{\ln\frac{\|s_1\|}{\|s_0\|},0\})},\\
\varphi_{g\circ\xi_{x_1}}'=-\indic_{[0,\max\{\ln\frac{\|s_0\|}{\|s_1\|},0\})},
\end{gather*}
which leads to
\[-\sum_{x\in\mathbb P(E)^{(1)}}[\kappa(x):k]\int_0^{+\infty}\varphi_{g_0\circ\xi_x}'(t)^2\,\mathrm{d}t=-\big|\ln\|s_0\|-\ln\|s_1\|\big|.\]
Moreover, \eqref{Equ: varphi g xi x} also implies that 
\[g_0(\eta_0)=g_1(\eta_0)=-\min\{\ln\|s_0\|,\ln\|s_1\|\}. \]
Hence the right hand side of \eqref{Equ: equality of two definitions} is equal to 
\[-2\min\{\ln\|s_0\|,\ln\|s_1\|\}-\big|\ln\|s_0\|-\ln\|s_1\|\big|=-\ln\|s_0\|-\ln\|s_1\|,\]
as desired.

\bigskip
(2) (the case where $\dim E = 2$)
Then one has a finite surjective morphism $f : X \to \PP(E)$ such that $f^*(\OO_{\PP(E)}(1)) = L$.
Let $\psi$ be the Fubini-Study metric induced by $E \otimes \OO_{\PP(E)} \to \OO_{\PP(E)}(1)$ and $\|\ndot\|$.
Then $f^*(\psi) = \varphi$, that is, $(L, \varphi) = f^*(\OO_{\PP(E)}(1), \psi)$. Therefore one can see
\[
\begin{cases}
((L, \varphi) \cdot (L,\varphi))_v = \deg(f) ((\OO_{\PP(E)}(1), \psi) \cdot (\OO_{\PP(E)}(1), \psi))_v, \\
((L, \varphi) \cdot (L,\varphi))_v' = \deg(f) ((\OO_{\PP(E)}(1), \psi) \cdot (\OO_{\PP(E)}(1), \psi))_v'.
\end{cases}
\]
Thus the assertion follows from (1).
\end{proof}

\begin{rema}
The above proposition suggests that \eqref{Equ: alternative intersection product} should be equal to the intersection number introduced in Definition \ref{Def: intersection in trivial valuation case}. We expect that an explicit computation of the resultant in the projective curve case would establish such an equality by using Theorem \ref{Thm: equality of local height}. 
\end{rema}


\chapter{Global intersection number}

Let $K$ be a field and $S = (K, (\Omega, {\mathcal A}, \nu), \phi)$ be an adelic curve the underlying field of which is $K$. For any $\omega\in\Omega$, we denote by $K_\omega$ the completion of $K$ with respect to $|\ndot|_\omega$.
We assume that, either the $\sigma$-algebra $\mathcal A$ is discrete, or there exists a countable subfield $K_0$ of $K$ which is dense in each $K_\omega$, $\omega\in\Omega$.
Let $X$ be a $d$-dimensional projective scheme over $K$. For any $\omega\in\Omega$,  let $X_\omega$ be the fiber product $X\times_{\Spec K}\Spec K_\omega$. Note that the morphism $\Spec K_\omega\rightarrow\Spec K$ is flat. Hence the morphism of projection $X_\omega\rightarrow X$ is also flat (see \cite[$\text{IV}_1.(2.1.4)$]{EGA}).

\section{Reminder on adelic vector bundles}
\label{Sec: Reminder on adelic vector bundles}

\begin{defi}
Let $E$ be a finite-dimensional vector space over $K$. We call \emph{norm family} of $E$ any family $\xi=(\|\ndot\|_\omega)_{\omega\in\Omega}$, where each $\|\ndot\|_\omega$ is a norm on $E_{K_\omega}:=E\otimes_KK_\omega$. If for any $\omega\in\Omega$, the norm $\|\ndot\|_\omega$ is either ultrametric (when $\omega$ is non-Archimedean) or induced by an inner product (when $\omega$ is Archimedean), we say that the norm family $\xi$ is \emph{Hermitian}.

If $\xi=(\|\ndot\|_\omega)_{\omega\in\Omega}$ and $\xi'=(\|\ndot\|_\omega)_{\omega\in\Omega}$ are two norm families of $E$, we define the \emph{local distance function} of $\xi$ and $\xi'$ as the function
\[(\omega\in\Omega)\longmapsto d_\omega(\xi,\xi'):=\sup_{s\in E_{K_\omega}\setminus\{0\}}\Big|\ln\|s\|_\omega-\ln\|s\|_\omega'\Big|.\] 
\end{defi}

\begin{exem}\label{Exe: norm family of a basis}
Let $\boldsymbol{e}=(e_i)_{i=1}^r$ be a basis of $E$ over $K$. For any $\omega\in\Omega$, we let $\|\ndot\|_{\boldsymbol{e},\omega}$ be the norm on $E_{K_\omega}$ such that, for $(\lambda_1,\ldots,\lambda_r)\in K_\omega^r$,  
\[ \|\lambda_1e_1+\cdots+\lambda_re_r\|_{\boldsymbol{e},\omega}=\begin{cases}
\max\{|\lambda_1|_\omega,\ldots,|\lambda_r|_\omega\},&\text{$\omega$ is non-Archimedean},\\
(|\lambda_1|_\omega^2+\cdots+|\lambda_r|_\omega^2)^{1/2}, &\text{$\omega$ is Archimedean}.
\end{cases}\]
Then $(\|\ndot\|_{\boldsymbol{e},\omega})_{\omega\in \Omega} $ forms a Hermitian norm family of $E$, which we denote by $\xi_{\boldsymbol{e}}$.
\end{exem}

\begin{defi}
Let $E$ be a finite-dimensional vector space over $K$ and $\xi$ be a norm family on $E$. We say that $\xi$ is \emph{measurable} if for any $s\in E$ the function 
\[(\omega\in\Omega)\longrightarrow\|s\|_\omega\]
is $\mathcal A$-measurable. We say that the family $\xi$ is \emph{strongly dominated} if there exists a basis $\boldsymbol{e}$ of $E$ such that the local distance function
\[(\omega\in\Omega)\longmapsto d_\omega(\xi,\xi_{\boldsymbol{e}})\]
is bounded from above by an integrable function. If $\xi$ is measurable and strongly dominated, we say that $(E,\xi)$ is a \emph{strongly adelic vector bundle}.
We refer the readers to \cite[\S4.1.4]{CMArakelovAdelic} for more details about this definition, and also to the Proposition 4.1.24 (1.b) for the measurability of the dual norm family of $\xi$ under our assumption on the adelic curve. 

If $(E,\xi)$ is a strongly adelic vector bundle, for any non-zero element $s$ of $E$, the function
\[(\omega\in\Omega)\longmapsto \ln\|s\|_{\omega}\]  
is integrable. We denote by $\widehat{\deg}_\xi(s)$, or simply by $\widehat{\deg}(s)$ if there is no ambiguity on the norm family, the integral
\[-\int_{\Omega}\ln\|s\|_\omega\,\nu(\mathrm{d}\omega),\]
called \emph{Arakelov degree} of $s$ (with respect to $\xi$).
\end{defi}

\begin{defi}
Let $X$ be a projective $K$-scheme and $L$ be an invertible $\mathcal O_X$-module. For any $\omega\in\Omega$, we denote by $L_\omega$ the pull-back of $L$ by the morphism of projection $X_\omega\rightarrow X$. We call \emph{metric family} of $L$ and family $\varphi=(\varphi_\omega)_{\omega\in\Omega}$, where each $\varphi_\omega$ is a continuous metric on $L_\omega$ (see Definition \ref{Defi: continuous metric}). Note that the dual metrics $(\varphi_\omega^\vee)_{\omega\in\Omega}$ form a metric family on the dual invertible $\mathcal O_X$-module $L^\vee$, which we denote by $\varphi^\vee$. If $L_1$ and $L_2$ are invertible $\mathcal O_X$-modules, and $\varphi_1$ and $\varphi_2$ are metric families on $L_1$ and $L_2$, respectively, then the metrics $(\varphi_{1,\omega}\otimes\varphi_{2,\omega})_{\omega\in\Omega}$ form a metric family of $L_1\otimes L_2$, which we denote by $\varphi_1\otimes\varphi_2$. 

If $\varphi$ and $\varphi'$ are two metric metrics of the same invertible $\mathcal O_X$-module $L$, we define the \emph{local distance function} between $\varphi$ and $\varphi'$ as the function
\[(\omega\in\Omega) \longmapsto d_{\omega}(\varphi, \varphi') := \sup_{x\in X_{\omega}^{\operatorname{an}}}\bigg|\ln\frac{|\ndot|_{\varphi_\omega}(x)}{|\ndot|_{\varphi_\omega'}(x)}\bigg|\]
\end{defi}

\begin{rema}
In the case where $X$ is the spectrum of a finite extension $K'$ of $K$, an invertible $\mathcal O_X$-module $L$ can be considered as a one-dimensional vector space over $K'$, and a metric family on $L$ identifies with a norm family of $L$ if we consider the adelic curve $S\otimes_KK'$.
\end{rema}

\begin{defi}
Let $f:Y\rightarrow X$ be a projective $K$-morphism of projective $K$-schemes. Let $L$ be an invertible $\mathcal O_X$-module, equipped with a metric family $\varphi=(\varphi_\omega)_{\omega\in\Omega}$. For any $\omega\in\Omega$, let $f_\omega:Y_\omega\rightarrow X_\omega$ be the $K_\omega$-morphism induce by $f$ by extension of scalars. Then, for any $\omega\in\Omega$, the metric $\varphi_\omega$ induces by pull-back a continuous metric $f_\omega^*(\varphi_\omega)$ on $f_\omega^*(L_\omega)$ such that, for any $y\in Y_\omega^{\operatorname{an}}$ and any $\ell\in L_\omega(f^{\operatorname{an}}(y))$, one has
\[|f_\omega^*(\ell)|_{f_\omega^*(\varphi_\omega)}(y)=|\ell|_{\varphi_\omega}(f^{\operatorname{an}}(y)).\]
We denote by $f^*(\varphi)$ the metric family $(f_\omega^*(\varphi_\omega))_{\omega\in\Omega}$ and call it the \emph{pull-back of $\varphi$ by $f$}. In the case where $f$ is an immersion, $f^*(\varphi)$ is also called \emph{restriction of $\varphi$}.
\end{defi}

\begin{exem}\label{Exe: quotient metric family}
A natural example of metric family is the \emph{quotient metric family} induced by a norm family. Denote by $\pi:X\rightarrow\Spec K$ the structural morphism. Let $E$ be a finite-dimensional vector space over $K$ and $f:\pi^*(E)\rightarrow L^{\otimes n}$ be a surjective homomorphism of $\mathcal O_X$-modules, where $n$ is a positive integer. For any $\omega\in\Omega$, the homomorphism $f$ induces by pull-back a surjective homomorphism of $\mathcal O_{X_\omega}$-modules $f_\omega:\pi_{K_\omega}^*(E)\rightarrow L_\omega$. 
Assume given a norm family $\xi=(\|\ndot\|_\omega)_{\omega\in\Omega}$ of $E$. We denote by $\varphi_{\xi}$ the metric family of $L$ consisting of quotient metrics associated with $\|\ndot\|_{\omega}$ (see Example \ref{Exa: Fubini-Study} \ref{Item: quotient metrics}), and call it the \emph{quotient metric family} induced by $\xi$.

Assume that the norm family $\xi$ is Hermitian. For each $\omega\in\Omega$, let $\varphi_{\xi,\omega}^{\operatorname{ort}}$ be the orthogonal quotient  metric  induced by $\|\ndot\|_\omega$ (see Definition \ref{Def: Fubini-Study metric}). Note that this metric coincides with $\varphi_{\xi,\omega}$ when $|\ndot|_\omega$ is non-Archimedean or $K_\omega$ is complex. The metric family $\varphi_{\xi}^{\operatorname{ort}}$ is called \emph{orthogonal quotient metric family} induced by $\xi$.
\end{exem}

\begin{exem}
Let $X$ be a projective $K$-scheme, $L$ be an invertible $\mathcal O_X$-module, and $\varphi=(\varphi_\omega)_{\omega\in\Omega}$ be a metric family on $L$. Let $K'/K$ be an algebraic extension of the field $K$, and \[S\otimes K'=(K',(\Omega',\mathcal A',\nu'),\phi')\] be the corresponding algebraic covering of the adelic curve $S$ (see \S\ref{Sec: algebraic coverings}). Recall that $\Omega'$ is defined as $\Omega\times_{M_K,\phi}M_{K'}$, where $M_K$ and $M_{K'}$ are the sets of all absolute values of $K$ and of $K'$, respectively.  

Let $X'$ be the fiber product $X\times_{\Spec K}\Spec K'$ and $L'$ be the pull-back of $L$ on $X'$. If $\omega'$ is an element of $\Omega'$ and $\omega$ is the image of $\omega'$ in $\Omega$ by the projection map 
\[\Omega'=\Omega\times_{M_{K},\phi}M_{K'}\longrightarrow\Omega,\]
then one has
\[X'_{\omega'}:=X'\times_{\Spec K'}\Spec{K'_{\omega'}}\cong (X\times_KK_\omega)\times_{K_\omega}K'_{\omega'}.\]
Moreover, the pull-back of $L_\omega$ on $X'_{\omega'}$ identifies with $L'_{\omega'}$. We denote by $p_{\omega'}$ the morphism of projection from $X'_{\omega'}$ to $X_\omega$. Then the map \[p_{\omega'}^{\natural}:(X'_{\omega'})^{\operatorname{an}}\longrightarrow X_\omega^{\mathrm{an}},\]
sending any point $x'=(j(x'),|\ndot|_{x'})$ to the pair consisting of the scheme point $p_{\omega'}(j(x'))$ of $X_\omega$ and the restriction of $|\ndot|_{x'}$ on the residue field of $p_{\omega'}(j(x'))$, is continuous (see \cite[Proposition 2.1.17]{CMArakelovAdelic}), where $j:(X_{\omega'}')^{\mathrm{an}}\rightarrow X'_{\omega'}$ denotes the map sending a point in the analytic space to its underlying scheme point. Therefore, the continuous metric $\varphi_\omega$ induces by composition with $p^{\natural}$ a continuous metric $\varphi_{\omega'}$ such that, for any $x'\in (X'_{\omega'})^{\operatorname{an}}$ and any $\ell\in L_\omega(p^{\natural}(x'))$, one has
\[\forall\,a\in\widehat{\kappa}(x'),\quad |a\otimes\ell|_{\varphi_{\omega'}}(x')=|a|_{x'}\cdot|\ell|_{\varphi_\omega}(x). \]
Therefore, $(\varphi_{\omega'})_{\omega'\in\Omega'}$ forms a metric family of $L'$ which we denote by $\varphi_{K'}$.
\end{exem}

\begin{defi}\label{Def: adelic line bundle}
Let $L$ be an invertible $\mathcal O_X$-module and $\varphi=(\varphi_\omega)_{\omega\in\Omega}$ be a metric family of $L$.
\begin{enumerate}[label=\rm(\arabic*)]
\item We say that $\varphi$ is \emph{dominated} if there exist invertible $\mathcal O_X$-modules $L_1$ and $L_2$, respectively equipped with metric families $\varphi_1$ and $\varphi_2$, which are quotient metric families associated with dominated norms families, such that $L\cong L_1\otimes L_2^\vee$ and that the local distance function \[(\omega\in\Omega)\longmapsto d_\omega(\varphi,\varphi_1\otimes\varphi_2^\vee)\]
is bounded from above by a $\nu$-integrable function (see \cite[\S6.1.1]{CMArakelovAdelic});
\item We say that $\varphi$ is \emph{measurable} if the following conditions are satisfied (see \cite[\S6.1.4]{CMArakelovAdelic}):
\begin{enumerate}[label=\rm(2.\roman*)] 
\item for any closed point $P$ of $X$, the norm family $P^*(\varphi)$ of $P^*(L)$ is measurable,
\item for any $\xi\in X^{\mathrm{an}}$ (where we consider the trivial absolute value on $K$ in the construction of $X^{\mathrm{an}}$) whose associated scheme is of dimension $1$ and such that the exponent\footnote{Since the schematic point associated with $\xi$ is of dimension $1$, the absolute value $|\ndot|_\xi$ is discrete and hence is of the form $|\ndot|_\xi=\exp(- t\operatorname{ord}_\xi(\ndot))$, where the (surjective) map $\operatorname{ord}_\xi(\ndot):\widehat{\kappa}(\xi)\rightarrow\mathbb Z\cup\{+\infty\}$ is the discrete valuation corresponding to the absolute value $|\ndot|_\xi$. The non-negative real number $t$ is called the \emph{exponent} of the absolute value $|\ndot|_\xi$.} of the absolute value $|\ndot|_\xi$ is rational, and for any $\ell\in L\otimes_{\mathcal O_X}\widehat{\kappa}(\xi)$, the function 
\[(\omega\in\Omega_0)\longmapsto |\ell|_{\varphi_\omega}(\xi)\]
is measurable, where $\Omega_0$ is the subset of $\omega\in\Omega$ such that $|\ndot|_\omega$ is trivial, and we consider the restriction of the $\sigma$-algebra $\mathcal A$ to $\Omega_0$.
\end{enumerate}
\end{enumerate}
If $\varphi$ is both dominated and measurable, we say that the pair $\overline L=(L,\varphi)$ is an \emph{adelic line bundle}.
\end{defi}

\begin{prop}\label{Pro: measurable direct image}
Let $\pi:X\rightarrow\Spec K$ be a projective scheme over $\Spec K$, $L$ be an invertible $\mathcal O_X$-module, $\varphi$ be a metric family of $L$, and $E=H^0(X,L)$. We equip $E$ with a norm family $\xi=(\|\ndot\|_{\omega})_{\omega\in\Omega}$. Consider the following norm family $\xi'=(\|\ndot\|_{\omega}')_{\omega\in\Omega}$ defined as
\[\forall\,s\in H^0(X_\omega,L_\omega),\quad\|s\|_{\omega}':=\max\bigg\{\sup_{x\in X_\omega^{\operatorname{an}}}|s|_{\varphi_\omega}(x),\|s\|_\omega\bigg\}.\]
Then one has the following:
\begin{enumerate}[label=\rm(\arabic*)]
\item\label{Item: xi prime measurable} If $\varphi$ and $\xi$ are both measurable, then $\xi'$ is also measurable. 
\item\label{Item: xi prime dominated} If $\varphi$ is dominated and $\xi$ is  strongly dominated, then $\xi'$ is strongly dominated.
\end{enumerate}
\end{prop}
\begin{proof}
\ref{Item: xi prime measurable} For any $\omega\in\Omega$, we let $\|\ndot\|_{\varphi_\omega}$ be the seminorm on $E\otimes_KK_\omega=H^0(X_\omega,L_\omega)$ defined as 
\[\forall\,s\in H^0(X_\omega,L_\omega),\quad\|s\|_{\varphi_\omega}:=\sup_{x\in X_\omega^{\operatorname{an}}}|s|_{\varphi_\omega}(x).\]
By \cite[Propositions 6.1.20 and 6.1.26]{CMArakelovAdelic}, for any $s\in H^0(X,L)$, the function 
\[(\omega\in\Omega)\longmapsto\|s\|_{\varphi_\omega}\]
is measurable. Therefore the function
\[(\omega\in\Omega)\longmapsto\|s\|_{\omega}'=\max\{\|s\|_{\varphi_\omega},\|s\|_\omega\}\]
is also measurable once the norm family $\xi$ is measurable. 

\ref{Item: xi prime dominated} We may assume without loss of generality that there exists a basis $\boldsymbol{e}=(e_i)_{i=1}^r$ of $E$ such that, for any $\omega\in\Omega$
\[\forall\,(\lambda_1,\ldots,\lambda_r)\in K_\omega^r,\quad \|\lambda_1e_1+\cdots+\lambda_re_r\|_{\omega}=\max_{i\in\{1,\ldots,r\}}|\lambda_i|_\omega. \]
By \cite[Remark 6.1.17]{CMArakelovAdelic}, for any $s\in H^0(X,L)$, the function 
\[(\omega\in\Omega)\longmapsto \ln\|s\|_{\varphi_\omega}\]
is bounded from above by an integrable function. Let $A:\Omega\rightarrow\mathbb R_{\geqslant 0}$ be a positive integrable function on $\Omega$ such that 
\[\forall\,\omega\in\Omega,\quad \max_{i\in\{1,\ldots,r\}}\ln\|e_i\|_{\varphi_\omega}\leqslant A(\omega).\]
For any $\omega\in\Omega\setminus\Omega_\infty$ and any $(\lambda_1,\ldots,\lambda_r)\in K_\omega^r$, one has
\[\ln\|\lambda_1e_1+\cdots+\lambda_re_r\|_\omega\leqslant \ln\|\lambda_1e_1+\cdots+\lambda_re_r\|_\omega'\leqslant\max_{i\in\{1,\ldots,r\}}(\ln|\lambda_i|_\omega+\ln\|e_i\|_\omega').
\]
Note that $\|e_i\|=1$ and hence
\[\ln\|e_i\|_\omega'=\max\{\ln\|e_i\|_{\varphi_\omega},\ln(1)\}\leqslant A(\omega).\]
Therefore one has 
\[d(\|\ndot\|_\omega,\|\ndot\|_{\omega}')\leqslant A(\omega).\]
In the case where $\omega\in\Omega_\infty$, for any $(\lambda_1,\ldots,\lambda_r)\in K_\omega^r$ one has
\[\ln\|\lambda_1e_1+\cdots+\lambda_re_r\|_\omega\leqslant \ln\|\lambda_1e_1+\cdots+\lambda_re_r\|_\omega'\leqslant\max_{i\in\{1,\ldots,r\}}\ln|\lambda_i|_\omega+A(\omega)+\ln(r).\]
Finally we obtain that 
\[\forall\,\omega\in\Omega,\quad d_{\omega}(\xi,\xi')\leqslant A(\omega)+\ln(r)\indic_{\Omega_\infty}(\omega).\]
Hence the norm family $\xi'$ is strongly dominated (see \cite[Proposition 3.1.2]{CMArakelovAdelic} for the fact that $\nu(\Omega_\infty)$ is finite).
\end{proof}

\begin{lemm}\label{lem:measurable:integrable:adelic:algebraic:cover}
Let $S = (K, (\Omega, \mathcal A, \nu), \phi)$ be an adelic curve, $K'$ be an algebraic extension of $K$ and
$S_{K'} = S \otimes_K K' = (K', (\Omega', \mathcal A', \nu'), \phi')$. Let $f$ be a function on $\Omega$.
Then one has the following:
\begin{enumerate}[label=\rm(\arabic*)]
\item $f$ is measurable if and only if $f \circ \pi_{K'/K}$ is measurable.
\item $f$ is integrable if and only if $f \circ \pi_{K'/K}$ is integrable.
\end{enumerate}
\end{lemm}

\begin{proof}
Clearly we may assume that $f$ is non-negative, so that
it is a consequence of \cite[Proposition~3.4.8 and Proposition~3.4.9]{CMArakelovAdelic}.
\end{proof}

\section{Integrability of local intersection numbers}
\label{Sec: Integrability of local intersection numbers}

In this section, we fix a projective $K$-scheme $X$. Let $d$ be the  dimension of $X$.

\begin{defi}\label{def:adelic:Cartier:divisor}
Let $D$ be a Cartier divisor on $X$. For any $\omega\in\Omega$, let $D_\omega$ be the pull-back of $D$ by the morphism of projection $X_\omega\rightarrow X$, which is well defined since the morphism of projection $X_\omega\rightarrow X$ is flat (see Remark \ref{Rem: pull back} and Definition \ref{Def: pull back}). We call \emph{Green function family} of $D$ any family $(g_\omega)_{\omega\in\Omega}$ parametrized by $\Omega$, where each $g_\omega$ is a Green function of $D_\omega$. We denote by $\varphi_g$ the metric family $(|\ndot|_{g_\omega})_{\omega\in\Omega}$ of $\mathcal O_X(D)$, where $|\ndot|_{g_\omega}$ is the continuous metric on $\mathcal O_{X_\omega}(D_\omega)$ induced by the Green function $g_\omega$ (see Remark \ref{Rem: metric induced by a green function}). If the metric family $\varphi_g$ is measurable, we say that the Green function family $g$ is \emph{measurable}. If the metric $\varphi_g$ is dominated, we say that the Green function family $g$ is \emph{dominated}. We refer to Definition \ref{Def: adelic line bundle} for the dominancy and measurability of metrics. If $g$ is both dominated and measurable, we say that $(D,g)$ is an \emph{adelic Cartier divisor}.

Let $D$ be an invertible $\mathcal O_X$-module and $g$ be a Green function family of $D$. If $D$ is ample and all metrics in the family $\varphi_g$ are semi-positive, we say that the Green function family $g$ is \emph{semi-positive}. We say that  $(D, g)$ is \emph{integrable} if there exist ample Cartier divisors $D_1$ and $D_2$, together with semi-positive Green function families $g_1$ and $g_2$ of $D_1$ and $D_2$ respectively, such that $D=D_1-D_2$ and $g=g_1-g_2$. Similarly, we say that an adelic line bundle $(L,\varphi)$ is \emph{integrable} if there exists ample invertible $\mathcal O_X$-modules $L_1$ and $L_2$, and metric families consisting of semi-positive metrics $\varphi_1$ and $\varphi_2$ on $L_1$ and $L_2$, respectively, such that $L=L_2\otimes L_1^\vee$ and $\varphi=\varphi_2\otimes\varphi_1^\vee$.

Let $D_0,\ldots,D_d$ be a family of Cartier divisors, which intersect properly. For any $i\in\{0,\ldots,d\}$, let $g_i$ be a Green function family of $D_i$ such that $(D_i,g_i)$ is integrable. Then, for any $\omega\in\Omega$, a local intersection number
\[((D_{0,\omega},g_{0,\omega}),\ldots,(D_{d,\omega},g_{d,\omega}))_{(K_\omega,|\cdot|_\omega)}\]
has been introduced in Definition \ref{def:local:intersection}, which we denote by 
\[(\overline D_0\cdots\overline D_d)_\omega\]
for simplicity. Thus the local intersection numbers define a function 
\[(\omega\in\Omega)\longrightarrow (\overline D_0\cdots\overline D_d)_\omega.\]
\end{defi}

\begin{defi}
Let $D_1$ and $D_2$ be Cartier divisors on $X$, and $g_1$ and $g_2$ be Green function families of $D_1$ and $D_2$, respectively. We say that $(D_1,g_1)$ and $(D_2,g_2)$ are \emph{linearly equivalent} and we note 
\[(D_1,g_1)\sim(D_2,g_2)\] if $\mathcal O_X(D_1)$ is isomorphic to $\mathcal O_X(D_2)$ and if there exists an isomorphism of $\mathcal O_X$-modules $\mathcal O_X(D_1)\rightarrow\mathcal O_X(D)$ which identifies the metric $\varphi_{g_1}$ to $\varphi_{g_2}$.
\end{defi}

\begin{prop}\label{prop:reduction:same:class} Assume that, for all Cartier divisors $E_0,\ldots,E_d$ which intersect properly, and measurable (resp. dominated) Green function families $h_0,\ldots,h_d$ of $E_0,\ldots,E_d$ respectively, such that all $(E_i,h_i)$ are integrable and linearly equivalent, the function of local intersection number
\[(\omega\in\Omega)\longmapsto(\overline E_0\cdots\overline E_d)_\omega\]
is measurable (resp. dominated). Then, for all Cartier divisors $D_0,\ldots,D_d$ which intersect properly and measurable (resp. dominated) Green function families $g_0,\ldots,g_d$ of $D_0,\ldots,D_d$ respectively, such that all $(D_i,g_i)$ are integrable (but not necessarily linearly equivalent), the function of local intersection number
\[(\omega\in\Omega)\longmapsto (\overline D_0\cdots\overline D_d)_\omega\]
is measurable (resp. dominated).
\end{prop}

\begin{proof}
First of all, by Lemma~\ref{lem:measurable:integrable:adelic:algebraic:cover}, we may assume that
$K$ is algebraically closed.
By Lemma~\ref{lem:IP:projective:ample},
we can choose a matrix \[( D_{i,j} )_{(i,j)\in\{0,\ldots,d\}^2}\]
consisting of Cartier divisors on $X$ such that $(D_{i_0,0}, \ldots, D_{i_d,d}) \in \IP_X$ for any
$(i_0, \ldots, i_d) \in \{ 0, \ldots, d \}^{d+1}$, and that $D_{i,j} \sim D_i$ for all $(i, j) \in \{0, \ldots, d \}^2$.
Let $g_{i,j}$ be a family of integrable Green functions of $D_{i,j}$ such that $(D_{i,j}, g_{i,j}) \sim (D_i, g_i)$.
By Proposition~\ref{pro:domain},
\begin{multline*}
\sum_{\sigma \in \mathfrak{S}(\{0, \ldots, d\})} \big(\overline{D}_{0,\sigma(0)} \cdots \overline{D}_{d,\sigma(d)}\big)_\omega =
\sum_{\sigma \in \mathfrak{S}(\{0, \ldots, d\})} \big(\overline{D}_{\sigma(0),0} \cdots \overline{D}_{\sigma(d),d}\big)_\omega\\
= \sum_{\emptyset \not= I \subseteq \{ 0, \ldots, d \}} (-1)^{(d+1) - \operatorname{card}(I)} \Bigg(\Big(\sum\nolimits_{i \in I} \overline{D}_{i,0} \Big) \cdots \Big(\sum\nolimits_{i \in I} \overline{D}_{i,d}\Big) \Bigg)_\omega,
\end{multline*}
where $\overline{D}_{i,j} = (D_{i,j}, g_{i,j})$.
Note that $\sum_{i \in I} \overline{D}_{i,a} \sim \sum_{i \in I} \overline{D}_{i,b}$, so that
\[
(\omega \in \Omega) \mapsto \sum_{\emptyset \not= I \subseteq \{ 0, \ldots, d \}} (-1)^{(d+1) - \operatorname{card}(I)} \Bigg(\Big(\sum\nolimits_{i \in I} \overline{D}_{i,0} \Big) \cdots \Big(\sum\nolimits_{i \in I} \overline{D}_{i,d}\Big) \Bigg)_\omega
\]
is measurable (resp. dominant) by our assumption. Moreover, by Proposition~\ref{prop:intersection:principal:div},
for each $\sigma \in \mathfrak{S}(\{0, \ldots, d\})$, there is an integrable function $A_{\sigma}$ on $\Omega$ such that
\[
(\overline{D}_{0,\sigma(0)} \cdots \overline{D}_{d,\sigma(d)})_\omega = 
(\overline{D}_0 \cdots \overline{D}_d)_\omega + A_{\sigma}(\omega).
\]
Thus the assertion follows. Note that $\int_{\Omega} A_{\sigma}(\omega)\, \nu(d\omega) = 0$
if $S$ is proper.
\end{proof}

\begin{theo}\label{Thm: measurability non archimedean}Assume that $\Omega_\infty=\emptyset$. 
Let $(L_i)_{i=0}^d$ be a family of invertible $\mathcal O_X$-modules. For each $i\in\{0,\ldots,d\}$, let $s_i$ be a regular meromorphic section of $L_i$ and $D_i=\operatorname{div}(s_i)$. We suppose that $D_0,\ldots,D_d$ intersect properly. For any $i\in\{0,\ldots,d\}$, let $\varphi_i=(\varphi_{i,\omega})_{\omega\in\Omega}$ be a measurable metric family on $L_i$ such that $(L_{i,\omega},\varphi_{i,\omega})$ is integrable, and let $g_i=(g_{i,\omega})_{\omega\in\Omega}$ be the family of Green functions of $D_i$ corresponding to $\varphi_i$. 
Then the 
function of local intersection numbers
\begin{equation}\label{Equ: local intersection measurable}(\omega\in\Omega)\longrightarrow (\overline D_0\cdots\overline D_d)_\omega\end{equation}
is $\mathcal A$-measurable.
\end{theo}
\begin{proof}  
By Lemma~\ref{lem:measurable:integrable:adelic:algebraic:cover}, we may assume that
$K$ is algebraically closed.
By using Proposition~\ref{prop:multilinear:symmetric:semiample:case}, we may further assume that $L_0,\ldots,L_d$ are very ample.
For any $i\in\{0,\ldots,d\}$, we denote by $\delta_i$ the intersection number
\[\deg(c_1(L_0)\cdots c_1(L_{i-1})c_1(L_{i+1})\cdots c_1(L_d)\cap[X]).\]
We introduce, for each $r\in\{-1,\ldots,d\}$, then following condition $(C_r)$:
\begin{quote}
\it\hskip\parindent For each $i\in\{0,\ldots,d\}$ such that $0\leqslant i\leqslant r$, there exist a positive integer $m_i$ and a measurable Hermitian norm family $\xi_i=(\|\ndot\|_{i,\omega})_{\omega\in\Omega}$ on $H^0(X,L_i^{\otimes m_i})$, such that $\varphi_i$ identifies with the quotient metric family induced by $\xi_i$.
\end{quote}
We will prove by inverted induction on $r$ that, under the condition $(C_r)$, the function \eqref{Equ: local intersection measurable} is $\mathcal A$-measurable. Note that the condition $(C_{-1})$ is always true and hence the measurability of  \eqref{Equ: local intersection measurable} under $(C_{-1})$ is just the statement of the theorem. We begin with the case where $r=d$. For any $i\in\{1,\ldots,d\}$, let $E_i=H^0(X,L_i^{\otimes m_i})$ and $f_i:X\rightarrow\mathbb P(E_i)$ be the canonical closed embedding. Note that $L_i^{\otimes m_i}$ is isomorphic to $f_i^*(\mathcal O_{E_i}(1))$. We denote by $R$ the resultant
\[R_{f_0,\ldots,f_d}^{X,s_0^{\otimes m_0},\ldots,s_d^{\otimes m_d}},\]
which is an element of  
\[S^{\delta_0 N_0}(E_0^\vee)\otimes_K\cdots\otimes_KS^{\delta_dN_d}(E_d^\vee),\]
where 
\[N_i=\frac{m_0\cdots m_d}{m_i}.\]
We equip this vector space with the family of $\varepsilon$-tensor product of $\varepsilon$-symmetric power norms of $\|\ndot\|_{i,\omega,*}$ (see Definition \ref{Def: epsilon tensor product}), which we denote by $\xi=(\|\ndot\|_\omega)_{\omega\in\Omega}$. By \cite[Proposition 4.1.24]{CMArakelovAdelic}, the norm family $\xi$ is measurable. By Theorem \ref{Thm: equality of local height}, one has
\[m_0\cdots m_d(\overline D_0\cdots\overline D_d)_{\omega}=(m_0\overline D_0\cdots m_d\overline D_d)_\omega=\ln\|R\|_\omega.\]
Hence the function 
\[(\omega\in\Omega)\longmapsto (\overline D_0\cdots\overline D_d)_\omega\] is measurable. 

We prove the measurability of \eqref{Equ: local intersection measurable} under $(C_{r-1})$ in assuming that the measurability of \eqref{Equ: local intersection measurable} is true under $(C_{r})$, where $r\in\{0,\ldots,d\}$. For any positive integer $m$, we let $g_r^{(m)}$ be the Green function family of $D_r$ corresponding to the metric family $\varphi_r^{(m)}=(\varphi_{r,\omega}^{(m)})_{\omega\in\Omega}$ (see Definition \ref{Def: Fubini-Study}). We first show that the function
\[(\omega\in\Omega)\longmapsto (\overline D_0\cdots\overline D_{r-1}(D_r,g_r^{(m)})\overline D_{r+1}\cdots\overline D_d)_\omega \] 
is measurable. For this purpose, we choose arbitrarily a measurable norm family $\xi_r=(\|\ndot\|_{\omega})_{\omega\in\Omega}$ on the vector space $H^0(X,L^{\otimes m})$ (one can choose $\xi_r=\xi_{\boldsymbol{e}}$, where $\boldsymbol{e}$ is a basis of $H^0(X,L^{\otimes m})$, see Example \ref{Exe: norm family of a basis}). For any $a>0$ and any $\omega\in\Omega$, we let $\varphi_{r,a,\omega}^{(m)}$ be the quotient metric on $L_r$ induced by the norm
\[\|\ndot\|_{a,\omega}:=\max\{\|\ndot\|_{\varphi_r^{\otimes m}},a\|\ndot\|_{\omega}\}\]
on $H^0(X_\omega,L_\omega^{\otimes m})$, and let $g_{r,a}^{(m)}$ be the Green function of $D_r$ corresponding to the metric $\varphi_{r,a,\omega}^{(m)}$. 
By Proposition \ref{Pro: measurable direct image}, the norm family $\xi_{r,a}:=(\|\ndot\|_{a,\omega})_{\omega\in\Omega}$ is measurable. Therefore $\overline D_0,\ldots,\overline D_{r-1},(D_r,g_{r,a}^{(m)}),\overline D_{r+1}\cdots\overline D_d$ satisfy the condition $(C_r)$. By the induction hypothesis, we obtain that 
the function
\[(\omega\in\Omega)\longmapsto (\overline D_0\cdots\overline D_{r-1}(D_r,g_{r,a}^{(m)})\overline D_{r+1}\cdots\overline D_d)_\omega \] 
is measurable. Moreover, by  Proposition \ref{Pro: varphi 1 as a continuous metric}, we obtain that, for any $\omega\in\Omega$, there exists $a_{\omega}>0$ such that $g_{r,a}^{(m)}=g_r^{(m)}$ when $0<a<a_\omega$. Therefore one has
\[(\overline D_0\cdots\overline D_{r-1}(D_r,g_r^{(m)})\overline D_{r+1}\cdots\overline D_d)_\omega=\lim_{a\in\mathbb Q,\,a\rightarrow 0+}(\overline D_0\cdots\overline D_{r-1}(D_r,g_{r,a}^{(m)})\overline D_{r+1}\cdots\overline D_d)_\omega\]
and hence the function
\[(\omega\in\Omega)\longmapsto (\overline D_0\cdots\overline D_{r-1}(D_r,g_r^{(m)})\overline D_{r+1}\cdots\overline D_d)_\omega \] 
is measurable. Finally, by Proposition \ref{Pro: convergence of varphin} and Corollary \ref{cor:intersection:trivial:divisor}, one has
\[(\overline D_0\cdots\overline D_d)_\omega=\lim_{m\rightarrow+\infty}(\overline D_0\cdots\overline D_{r-1}(D_r,g_r^{(m)})\overline D_{r+1}\cdots\overline D_d)_\omega\]
and therefore the function
\[(\omega\in\Omega)\longmapsto (\overline D_0\cdots\overline D_d)_\omega\]
is measurable.  
\end{proof}

In the following, we study the measurability of the function of local intersection number over Archimedean places. We assume that $\Omega_{\infty}=\Omega$.
If $K$ contains a square root $\sqrt{-1}$ of $-1$, then, by Lemma~\ref{lemma:measurable:family:embeddings},
for each $\omega \in \Omega$, there is an embedding 
$\sigma_{\omega} : K \hookrightarrow \CC$
with the following properties:
\begin{enumerate}[label=\rm(\arabic*)]
\item $|\ndot|_{\omega} = |\sigma_{\omega}(\ndot)|$ for all $\omega \in \Omega$.

\item $\sigma_{\omega}(\sqrt{-1}) = i$, so that $\sigma_{\omega}(a + \sqrt{-1}b) = a + ib$ for all $a, b \in \QQ$, where $i$ is the usual imaginary unit in $\mathbb C$.

\item For $a \in K$, $(\omega \in \Omega) \mapsto \sigma_{\omega}(a)$ is measurable.
\end{enumerate}

\begin{prop}\label{prop:measurable:partial:X}
We assume that $\Omega = \Omega_{\infty}$ and $\sqrt{-1} \in K$. 
Let $n$ and $d$ be non-negative integers with $n \geq d$ and
$\pi: \A^n_K \to \A^d_K$ be the projection given by $(x_1, \ldots, x_n) \mapsto (x_1, \ldots, x_d)$.
Let $U$ be a non-empty Zariski open set of $\A^d_K$ and $X$ be a reduced closed subscheme of $\pi^{-1}(U)$ such that $\rest{\pi}{X} : X \to U$ is finite, surjective and \'etale.

We assume that either (i) $n=d$ 
and $X = \pi^{-1}(U)$, or (ii) $K$ is algebraically closed field.
Let $f= \{ f_{\omega} \}_{\omega \in \Omega}$ be a family of functions indexed by $\Omega$ such that
$f_{\omega}$ is a $C^{\infty}$-function on $\pi_{\omega}^{-1}(U_{\omega})$ and that, for any $K$-rational point $P \in \pi^{-1}(U)(K)$,
the function given by $(\omega \in \Omega) \mapsto f_{\omega}(P_{\omega})$ is measurable.
If we set $g_{\omega} = \rest{f_{\omega}}{X_{\omega}}$ for $\omega \in \Omega$, then, for any $P \in X(K)$ and $l \in \{ 1, \ldots, d\}$,
\[
(\omega \in \Omega) \mapsto \frac{\partial g_{\omega}}{\partial z_{l\omega}}(P_{\omega})\quad\text{and}\quad
(\omega \in \Omega) \mapsto \frac{\partial g_{\omega}}{\partial \overline{z}_{l\omega}}(P_{\omega})
\]
are measurable, where $(z_{1\omega}, \ldots, z_{d\omega})$ be the coordinate of $\A^d \times_{\sigma_{\omega}} \CC$. 
\end{prop}

\begin{proof}
{\bf Case (i)}: $n=d$ (so that $\pi = \operatorname{id}$) 
and $X = \pi^{-1}(U)$.

Let $x_{l\omega}$ (resp. $y_{l\omega}$) be the real part (resp. the imaginary part) of $z_{l\omega}$.
It is sufficient to show that
\[
(\omega \in \Omega) \mapsto \frac{\partial f_{\omega}}{\partial x_{l\omega}}(P_{\omega})\quad\text{and}\quad
(\omega \in \Omega) \mapsto \frac{\partial f_{\omega}}{\partial y_{l\omega}}(P_{\omega})
\]
are measurable.
We set $P_{\omega} = \sigma_{\omega}(P) = (a_{1\omega} + i b_{1\omega}, \ldots, a_{n\omega} + ib_{n\omega})$.
Then, for $\varepsilon \in \QQ^{\times}$,
\[
\begin{cases}
(P + \varepsilon e_l)_{\omega} = \sigma_{\omega}(P + \varepsilon e_l) = (a_{1\omega} + i b_{1\omega}, \ldots, (a_{l\omega} + \varepsilon) + ib_{l\omega}, \ldots, a_{n\omega} + ib_{n\omega}),\\
(P + \varepsilon ie_l)_{\omega} = \sigma_{\omega}(P + \varepsilon i e_l) = (a_{1\omega} + i b_{1\omega}, \ldots, a_{l\omega} + i(b_{l\omega}+\varepsilon), \ldots, a_{n\omega} + ib_{n\omega}),
\end{cases}
\]
where $\{ e_1, \ldots, e_n \}$ is the standard basis of $K^n$, so that
\[
\begin{cases}
{\displaystyle \lim_{\substack{\varepsilon \in \QQ^{\times}\\ \varepsilon \to 0}} \frac{f_{\omega}((P + \varepsilon e_l)_{\omega}) - f_{\omega}(P_{\omega})}{\varepsilon} = \frac{\partial f_{\omega}}{\partial x_{l\omega}}(P_{\omega}),}\\
{\displaystyle \lim_{\substack{\varepsilon \in \QQ^{\times}\\ \varepsilon \to 0}} \frac{f_{\omega}((P + \varepsilon ie_l)_{\omega}) - f_{\omega}(P_{\omega})}{\varepsilon} = \frac{\partial f_{\omega}}{\partial y_{l\omega}}(P_{\omega}).}
\end{cases}
\]
Note that
\[
\begin{cases}
{\displaystyle (\omega \in \Omega) \mapsto \frac{f_{\omega}((P + \varepsilon e_l)_{\omega}) - f_{\omega}(P_{\omega})}{\varepsilon},}\\[1ex]
{\displaystyle (\omega \in \Omega) \mapsto \frac{f_{\omega}((P + \varepsilon ie_l)_{\omega}) - f_{\omega}(P_{\omega})}{\varepsilon}.}
\end{cases}
\]
are measurable. Thus the assertion follows.

\bigskip
{\bf Case (ii)}: $K$ is algebraically closed field.

By replacing $U$ and $X$ by $U \setminus \pi(P)$ and $X \setminus P$, we may assume that $P=(0, \ldots, 0)$.
If we set $Q = \pi(P)$, then $(\rest{\pi}{X})^* : \OO_{U, Q}^h \overset{\sim}{\longrightarrow} \OO_{X, P}^h$, where
$\OO_{U, Q}^h$ and $\OO_{X, P}^h$ are the Henselizations of $\OO_{U, Q}$ and $\OO_{X, P}$, respectively.
Thus there are $\varphi_{d+1}, \ldots, \varphi_{n} \in \OO_{U, Q}^h$ such that $(\rest{\pi}{X})^*(\varphi_{j}) = \rest{x_j}{X}$ for
$j \in \{ d+1, \ldots, n \}$.
We set 
\[
\varphi_j = \sum_{e_1\cdots e_d \in \ZZ_{\geq 0} } a_{j, e_1\cdots e_d} X_1^{e_1} \cdots X_d^{e_d}
\]
as an element of $K[\![ X_1, \ldots, X_{d} ]\!]$. Note that if we set
\[
\varphi_{j\omega} = \sum_{e_1\cdots e_d \in \ZZ_{\geq 0} } \sigma_{\omega}(a_{j, e_1\cdots e_d}) X_1^{e_1} \cdots X_d^{e_d},
\]
then \[g_{\omega} = f_{\omega}(z_{1\omega},\ldots,z_{d\omega},\varphi_{d+1 \omega}(z_{1\omega},\ldots,z_{d\omega}), \ldots, \varphi_{n \omega}(z_{1\omega},\ldots,z_{d\omega}))
\]
as a function on $U$ around $Q$. Then, for $l \in \{ 1, \ldots, d \}$,
\[
\begin{cases}
{\displaystyle \frac{\partial g_{\omega}}{\partial z_{l\omega}}(P_{\omega}) = \frac{\partial f_{\omega}}{\partial z_{l\omega}}(0, \ldots, 0) + \sum_{j=d+1}^n \frac{\partial f_{\omega}}{\partial z_{j \omega}}(0, \ldots, 0)\frac{\partial \varphi_{j\omega}}{\partial z_{l\omega}}(0, \ldots, 0)}, \\[2ex]
{\displaystyle \frac{\partial g_{\omega}}{\partial \overline{z}_{l\omega}}(P_{\omega}) = \frac{\partial f_{\omega}}{\partial \overline{z}_{l\omega}}(0, \ldots, 0)}.
\end{cases}
\]
If we denote $a_{j, e_1, \ldots, e_d}$ by $a_{j, l}$ in the case where $e_1 = 0, \ldots, e_l = 1, \ldots, e_d  = 0$, then
\[\begin{cases}
{\displaystyle \frac{\partial g_{\omega}}{\partial z_{l\omega}}(P_{\omega}) = \frac{\partial f_{\omega}}{\partial z_{l\omega}}(0, \ldots, 0) + \sum_{j=d+1}^n \frac{\partial f_{\omega}}{\partial z_{j \omega}}(0, \ldots, 0)\sigma_{\omega}(a_{j, l})}, \\[2ex]
{\displaystyle \frac{\partial g_{\omega}}{\partial \overline{z}_{l\omega}}(P_{\omega}) = \frac{\partial f_{\omega}}{\partial \overline{z}_{l\omega}}(0, \ldots, 0)},
\end{cases}
\]
so that the assertions follow from the case (i).
\end{proof}

\begin{prop}\label{prop:measurable:vol:int:A}
We assume that $\Omega = \Omega_{\infty}$ and $\sqrt{-1} \in K$.
Let $U$ be a non-empty Zariski open set of $\A^n_K$.
Let $h= \{ h_{\omega} \}_{\omega \in \Omega}$ be a family of functions indexed by $\Omega$ such that
$h_{\omega}$ is a $C^{\infty}$-function on $U_{\omega}$ and that, for any $K$-rational point $P \in U(K)$,
the function given by $(\omega \in \Omega) \mapsto h_{\omega}(P_{\omega})$ is measurable.
For each $\omega \in \Omega$, let $(z_{1\omega}, \ldots, z_{n\omega})$ is the coordinate of $\A^n \otimes_{\sigma_{\omega}} \CC$. If
\[\int_{U_{\omega}} \Big(\frac i2\Big)^nh_{\omega} (z_{1\omega},\ldots,z_{n\omega})\, \mathrm{d}z_{1\omega} \wedge \mathrm{d}\bar{z}_{1\omega} \wedge \cdots \wedge \mathrm{d}z_{n\omega} \wedge \mathrm{d}\bar{z}_{n\omega} \] exists for any $\omega \in \Omega$, then
\[
(\omega \in \Omega) \mapsto \int_{U_{\omega}} \Big(\frac i2\Big)^nh_{\omega} (z_{1\omega},\ldots,z_{n\omega}) \,\mathrm{d}z_{1\omega} \wedge \mathrm{d}\bar{z}_{1\omega} \wedge \cdots \wedge \mathrm{d}z_{n\omega} \wedge \mathrm{d}\bar{z}_{n\omega}
\]
is measurable.
\end{prop}

\begin{proof} 
Shrinking $U$ if necessarily, we may assume that $\A^n_K \setminus U$ is defined by $\{ F =  0 \}$ for some $F \in K[X_1, \ldots, X_n] \setminus \{ 0 \}$.
We set 
\[
U_{\omega, N} = \Big\{ (z_{1\omega},\ldots,z_{n\omega}) \in \CC^n \,\Big|\, 
\text{$\max_{j\in\{1, \ldots, n\}} |z_{j\omega}| \leq N$ and $|F(z_{1\omega}, \ldots, z_{n\omega})| \geqslant 1/N$} \Big\}.
\]
Let $x_{i\omega}$ (resp. $y_{i\omega}$)
be the real part (resp. imaginary part) of $z_{i\omega}$.
Then
\begin{multline*}
 \Big(\frac i2\Big)^n h_{\omega}  \,\mathrm{d}z_{1\omega} \wedge \mathrm{d}\bar{z}_{1\omega} \wedge \cdots \wedge \mathrm{d}z_{n\omega} \wedge \mathrm{d}\bar{z}_{n\omega} \\
 =
 h_{\omega}\, \mathrm{d}x_{1\omega} \wedge \mathrm{d}y_{1\omega} \wedge \cdots
\wedge \mathrm{d}x_{n\omega} \wedge \mathrm{d}y_{n\omega}.
\end{multline*} 
Moreover,
\begin{multline}\label{eqn:prop:measurable:projective:space:01}
\int_{U_{\omega, N}} h_{\omega}\, \mathrm{d}x_{1\omega} \wedge \mathrm{d}y_{1\omega} \wedge \cdots
\wedge \mathrm{d}x_{n\omega} \wedge \mathrm{d}y_{n\omega} \\
\kern-12em = \lim_{m\to\infty} \sum_{\substack{a_1, b_1, \ldots, a_n, b_n \in \ZZ \\
\big(\frac{a_1 + i b_1}{m}, \ldots, \frac{a_n + i b_n}{m}\big) \in U_{\omega, N}}} \frac{1}{m^{2n}} h_{\omega}\Big(\frac{a_1 + i b_1}{m}, \ldots, \frac{a_n + i b_n}{m}\Big).
\end{multline}
Note that 
\[
(\omega \in \Omega_{\infty}) \longmapsto h_{\omega}\Big(\frac{a_1 + i b_1}{m}, \ldots, \frac{a_n + i b_n}{m}\Big)
\]
is measurable, so that \eqref{eqn:prop:measurable:projective:space:01} means that
\[
(\omega \in \Omega_{\infty}) \longmapsto \int_{U_{\omega, N}} h_{\omega}\, \mathrm{d}x_{1\omega} \wedge \mathrm{d}y_{1\omega} \wedge \cdots
\wedge \mathrm{d}x_{n\omega} \wedge \mathrm{d}y_{n\omega}
\]
is measurable. Therefore, one has the assertion because
\begin{multline*}
\lim_{N\to\infty} \int_{U_{\omega, N}} h_{\omega}\, \mathrm{d}x_{1\omega} \wedge \mathrm{d}y_{1\omega} \wedge \cdots
\wedge \mathrm{d}x_{n\omega} \wedge \mathrm{d}y_{n\omega}
\\ = \int_{U_{\omega}} h_{\omega}\, \mathrm{d}x_{1\omega} \wedge \mathrm{d}y_{1\omega} \wedge \cdots
\wedge \mathrm{d}x_{n\omega} \wedge \mathrm{d}y_{n\omega}. \end{multline*}
\end{proof}

\begin{theo}\label{thm:L2:norm:very:ample}
We assume that $\Omega = \Omega_{\infty}$ and $K$ is algebraically closed.
Let $X$ be a $d$-dimensional 
projective and integral variety over $K$ and $L$ be a very ample invertible $\OO_X$-module.
Let $\{ \|\ndot\|_{\omega} \}_{\omega \in \Omega}$ be a measurable family of hermitian norms on $H^0(X, L)$.
Let $\varphi = \{ \varphi_{\omega} \}_{\omega \in \Omega}$ be a family of metrics on $L$ induced by the surjective homomorphism $H^0(X, L) \otimes \OO_X \to L$ and
$\{ \|\ndot\|_{\omega} \}_{\omega \in \Omega}$. For $s \in H^0(X, L) \setminus \{ 0 \}$,
\[(\omega \in \Omega) \mapsto \int_{X_{\omega}} \log |s|_{\varphi_{\omega}} c_1(L_{\omega}, \varphi_{\omega})^{\wedge d}\]
is measurable.
\end{theo}

\begin{proof}
Let $n = \dim_K H^0(X, L) -1$ and $X \hookrightarrow \PP^n_K$ be the embedding by $L$. Note that $L = \rest{\OO_{\PP^n_K}(1)}{X}$.
Since $H^0(\PP^n_K, \OO_{\PP^n_K}(1)) \simeq H^0(X, L)$, one has $t \in H^0(\PP^n_K, \OO_{\PP^n_K}(1))$ with $\rest{t}{X} = s$.
Let $\psi = \{ \psi_{\omega} \}_{\omega \in \Omega}$ be a family of metics of $\OO_{\PP^n_K}(1)$ induced by the surjective homomorphism
$H^0(\PP^n_K, \OO_{\PP^n_K}(1)) \otimes \OO_{\PP^n_K} \to \OO_{\PP^n_K}(1)$ and
$\{ \|\ndot\|_{\omega} \}_{\omega \in \Omega}$. Note that $\rest{\psi}{X} = \varphi$.
By Proposition~\ref{prop:projection:to:projective:space},
we can choose a linear subspace $M$ in $\PP^n_K$ such that $\codim M = d+1$, $M \cap X = \emptyset$ and
$M \subseteq \{ t = 0 \}$, so that, by Proposition~\ref{prop:projection:to:projective:space} again,
the morphism $\pi : X \to \PP^d_K$ induced by 
the projection $\pi_M : \PP^n_K \setminus M \to \PP^d_K$ with the center $M$ is finite and surjective.
We choose a homogenous coordinate $(T_0:\ldots:T_n)$ on $\PP^n_K$ such that \[t = T_0\quad\text{and}\quad
M = \{ T_0 = \cdots = T_d = 0 \}.\] Then
$\pi_M$ is given by $(T_0 : \cdots : T_n) \mapsto (T_0 : \cdots : T_d)$.
Let $U$ be a non-empty open of $\PP^d_K$ such that $\pi : X \to \PP^d_K$ is \'{e}tale over $U$.
We may assume that $U \subseteq \{ T_0 \not = 0 \}$. We set $X_j = T_j/T_0$ ($j=1, \ldots, n$).
Then \[
\begin{cases} \PP^n_K \setminus \{ T_0 = 0 \} = \Spec(K[X_1, \ldots, X_n]) = \A^n_K, \\
\PP^d_K \setminus \{ T_0 = 0 \} = \Spec(K[X_1, \ldots, X_d]) = \A^d_K
\end{cases}
\] and $\pi_M$ on $\PP^n_K \setminus \{ T_0 = 0 \}$
is given by $(X_1, \ldots, X_n) \mapsto (X_1, \ldots, X_d)$.
Let \[(z_{1\omega}, \ldots, z_{n\omega})\quad\text{and}\quad (z_{1\omega}, \ldots, z_{d\omega})\] be
the coordinates of $\A^n_K \otimes_{\sigma_{\omega}} \CC$ and $\A^d_K \otimes_{\sigma_{\omega}} \CC$,
respectively. Note that $f_{\omega} := \log |t|_{\psi_{\omega}}$ is $C^{\infty}$ on $\A^n_K \otimes_{\sigma_{\omega}} \CC$.
Then, by Proposition~\ref{prop:measurable:partial:X},
if we set
\[
\rest{f_{\omega}}{X_{\omega}} c_1(L_{\omega}, \varphi_{\omega})^{\wedge d} = i^d h_{\omega} (dz_{1\omega} \wedge d\bar{z}_{1\omega}) \wedge
\cdots \wedge (dz_{d\omega} \wedge d\bar{z}_{d\omega})
\]
on $\pi_{\omega}^{-1}(U_{\omega})$, then, for $P \in \pi^{-1}(U)$,
$(\omega \in \Omega) \mapsto h(P_{\omega})$ is measurable.
Note that
\begin{align*}
\int_{X_{\omega}} \log |s|_{\varphi_{\omega}} c_1(L_{\omega}, \varphi_{\omega})^{\wedge d} & = \int_{\pi_{\omega}^{-1}(U_{\omega})} \rest{f_{\omega}}{X_{\omega}} c_1(L_{\omega}, \varphi_{\omega})^{\wedge d} \\
& = \int_{\pi_{\omega}^{-1}(U_{\omega})} i^d h_{\omega} (dz_{1\omega} \wedge d\bar{z}_{1\omega}) \wedge
\cdots \wedge (dz_{d\omega} \wedge d\bar{z}_{d\omega}) \\
& = \int_{U_{\omega}} i^d (\pi_{\omega})_*(h_{\omega}) (z_{1\omega} \wedge d\bar{z}_{1\omega}) \wedge
\cdots \wedge (dz_{d\omega} \wedge d\bar{z}_{d\omega}).
\end{align*}
Moreover, $(\pi_{\omega})_*(h_{\omega})$ is $C^{\infty}$ over $U_{\omega}$. Further,
for $P \in U(K)$, if we set $\pi^{-1}(P) = \{ Q_1, \ldots, Q_r \}$, then
\[
(\pi_{\omega})_*(h_{\omega})(P_{\omega}) = \sum_{i=1}^r h_{\omega}(Q_{i\omega}),
\]
so that $(\omega \in \Omega) \mapsto (\pi_{\omega})_*(h_{\omega})(P_{\omega})$ is measurable.
Therefore, by Proposition~\ref{prop:measurable:vol:int:A}, 
\[
(\omega \in \Omega) \mapsto \int_{U_{\omega}} i^d (\pi_{\omega})_*(h_{\omega}) (dz_{1\omega} \wedge d\bar{z}_{1\omega}) \wedge
\cdots \wedge (dz_{d\omega} \wedge d\bar{z}_{d\omega})
\]
is measurable. Thus the assertion follows.
\end{proof}

\begin{theo}\label{thm:measurable:archimedean}
We assume that $\Omega = \Omega_{\infty}$.
Let $X$ be a projective scheme over $K$ and $L$ be an ample invertible $\OO_X$-module.
Let $\varphi = \{ \varphi_{\omega} \}_{\omega \in \Omega}$ be a measurable family of semipositive metrics.
Then, for $s \in H^0(X, L) \setminus \{ 0 \}$, 
\[ (\omega \in \Omega) \longmapsto \int_{X_{\omega}} \log |s|_{\varphi_{\omega}} c_1(L_{\omega},\varphi_{\omega})^d \]
is measurable.
\end{theo}

\begin{proof}
By Lemma~\ref{lem:measurable:integrable:adelic:algebraic:cover}, 
we may assume that $K$ is algebraically closed. 
We choose a positive integer $N$ such that $L^{\otimes n}$ is very ample for 
for all $n \geqslant N$.
Let $\varphi_n = \{ \varphi_{n,\omega} \}_{\omega \in \Omega}$ be the Fubini-Study metric of $L^{\otimes n}$ induced by $H^0(X, L^{\otimes n}) \otimes \OO_X \to L^{\otimes n}$ and $\|\ndot\|_{n\varphi} = \{ \|\ndot\|_{n\varphi_{\omega}} \}_{\omega \in \Omega}$.
Moreover, by \cite[Theorem~4.1.26]{CMArakelovAdelic}, there is a measurable Hermitian norm family
$\{ \|\ndot\|_{n, \omega}^H \}_{\omega \in \Omega}$ on $H^0(X, L^{\otimes n})$ such that
\[
\|\ndot\|_{n\varphi_\omega} \leqslant \|\ndot\|_{n, \omega}^H \leqslant (h^0(L^{\otimes n}) + 1)^{1/2} \|\ndot\|_{n\varphi_\omega}
\]
for $\omega \in \Omega$.
Let $\varphi^H_{n,\omega}$ be the Fubini-Study metric of $L^{\otimes n}$ induced by $H^0(X, L^{\otimes n}) \otimes \OO_X \to L^{\otimes n}$ and $\|\ndot\|^H_{n, \omega}$.
Note that
\[
d_{\omega}\left(\frac{1}{n}\varphi_{n}, \frac{1}{n}\varphi^H_{n}\right) \leqslant \frac{ d_{\omega}(\|\ndot\|_{n\varphi},
\|\ndot\|_{n}^H)}{n} \leqslant \frac{\ln (h^0(L^{\otimes n}) + 1)}{2n}.
\]
Therefore, if we set
$\psi_{n,\omega} = (1/n)\varphi^H_{n,\omega}$,
then $\lim_{n\to\infty} d_{\omega}(\varphi, \psi_n) = 0$ for all $\omega \in \Omega$
because $\lim_{n\to\infty} d_{\omega}(\varphi, (1/n)\varphi_n) = 0$.
By Theorem~\ref{thm:L2:norm:very:ample},
\[
(\omega \in \Omega_{\infty}) \mapsto \int_{X_{\omega}} \log |s|_{\psi_{n,\omega}} c_1(L_{\omega}, \psi_{n, \omega})^d
= \frac{1}{n^{d+1}} \int_{X_{\omega}} \log |s^n|_{\varphi^H_{n,\omega}} c_1(nL_{\omega}, \varphi^H_{n, \omega})^d
\]
is measurable. Further, by \cite[Corollary~3.6]{Demagbook},
\[
\lim_{n\to\infty} \int_{X_{\omega}} \log |s|_{\psi_{n,\omega}} c_1(L_{\omega}, \psi_{n, \omega})^d = \int_{X_{\omega}} \log |s|_{\varphi_{n}} c_1(L_{\omega}, \varphi_{\omega})^d.
\]
Therefore, the assertion follows.
\end{proof}

Combining Theorems \ref{Thm: measurability non archimedean} and \ref{thm:measurable:archimedean}, we obtain the following result. 

\begin{theo}\label{Thm: measurability}
Let $X\rightarrow\Spec K$ be a projective scheme over $K$ and $d$ be the dimension of $X$. Let $D_0,\ldots,D_d$ be Cartier divisors on $X$, which intersect properly. We equip each $D_i$ with a measurable Green function $g_i$ such that $(D_{i},g_{i})$ is integrable. Then the local intersection function
\[(\omega\in\Omega)\longmapsto ((D_0,g_0)\cdots(D_d,g_d))_{\omega}\] 
is $\mathcal A$-measurable.
\end{theo}
\begin{proof}The measurability over $\Omega\setminus\Omega_\infty$ follows directly from Theorem \ref{Thm: measurability non archimedean}. Moreover, in view of Theorem \ref{thm:measurable:archimedean}, the measurability over $\Omega_\infty$ follows from Proposition~\ref{prop:intersection:trivial:divisor} and the multi-linearity of 
the local intersection measure.
\end{proof}

\begin{theo}\label{thm:integrability:local:intersection}
Let $X\rightarrow\Spec K$ be a projective scheme over $K$ and $d$ be the dimension of $X$. Let $D_0,\ldots,D_d$ be Cartier divisors on $X$, which intersect properly. We equip each $D_i$ with a dominated Green function $g_i$ such that $(D_{i},g_{i})$ is integrable. Then the local intersection function
\begin{equation}\label{eqn:thm:integrability:local:intersection:01}(\omega\in\Omega)\longmapsto ((D_0,g_0)\cdots(D_d,g_d))_{\omega}\end{equation}
is dominated.
\end{theo}

\begin{proof}
By Lemma~\ref{lem:measurable:integrable:adelic:algebraic:cover}, we may assume that
$K$ is algebraically closed.
By using Proposition~\ref{prop:multilinear:symmetric:semiample:case},
we may further assume that $D_0, \ldots, D_d$ are very ample.
Moreover, by Proposition~\ref{prop:reduction:same:class}, we may assume without loss of generality that
there are an integrable adelic line bundle $(L, \varphi)$ and non-zero rational sections $s_0, \ldots, s_d$ of $L$
such that $\OO_X(D_i) = L$ and $g_i = -\log |s_i|_{\varphi}$ for $i \in \{0,\ldots, d\}$.
Note that $L$ is very ample. Thus, by Proposition~\ref{prop:projection:to:projective:space},
there is a finite and surjective morphism $\pi : X \to \PP^d_K$ such that
$L = \pi^*(\OO_{\PP^d}(1))$. Let $(T_0 : \cdots : T_d)$ be a homogeneous coordinate of $\PP^d_K$.
Let $\varphi_{\rm{FS}}$ be a Fubini-Study metric of $\PP^d_K$ and $H_i = \{ T_i = 0 \}$ as in Proposition~\ref{prop:local:intersection:PP}. Moreover, we set $h_i = -\log |T_i|_{\varphi_{\rm{FS}}}$.

First we assume that $\varphi = \pi^*(\varphi_{\rm{FS}})$.
If $D_i = \pi^*(H_i)$ for $i \in \{ 0, \ldots, d \}$, then the dominancy of \eqref{eqn:thm:integrability:local:intersection:01} follows from
Proposition~\ref{prop:intersection:finite:morphism} and Proposition~\ref{prop:local:intersection:PP}. 
In general, there are non-zero rational functions $f_0, \ldots, f_d$ on $X$
such that $D_i = \pi^*(H_i) + (f_i)$ for $i \in \{ 0, \ldots, d \}$. Then, by Proposition~\ref{prop:intersection:principal:div},
there is an integrable function $\theta$ on $\Omega$ such that
\[
((D_0, g_0) \cdots (D_d, g_d))_{\omega} = ((\pi^*(H_0), \pi^*(h_0)) \cdots (\pi^*(H_d), \pi^*(h_d)))_{\omega} + \theta(\omega).
\]
Thus one has the dominancy of \eqref{eqn:thm:integrability:local:intersection:01}.

In general, there is a family $g$ of integrable continuous functions such that $\varphi = \exp(g) \pi^*(\varphi_{\rm{FS}})$.
In this case, the dominancy of \eqref{eqn:thm:integrability:local:intersection:01} follows from 
Corollary~\ref{cor:intersection:trivial:divisor}. \end{proof}

Finally, we obtain the following integrability theorem.

\begin{theo}\label{Thm: integrability of local intersection}
Let $X$ be a  projective $K$-scheme of dimension $d$, and $\overline D_0,\ldots,\overline D_d$ be a family of integrable  adelic Cartier divisors. Assume that the underlying Cartier divisors $D_0,\ldots,D_d$ intersect properly. Then the function of local intersection numbers
\begin{equation}\label{Equ: local intersection function}(\omega\in\Omega)\longmapsto (\overline D_0\cdots\overline D_d)_\omega\end{equation}
is integrable on the measure space $(\Omega,\mathcal A,\nu)$.
\end{theo}

\begin{defi}
Let $X$ be a projective $K$-scheme of dimension $d$, and $\overline D_0,\ldots,\overline D_d$ be a family of integrable adelic Cartier divisors, such that $D_0,\ldots,D_d$ intersect properly. We define the global intersection number of $\overline D_0,\ldots,\overline D_d$ as
\[(\overline D_0\cdots\overline D_d)_{S}:=\int_{\omega\in\Omega}(\overline D_0\cdots\overline D_d)_\omega\,\nu(\mathrm{d}\omega).\]
\end{defi}

\begin{rema}\label{Rem: global heights equality}
Let $X$ be a projective $K$-scheme of dimension $d$. For any $i\in\{0,\ldots,d\}$, let \[(E_i,\xi_i=(\|\ndot\|_{i,\omega})_{\omega\in\Omega})\] be a Hermitian adelic vector bundle on $S$, and $f_i:X\rightarrow\mathbb P(E_i)$ be a closed embedding. Let $L_i$ be the restriction of $\mathcal O_{E_i}(1)$ to $X$, which is equipped with the orthogonal quotient metric family $\varphi_i$ induced by $\xi_i$. We choose a global section $s_i$ of $L_i$ such that $s_0,\ldots,s_d$ intersect properly. For each $i\in\{0,\ldots,d\}$, let $D_i$ be the Cartier divisor $\operatorname{div}(s_i)$ and $g_i$ be the Green function family of $D_i$ corresponding to $\varphi_i$. By Theorem \ref{Thm: equality of local height}, if we denote by $R$ the resultant
\[R^{X,s_0,\ldots,s_d}_{f_0,\ldots,f_d}\in S^{\delta_0}(E_0^\vee)\otimes_K\cdots\otimes_KS^{\delta_d}(E_d^\vee),\]
where $\delta_i=(D_0\cdots D_{i-1}D_{i+1}\cdots D_d)$,
then the following equality holds
\[\begin{split}&(\overline D_0\cdots\overline D_d)=\int_{\omega\in\Omega\setminus\Omega_\infty}\ln\|R\|_\omega\,\nu(\mathrm{d}\omega)\\&\quad+
\int_{\sigma\in\Omega_\infty}\nu(\mathrm{d}\sigma)\int_{\mathbb S_{0,\sigma}\times\cdots\times \mathbb S_{d,\sigma}}\ln|R_\sigma(z_0,\ldots,z_d)|\,\eta_{\mathbb S_{0,\sigma}}(\mathrm{d}z_0)\otimes\cdots\otimes\eta_{\mathbb S_{d,\sigma}}(d\mathrm{z}_d)\\
&\qquad+\nu(\Omega_\infty)\frac 12\sum_{i=0}^{d}\delta_i\sum_{\ell=1}^{r_i}\frac{1}{\ell},
\end{split}\]
where 
\begin{enumerate}[label=\rm(\arabic*)]
\item $\|\ndot\|_\omega$ is the $\varepsilon$-tensor product   of $\delta_i$-th $\varepsilon$-symmetric tensor power of $\|\ndot\|_{i,\omega,*}$,
\item $R_\sigma$ is the element of
\[S^{\delta_0}(E_{0,\mathbb C_\sigma}^\vee)\otimes_{\mathbb C_\sigma}\cdots\otimes_{\mathbb C_\sigma}S^{\delta_{d,\mathbb C_\sigma}}(E_{d,\mathbb C_\sigma}^\vee)\]
indued by $R$,
\item $\mathbb S_{i,\sigma}$ is the unique sphere of $(E_{i,\mathbb C_\sigma},\|\ndot\|_{i,\sigma,\mathbb C_\sigma})$,
\item $\eta_{\mathbb S_{i,\sigma}}$ is the $U(E_{i,\mathbb C_\sigma},\|\ndot\|_{i,\mathbb C_\sigma})$-invariant Borel probaility measure on $\mathbb S_{i,\sigma}$.
\end{enumerate}

\end{rema}

\section{Invariance of intersection number by coverings}

Let $S=(K,(\Omega,\mathcal A,\nu),\phi)$ be an adelic curve. Consider a covering \[\alpha=(\alpha^{\#},\alpha_{\#},I_\alpha)\] from another adelic curve $S'=(K',(\Omega',\mathcal A',\nu'),\phi')$ to $S$ (see Definition \ref{Def: covering of adelic curves}). We assume that, either both $\sigma$-algebra $\mathcal A$ and $\mathcal A'$ is discrete, or there exist countable subfields $K_0$ and $K_0'$ of $K$ and $K'$ respectively, such that $K_0$ is dense in each $K_\omega$ with $\omega\in\Omega$, and $K_0'$ is dense in each $K'_{\omega'}$ with $\omega'\in\Omega'$. Recall that $\alpha^{\#}:K\longrightarrow K'$ is a field homomorphism, \[\alpha_{\#}:(\Omega',\mathcal A')\rightarrow (\Omega,\mathcal A)\] is a measurable map, and \[I_\alpha:\mathscr L^1(\Omega',\mathcal A',\nu')\longrightarrow\mathscr L^1(\Omega,\mathcal A,\nu)\] is a disintegration kernel of $\nu'$ over $\nu$ such that, for any $g\in\mathscr L^1(\Omega,\mathcal A,\nu)$, one has $g\circ\alpha_{\#}\in\mathscr L^1(\Omega',\mathcal A',\nu')$ and 
$I_\alpha(g\circ\alpha_{\#})=g$.
In this section, we consider a projective scheme $X$ of dimension $d$ over $\Spec K$ and a family \[\overline D_0=(D_0,g_0),\ldots,\overline{D}_d=(D_d,g_d)\] of adelic Cartier divisors, such that $D_0,\ldots,D_d$ intersect properly. The purpose of this section is to define the extension of scalars $\overline D_{i,\alpha}$ of each adelic Cartier divisor $\overline D_i$ by $\alpha$ and show the following equality
\[(\overline D_{0,\alpha}\cdots\overline D_{d,\alpha})_{S'}=(\overline D_0\cdots\overline D_d)_S.\]

\begin{defi}
Let $D$ be a Cartier divisor on $X$ and $g=(g_\omega)_{\omega\in\Omega}$ be a Green function family of $D$ (see Definition \ref{def:adelic:Cartier:divisor}). Let $X_{\alpha}$ be the fiber product \[X\times_{\Spec K,\alpha^{\#}}\Spec K'\] and $D_\alpha$ be the pull-back of $D$ by the morphism of projection $X_{\alpha}\rightarrow X$. If $\omega'$ is an element of $\Omega'$ and $\omega=\alpha_{\#}(\omega')$, then the Cartier divisor $D_{\alpha,\omega'}$ identifies with the pull-back of $D_\omega$ by the morphism of projection
\[X_{\alpha,\omega'}=X_{\alpha}\times_{\Spec K'}\Spec K'_{\omega'}\cong X_\omega\times_{\Spec K_\omega}\Spec K'_{\omega'}\longrightarrow X_\omega.\]
We denote by $g_{\alpha,\omega'}$ the Green function $g_{\omega,K'_{\omega'}}$ (see Remark \ref{Rem: extension scalar Green function}). Then the family $g_\alpha:=(g_{\alpha,\omega'})_{\omega'\in\Omega'}$ forms a Green function family of the Cartier divisor $D_\alpha$.

Let $L$ be an invertible $\mathcal O_X$-module and $\varphi=(\varphi_\omega)_{\omega\in\Omega}$ be a metric family on $L$. We denote by $L_\alpha$ the pull-back of $L$ by the morphism of projection $X_{\alpha}\rightarrow X$. If $\omega'$ is an element of $\Omega'$ and $\omega=\alpha_{\#}(\omega')$, then the invertible sheaf $L_{\alpha,\omega'}$ identifies with the pull-back of $L_\omega$ by the morphism of projection $X_{K',\omega'}\rightarrow X_\omega$. We denote by $\varphi_{\alpha,\omega'}$ the continuous metric $\varphi_{\omega,K'_{\omega'}}$ (see Example \ref{Exa: Fubini-Study} \ref{Item: extension of scalars metric}) on $L_{\alpha,\omega'}$. Then the family $\varphi_\alpha:=(\varphi_{\alpha,\omega'})_{\omega'\in\Omega'}$ forms a metric family of $L_\alpha$. Note that, if $s$ is a regular meromorphic section of $L$, $D=\operatorname{div}(s)$ and $g=(g_\omega)_{\omega\in\Omega}$ is the Green function family of $D$ corresponding to the metric family $\varphi$, then $g_\alpha$ is the Green function family of $D_\alpha$ corresponding to $\varphi_{\alpha}$. 
\end{defi}

\begin{prop}Let $\pi:X\rightarrow\Spec K$ be a projective $K$-scheme.
\begin{enumerate}[label=\rm(\arabic*)]
\item\label{Item: pull-back adelic line bundle} Let $L$ be an invertible $\mathcal O_X$-module and $\varphi$ be a metric family on $L$. If $\varphi$ is dominated, then $\varphi_\alpha$ is also dominated.
\item\label{Item: pull-back dominated} Let $D$ be a Cartier divisor on $X$ and $g$ is a Green function family of $g$. If $g$ is dominated, then $g_\alpha$ is also dominated.
\end{enumerate}
\end{prop} 
\begin{proof}
It suffices to prove the first statement. Assume that $\psi$ is another metric family on $L$. If $\omega'$ is an element of $\Omega'$ and if $\omega=\alpha_{\#}(\omega')$, then by \eqref{Equ: distance extension of scalars} one has
\[d_{\omega'}(\varphi_\alpha,\psi_\alpha)=d_\omega(\varphi,\psi).\]
Therefore, if the function $(\omega\in\Omega)\mapsto d_\omega(\varphi,\psi)$ is dominated, then also is the function $(\omega'\in\Omega')\mapsto d_{\omega'}(\varphi_\alpha,\psi_\alpha)$. To prove that the metric family $\varphi$ is dominated, we can assume without loss of generality that there exist a finite-dimensional vector space over $K$, a strongly dominated norm family $\xi=(\|\ndot\|_{\omega})_{\omega\in\Omega}$ on $E$, a positive integer $n$ and a surjective homomorphism $f:\pi^*(E)\rightarrow L^{\otimes n}$ such that $\varphi$ identifies with the orthogonal quotient metric family induced by $\xi$ (see Definition \ref{Def: Fubini-Study metric}). We may assume further that $\xi$ is Hermitian and $E$ admits a basis $\boldsymbol{e}$ which is orthonormal with respect to all norms $\|\ndot\|_\omega$. 

For any $\omega'\in\Omega'$, let $\|\ndot\|_{\omega'}$ be the norm $\|\ndot\|_{\omega,K'_{\omega'}}$, where $\omega=\alpha_{\#}(\omega')$. Then $\xi_\alpha^{H}=(\|\ndot\|_{\omega'})_{\omega'\in\Omega'}$ is a norm family on $E_{K'}$. Moreover, if we view $\boldsymbol{e}$ as a basis of $E_{K'}$ over $K'$, then it is orthonormal with respect to all norms $\|\ndot\|_{\omega'}$. In particular, the norm family $\xi_\alpha^{H}$ is strongly dominated. Since $\varphi_\alpha$ coincides with the orthogonal quotient metric family induced by $\xi_\alpha^{H}$, we deduce that the metric family $\varphi_\alpha$ is also dominated.
\end{proof}

\begin{defi}
Let $E$ be a finite-dimensional vector space over $K$ and $\xi=(\|\ndot\|_\omega)_{\omega\in\Omega}$ be a norm family on $E$. We denote by $\xi^\vee=(\|\ndot\|_{\omega,*})_{\omega\in\Omega}$ the \emph{dual norm family} on $E^\vee$, which is defined as 
\[\forall\,f\in E^\vee_{K_\omega},\quad \|f\|_{\omega,*}:=\sup_{s\in E_{K_\omega}\setminus\{0\}}\frac{|f(s)|_\omega}{\|s\|_\omega}.\] 
By \cite[Proposition 4.1.24 (1.b)]{CMArakelovAdelic}, the norm family $\xi^\vee$ is measurable.

We define $\xi_{\alpha}=(\|\ndot\|_{\omega'})_{\omega'\in\Omega'}$ the following norm family on $E_\alpha:=E\otimes_{K,\alpha^{\#}}K'$. In the case where $\omega'$ is non-Archimedean, the norm $\|\ndot\|_{\omega'}$ is the $\varepsilon$-extension of scalars of $\|\ndot\|_\omega$, where $\omega=\alpha_{\#}(\omega')$; in other words, one has
\[\forall\,s\in E_{\alpha,K'_{\omega'}},\quad \|s\|_{\omega'}=\sup_{f\in E_{K_\omega}^\vee\setminus\{0\}}\frac{|f(s)|_{\omega'}}{\|f\|_{\omega,*}}.\]
In the case where $\omega'$ is Archimedean, the norm $\|\ndot\|_{\omega'}$ is the $\pi$-extension of scalars of $\|\ndot\|_\omega$, in other words, one has
\[\forall\,s\in E_{\alpha,K_{\omega'}'},\quad \|s\|_{\omega'}=\inf\left\{|\lambda_1|_{\omega'}\cdot\|s_1\|_\omega+\cdots+|\lambda_N|_{\omega'}\cdot\|s_N\|_\omega\,\left|\,\begin{subarray}{l}
N\in\mathbb N,\,N\geqslant 1\\
(\lambda_1,\ldots,\lambda_N)\in (K'_{\omega'})^N\\
(s_1,\ldots,s_N)\in E_\omega^N\\
s=\lambda_1s_1+\cdots+\lambda_Ns_N
\end{subarray}\right.\right\}.\]

Similarly, we define $\xi_{\alpha,\varepsilon}$ the norm family on $E_\alpha$ consisting of $\varepsilon$-extension of scalars (for both non-Archimedean and Archimedean absolute values). 
\end{defi}

\begin{lemm}\label{Lem: norm family extension of scalars}
Let $E$ be a finite-dimension vector space over $K$ and $\xi=(\|\ndot\|_{\omega})_{\omega\in\Omega}$ be a measurable norm family on $E$. Then the norm families $\xi_{\alpha,\varepsilon}$ and $\xi_{\alpha}$ defined above are also measurable.
\end{lemm}
\begin{proof}
The proof is very similar to that of \cite[Proposition 4.1.24 (1.c)]{CMArakelovAdelic}. The case where $\mathcal A$ and $\mathcal A'$ are discrete is trivial. In the following, we will treat the case where $K$ and $K'$ admit countable subfields $K_0$ and $K_0'$   such that $K_0$ is dense in each $K_\omega$ with $\omega\in\Omega$, and $K_0'$ is dense in each $K'_{\omega'}$ with $\omega'\in\Omega'$, respectively. We first check the measurability of $\xi_{\alpha,\varepsilon}$. For any $\omega'\in\Omega'$, let $\|\ndot\|_{\omega',\varepsilon}$ be the norm indexed by $\omega'$ in the family $\xi_{\alpha,\varepsilon}$. Let $H_0$ be a finite-dimensional $K_0$-vector subspace of $E^\vee$ which generates $E^\vee$ as a vector space over $K$. Then $H_0\setminus\{0\}$ is dense in $E_{K_\omega}^\vee\setminus\{0\}$ for any $\omega\in\Omega$. If $s$ is an element of $E_\alpha$, then for any $\omega'\in\Omega'$,
\[\|s\|_{\omega',\varepsilon}=\sup_{f\in H_0\setminus\{0\}}\frac{|f(s)|_{\omega'}}{\|f\|_{\omega,*}}.\]
Hence it is the supremum of a countable family of $\mathcal A'$-measurable function in $\omega'$. As for the second statement, it suffices to apply the first statement to $\xi^\vee$ to obtain the measurability of $(\xi^{\vee})_{\alpha,\varepsilon}$. Since $\xi_{\alpha}$ is the dual norm family of $(\xi^{\vee})_{\alpha,\varepsilon}$ (see \cite[Proposition 1.3.20]{CMArakelovAdelic}), by \cite[Proposition 4.1.24 (1.c)]{CMArakelovAdelic} we obtain the measurability of $\xi_\alpha$.
\end{proof}

\begin{prop}
Let $X$ be a projective scheme over $\Spec K$. \begin{enumerate}[label=\rm(\arabic*)]
\item\label{Item: pull-back adelic line bundle} Let $L$ be an invertible $\mathcal O_X$-module and $\varphi$ be a metric family on $L$. We assume that $L$ is ample and all metrics in the family $\varphi$ is semi-positive If $\varphi$ is measurable, then $\varphi_\alpha$ is also measurable.
\item\label{Item: pull-back dominated} Let $D$ be a Cartier divisor on $X$ and $g$ be a Green function family of $g$. Assume that $D$ is ample and $g$ is semi-positive.  If $g$ is measurable, then $g_\alpha$ is also measurable.
\end{enumerate}
\end{prop}
\begin{proof}
It suffices to prove the first statement. Similarly  to the proof of Theorem \ref{Thm: measurability non archimedean}, for any $m\in\mathbb N_{\geqslant 1}$ such that $L^{\otimes m}$ is very ample we choose a norm family $\xi_m=(\|\ndot\|_{\omega}^{(m)})_{\omega\in\Omega}$ on $H^0(X,L^{\otimes m})$ such that $H^0(X,L^{\otimes m})$ admet a basis which is orthonormal with respect to each norm $\|\ndot\|_{\omega}^{(m)}$. This norm family is clearly measurable. For any $b>0$ and any $\omega\in\Omega$, let $\varphi_{b,\omega}^{(m)}$ the quotient metric on $L$ induced by the norm
\[\|\ndot\|^{(m)}_{b,\omega}=\max\{\|\ndot\|_{\varphi^{\otimes m}},b\|\ndot\|^{(m)}_\omega\}\]
on $H^0(X_\omega,L_\omega^{\otimes m})$. By Proposition \ref{Pro: measurable direct image}, the norm family $\xi_{b}^{(m)}:=(\|\ndot\|_{b,\omega}^{(m)})_{\omega\in\Omega}$ is measurable. By Lemma \ref{Lem: norm family extension of scalars}, we deduce that the norm family $\xi_{b,\alpha}^{(m)}$ of $H^0(X_\alpha,L_\alpha^{\otimes m})$ is $\mathcal A'$-measurable.

Let $\varphi_{b}^{(m)}$ be the quotient metric family on $L$ induced by $\xi_{b}^{(m)}$. By \cite[Remark 2.2.14]{CMArakelovAdelic}, the metric $\varphi_{b,\alpha}^{(m)}$ identifies with the quotient metric family on $L_\alpha$ induced by $\xi_{b,\alpha}^{(m)}$. Since the norm family $\xi_{b,\alpha}^{(m)}$ is measurable, by \cite[Proposition 6.1.30]{CMArakelovAdelic}, the metric family $\varphi_{b,\alpha}^{(m)}$ is measurable. By Proposition \ref{Pro: varphi 1 as a continuous metric}, for any fixed $\omega'\in\Omega'$ and $\omega=\alpha_{\#}(\omega')$, for sufficiently small $b$ one has
$\varphi_{b,\omega}^{(m)}=\varphi^{(m)}_{\omega}$ and hence $\varphi_{b,\alpha,\omega'}^{(m)}=\varphi^{(m)}_{\alpha,\omega'}$. Therefore, by \cite[Proposition 6.1.29]{CMArakelovAdelic} we obtain that $\varphi^{(m)}_\alpha$ is measurable.
By \eqref{Equ: distance extension of scalars}, for any $\omega'\in\Omega'$ and $\omega=\alpha_{\#}(\omega)$, one has
\[d_{\omega'}(\varphi^{(m)}_\alpha,\varphi_\alpha)\leqslant d_{\omega}(\varphi^{(m)},\varphi).\]
Since the metric family $\varphi$ is semi-positive, by Proposition \ref{Pro: convergence of varphin}, we deduce that, for any $\omega'\in\Omega'$, one has
\[\lim_{m\rightarrow+\infty}d_{\omega'}(\varphi_\alpha^{(m)},\varphi_\alpha)=0.\]
Still by \cite[Proposition 6.1.29]{CMArakelovAdelic}, we obtain that the metric family $\varphi$ is measurable.
\end{proof}

\begin{theo}
Let $X$ be a projective scheme over $\Spec K$ and $d$ be the dimension of $X$. Let $D_0,\ldots,D_d$ be Cartier divisors on $X$ which intersects properly. We assume that each Cartier divisor $D_i$ is equipped with an integrable Green function family $g_i$. The the following equality holds
\[((D_{0,\alpha},g_{0,\alpha})\cdots(D_{d,\alpha},g_{d,\alpha}))_{\omega'}=((D_0,g_0)\cdots(D_d,g_d))_{\alpha_{\#}(\omega')}.\]
In particular, if all Green function family $g_i$ are dominated (resp. measurable), then the function \[(\omega'\in\Omega')\longmapsto ((D_{0,\alpha},g_{0,\alpha})\cdots(D_{d,\alpha},g_{d,\alpha}))_{\omega'}\]is dominated (resp. measurable). If all $(D_i,g_i)$ are adelic Cartier divisors, then the following equality holds
\[((D_0,g_0)\cdots(D_d,g_d))_S=((D_{0,\alpha},g_{d,\alpha})\cdots(D_{d,\alpha},g_{d,\alpha}))_{S'}.\]
\end{theo}
\begin{proof}
For any $\omega'\in\Omega'$ and $\omega=\alpha_{\#}(\omega')$, the equality 
\[((D_{0,\alpha},g_{0,\alpha})\cdots(D_{d,\alpha},g_{d,\alpha}))_{\omega'}=((D_0,g_0)\cdots(D_d,g_d))_\omega\]
follows from \ref{prop:invariance:local:intersection:field:extension} (see also Remark \ref{Rem: extension of scalars archimedean}).  

If $g_0,\ldots,g_d$ are measurable, by Theorem \ref{Thm: measurability}, the function 
\[(\omega\in\Omega)\longmapsto ((D_0,g_0)\cdots(D_d,g_d))_\omega\]
is $\mathcal A$-measurable. Since $\alpha_{\#}$ is a measurable map, we deduce that the function 
\[(\omega'\in\Omega')\longmapsto((D_{0,\alpha},g_{0,\alpha})\cdots(D_{d,\alpha},g_{d,\alpha})_{\omega'}\]
is $\mathcal A'$-measurable.

Assume that the Green function families $g_0,\ldots,g_d$ are dominated. By Theorem \ref{thm:integrability:local:intersection}, there exists an integrable function $F$ on the measure space $(\Omega,\mathcal A,\nu)$ such that 
\[\forall\,\omega\in\Omega,\quad |((D_0,g_0)\cdots(D_d,g_d))_{\omega}|\leqslant F(\omega).\]
Hence
\[\forall\,\omega'\in\Omega,\quad |((D_{0,\alpha},g_{0,\alpha})\cdots(D_{d,\alpha},g_{d,\alpha}))_{\omega'}|\leqslant F(\alpha_{\#}(\omega')).\]
Hence the function
\[(\omega'\in\Omega')\longmapsto ((D_{0,\alpha},g_{0,\alpha}),\cdots,(D_{d,\alpha},g_{d,\alpha}))_{\omega'}\]
is dominated. Finally, if the function
\[(\omega\in\Omega)\longmapsto ((D_0,g_0)\cdots(D_d,g_d))_\omega\]
is integrable, then also is the function 
\[(\omega'\in\Omega')\longmapsto ((D_{0,\alpha},g_{0,\alpha})\cdots(D_{d,\alpha},g_{d,\alpha}))_{\omega'}=((D_0,g_0)\cdots(D_d,g_d))_{\alpha_{\#}(\omega')}\]
is also integrable, and one has
\[\begin{split}&\quad\;((D_{0,\alpha},g_{d,\alpha})\cdots(D_{d,\alpha},g_{d,\alpha}))_{S'}=\int_{\Omega'}((D_{0,\alpha},g_{0,\alpha})\cdots(D_{d,\alpha},g_{d,\alpha}))_{\omega'}\,\nu'(\mathrm{d}\omega')\\
&=\int_{\Omega}I_\alpha(\omega'\longmapsto((D_{0,\alpha},g_{0,\alpha})\cdots(D_{d,\alpha},g_{d,\alpha}))_{\omega'})\,\nu(\mathrm{d}\omega)\\
&=\int_{\Omega}((D_0,g_0)\cdots(D_d,g_d))_\omega\,\nu(\mathrm{d}\omega)=((D_0,g_0)\cdots(D_d,g_d))_S.
\end{split}\]
\end{proof}

\section{Multi-heights}

From now on, we assume that the adelic curve $S$ is proper.

\begin{defi}
Let $X$ be a projective scheme over $\Spec K$. If $f$ is a regular meromorphic function on $X$, we denote by $\widehat{\operatorname{div}}(f)$ the following adelic Cartier divisor
\[(\operatorname{div}(f),(-\ln|f|_\omega)_{\omega\in\Omega}).\]
If $\overline L=(L,\varphi)$ is an adelic line bundle on $X$ and if $s$ is a regular meromorphic section of $L$ on $X$, we denote by $\widehat{\operatorname{div}}(s)$ the following adelic Cartier divisor
\[(\operatorname{div}(s),(-\ln|s|_{\varphi_\omega})_{\omega\in\Omega}).\]
\end{defi}

\begin{prop}\label{prop:global:intersection:proper:adelic}
Let $X$ be a projective $K$-scheme of dimension $d$, and $\overline D_0,\ldots,\overline D_d$ and
$\overline D'_0,\ldots,\overline D'_d$ be families of integrable adelic Cartier divisors, such that $D_0,\ldots,D_d$ and $D'_0,\ldots,D'_d$ intersect properly. 
If there is a family of regular meromorphic functions $f_0, \ldots, f_d$ on $X$
such that $\overline D_i = \overline D'_i + \widehat{\operatorname{div}}(f_i)$ for $i \in \{ 0, \ldots, d \}$.
Then \[(\overline D_0\cdots\overline D_d)_S =(\overline D'_0\cdots\overline D'_d)_S.\]
\end{prop}

\begin{proof}
It is sufficient to prove that if $f$ is a regular meromorphic function on $X$ and
$\overline D_1,\ldots,\overline D_d$ are integrable adelic Cartier divisors such that
$\operatorname{div}(f), D_1,\ldots,D_d$ intersect properly,
then $(\widehat{\operatorname{div}}(f) \cdot \overline D_1\cdots\overline D_d)_S = 0$.
Clearly we may assume that $K$ is algebraically closed, so that the assertion follows from
Proposition~\ref{prop:intersection:principal:div} and the product formula.
\end{proof}

\begin{defi}
Let $\overline{L}_0=(L_0,\varphi_0), \ldots, \overline{L}_d=(L_d,\varphi_d)$ be a family of integrable adelic line bundles.
Let $s_0, \ldots, s_d$ be regular meromorphic sections of $L_0, \ldots, L_d$, respectively such that
$\operatorname{div}(s_0), \ldots, \operatorname{div}(s_d)$ intersect properly.
Then, by Proposition~\ref{prop:global:intersection:proper:adelic}, the global intersection number 
\[
(\widehat{\operatorname{div}}(s_0) \cdots \widehat{\operatorname{div}}(s_d))_S
\] 
does not depend on the choice of $s_0, \ldots, s_d$. The \emph{global intersection number}
\[
(\overline{L}_0 \cdots \overline{L}_d )_S
\]
of $\overline{L}_0 \cdots \overline{L}_d$ over $S$ is then defined as \[(\widehat{\operatorname{div}}(s_0) \cdots \widehat{\operatorname{div}}(s_d))_S.\]
This number is also called the \emph{multi-height} of $X$ with respect to $\overline L_0,\ldots,\overline L_d$ and is denoted by 
\[h_{\overline L_0,\ldots\overline L_d}(X).\]
In the particular case where $\overline L_0,\ldots\overline L_d$ are all equal to the same integrable adelic line bundle $\overline L$, the number $h_{\overline L,\ldots,\overline L}(X)$ is denoted by $h_{\overline L}(X)$ in abbreviation, and is called the \emph{height} of $X$ with respect to $\overline L$. 
\end{defi}

\begin{prop}
\begin{enumerate}[label=\rm(\arabic*)]
\item The global intersection pairing is a symmetric bilinear form on the group consisting of integrable adelic line bundle.

\item
Let $X_1, \ldots, X_\ell$ be irreducible components of $X$ and $\eta_1, \ldots, \eta_\ell$ be the generic points of $X_1, \ldots, X_\ell$, respectively.
Then
\[
\big(\overline{L}_0 \cdots \overline{L}_d\big)_S = \sum_{j=1}^{\ell} \operatorname{length}_{\mathcal O_{X, \eta_j}}(\mathcal O_{X, \eta_j}) \big(\rest{\overline{L}_0}{X_j} \cdots \rest{\overline{L}_d}{X_j}\big)_S.
\]

\item Let $s_d$ be a regular meromorphic section of $L_d$ and $\operatorname{div}(s_d) = a_1 Z_1 + \cdots + a_n Z_n$ be the decomposition
as cycles. Then
\[\begin{split}
(\overline{L}_0 \cdots \overline{L}_{d})_S =\int_{\Omega} & {\left(\int_{X_{\omega}^{\operatorname{an}}} -\log |s_d|_{\varphi_\omega}(x)\, \mu_{(L_{0,\omega},\varphi_{0,\omega}), \cdots (L_{d-1,\omega},\varphi_{d-1,\omega})}(\mathrm{d}x)\right)}\nu(\mathrm{d}\omega)\\
&+ 
\sum_{i=1}^n a_i (\rest{\overline{L}_0}{Z_i} \cdots \rest{\overline{L}_{d-1}}{Z_i})_S.
\end{split}
\]
\end{enumerate}
\end{prop}

\begin{proof}
They follows from \eqref{eqn:def:local:intersection} and Proposition~\ref{prop:multilinear:symmetric:semiample:case}.
\end{proof}

Finally let us consider the projection formula for our intersection theory.
For this purpose, we need three lemmas.

\begin{lemm}\label{lem:Artin:local:algebra:formula}
Let $(A,\mathfrak m)$ be a local Artinian ring and $B$ be an $A$-algebra such that $B$ is finitely generated as an $A$-module.
Let $M$ be a finitely generated $B$-module. Then
\[
\length_A(M) = \sum_{\mathfrak n \in \Spec(B)} [B/\mathfrak n : A/\mathfrak m] \length_{B_{\mathfrak n}}(M_{\mathfrak n}).
\]
In particular, if $B$ is flat over $A$, then
\[\rank_A(B) \length_A(A) = \sum_{\mathfrak n \in \Spec(B)} [B/\mathfrak n : A/\mathfrak m] \length_{B_{\mathfrak n}}(B_{\mathfrak n}).\]
\end{lemm}

\begin{proof}
Let $0 \to M' \to M \to M'' \to 0$ be an exact sequence of finitely generated $B$-modules. Then,
both sides of the above first equation are additive with respect to the exact sequence.
Therefore, we may assume that $M = B/\mathfrak{n}$ for some $\mathfrak n \in \Spec(B)$. In this case, it is obvious.
\end{proof}

\begin{lemm}\label{lem:associated:prime:maps:zero:ideal}
Let $A$ be an integral domain and $B$ be a flat $A$-algebra.
If we denote the structure homomorphism $A \to B$ by $\phi$, then
$\phi^{-1}(P) = \{ 0 \}$ for any $P \in \mathrm{Ass}_B(B)$.
\end{lemm}

\begin{proof}
We set $P = \mathrm{ann}(b)$ for some $b \in B \setminus \{ 0 \}$.
If there is $a \in \phi^{-1}(P) \setminus \{ 0 \}$,
then $\phi(a) b = 0$.
Since $B$ is flat over $A$, $\phi(a)$ is regular, so that $b = 0$.
This is a contradiction.
\end{proof}

\begin{lemm}\label{lem:cycle:push:forward:base:change}
Let $f : Y \to X$ be a proper and surjective morphism of integral scheme of finite type over a field $k$ such that
$\dim X = \dim Y$.
For an extension filed $k'$ of $k$, if $X' := X \times_{\Spec(k)}\Spec(k')$, $Y' := Y \times_{\Spec(k)}\Spec(k')$ and
$f' : X' \to Y'$ is the induced morphism, then
\[f'_*([X']) = [k(Y) : k(X)][Y'].\]
\end{lemm}

\begin{proof}
By Lemma~\ref{lem:associated:prime:maps:zero:ideal}, any irreducible component of $X'$ (resp. $Y'$) maps surjectively to $X$ (resp. $Y$) by $X' \to X$ (resp. $Y' \to Y$).
Moreover, we can find a non-empty Zariski open set $U$ of $X$ such that $f^{-1}(U) \to U$ is finite and flat. 
Note that if we set $U' := U \times_{\Spec(k)}\Spec(k')$, then
${f'}^{-1}(U') = f^{-1}(U) \times_{\Spec(k)}\Spec(k')$ and ${f'}^{-1}(U') \to U'$ is finite and flat.
Therefore, we may assume that $f$ is finite and flat, so that
the assertion
is a consequence of the second formula in Lemma~\ref{lem:Artin:local:algebra:formula}.
\end{proof}

\begin{defi}
Let $Z = a_1 Z_1 + \cdots + a_r Z_r$  be an $l$-dimensional cycle on $X$ and
$\overline{L}_0, \ldots, \overline{L}_l$ be integrable adelic line bundles.
Then $(\overline{L}_0 \cdots \overline{L}_l \mid Z)_S$ is defined to be
\[
(\overline{L}_0 \cdots \overline{L}_l \mid Z)_S := \sum_{j=1}^r a_j \Big(\rest{\overline{L}_0}{Z_j} \cdots \rest{\overline{L}_l}{Z_j}\Big)_S.
\]
In the case where $\overline L_0,\ldots,\overline{L}_l$ are all equal to the same adelic line bundle $\overline L$, we call it the \emph{height} of the cycle $Z$ with respect to $\overline L$, and denote it by $h_{\overline L}(Z)$.
\end{defi}

\begin{theo}[Projection formula]\label{prop:projection:formula:intersection}
Let $f : Y \to X$ be a morphism of projective schemes over $K$ and $\overline{L}_0, \ldots, \overline{L}_l$
be integrable adelic line bundles on $X$.
For an $l$-cycle $Z$ on $Y$,
\[(f^*(\overline{L}_0) \cdots f^*(\overline{L}_l) \mid Z)_S = (\overline{L}_0 \cdots \overline{L}_l \mid f_*(Z))_S.\]
\end{theo}

\begin{proof}
First let us see the following:

\begin{enonce}{Claim}\label{claim:prop:projection:formula:intersection:01}
If $f$ is a surjective morphism of projective integral schemes over $K$, then
\[
(f^*(\overline{L}_0) \cdots f^*(\overline{L}_l))_S = \begin{cases}
\deg(f) (\overline{L}_0 \cdots \overline{L}_l)_S & \text{if $\dim X = \dim Y$},\\
0 & \text{if $\dim X < \dim Y$}.
\end{cases}
\]
In other words,
\[ (f^*(\overline{L}_0) \cdots f^*(\overline{L}_l) \mid Y)_S = (\overline{L}_0 \cdots \overline{L}_l \mid f_*(Y))_S.
\]
\end{enonce}

\begin{proof}
We choose rational sections $s_0, \ldots, s_d$ of $L_0,\ldots, L_d$, respectively such that
$\operatorname{div}(s_0), \ldots, \operatorname{div}(s_d)$ intersect properly on $X$ and
$f^*(\operatorname{div}(s_0)), \ldots, f^*(\operatorname{div}(s_d))$ intersect properly on $Y$.
Let $K_{\omega}$ be the completion of $K$ with respect to $\omega \in \Omega$,
$X_{\omega} := X \times_{\Spec(K)} \Spec(K_{\omega})$, $Y_{\omega} := Y \times_{\Spec(K)} \Spec(K_{\omega})$
and $f_{\omega} : Y_{\omega} \to X_{\omega}$  be the induced morphism.
Further let $\pi_{X,\omega} : X_{\omega} \to X$ and $\pi_{Y,\omega} : Y_{\omega} \to Y$ be the projections.
Then the following diagram is commutative.
\[
\xymatrix{
Y_{\omega} \ar[r]^-{f_{\omega}}\ar[d]_-{\pi_{Y,\omega}}& X_{\omega}\ar[d]^-{\pi_{X,\omega}} \\
Y \ar[r]_-{f}& X
}
\]
Since $X$ and $Y$ are integral, $f^*(\operatorname{div}(s_i))$ is well defined as a Cartier divisor.
Moreover, $\pi_{Y,\omega}^*(f^*(\operatorname{div}(s_i)))$ and $\operatorname{div}(s_i)_{\omega} := \pi_{X,\omega}^*(\operatorname{div}(s_i))$ are defined because $\pi_{Y,\omega}$ and $\pi_{X,\omega}$ are flat.
Therefore, $f_{\omega}^*(\operatorname{div}(s_i)_{\omega})$ is defined as a Cartier divisor on $Y_{\omega}$ for each
$i = 0, \ldots, d$.
Let $Y_{\omega, 1}, \ldots, Y_{\omega, m_{\omega}}$ (resp. $X_{\omega, 1}, \ldots, X_{\omega, n_{\omega}}$)
be irreducible components of $Y_{\omega}$ (resp. $X_{\omega}$).

First we assume that $\dim X < \dim Y$. Then, by Proposition~\ref{prop:intersection:finite:morphism},
\[
\Big(\rest{f_{\omega}^*(\operatorname{div}(s_0)_{\omega}, -\log |s_0|_{\varphi_\omega})}{Y_{\omega,j}} \cdots 
\rest{f_{\omega}^*(\operatorname{div}(s_d)_{\omega}, -\log |s_d|_{\varphi_\omega})}{Y_{\omega,j}}\Big)_{\omega} =
0\]
for all $j = 1, \ldots, m_{\omega}$.
Therefore, \[
\Big(f_{\omega}^*(\operatorname{div}(s_0)_{\omega}, -\log |s_0|_{\varphi_\omega}) \cdots 
f_{\omega}^*(\operatorname{div}(s_d)_{\omega}, -\log |s_d|_{\varphi_\omega})\Big)_{\omega} = 0,
\] 
and hence
the assertion follows.

Next we assume that $\dim X = \dim Y$.
For each $i \in\{ 1, \ldots, n_{\omega}\}$, let \[J_{\omega, i} := \{ j \in \{ 1,\ldots, m_{\omega} \} \mid f_{\omega}(Y_{\omega, j}) = X_{\omega, i} \}\] and \[J_{\omega, 0} := \{ 1, \ldots, n_{\omega} \} \setminus (J_{\omega, 1} \cup \cdots \cup J_{\omega, n_{\omega}}).\] By Proposition~\ref{prop:intersection:finite:morphism}, if $j \in J_{\omega, i}$ ($i\in\{1, \ldots, n_{\omega}\}$), then
\begin{multline*}
\Big(\rest{f_{\omega}^*(\operatorname{div}(s_0)_{\omega}, -\log |s_0|_{\varphi_\omega})}{Y_{\omega,j}} \cdots 
\rest{f_{\omega}^*(\operatorname{div}(s_d)_{\omega}, -\log |s_d|_{\varphi_\omega})}{Y_{\omega,j}}\Big)_{\omega} \\
=
\deg(\rest{f_{\omega}}{Y_{\omega, j}})\Big(\rest{(\operatorname{div}(s_0)_{\omega}, -\log |s_0|_{\varphi_\omega})}{X_{\omega,i}} \cdots 
\rest{(\operatorname{div}(s_d)_{\omega}, -\log |s_d|_{\varphi_\omega})}{X_{\omega,i}}\Big)_{\omega}.
\end{multline*}
Moreover, if $j \in J_{\omega, 0}$, then
\[
\Big(\rest{f_{\omega}^*(\operatorname{div}(s_0)_{\omega}, -\log |s_0|_{\varphi_\omega})}{Y_{\omega,j}} \cdots 
\rest{f_{\omega}^*(\operatorname{div}(s_d)_{\omega}, -\log |s_d|_{\varphi_\omega})}{Y_{\omega,j}}\Big)_{\omega} \\
= 0.
\]
Thus, by Lemma~\ref{lem:cycle:push:forward:base:change}, one has
\begin{multline*}
\Big(f_{\omega}^*(\operatorname{div}(s_0)_{\omega}, -\log |s_0|_{\varphi_\omega}) \cdots 
f_{\omega}^*(\operatorname{div}(s_d)_{\omega}, -\log |s_d|_{\varphi_\omega})\Big)_{\omega} \\
= \deg(f)\Big((\operatorname{div}(s_0)_{\omega}, -\log |s_0|_{\varphi_\omega}) \cdots 
(\operatorname{div}(s_d)_{\omega}, -\log |s_d|_{\varphi_\omega})\Big)_{\omega}.
\end{multline*}
Therefore,
\begin{align*}
(f^*(\overline{L}_0) \cdots f^*(\overline{L}_l))_S & \\
& \kern-5em = \int_{\Omega} \Big(f_{\omega}^*(\operatorname{div}(s_0)_{\omega}, -\log |s_0|_{\varphi_\omega}) \cdots 
f_{\omega}^*(\operatorname{div}(s_d)_{\omega}, -\log |s_d|_{\varphi_\omega})\Big)_{\omega} \nu(\mathrm{d}\omega) \\
& \kern-5em = \deg(f) \int_{\Omega} \Big((\operatorname{div}(s_0)_{\omega}, -\log |s_0|_{\varphi_\omega}) \cdots 
(\operatorname{div}(s_d)_{\omega}, -\log |s_d|_{\varphi_\omega})\Big)_{\omega} \nu(\mathrm{d}\omega) \\
& \kern-5em = \deg(f) (\overline{L}_0 \cdots \overline{L}_l)_S.
\end{align*}
as required.
\end{proof}

In general, if we set $Z = a_1 Z_1 + \cdots + a_r Z_r$, then, by Claim~\ref{claim:prop:projection:formula:intersection:01},
\begin{multline*}
(f^*(\overline{L}_0) \cdots f^*(\overline{L}_l) \mid Z)_S =
\sum_{j=1}^r a_j (f^*(\overline{L}_0) \cdots f^*(\overline{L}_l) \mid Z_j)_S \\
=
\sum_{j=1}^r a_j (\overline{L}_0 \cdots \overline{L}_l \mid f_*(Z_j))_S =
(\overline{L}_0 \cdots \overline{L}_l \mid f_*(Z))_S.
\end{multline*}
\end{proof}

\section{Polarized adelic structure case}
Let $K$ be a finitely generated field over $\mathbb Q$ and $n$ be the transcendental degree of $K$ over $\QQ$.
Let $(\mathscr B; \overline{\mathscr H}_1, \ldots, \overline{\mathscr H}_n)$ be a polarization of $K$ and
$S = (K, (\Omega, \mathcal A, \nu), \phi)$ be the polarized adelic structure by 
$(\mathscr B; \overline{\mathscr H}_1, \ldots, \overline{\mathscr H}_n)$ (for details, see Section~\ref{sec:Polarized adelic structure}).

Let $X$ be a $d$-dimensional projective and integral scheme over $K$.
We choose a projective arithmetic variety $\mathscr X$ and a morphism $\pi : \mathscr X \to \mathscr B$
such that the generic fiber of $\mathscr X \to \mathscr B$ is $X$.
Let $L_0, \ldots, L_d$ be invertible $\OO_X$-modules.
We assume that there are $C^{\infty}$-metrized invertible $\OO_{\mathscr X}$-modules
$\overline{\mathscr L}_0 = (\mathscr L, h_0), \ldots, \overline{\mathscr L}_d = (\mathscr L_d, h_d)$ in the usual sense on arithmetic varieties such that
$\mathscr L_0, \ldots, \mathscr L_d$ coincides with $L_0, \ldots, L_d$ on $X$.
Note that, for each $\omega \in \Omega$, $\overline{\mathscr L}_i$ yields a smooth metric $\varphi_{i,\omega}$
of $L_{i,\omega}$, that is,
if $\omega \in \Omega_{\infty}$, then $\varphi_{i,\omega} = \rest{h_i}{\pi^{-1}(\omega)}$;
if $\omega \in \Omega \setminus \Omega_{\infty}$, then $\varphi_{i,\omega}$ is the model metric induced by
the model $(\mathscr X, \mathscr L_i)$.
We denote $\{(L_{i,\omega}, \varphi_{i,\omega})\}_{\omega \in \Omega}$ by $\overline{L}_i$.

\begin{prop}
$(\overline{L}_0\cdots \overline{L}_d)_S = (\overline{\mathscr L}_0 \cdots \overline{\mathscr L}_d \cdot
\pi^*(\overline{\mathscr H}_1) \cdots \pi^*(\overline{\mathscr H}_n))$
\end{prop}

\begin{proof}
We prove the assertion by induction on $d$.
Clearly we may assume that $\mathscr X$ is normal.
If $d=0$, that is, $\dim \mathscr X = n+1$, then it is an easy consequence of \cite[Lemma~1.12, Lemma~1.15, Proposition~5.3, Lemma~5.15 and Theorem~5.20]{MArakelov}.

We assume that $d > 0$. Let us choose a non-zero rational section $s_0$ of $\mathscr L_0$.
Let $\operatorname{div}(s_0) = a_1 \mathcal Z_1 + \cdots + a_r \mathcal Z_r$ be the decomposition as a cycle.
Then one has
\begin{multline*}
 (\overline{\mathscr L}_0 \cdots \overline{\mathscr L}_d \cdot
\pi^*(\overline{\mathscr H}_1) \cdots \pi^*(\overline{\mathscr H}_n)) \\
\kern-10em = \sum_{i=1}^r a_i (\overline{\mathscr L}_1 \cdots \overline{\mathscr L}_d \cdot
\pi^*(\overline{\mathscr H}_1) \cdots \pi^*(\overline{\mathscr H}_n) \cdot (\mathcal Z_i, 0)) \\
+ \int_{\mathscr X(\CC)} -\log |s_0|_{h_0} c_1(\overline{\mathscr L_1}) \wedge \cdots \wedge c_1(\overline{\mathscr L_d})
\wedge c_1(\pi^*\overline{\mathscr H_1}) \wedge \cdots \wedge c_1(\pi^*\overline{\mathscr H_n}).
\end{multline*}
Note that
\begin{multline*}
\int_{\mathscr X(\CC)} -\log |s_0|_{h_0} c_1(\overline{\mathscr L_1}) \wedge \cdots \wedge c_1(\overline{\mathscr L_d})
\wedge c_1(\pi^*\overline{\mathscr H_1}) \wedge \cdots \wedge c_1(\pi^*\overline{\mathscr H_n}) \\
= \int_{\mathscr B(\CC)} \left(\int_{\mathscr X(\CC)/\mathscr B(\CC)} -\log |s_0|_{h_0} c_1(\overline{\mathscr L_1}) \wedge \cdots \wedge c_1(\overline{\mathscr L_d}) \right)
c_1(\overline{\mathscr H_1}) \wedge \cdots \wedge c_1(\overline{\mathscr H_n}).
\end{multline*}
Here we consider the following claim:

\begin{enonce}{Claim}
Let $\psi : \mathscr Y \to \mathscr C$ be a surjective morphism of projective arithmetic varieties.
Let $\overline{\mathscr M}_1, \ldots, \overline{\mathscr M}_d$ 
(resp. $\overline{\mathscr D}_1, \ldots, \overline{\mathscr D}_n$) be metrized integrable 
invertible $\OO_{\mathscr Y}$-modules
(resp. $\OO_{\mathscr C}$-modules) such that $d + n = \dim \mathscr Y$.
Let $\mathscr Y_{\eta}$ be the generic fiber of $\psi : \mathscr Y \to \mathscr C$.
Then
\begin{multline*}
(\overline{\mathscr M}_1 \cdots \overline{\mathscr M}_d \cdot \pi^*\overline{\mathscr D}_1 \cdots 
\pi^*\overline{\mathscr D}_n) \\
= \begin{cases}
\big(\rest{\mathscr M_1}{\mathscr Y_{\eta}} \cdots \rest{\mathscr M_d}{\mathscr Y_{\eta}}\big) (\overline{\mathscr D}_1 \cdots \overline{\mathscr D}_n), & \text{if $d = \dim \mathscr Y_{\eta}$}, \\
0, & \text{if $d < \dim \mathscr Y_{\eta}$}.
\end{cases}
\end{multline*}
\end{enonce}

\begin{proof}
This is a consequence of the projection formula (cf. \cite[Theorem~5.20]{MArakelov}).
\end{proof}

By the above claim, if $\mathcal Z$ is a prime divisor on $\mathscr X$
with $\pi(\mathcal Z)\not=\mathscr B$, then
\begin{multline*}
\big(\overline{\mathscr L}_1 \cdots \overline{\mathscr L}_n \cdot
\pi^*(\overline{\mathscr H}_1) \cdots \pi^*(\overline{\mathscr H}_d) \cdot (\mathcal Z, 0)\big) \\
= \begin{cases}
\big(\rest{\mathscr L_1}{\mathcal Z_{\eta}} \cdots \rest{\mathscr L_n}{\mathcal Z_{\eta}}\big) (\overline{\mathscr H}_1 \cdots \overline{\mathscr H}_d \cdot (\pi(\mathcal Z), 0)), & \text{if $\codim(\pi(\mathcal Z); \mathscr B) = 1$},\\
0, & \text{if $\codim(\pi(\mathcal Z); \mathscr B) \geqslant 2$},
\end{cases}
\end{multline*}
where $\mathcal Z_{\eta}$ is the generic fiber of $\mathcal Z \to \pi(\mathcal Z)$. Therefore, if we set 
\[
\begin{cases}
I_{h} := \{ i \in \{ 1, \ldots, r \} \mid \pi(\mathcal Z_i) = \mathscr B \}, \\
I_{\Gamma} := \{ i \in \{ 1, \ldots, r \} \mid \pi(\mathcal Z_i) = \Gamma \}
\end{cases}
\]
for $\Gamma \in \Omega \setminus \Omega_{\infty}$, and denote 
\[
\sum_{i=1}^r a_i (\overline{\mathscr L}_1 \cdots \overline{\mathscr L}_d \cdot
\pi^*(\overline{\mathscr H}_1) \cdots \pi^*(\overline{\mathscr H}_n) \cdot (\mathcal Z_i, 0)) 
\]
by $T$, then,
by Example~\ref{Exe:chambert-loir measure} and hypothesis of induction on $d$, one has
{\allowdisplaybreaks\begin{align*}
T & = \sum_{i \in I_{h}} a_i (\overline{\mathscr L}_1 \cdots \overline{\mathscr L}_d \cdot
\pi^*(\overline{\mathscr H}_1) \cdots \pi^*(\overline{\mathscr H}_n) \cdot (\mathcal Z_i, 0)) \\
& \kern4em +
\sum_{\Gamma \in \Omega \setminus \Omega_{\infty}}
\sum_{i \in I_{\Gamma}} a_i (\overline{\mathscr L}_1 \cdots \overline{\mathscr L}_d \cdot
\pi^*(\overline{\mathscr H}_1) \cdots \pi^*(\overline{\mathscr H}_n) \cdot (\mathcal Z_i, 0)) \\
& = \sum_{i \in I_{h}} a_i \big(\rest{\overline{L}_1}{Z_i} \cdots \rest{\overline{L}_d}{Z_i}\big)_S \\
& \kern2em + \sum_{\Gamma \in \Omega \setminus \Omega_{\infty}} (\overline{\mathscr H}_1 \cdots \overline{\mathscr H}_d \cdot
(\Gamma, 0))
\int_{X_{\Gamma}^{\an}} -\log | s_0|_{\varphi_{0,\Gamma}} c_1(L_1, \varphi_{1,\Gamma}) \cdots c_1(L_d, \varphi_{d,\Gamma}),
\end{align*}}
where $Z_i$ is the generic fiber of $\mathcal Z_i \to \mathscr B$ for $i \in I_h$.
Thus, by \eqref{eqn:def:local:intersection},
\begin{multline*}
 (\overline{\mathscr L}_0 \cdots \overline{\mathscr L}_d \cdot
\pi^*(\overline{\mathscr H}_1) \cdots \pi^*(\overline{\mathscr H}_n)) = \sum_{i \in I_{h}} a_i \big(\rest{\overline{L}_1}{Z_i} \cdots \rest{\overline{L}_d}{Z_i}\big)_S \\
\kern3em + \int_{\Omega}
\left(\int_{X_{\omega}^{\an}} -\log | s_0|_{\varphi_{0,\omega}} c_1(L_1, \varphi_{1,\omega}) \cdots c_1(L_d, \varphi_{d,\omega})\right) \nu(\mathrm{d}\omega) = (\overline{L}_1 \cdots \overline{L}_d)_S,
\end{multline*}
as required.
\end{proof}


\appendix


\chapter{}
\section{Measurable family of embeddings}


In this appendix, we consider the following theorem, which was proved in \cite[Step~1 in Theorem~4.1.26]{CMArakelovAdelic}. We present here an alternative proof.

\begin{theo}\label{thm:measurable:family:embeddings}
Let $S=(K,(\Omega,\mathcal A,\nu),\phi)$ be an adelic curve such that $K$ is countable and $\Omega_\infty$ is not empty. There exists
 a family $( \iota_{\omega} )_{\omega \in \Omega_{\infty}}$ 
of embeddings $K \to \CC$ such that $|\ndot|_{\omega} = |\iota_{\omega}(\ndot)|$ for all
$\omega \in \Omega_{\infty}$ and that the map $(\omega \in \Omega_{\infty}) \mapsto
\iota_{\omega}(a)$ is measurable for each $a \in K$.
\end{theo}

Let 
$\CC^{\NN}$ be the set of all sequences $(x_n)_{n=0}^{\infty}$
consisting of complex numbers. For $x = (x_n)_{n=0}^{\infty} \in \CC^{\NN}$,
the $n$-th entry of $x_n$ of $x$ is often denoted by $x(n)$.
One can easily check the following proposition (see \cite[\S3.5]{Ovch} for (1), and (2) is straightforward). 

\begin{prop}\label{prop:measurability:C:N}
\begin{enumerate}
\renewcommand{\labelenumi}{\textup{(\arabic{enumi})}}
\item
If we define $d : \CC^{\NN} \times \CC^{\NN} \to \RR_{\geqslant 0}$ to be
\[
d(x, y) = \sum_{n=0}^{\infty} 2^{-n}\frac{|x(n) - y(n)|}{1 + |x(n) - y(n)|}
\quad (x, y \in \CC^{\NN}),
\]
then $d$ yields a distance function on $\CC^{\NN}$. Moreover,
the topology determined by $d$ coincides with the product topology, and
$(\CC^{\NN}, d)$ is a complete and second-countable space.

\item
Let $\mathcal{B}_{\CC^{\NN}}$ be the Borel $\sigma$-algebra by the product topology
on $\CC^{\NN}$ and $\mathcal{B}_{\CC}$ be the Borel $\sigma$-algebra  by
the standard topology on $\CC$.
Let $(\Omega, \mathcal{A})$ be a measurable space and
$( f_n )_{n=0}^{\infty}$ be a family of maps $\Omega \to \CC$. Then
$f : (\Omega, \mathcal{A}) \to (\CC^{\NN}, \mathcal{B}_{\CC^{\NN}})$ given by $f(\omega) = (f_n(\omega))_{n=0}^{\infty}$ is measurable
if and only if $f_n : (\Omega, \mathcal{A}) \to (\CC, \mathcal{B}_{\CC})$ is measurable for all $n \in \NN$.
\end{enumerate}
\end{prop}


\begin{lemm}\label{lemma:measurable:family:embeddings}
Let $S = (K, (\Omega, \mathcal{A}, \nu), \phi)$ be an adelic curve such that $\Omega_\infty$ is not empty. Suppose that $-1\in K$ admits a square root $\xi$ in $K$. Then there is a family $( \iota_{\omega} )_{\omega \in \Omega_{\infty}}$ of embeddings $K \to \CC$ which satisfy the following conditions:
\begin{enumerate}[label=\rm(\arabic*)]
\item for any $\omega\in\Omega_\infty$,
$\iota_{\omega}(\xi) = \sqrt{-1}$, 
\item for any $\omega \in \Omega_{\infty}$, $|\ndot|_{\omega} = |\iota_{\omega}(\ndot)|$,
\item for any $a \in K$, the function $(\omega \in \Omega_{\infty}) \mapsto
\iota_{\omega}(a)$ is measurable. 
\end{enumerate}
\end{lemm}

\begin{proof}
Fix a family $( \sigma_{\omega} )_{\omega \in \Omega_{\infty}}$ of embeddings $K \to \CC$ such that $|\ndot|_{\omega} = |\sigma_{\omega}(\ndot)|$ for all $\omega \in \Omega_{\infty}$.
Note that $\sigma_{\omega}(\xi) = \pm \sqrt{-1}$ because
\[(\sigma_{\omega}(\xi))^2 = \sigma_{\omega}(\xi^2) = \sigma_{\omega}(-1) = -1.\]
We define a family $( \iota_{\omega} )_{\omega \in \Omega_{\infty}}$ of embeddings 
by
\[
\iota_{\omega} = \begin{cases}
\sigma_{\omega} & \text{if $\sigma_{\omega}(\xi) = \sqrt{-1}$}, \\[1ex]
\overline{\sigma_{\omega}} & \text{if $\sigma_{\omega}(\xi) = -\sqrt{-1}$},
\end{cases}
\]
where $\overline{\phantom{x}}$ means the complex conjugation.
Then $\iota_{\omega}(\xi) = \sqrt{-1}$ for all $\omega \in \Omega_{\infty}$.
Thus one can see
\[
\iota_{\omega}(a) = \left(|a + (1/2)|^2_{\omega} - |a|^2_{\omega} - |1/2|_{\omega}^2\right)
+ \sqrt{-1} \left(|a + (\xi/2)|^2_{\omega} - |a|^2_{\omega} - |\xi/2|_{\omega}^2\right),
\]
as required.
\end{proof}

\begin{proof}[Proof of Theorem~\ref{thm:measurable:family:embeddings}]
By Lemma~\ref{lemma:measurable:family:embeddings}, we may assume that $-1\in K$ does not have any square root in $K$.
We set \[K = (a_n)_{n=1}^{\infty}, \quad L := K(\sqrt{-1})\quad \text{and}\quad
(L, (\Omega_L, \mathcal{A}_L, \nu_L), \phi_L) = S \otimes_K L.\]
Then, by Lemma~\ref{lemma:measurable:family:embeddings} again, there is a family $( \iota_{\chi} )_{\chi \in \Omega_{L,\infty}}$
of embeddings $L \to \CC$ such that $|\ndot|_{\chi} = |\iota_{\chi}(\ndot)|$ for all
$\chi \in \Omega_{L,\infty}$ and $(\chi \in \Omega_{L,\infty}) \mapsto \iota_{\chi}(b)$ is
measurable for all $b \in L$. Thus, if we define $h : \Omega_{L,\infty} \to \CC^{\NN}$
by $h(\chi) = (\iota_{\chi}(a_n))_{n=1}^{\infty}$, then $h$ is also measurable by (2) in Proposition~\ref{prop:measurability:C:N}, so that,
for an open set $U$ in $\CC^{\NN}$, $h^{-1}(U)$ is a measurable subset in  $\Omega_{L,\infty}$.
Let $\pi_{L/K} : \Omega_L \to \Omega$ be the canonical map. We consider a map
$F : \Omega_{\infty} \to \mathcal{P}(\CC^{\NN})$ given by
$F(\omega) = h(\pi_{L/K}^{-1}(\omega))$.
Then one can see that
\begin{align*}
\{ \omega \in \Omega_{\infty} \mid F(\omega) \cap U \not= \emptyset \} & =
\{ \omega \in \Omega_{\infty} \mid \pi^{-1}_{L/K}(\omega) \cap h^{-1}(U) \not= \emptyset \} \\
& = \{ \omega \in \Omega_{\infty} \mid
I_{L/K}(\indic_{h^{-1}(U)})(\omega) > 0 \}.
\end{align*}
By \cite[Theorem~3.3.4]{CMArakelovAdelic}, $I_{L/K}(\indic_{h^{-1}(U)})$ is measurable, so that
\[\{ \omega \in \Omega_{\infty} \mid F(\omega) \cap U \not= \emptyset \} \in \mathcal A.\]
Thus, by Kuratowski and Ryll-Nardzewski measurable selection theorem (cf. \cite[Theorem~A.2.1]{CMArakelovAdelic}) together with (1) in Proposition~\ref{prop:measurability:C:N}, there is a measurable map
$f : \Omega_{\infty} \to \CC^{\NN}$ such that $f(\omega) \in F(\omega)$ for all $\omega \in \Omega_{\infty}$.
For each $\omega \in \Omega_{\infty}$, we choose $\chi_{\omega} \in \pi_{L/K}^{-1}(\omega)$ such that $f(\omega) = h(\chi_{\omega})$.
If we set $\iota_{\omega} = \rest{\iota_{\chi_{\omega}}}{K}$, then $\iota_{\omega}$ yields an embedding $K \to \CC$
such that $|a|_{\omega} = |\iota_{\omega}(a)|$ for all $a \in K$.
Moreover, for all $n\in \NN$, $(\omega \in \Omega_{\infty}) \mapsto \iota_{\omega}(a_n)$
is measurable by (2) in Proposition~\ref{prop:measurability:C:N}. Thus the assertion follows.
\end{proof}

\if01
\begin{theo}\label{thm:measurable:family:embeddings}
If $\#(K) = \aleph_0$, then there is a family $\{ \iota_{\omega} \}_{\omega \in \Omega_{\infty}}$ 
of embeddings $K \to \CC$ such that $|\ndot|_{\omega} = |\iota_{\omega}(\ndot)|$ for all
$\omega \in \Omega_{\infty}$ and a map $(\omega \in \Omega_{\infty}) \mapsto
\iota_{\omega}(a)$ is measurable for each $a \in K$.
\end{theo}

\subsubsection{Topology of product space}
Let $\{ (X_n, d_n) \}_{n=1}^{\infty}$ be a family of metric spaces. 
Let \[\XXX := \prod_{n=1}^{\infty} X_n\] and 
$\mathscr{T}$ be the product topology of $\XXX$, that is,
the coarsest topology such that the projection $p_n : \XXX \to X_n$ to the $n$-th factor is continuous for all $n \in \NN$, where $\NN$ is the set of all positive integers.
For $x = (x_n)_{n=1}^{\infty} \in \XXX$, the $n$-th entry $x_n$ of $x$ is often denoted by $x(n)$.
Let $\OO_n$ be the set of all open sets in $X_n$ and
\begin{multline*}
\left(\prod\nolimits_{n=1}^{\infty} \OO_n\right)_* := \\
\left\{ (U_n)_{n=1}^{\infty} \in \prod\nolimits_{n=1}^{\infty} \OO_n \mid \text{$U_n = X_n$ except finitely many $n \in \NN$} \right\}.
\end{multline*}
For $(U_n)_{n=1}^{\infty} \in \left(\prod\nolimits_{n=1}^{\infty} \OO_n\right)_*$,
we set
\[
\UUU((U_n)_{n=1}^{\infty}) := \{ x \in X \mid \text{$x(n) \in U_n$ for all $n \in \NN$} \}
= \bigcap_{n=1}^{\infty} p_n^{-1}(U_n).
\]
Note that $\{ \UUU((U_n)_{n=1}^{\infty}) \}_{(U_n)_{n=1}^{\infty} \in \left(\prod\nolimits_{n=1}^{\infty} \OO_n\right)_*}$ forms an open basis of $\mathscr{T}$.
We define $d : \XXX \times \XXX \to \RR_{\geqslant 0}$ to be
\[
d(x, y) := \sum_{n=1}^{\infty} 2^{-n} \frac{d_n(x(n),y(n))}{1 + d_n(x(n),y(n))}\quad(x, y \in \XXX).
\]
Note that the function $t/(1+t)$ on $\RR_{\geqslant 0}$
is strictly increasing and subadditive, so that one can see that
the above $d$ yields a distance on $\XXX$.
For $r \in \RR_{> 0}$, $x \in \XXX$ and $x_n \in X_n$, we set
\[
\UUU(x;r) := \{ y \in X \mid d(x, y) < r\}\quad\text{and}\quad
U_n(x_n; r) := \{ y_n \in X_n \mid d_n(x_n, y_n) < r \}.
\]

\begin{prop}\label{prop:topology:product}
\begin{enumerate}
\renewcommand{\labelenumi}{\textup{(\arabic{enumi})}}
\item
The topology determined by the distance $d$ coincides with the product topology $\mathscr{T}$.

\item
If $(X_n, d_n)$ is complete for all $n \in \NN$, then
$(\XXX, d)$ is also complete.

\item We assume that $(X_n, d_n)$ is a second-countable space for all $n \in \NN$, that is,
there is $\OO'_n \subseteq \OO_n$ such that
$\OO'_n$ is a basis and $\#(\OO'_n) = \aleph_0$. 
If we set
\[
\kern2.2em \left(\prod\nolimits_{n=1}^{\infty} \OO_n\right)'_* = \left\{
(U_n)_{n=1}^{\infty} \in \left(\prod\nolimits_{n=1}^{\infty} \OO_n\right)_* \mid \text{$U_n \in \OO'_n$ or $U_n = X_n$} \right\},
\]
then $\#(\left(\prod\nolimits_{n=1}^{\infty} \OO_n\right)'_*) = \aleph_0$ and
$\{ \UUU((U_n)_{n=1}^{\infty}) \}_{(U_n)_{n=1}^{\infty} \in \left(\prod\nolimits_{n=1}^{\infty} \OO_n\right)'_*}$ forms a basis of $\mathscr{T}$.
In particular, $\XXX$ is also a second-countable space.
\end{enumerate}
\end{prop}

\begin{proof}
(1) Let $\mathscr{T}'$ be the topology determined by $d$.
It is sufficient to see the following (i) and (ii):
\begin{enumerate}
\renewcommand{\labelenumi}{\textup{(\roman{enumi})}}
\item The projection $p_n : \XXX \to X_n$ is continuous with respect to $\mathscr{T}'$ for
all $n \in \NN$.
In particular, $\mathscr{T} \subseteq \mathscr{T}'$.

\item For $\UUU(x;r)$ and $y \in \UUU(x;r)$,
there is $(U_n)_{n=1}^{\infty} \in \left(\prod\nolimits_{n=1}^{\infty} \OO_n\right)_*$
such that $y \in \UUU((U_n)_{n=1}^{\infty}) \subseteq \UUU(x;r)$. In particular,
$\mathscr{T}' \subseteq \mathscr{T}$.
\end{enumerate}

(i) follows from
\begin{equation}\label{eqn:prop:topology:product:01}
\frac{d_n(x(n), y(n))}{1 + d_n(x(n), y(n))} \leqslant 2^n d(x, y)
\quad (\forall x, y \in X).
\end{equation}

(ii) As $\UUU(y; r - d(x,y)) \subseteq \UUU(x;r)$, we may assume that $y=x$.
We choose $N \in \NN$ and $r' \in \RR_{>0}$ such that $\sum_{n=N+1}^{\infty} 2^{-n} < r/2$ and $\sum_{n=1}^{N} 2^{-n} r'/(1 + r') < r/2$. 
If we set
\[
U_n := \begin{cases} U_n(x(n);r') & \text{if $n \in \{ 1,\ldots,N\}$},\\
X_n & \text{if $n > N$},
\end{cases}
\]
then $x \in \UUU((U_n)_{n=1}^{\infty})\subseteq \UUU(x;r)$. Indeed, if $y \in \UUU((U_n)_{n=1}^{\infty})$, then
\[
d(y, x) = \sum_{n=1}^{\infty} 2^{-n}\frac{d_n(y(n),x(n))}{1+d_n(y(n),x(n))} <
\sum_{n=1}^N 2^{-n} \frac{r'}{1+ r'} + \sum_{n=N+1}^{\infty} 2^{-n} < r.
\]

\medskip
(2) Let $\{ x_m \}_{m=1}^{\infty}$ be a Cauchy sequence in $(\XXX, d)$.
Then, by \eqref{eqn:prop:topology:product:01}, 
$\{ x_m(n) \}_{m=1}^{\infty}$ is also a Cauchy sequence in $X_n$,
so that there is $x_n \in X_n$ such that \[\lim_{m\to\infty} d_n(x_m(n), x_n) = 0.\]
Let $x \in \XXX$ such that $x(n) = x_n$ for all $n \in \NN$. For any $\varepsilon > 0$,
by the above (ii), there is $(U_n)_{n=1}^{\infty} \in \left(\prod\nolimits_{n=1}^{\infty} \OO_n\right)_*$
such that $x \in \UUU((U_n)_{n=1}^{\infty}) \subseteq \UUU(x;\varepsilon)$.
If we set $I = \{ n \in \NN \mid U_n \not= X_n \}$,
then there is $m_0 > 0$ such that
$x_m(n) \in U_n$ for all $n \in I$ and $m \geqslant m_0$.
Thus $d(x_m, x) < \varepsilon$,
so that $\lim_{m\to\infty} d(x_m, x) = 0$, as required.

\medskip
(3) is obvious.
\end{proof}

Let $\mathcal{B}_{X_n}$ and $\mathcal{B}_{\XXX}$ be the Borel $\sigma$-algebras
of $X_n$ and $\XXX$, respectively.

\begin{prop}\label{prop:measurability:product}
We assume that $(X_n, d_n)$ is a second-countable space for all $n \in \NN$.
\begin{enumerate}
\renewcommand{\labelenumi}{\textup{(\arabic{enumi})}}
\item
$\mathcal{B}_{\XXX}$ is the smallest $\sigma$-algebra such that $p_n : (\XXX, \mathcal{B}_{\XXX}) \to (X_n,  \mathcal{B}_{X_n})$
is measurable for all $n \in \NN$.

\item
Let $(\Omega, \mathcal{A})$ be a measurable space and
$f_n : \Omega \to X_n$ be a map for each $n \in \NN$. Then
$f : (\Omega, \mathcal{A}) \to (\XXX, \mathcal{B}_{\XXX})$ given by $f(\omega) = 
(f_n(\omega))_{n=1}^{\infty}$ is measurable
if and only if $f_n : (\Omega, \mathcal{A}) \to (X_n, \mathcal{B}_{X_n})$ is measurable for all $n \in \NN$.
\end{enumerate}
\end{prop}

\begin{proof}
(1)
Let $\mathcal{B}'_{\XXX}$ be the smallest $\sigma$-algebra such that $p_n : (\XXX, \mathcal{B}'_{\XXX}) \to (X_n,  \mathcal{B}_{X_n})$
is measurable for all $n \in \NN$. Obviously $\mathcal{B}'_{\XXX} \subseteq \mathcal{B}_{\XXX}$ because
$p_n : \XXX \to X_n$
is continuous for all $n \in \NN$. For its converse, we need to show that $\UUU \in \mathcal{B}'_{\XXX}$ for all open sets $\UUU$ of $\XXX$.
By (3) in Proposition~\ref{prop:topology:product}, we may assume that $\UUU = \UUU((U_n)_{n=1}^{\infty})$ for some $(U_n)_{n=1}^{\infty} \in 
(\prod\nolimits_{n=1}^{\infty} \OO_n)'_*$.
Then \[ \UUU = \bigcap_{n=1}^{\infty} p_n^{-1}(U_n) \in \mathcal{B}'_{\XXX},\] as required.

\medskip
(2) is a consequence of (1).
\end{proof}

As a corollary of Proposition~\ref{prop:topology:product} and Proposition~\ref{prop:measurability:product}. one has the following.

\begin{Corollary}\label{cor:measurability:C:N}
\begin{enumerate}
\renewcommand{\labelenumi}{\textup{(\arabic{enumi})}}
\item
If we define $d : \CC^{\NN} \times \CC^{\NN} \to \RR_{\geqslant 0}$ to be
\[
d(x, y) = \sum_{n=1}^{\infty} 2^{-n}\frac{|x(n) - y(n)|}{1 + |x(n) - y(n)|}
\quad (x, y \in \CC^{\NN}),
\]
then the topology determined by $d$ coincides with the product topology, and
$(\CC^{\NN}, d)$ is a complete and second-countable space.
In particular, $(\CC^{\NN}, d)$ is separable.

\item
Let $\mathcal{B}_{\CC^{\NN}}$ be the Borel $\sigma$-algebra by the product topology
on $\CC^{\NN}$ and $\mathcal{B}_{\CC}$ be the Borel $\sigma$-algebra  by
the standard topology on $\CC$.
Let $(\Omega, \mathcal{A})$ be a measurable space and
$\{ f_n \}_{n \in \ZZ_{\geqslant 1}}$ be a family of maps $\Omega \to \CC$. Then
$f : (\Omega, \mathcal{A}) \to (\CC^{\NN}, \mathcal{B}_{\CC^{\NN}})$ given by $f(\omega) = (f_n(\omega))_{n=1}^{\infty}$ is measurable
if and only if $f_n : (\Omega, \mathcal{A}) \to (\CC, \mathcal{B}_{\CC})$ is measurable for all $n \in \NN$.
\end{enumerate}
\end{Corollary}

\subsubsection{Proof of Theorem~\ref{thm:measurable:family:embeddings}}

Let us begin with the following lemma.

\beginlemm\label{lemma:measurable:family:embeddings}
Let $S = (K, (\Omega, \mathcal{A}, \nu), \phi)$ be an adelic structure of a field $K$.
If $\sqrt{-1} \in K$, then there is a family $\{ \iota_{\omega} \}_{\omega \in \Omega_{\infty}}$ of embeddings $K \to \CC$ such that
$\iota_{\omega}(\sqrt{-1}) = \sqrt{-1}$ and $|\ndot|_{\omega} = |\iota_{\omega}(\ndot)|$
for all $\omega \in \Omega_{\infty}$, and that
a map $(\omega \in \Omega_{\infty}) \mapsto
\iota_{\omega}(a)$ is measurable for all $a \in K$.
\endlemm

\begin{proof}
Fix a family $\{ \sigma_{\omega} \}_{\omega \in \Omega_{\infty}}$ of embeddings $K \to \CC$ such that $|\ndot|_{\omega} = |\sigma_{\omega}(\ndot)|$ for all $\omega \in \Omega_{\infty}$.
Note that $\sigma_{\omega}(\sqrt{-1}) = \pm \sqrt{-1}$ because
$(\sigma_{\omega}(\sqrt{-1}))^2 = \sigma_{\omega}(\sqrt{-1}^2) = \sigma_{\omega}(-1) = -1$.
We define a family $\{ \iota_{\omega} \}_{\omega \in \Omega_{\infty}}$ of embeddings 
by
\[
\iota_{\omega} = \begin{cases}
\sigma_{\omega} & \text{if $\sigma_{\omega}(\sqrt{-1}) = \sqrt{-1}$}, \\[1ex]
\overline{\sigma_{\omega}} & \text{if $\sigma_{\omega}(\sqrt{-1}) = -\sqrt{-1}$},
\end{cases}
\]
where $\overline{\phantom{x}}$ means the complex conjugation.
Then $\iota_{\omega}(\sqrt{-1}) = \sqrt{-1}$ for all $\omega \in \Omega_{\infty}$.
Thus one can see
\begin{multline*}
\iota_{\omega}(a) = \left(|a + (1/2)|^2_{\omega} - |a|^2_{\omega} - |1/2|_{\omega}^2\right) \\
+ \sqrt{-1} \left(|a + (\sqrt{-1}/2)|^2_{\omega} - |a|^2_{\omega} - |\sqrt{-1}/2|_{\omega}^2\right),
\end{multline*}
as required.
\end{proof}

\begin{proof}[Proof of Theorem~\ref{thm:measurable:family:embeddings}]
By Lemma~\ref{lemma:measurable:family:embeddings}, we may assume that $\sqrt{-1} \not\in K$.
We set \[K = \{ a_n \}_{n=1}^{\infty}, \quad L := K(\sqrt{-1})\quad \text{and}\quad
(L, (\Omega_L, \mathcal{A}_L, \nu_L), \phi_L) = S \otimes_K L.\]
Then, by Lemma~\ref{lemma:measurable:family:embeddings} again, there is a family $\{ \iota_{\chi} \}_{\chi \in \Omega_{L,\infty}}$
of embeddings $L \to \CC$ such that $|\ndot|_{\chi} = |\iota_{\chi}(\ndot)|$ for all
$\chi \in \Omega_{L,\infty}$ and $(\chi \in \Omega_{L,\infty}) \mapsto \iota_{\chi}(b)$ is
measurable for all $b \in L$. Thus, if we define $h : \Omega_{L,\infty} \to \CC^{\NN}$
by $h(\chi) = (\iota_{\chi}(a_n))_{n=1}^{\infty}$, then $h$ is also measurable by (2) in Corollary~\ref{cor:measurability:C:N}, so that,
for an open set $U$ in $\CC^{\NN}$, $h^{-1}(U)$ is a measurable subset in  $\Omega_{L,\infty}$.
Let $\pi_{L/K} : \Omega_L \to \Omega$ be the canonical map. We consider a map
$F : \Omega_{\infty} \to \mathcal{P}(\CC^{\NN})$ given by
\[
F(\omega) = h(\pi_{L/K}^{-1}(\omega)).
\]
Then one can see that
\begin{align*}
\{ \omega \in \Omega_{\infty} \mid F(\omega) \cap U \not= \emptyset \} & =
\{ \omega \in \Omega_{\infty} \mid \pi^{-1}_{L/K}(\omega) \cap h^{-1}(U) \not= \emptyset \} \\
& = \{ \omega \in \Omega_{\infty} \mid
I_{L/K}(\indic_{h^{-1}(U)})(\omega) > 0 \}.
\end{align*}
By \cite[Theorem~3.3.4]{CMArakelovAdelic}, $I_{L/K}(\indic_{h^{-1}(U)})$ is measurable, so that
\[\{ \omega \in \Omega_{\infty} \mid F(\omega) \cap U \not= \emptyset \} \in \mathcal A.\]
Thus, by Kuratowski and Ryll-Nardzewski measurable selection theorem (cf. \cite[Theorem~A.2.1]{CMArakelovAdelic}) together with (1) in Corollary~\ref{cor:measurability:C:N}, there is a measurable map
$f : \Omega_{\infty} \to \CC^{\NN}$ such that $f(\omega) \in F(\omega)$ for all $\omega \in \Omega_{\infty}$.
For each $\omega \in \Omega_{\infty}$, we choose $\chi_{\omega} \in \pi_{L/K}^{-1}(\omega)$ such that $f(\omega) = h(\chi_{\omega})$.
If we set $\iota_{\omega} = \rest{\iota_{\chi_{\omega}}}{K}$, then $\iota_{\omega}$ yields an embedding $K \to \CC$
such that $|a|_{\omega} = |\iota_{\omega}(a)|$ for all $a \in K$.
Moreover, for all $n\in \NN$, $(\omega \in \Omega_{\infty}) \mapsto \iota_{\omega}(a_n)$
is measurable by (2) in Corollary~\ref{cor:measurability:C:N}. Thus the assertion follows.
\end{proof}
\fi

\backmatter
\bibliography{intersection}
\bibliographystyle{smfplain_nobysame}

\printindex

\end{document}